\documentclass[a4paper,10.5pt, reqno, oneside]{amsart}

\usepackage{amssymb}
\usepackage{amstext}
\usepackage{amsmath}
\usepackage{amscd}
\usepackage{amsthm}
\usepackage{amsfonts}
\usepackage{enumerate}
\usepackage{graphicx}
\usepackage{latexsym}
\usepackage{mathrsfs}
\usepackage{mathtools}
\usepackage{tikz}
\usetikzlibrary{arrows}
\usetikzlibrary{snakes}
\usetikzlibrary{shapes}
\usetikzlibrary{calc,decorations.markings}

\usepackage{tikz-cd}

\usepackage{hyperref}

\usepackage{lscape}
\usepackage{comment}

\makeatletter
  \renewcommand{\theequation}{%
  \thesection.\arabic{equation}}
  \@addtoreset{equation}{section}
\makeatother

\newtheorem{theorem}{Theorem}[section]
\newtheorem{lemma}[theorem]{Lemma}
\newtheorem{proposition}[theorem]{Proposition}
\newtheorem{corollary}[theorem]{Corollary}

\newtheorem{question}[theorem]{Question}

\renewcommand{\theac}{ }

\theoremstyle{definition}
\newtheorem{definition}[theorem]{Definition}
\newtheorem{definiton-theorem}[theorem]{Definition-Theorem}
\newtheorem{lemma-definition}[theorem]{Lemma-Definition}
\newtheorem{example}[theorem]{Example}

\theoremstyle{remark}
\newtheorem{remark}[theorem]{Remark}

\newcommand{\rma}{\mathrm{a}}
\newcommand{\rmb}{\mathrm{b}}
\newcommand{\rmc}{\mathrm{c}}
\newcommand{\rmd}{\mathrm{d}}
\newcommand{\rme}{\mathrm{e}}
\newcommand{\rmf}{\mathrm{f}}
\newcommand{\rmg}{\mathrm{g}}
\newcommand{\rmh}{\mathrm{h}}
\newcommand{\rmi}{\mathrm{i}}
\newcommand{\rmj}{\mathrm{j}}
\newcommand{\rmk}{\mathrm{k}}
\newcommand{\rml}{\mathrm{l}}
\newcommand{\rmm}{\mathrm{m}}
\newcommand{\rmn}{\mathrm{n}}
\newcommand{\rmo}{\mathrm{o}}
\newcommand{\rmp}{\mathrm{p}}
\newcommand{\rmq}{\mathrm{q}}
\newcommand{\rmr}{\mathrm{r}}
\newcommand{\rms}{\mathrm{s}}
\newcommand{\rmt}{\mathrm{t}}
\newcommand{\rmu}{\mathrm{u}}
\newcommand{\rmv}{\mathrm{v}}
\newcommand{\rmw}{\mathrm{w}}
\newcommand{\rmx}{\mathrm{x}}
\newcommand{\rmy}{\mathrm{y}}
\newcommand{\rmz}{\mathrm{z}}

\newcommand{\rmA}{\mathrm{A}}
\newcommand{\rmB}{\mathrm{B}}
\newcommand{\rmC}{\mathrm{C}}
\newcommand{\rmD}{\mathrm{D}}
\newcommand{\rmE}{\mathrm{E}}
\newcommand{\rmF}{\mathrm{F}}
\newcommand{\rmG}{\mathrm{G}}
\newcommand{\rmH}{\mathrm{H}}
\newcommand{\rmI}{\mathrm{I}}
\newcommand{\rmJ}{\mathrm{J}}
\newcommand{\rmK}{\mathrm{K}}
\newcommand{\rmL}{\mathrm{L}}
\newcommand{\rmM}{\mathrm{M}}
\newcommand{\rmN}{\mathrm{N}}
\newcommand{\rmO}{\mathrm{O}}
\newcommand{\rmP}{\mathrm{P}}
\newcommand{\rmQ}{\mathrm{Q}}
\newcommand{\rmR}{\mathrm{R}}
\newcommand{\rmS}{\mathrm{S}}
\newcommand{\rmT}{\mathrm{T}}
\newcommand{\rmU}{\mathrm{U}}
\newcommand{\rmV}{\mathrm{V}}
\newcommand{\rmW}{\mathrm{W}}
\newcommand{\rmX}{\mathrm{X}}
\newcommand{\rmY}{\mathrm{Y}}
\newcommand{\rmZ}{\mathrm{Z}}

\newcommand{\cala}{\mathcal{a}}
\newcommand{\calb}{\mathcal{b}}
\newcommand{\calc}{\mathcal{c}}
\newcommand{\cald}{\mathcal{d}}
\newcommand{\cale}{\mathcal{e}}
\newcommand{\calf}{\mathcal{f}}
\newcommand{\calg}{\mathcal{g}}
\newcommand{\calh}{\mathcal{h}}
\newcommand{\cali}{\mathcal{i}}
\newcommand{\calj}{\mathcal{j}}
\newcommand{\calk}{\mathcal{k}}
\newcommand{\call}{\mathcal{l}}
\newcommand{\calm}{\mathcal{m}}
\newcommand{\caln}{\mathcal{n}}
\newcommand{\calo}{\mathcal{o}}
\newcommand{\calp}{\mathcal{p}}
\newcommand{\calq}{\mathcal{q}}
\newcommand{\calr}{\mathcal{r}}
\newcommand{\cals}{\mathcal{s}}
\newcommand{\calt}{\mathcal{t}}
\newcommand{\calu}{\mathcal{u}}
\newcommand{\calv}{\mathcal{v}}
\newcommand{\calw}{\mathcal{w}}
\newcommand{\calx}{\mathcal{x}}
\newcommand{\caly}{\mathcal{y}}
\newcommand{\calz}{\mathcal{z}}

\newcommand{\calA}{\mathcal{A}}
\newcommand{\calB}{\mathcal{B}}
\newcommand{\calC}{\mathcal{C}}
\newcommand{\calD}{\mathcal{D}}
\newcommand{\calE}{\mathcal{E}}
\newcommand{\calF}{\mathcal{F}}
\newcommand{\calG}{\mathcal{G}}
\newcommand{\calH}{\mathcal{H}}
\newcommand{\calI}{\mathcal{I}}
\newcommand{\calJ}{\mathcal{J}}
\newcommand{\calK}{\mathcal{K}}
\newcommand{\calL}{\mathcal{L}}
\newcommand{\calM}{\mathcal{M}}
\newcommand{\calN}{\mathcal{N}}
\newcommand{\calO}{\mathcal{O}}
\newcommand{\calP}{\mathcal{P}}
\newcommand{\calQ}{\mathcal{Q}}
\newcommand{\calR}{\mathcal{R}}
\newcommand{\calS}{\mathcal{S}}
\newcommand{\calT}{\mathcal{T}}
\newcommand{\calU}{\mathcal{U}}
\newcommand{\calV}{\mathcal{V}}
\newcommand{\calW}{\mathcal{W}}
\newcommand{\calX}{\mathcal{X}}
\newcommand{\calY}{\mathcal{Y}}
\newcommand{\calZ}{\mathcal{Z}}

\newcommand{\fka}{\mathfrak{a}}
\newcommand{\fkb}{\mathfrak{b}}
\newcommand{\fkc}{\mathfrak{c}}
\newcommand{\fkd}{\mathfrak{d}}
\newcommand{\fke}{\mathfrak{e}}
\newcommand{\fkf}{\mathfrak{f}}
\newcommand{\fkg}{\mathfrak{g}}
\newcommand{\fkh}{\mathfrak{h}}
\newcommand{\fki}{\mathfrak{i}}
\newcommand{\fkj}{\mathfrak{j}}
\newcommand{\fkk}{\mathfrak{k}}
\newcommand{\fkl}{\mathfrak{l}}
\newcommand{\fkm}{\mathfrak{m}}
\newcommand{\fkn}{\mathfrak{n}}
\newcommand{\fko}{\mathfrak{o}}
\newcommand{\fkp}{\mathfrak{p}}
\newcommand{\fkq}{\mathfrak{q}}
\newcommand{\fkr}{\mathfrak{r}}
\newcommand{\fks}{\mathfrak{s}}
\newcommand{\fkt}{\mathfrak{t}}
\newcommand{\fku}{\mathfrak{u}}
\newcommand{\fkv}{\mathfrak{v}}
\newcommand{\fkw}{\mathfrak{w}}
\newcommand{\fkx}{\mathfrak{x}}
\newcommand{\fky}{\mathfrak{y}}
\newcommand{\fkz}{\mathfrak{z}}

\newcommand{\fkA}{\mathfrak{A}}
\newcommand{\fkB}{\mathfrak{B}}
\newcommand{\fkC}{\mathfrak{C}}
\newcommand{\fkD}{\mathfrak{D}}
\newcommand{\fkE}{\mathfrak{E}}
\newcommand{\fkF}{\mathfrak{F}}
\newcommand{\fkG}{\mathfrak{G}}
\newcommand{\fkH}{\mathfrak{H}}
\newcommand{\fkI}{\mathfrak{I}}
\newcommand{\fkJ}{\mathfrak{J}}
\newcommand{\fkK}{\mathfrak{K}}
\newcommand{\fkL}{\mathfrak{L}}
\newcommand{\fkM}{\mathfrak{M}}
\newcommand{\fkN}{\mathfrak{N}}
\newcommand{\fkO}{\mathfrak{O}}
\newcommand{\fkP}{\mathfrak{P}}
\newcommand{\fkQ}{\mathfrak{Q}}
\newcommand{\fkR}{\mathfrak{R}}
\newcommand{\fkS}{\mathfrak{S}}
\newcommand{\fkT}{\mathfrak{T}}
\newcommand{\fkU}{\mathfrak{U}}
\newcommand{\fkV}{\mathfrak{V}}
\newcommand{\fkW}{\mathfrak{W}}
\newcommand{\fkX}{\mathfrak{X}}
\newcommand{\fkY}{\mathfrak{Y}}
\newcommand{\fkZ}{\mathfrak{Z}}

\newcommand{\NN}{\mathbb{N}}
\newcommand{\ZZ}{\mathbb{Z}}
\newcommand{\QQ}{\mathbb{Q}}
\newcommand{\RR}{\mathbb{R}}
\newcommand{\CC}{\mathbb{C}}
\newcommand{\FF}{\mathbb{F}}

\newcommand{\sfA}{\mathsf{A}}
\newcommand{\sfB}{\mathsf{B}}
\newcommand{\sfC}{\mathsf{C}}
\newcommand{\sfD}{\mathsf{D}}
\newcommand{\sfE}{\mathsf{E}}

\newcommand{\sfU}{\mathsf{U}}
\newcommand{\ttI}{\mathtt{I}}
\newcommand{\sfM}{\mathsf{M}}
\newcommand{\sfN}{\mathsf{N}}
\newcommand{\sfP}{\mathsf{P}}
\newcommand{\sfR}{\mathsf{R}}
\newcommand{\sfT}{\mathsf{T}}

\newcommand{\Min}{\operatorname{Min}}
\newcommand{\height}{\operatorname{ht}}
\newcommand{\codim}{\operatorname{codim}}
\newcommand{\ind}{\operatorname{ind}}
\newcommand{\Spec}{\operatorname{Spec}}
\newcommand{\Soc}{\operatorname{Soc}}
\newcommand{\soc}{\operatorname{soc}}
\newcommand{\Ass}{\operatorname{Ass}}
\newcommand{\depth}{\operatorname{depth}}
\newcommand{\grade}{\operatorname{grade}}
\newcommand{\rad}{\operatorname{rad}}
\newcommand{\rank}{\operatorname{rank}}
\newcommand{\Ext}{\operatorname{Ext}}
\newcommand{\Tor}{\operatorname{Tor}}
\newcommand{\tor}{\operatorname{tor}}
\newcommand{\Hom}{\operatorname{Hom}}
\newcommand{\id}{\operatorname{id}}
\newcommand{\pd}{\operatorname{pd}}
\newcommand{\Image}{\operatorname{Im}}
\newcommand{\Tr}{\operatorname{Tr}}
\newcommand{\Ker}{\operatorname{Ker}}
\newcommand{\Coker}{\operatorname{Coker}}
\newcommand{\GL}{\operatorname{GL}}
\newcommand{\SL}{\operatorname{SL}}
\newcommand{\PSL}{\operatorname{PSL}}
\newcommand{\Sp}{\operatorname{Sp}}
\newcommand{\SO}{\operatorname{SO}}
\newcommand{\Sym}{\operatorname{Sym}}
\newcommand{\charac}{\operatorname{char}}
\newcommand{\pr}{\operatorname{pr}}
\newcommand{\Pic}{\operatorname{Pic}}
\newcommand{\Cl}{\operatorname{Cl}}
\newcommand{\Proj}{\operatorname{Proj}}
\newcommand{\Div}{\operatorname{Div}}
\newcommand{\Irr}{\operatorname{Irr}}
\newcommand{\End}{\operatorname{End}}
\newcommand{\gldim}{\mathrm{gl.dim}}
\newcommand{\projdim}{\mathrm{proj.dim}}
\newcommand{\twoheadlongrightarrow}{\relbar\joinrel\twoheadrightarrow}

\newcommand{\add}{\mathsf{add}}
\newcommand{\CM}{\mathsf{CM}}
\newcommand{\MCM}{\mathsf{MCM}}
\newcommand{\Mod}{\mathsf{Mod}}
\newcommand{\mc}{\mathsf{mod}}
\newcommand{\refl}{\mathsf{ref}}
\newcommand{\coh}{\mathsf{coh}}
\newcommand{\proj}{\mathsf{proj}}
\newcommand{\EG}{\mathsf{EG}}
\newcommand{\MM}{\operatorname{MM}}
\newcommand{\TMM}{\operatorname{MM}_1}
\newcommand{\MMG}{\operatorname{MMG}}
\newcommand{\TMMG}{\operatorname{MMG}_1}

\newcommand{\new}[1]{{\blue #1}}
\newcommand{\old}[1]{{\red #1}}
\renewcommand{\comment}[1]{{\green #1}}

\begin{document}
%\tikzset{auto}

%%%%%---title---%%%%%
\title[Mutations of splitting maximal modifying modules]{Mutations of splitting maximal modifying modules:\\ the case of reflexive polygons}
\author[Yusuke Nakajima]{Yusuke Nakajima}
\date{}


\address[Y. Nakajima]{Kavli Institute for the Physics and Mathematics of the Universe (WPI), UTIAS, The University of Tokyo, Kashiwa, Chiba 277-8583, Japan} 
\email{yusuke.nakajima@ipmu.jp}

%\address[Yusuke Nakajima]{Graduate School Of Mathematics, Nagoya University, Chikusa-Ku, Nagoya, 464-8602 Japan
% \footnote{Present address: Kavli Institute for the Physics and Mathematics of the Universe (WPI), UTIAS, The University of Tokyo, Kashiwa, Chiba 277-8583, Japan}} 
%\email{m06022z@math.nagoya-u.ac.jp
%\footnote{Present email: yusuke.nakajima@ipmu.jp}}

\maketitle

%%%---abstract---%%%
\begin{abstract}
It is known that every three dimensional Gorenstein toric singularity has a crepant resolution. 
Although it is not unique, all crepant resolutions are connected by repeating the operation ``flop". 
On the other hand, this singularity also has a non-commutative crepant resolution (= NCCR) which is constructed from a consistent dimer model. 
Such an NCCR is given as the endomorphism ring of a certain module which we call splitting maximal modifying module. 
In this paper, we show that all splitting maximal modifying modules are connected by repeating the operation ``mutation" of 
splitting maximal modifying modules for the case of toric singularities associated with reflexive polygons. 
\end{abstract}

\setcounter{tocdepth}{1}
\tableofcontents

%%%%%---text_start---%%%%%
\section{Introduction} 

%%%%%%%%%%%%%%%%%%%%%%%%%%%%%%%%%%%%%%%%%%%%%%%%%%%%%%%%%%%%%%%%%%%%%%

\subsection{Backgrounds}
Bondal and Orlov conjectured that any bounded derived categories of coherent sheaves over crepant resolutions of a certain singularity are derived equivalent \cite{BO}. 
Bridgeland proved this conjecture for the three dimensional case \cite{Bri}. 
Also, Van den Bergh proved it in more general settings containing the three dimensional case \cite{VdB1}. 
The idea in \cite{VdB1} is to consider a non-commutative ring $\Lambda$ derived equivalent to crepant resolutions. 
In \cite{VdB2}, Van den Bergh formulated such a ring $\Lambda$ as a non-commutative crepant resolution (= NCCR), 
which is obtained as the endomorphism ring of a reflexive module. (For the precise definition, see Definition~\ref{NCCR_def}.) 
Note that an NCCR does not necessarily exist, and even if it exists it is not unique (sometimes NCCRs are infinite families).  
On the other hand, it is known that the non-commutative version of Bondal-Orlov conjecture holds for some singularities \cite{IW1,IR}, 
that is, if we assume that $\End_R(M)$ and $\End_R(N)$ are NCCRs of three dimensional Cohen-Macaulay normal domain $R$, 
then they are derived equivalent. In this situation, is there a relationship between $M$ and $N$ ? 
The present paper is dedicated to investigate this question for the case of three dimensional Gorenstein toric singularities. 
In the rest of this paper, we call the ordinary crepant resolution commutative crepant resolution (= CCR) to emphasize the difference with an NCCR.    

%%%%%%%%%%%%%%%%%%%%%%%%%%%%%%%%%%%%%%%%%%%%%%%%%%%%%%%%%%%%%%%%%%%%%%
\subsection{The case of toric singularities} 

Let $R$ be a three dimensional Gorenstein toric singularity associated with a cone $\sigma\subset\RR^3$ generated by $v_1,\cdots,v_n\in\ZZ^3$. 
Since $R$ is Gorenstein, we can take a hyperplane so that any vectors $v_1,\cdots,v_n$ lie on this hyperplane. 
Thus, we obtain the lattice polygon $\Delta\subset\RR^2$ on this hyperplane for a given $R$ (see subsection~\ref{toric_pre} and \ref{subsec_dimer}). 
It is known that CCRs of $R$ correspond to triangulations of the associated polygon $\Delta$. 
Thus, a CCR of $R$ exists but it is not unique in general. More precisely, if there is a quadrangle consisting of two elementary triangles in a given triangulation, 
we obtain another triangulation by switching the diagonal as in the figure below, and it induces another CCR. 
This operation is called flop and any CCRs are connected by repeating this operation. 

\medskip

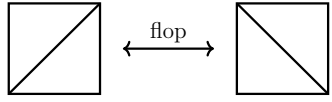
\begin{center}
\scalebox{0.6}{
\begin{tikzpicture}
\draw[line width=0.05cm]  (0,0) rectangle (2,2);
\draw[line width=0.05cm]  (5,0) rectangle (7,2);

\draw [line width=0.05cm] (0,0)--(2,2); \draw [line width=0.05cm] (5,2)--(7,0);

\draw[<->, line width=0.05cm] (2.5,1)--node[above] {{\LARGE{flop}}}(4.5,1);

\end{tikzpicture} }
\end{center}

\smallskip 

On the other hand, every three dimensional Gorenstein toric singularity has an NCCR, and it is obtained from a dimer model. 
A dimer model is a polygonal cell decomposition of the real two-torus whose vertices and edges form a finite bipartite graph (see subsection~\ref{subsec_dimer}). 
It is introduced in the field of statistical mechanics in 1960s. 
From 2000s, string theorists have been used it for studying quiver gauge theories. 
For more details, see e.g., \cite{Kenn,Keny} and references therein. 
Subsequently, a dimer model has been investigated actively, and recently relations with many branches of mathematics 
(for example, the McKay correspondence, crepant resolutions, non-commutative crepant resolutions, Calabi-Yau algebras, mirror symmetry, etc) have been discovered. 
The important point is that we can construct a three dimensional Gorenstein toric singularity $R$ by using a dimer model (see subsection~\ref{subsec_dimer}). 
On the other hand, we obtain a quiver with potential (= QP) as the dual of a dimer model (see subsection~\ref{sec_QP_from_dimer}). 
Under a certain condition, we can construct an NCCR of $R$ by using a quiver with potential associated with a dimer model, 
and for every three dimensional Gorenstein toric singularity, there exists a dimer model which gives an NCCR of it (see subsection~\ref{sec_consist}). 
Note that a module giving such an NCCR is called a splitting maximal modifying (= MM) module (see Definition~\ref{def_splitting}). 
Thus, an NCCR of $R$ always exists, but a dimer model which gives an NCCR of $R$ is not unique in general, 
hence a splitting MM module is also not unique. 
Although it is not unique, the number of NCCRs arising from dimer models is finite (see Section~\ref{sec_toric_mutation}). 
Therefore, it is natural to ask the question below. 
We will restate this question later in a detailed form (see Question~\ref{que_mutation_dimer} and \ref{que_mutation_MM}). 

\begin{question}
\label{motiv_que}
Let $R$ be a three dimensional Gorenstein toric singularity, 
and $\End_R(M)$, $\End_R(N)$ be NCCRs of $R$ arising from dimer models. 
Then, is there a relationship between $\End_R(M)$ and $\End_R(N)$? More precisely, is there a relationship between $M$ and $N$?
\end{question}

We recall that every CCR of $R$ is connected by repeating the operation flop. 
How about NCCRs? Namely, is there a good operation that connects all NCCRs? 
In this paper, we will consider the operation called ``mutation", and show that all splitting MM modules 
are connected by repeating a certain mutation for several cases. 

In the theory of cluster algebras due to Fomin and Zelevinsky (e.g., \cite{FZ1,FZ2}), the operation called the mutation of skew symmetric matrices 
(these are interpreted as quivers without loops and $2$-cycles) plays the crucial role. 
In addition, cluster categories were introduced in \cite{BMRRT,CCS} as a categorification of cluster algebras. 
They are defined as orbit categories of the derived category of quiver representations, and the mutation (it is also called ``exchange") of cluster tilting objects 
in such categories was also introduced. 
Note that this mutation corresponds to the mutation in the theory of cluster algebras, and the notion of cluster categories is generalized in \cite{Ami} (see also \cite{Guo}). 
Also, Geiss, Leclerc and Schr\"{o}er used such mutations to investigate rigid modules over preprojective algebras \cite{GLS}. 
Subsequently, several mutations are introduced, for example the mutations of cluster tilting objects in Hom-finite triangulated $2$-Calabi-Yau categories \cite{IY}, 
the mutations of tilting modules over $3$-Calabi-Yau algebras \cite{IR}, the mutations of MM modules \cite{IW2}. 
Also, Derksen-Weyman-Zelevinsky introduced the mutations of QPs \cite{DWZ}. 
We can also define the mutation of dimer models as the dual of mutation of QPs under some restrictions (see Section~\ref{sec_toric_mutation}). 
From a viewpoint of physics, dimer models and their mutations correspond to quiver gauge theories and Seiberg duality (see e.g., \cite{Eag,EF,HV,HS,Vit}). 

In this paper, we mainly investigate the mutation of MM modules (see subsection~\ref{subsec_MMmutation}), 
and especially apply it to splitting MM modules which give NCCRs arising from dimer models. 
By considering the mutation of an MM module, we have another MM module (see Proposition~\ref{prop_mutation_MM}), 
although it might be isomorphic to the original one. 
We further note that even if a given MM module is splitting, the mutated one is not necessarily splitting in general. 
Thus, we need a further assumption for making the mutated one splitting (see Lemma~\ref{mutatable_MM}), 
and we call the argued operation the mutation of splitting MM modules only when it makes a splitting one splitting again. 
Under these backgrounds, in this paper, we consider such mutations for the case of a three dimensional Gorenstein toric singularity whose associated lattice polygon 
$\Delta$ is a reflexive polygon (see Section~\ref{mutation_reflexive}), and show the theorem below. (For further details on terminologies, see later sections.) 
Note that reflexive polygons are also studied in the context of toric Fano varieties (see e.g.,\cite[Section 8.3]{CLS}) and mirror symmetry (see e.g., \cite{Bat}). 

\begin{theorem}[{see Theorem~\ref{main_reflexive} and Corollary~\ref{mutation_MMmod}}] 
\label{main_intro}
Let $R$ be a three dimensional complete local Gorenstein toric singularity associated with a reflexive polygon. 
Then, any two splitting MM $R$-modules are transformed into each other by repeating the mutation of splitting MM modules. 
\end{theorem}

In this theorem, we emphasize that all connections between splitting MM modules are given by only mutations of splitting MM modules. 
That is, they do not factor through non-splitting ones. 

In order to show this theorem, we discuss the relationship between the mutation of QPs associated with dimer models and 
that of splitting MM generators arising from dimer models in Section~\ref{sec_toric_mutation} along the idea as in \cite{Boc2}. 
After that we will consider the mutation of splitting MM generators for toric singularities associated with reflexive polygons. 
Since dimer models giving NCCRs are classified for those cases, 
we hence compute all splitting MM generators arising from such dimer models in Section~\ref{mutation_reflexive}, 
and prove the theorem for the case of splitting MM generators by a case-by-case check. 
Finally, we extend our discussion to splitting MM modules in Section~\ref{sec_toric_mutation_module}. 

\subsection*{Conventions and Notations}  
In this paper, all modules are left modules. 
For a ring $R$, we denote 
by $\mc R$ the category of finitely generated $R$-modules. 
We denote $\add_R M$ to be the full subcategory consisting of direct summands of finite direct sums of some copies of $M\in\mc R$. 
We call $\add_R M$ the \emph{additive closure} of $M$. 
We say that $M\in\mc R$ is a \emph{generator} if $R\in\add_R M$. 
When we consider a composition of morphism, $fg$ means we first apply $f$ then $g$. 
With this convention, $\Hom_R(M, X)$ is an $\End_R(M)$-module and $\Hom_R(X, M)$ is an $\End_R(M)^{\rm op}$-module. 
Similarly, in a quiver, a path $ab$ means $a$ then $b$. 

Also, we denote by $\Cl(R)$ the class group of $R$. 
When we consider a divisorial ideal (rank one reflexive $R$-module) $I$ as an element of $\Cl(R)$, we denote it by $[I]$. 

In this paper, we consider the real two torus $\sfT\coloneqq\RR^2/\ZZ^2$ and its homology group $\rmH_1(\sfT)\cong\ZZ^2$. 
We fix the fundamental domain of $\sfT$ and orientation of cycles as 

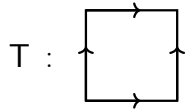
\begin{center}
\scalebox{0.6}{
\begin{tikzpicture}
\draw[line width=0.05cm]  (0,0) rectangle (2,2);

%edge
\draw[->, line width=0.05cm]  (0,0)--(0,1.2); \draw[->, line width=0.05cm]  (2,0)--(2,1.2);
\draw[->, line width=0.05cm]  (0,0)--(1.2,0); \draw[->, line width=0.05cm]  (0,2)--(1.2,2);

\node  at (-1.2,1) {\Huge{$\sfT \,:$}} ;
\end{tikzpicture} }
\end{center}
and the horizontal (resp. vertical) cycle corresponds to the first coordinate (resp. second) of $\ZZ^2$.

%%%%%%%%%%%%%%%%%%%%%%%%%%%%%%%%%%%%%%%%%%%%%%%%%%%%%%%%%%%%%%%%%%%%%%
\section{Preliminaries} 
\subsection{Maximal modifying modules}

Let $R$ be a commutative Noetherian ring. 
In this paper, we denote the $R$-dual functor by 
\[(-)^*:=\Hom_R(-, R) : \mc R\rightarrow \mc R.\] 
We say that $M\in\mc R$ is \emph{reflexive} if the natural morphism $M\rightarrow M^{**}$ is an isomorphism. 
We denote $\refl R$ to be the category of reflexive $R$-modules. 
Let $(R,\fkm)$ be a $d$-dimansional commutative Noetherian local ring. 
For $M\in\mc R$, we define the \emph{depth} of $M$ as 
\[
\depth_RM:=\mathrm{inf} \{ i\ge0 \mid \Ext^i_R(R/\fkm, M)\neq 0\}. 
\]
We say that $M$ is a \emph{maximal Cohen-Macaulay module} (= \emph{MCM module})  if $\depth_RM=d$. 
When $R$ is non-local, we say that $M$ is an MCM module if $M_\fkp$ is an MCM module for all $\fkp\in\Spec R$. 
Furthermore, we say that $R$ is a \emph{Cohen-Macaulay ring} (= \emph{CM ring}) if $R$ is an MCM $R$-module. 
We denote $\CM R$ to be the category of MCM $R$-modules. 
We say that a ring $R$ is \emph{equi-codimensional} if every maximal ideals have the same height. 
For example, all domains finitely generated over a field are equi-codimensional \cite[13.4]{Eis}. 
Also, we consider the notion of (singular) $d$-Calabi-Yau rings. 
Originally, these rings are defined by a functorial condition, but we adopt the characterization of these rings due to \cite[Proposition~3.10]{IR} as the definition in this paper. 

\begin{definition}
Let $R$ be a commutative Noetherian ring. 
We say that $R$ is \emph{$d$-Calabi-Yau} (= \emph{$d$-CY}) if $R$ is regular and equi-codimensional with $\dim R=d$. 
Also, we say that $R$ is \emph{singular $d$-Calabi-Yau} (= \emph{$d$-sCY}) if $R$ is Gorenstein and equi-codimensional with $\dim R=d$. 
\end{definition}

For example, $d$-dimensional Gorenstein toric singularities are $d$-sCY rings (see subsection~\ref{toric_pre}), 
and hence we will apply several results discussed in this subsection to these singularities. 

Next, we recall the definition of non-commutative crepant resolutions due to Van den Bergh \cite{VdB2}. 

\begin{definition}
Let $R$ be a CM ring and let $\Lambda$ be a module finite $R$-algebra (i.e., $\Lambda$ is finitely generated as an $R$-module). We say that 
 \begin{enumerate}
 \setlength{\itemsep}{-1pt}
 \item $\Lambda$ is an \emph{$R$-order} if $\Lambda$ is an MCM $R$-module, 
 \item an $R$-order $\Lambda$ is \emph{non-singular} if ${\rm gl.dim}\,\Lambda_\fkp={\rm dim}\,R_\fkp$ for all $\fkp\in\Spec R$. 
 \end{enumerate}
\end{definition}

\begin{definition}
\label{NCCR_def}
Let $R$ be a CM ring, and $0\neq M\in\refl R$. 
We say that $\Lambda\coloneqq\End_R(M)$ is a \emph{non-commutative crepant resolution} (= \emph{NCCR}) of $R$ 
or $M$ \emph{gives an NCCR} of $R$ if $\Lambda$ is a non-singular $R$-order.    
\end{definition}

An NCCR does not necessarily exist for a given singularity, thus we also introduce a weaker notion. 

\begin{definition}[{see \cite{IW2}}]
Let $R$ be a CM ring. 
\begin{enumerate}
\item We say that $M\in\refl R$ is a \emph{modifying module} if $\End_R(M)\in\CM R$. 
\item We say that $M\in\refl R$ is a \emph{maximal modifying module} (= \emph{MM module}) if it is a modifying module and if $M\oplus X$ is modifying for $X\in\refl R$ 
then $X\in\add_RM$. In other words, 
\[
\add_RM=\{X\in\refl R \mid \End_R(M\oplus X)\in\CM R\}. 
\]
For an MM module $M$, we call $\End_R(M)$ a \emph{maximal modifying algebra} (= \emph{MMA}). 
\end{enumerate}
\end{definition}

Also, we define the following nice class of modules. We will discuss these in Section~\ref{sec_toric_mutation}. 

\begin{definition}[{see \cite{IN}}]
\label{def_splitting}
We say that a reflexive module $M$ is \emph{splitting} if it is a finite direct sum of rank one reflexive modules. 
Also, we say that a (maximal) modifying module $M$ is a \emph{splitting (maximal) modifying module} if $M$ is splitting. 
Similarly, we say that an NCCR (resp. MMA) $\End_R(M)$ a \emph{splitting NCCR} (resp. \emph{splitting MMA}) if $M$ is splitting. 
\end{definition}

We then prepare some properties of MM modules used in later sections. 

\begin{proposition}[{see \cite[Proposition~4.5, 5.11 and Theorem~4.16]{IW2}}] 
\label{maximal_NCCR}
\begin{enumerate}[\rm(1)]
\item Let $R$ be a $d$-dimensional normal CM local ring with the canonical module $\omega_R$. 
Then, all reflexive modules giving NCCRs  are MM modules. 
\item If $R$ is a $3$-sCY normal domain admitting an NCCR, 
then MM modules are precisely reflexive modules giving NCCRs of $R$. 
Furthermore, suppose that $M$ and $N$ are MM $R$-modules, then $\End_R(M)$ and $\End_R(N)$ are derived equivalent.
\end{enumerate}
\end{proposition}

%%%%%%%%%%%%%%%%%%%%%%%%%%%%%%%%%%%%%%%%%%%%%%%%%%%%%%%%%%%%%%%%%%%%%%
\subsection{Mutations of MM modules and tilting modules}
\label{subsec_MMmutation}

Following \cite[Section~6]{IW2}, we introduce the mutation of MM modules to consider Question~\ref{motiv_que}. 
Let $R$ be a complete local normal $d$-sCY ring, and $M\coloneqq\oplus_{i\in I}M_i$ be a modifying $R$-module. 
Here, we assume that $M$ is \emph{basic}, that is, $M_i$'s are mutually non-isomorphic. 
For a subset $J$ of $I$, we define $M_J\coloneqq\oplus_{j\in J}M_j$ and $J^c\coloneqq I\setminus J$. 
In particular, we have that $M=M_J\oplus M_{J^c}$. 
We then consider a \emph{minimal right $(\add_R M_{J^c})$-approximation} of $M_J$. 
Namely, it is a morphism $\varphi:N\rightarrow M_J$ such that $N\in\add_R M_{J^c}$ and 
\[
\Hom_R(M_{J^c}, N)\xrightarrow{\cdot\varphi}\Hom_R(M_{J^c}, M_J)
\]
 is surjective, and if $\phi\in\End_R(N)$ satisfies $\phi\varphi=\varphi$, then $\phi$ is an automorphism. 
Since $R$ is complete, $\varphi$ exists and is unique up to isomorphism. Also, $\varphi$ is surjective if $N$ is a generator. 
Then, we consider the kernel 
\[ 0\rightarrow \calK \coloneqq\Ker\varphi\rightarrow N\xrightarrow{\varphi} M_J\]
and we call this sequence the \emph{exchange sequence}. 
We define the \emph{right mutation} of $M$ at $J$ as $\mu^+_J(M)\coloneqq M_{J^c}\oplus \calK$, 
and define the \emph{left mutation} of $M$ at $J$ as $\mu^-_J(M)\coloneqq (\mu^+_J(M^*))^*$. 

\medskip

Here, we collect some properties of the mutations of MM modules. 

\begin{proposition}[{see \cite[Section~6]{IW2}}] 
\label{prop_mutation_MM}
Let $R$ be a complete local normal $d$-sCY ring, $M=\oplus_{i\in I}M_i$ be a basic modifying $R$-module. 
Suppose that $J$ is a subset of $I$. Then, we have the following. 
\begin{itemize}
\item[(1)] $\mu^+_J$ and $\mu^-_J$ are mutually inverse operations. That is, 
\[ \mu^-_J(\mu^+_J(M))=M \text{\quad and \quad} \mu^+_J(\mu^-_J(M))=M\] 
up to additive closure. 
\item[(2)] $\mu^+_J(M)$ and $\mu^-_J(M)$ are modifying $R$-modules. 
\item[(3)] If $M$ gives an NCCR, then $\mu^+_J(M)$ and $\mu^-_J(M)$ also give NCCRs. 
\item[(4)] $\End_R(M)$, $\End_R(\mu^+_J(M))$ and $\End_R(\mu^-_J(M))$ are derived equivalent. 
\item[(5)] Assume that $d=3$, if $M$ is an MM module, then so are $\mu^+_J(M)$ and $\mu^-_J(M)$. 
\item[(6)] Assume that $d=3$, if $M$ is an MM module and $J=\{k\}$, then $\mu^+_J(M)\cong\mu^-_J(M)$, hence we denote it by $\mu_k(M)$. 
\end{itemize}
\end{proposition}

In the rest of this subsection, we keep the notations in Proposition~\ref{prop_mutation_MM} and assume $d=3$. 
We denote the set of isomorphism classes of basic MM modules (resp. generators) by $\MM(R)$ (resp. $\MMG(R)$). 
In this paper, we only consider the mutation at $J=\{k\}$, 
and hence for $M=\bigoplus_{i\in I}M_i\in\MM(R)$ we have the new element $\mu_k(M)\in\MM(R)$ by Proposition~\ref{prop_mutation_MM}(5)(6). 
Furthermore, since we are interested in splitting MM modules and their mutations, we restrict the above notation to splitting ones. 
Namely, we denote the set of isomorphism classes of basic splitting MM modules (resp. generators) by $\TMM(R)$ (resp. $\TMMG(R)$). 
%Then, we consider the mutation of $M=\bigoplus_{i\in I}M_i\in\TMM(R)$, but 
As we mentioned, even if $M=\bigoplus_{i\in I}M_i\in\TMM(R)$, we does not necessarily have $\mu_k(M)\in\TMM(R)$. 
(Note that we always have $\mu_k(M)\in\MM(R)$.) 
Therefore, we call the above operation $\mu_k$ \emph{the mutation of splitting MM modules} only when both of $M$ and $\mu_k(M)$ are splitting. 
Then, we define the \emph{exchange graph} of $\TMM(R)$ as follows, and denote it by $\EG(\TMM(R))$. 
That is, the vertices of $\EG(\TMM(R))$ are elements in $\TMM(R)$, and draw an edge between $M$ 
and $\mu_k(M)$ for each $M\in\TMM(R)$ and $k\in I$ with $\mu_k(M)\in\TMM(R)$. 
We also define the exchange graph of $\TMMG(R)$ in a similar way, and denote it by $\EG(\TMMG(R))$.

\medskip

Next, we recall the definition of a tilting module. 

\begin{definition}
For a ring $\Lambda$, we say that $T\in\mc\Lambda$ is a \emph{partial tilting module} if 
$\projdim_\Lambda T\le 1$ and $\Ext^1_\Lambda(T, T)=0$. In addition, if there exists an exact sequence:
$0\rightarrow \Lambda\rightarrow T_0\rightarrow T_1\rightarrow 0$ 
such that $T_i\in\add_\Lambda T \,(i=0, 1)$, we say that $T$ is a \emph{tilting module}. 
\end{definition}

For an MM $R$-module $M$, we can obtain reflexive tilting $\End_R(M)$-modules as follows. 
Here, for an $R$-algebra $\Lambda$, we say that $M\in\mc\Lambda$ is a \emph{reflexive $\Lambda$-module} if it is reflexive as an $R$-module. 

\begin{theorem}[{see \cite[Theorem~4.17]{IW2}}] 
\label{MM_tilting}
Let $R$ be a normal $3$-sCY ring and $M$ be an MM module. Then, we have the following bijection via 
the functor $\Hom_R(M, -):\mc R\rightarrow\mc\End_R(M)$. 
\[
\{\text{maximal modifying $R$-modules}\}\overset{1:1}{\longleftrightarrow}\{\text{reflexive tilting $\End_R(M)$-modules}\}. 
\]
\end{theorem}

For a module-finite $R$-algebra $\Lambda$, we say that a basic partial tilting $\Lambda$-module $T$ is 
an \emph{almost complete tilting module} if it has $r$ non-isomorphic indecomposable direct summands, where $(r+1)$ is the number of non-isomorphic 
indecomposable projective $\Lambda$-modules. 
Then $X\in\mc\Lambda$ is called a \emph{complement} of $T$ if $T\oplus X$ is a basic tilting $\Lambda$-module. 
A complement of an almost complete tilting module $T$ always exists \cite{Bon}, 
and there are at most two complements of $T$ \cite[Proposition~5.1]{IR}. 

Let $M$ be an MM $R$-module, and let $\Lambda=\End_R(M)$. 
For an MM $R$-module $N=\bigoplus_{i\in I}N_i$ where $I=\{0, 1, \cdots, r\}$, $\Hom_R(M, N)$ is a tilting $\Lambda$-module by Theorem~\ref{MM_tilting}.  
Therefore $\Hom_R(M, \bigoplus_{i\in I\setminus\{k\}}N_i)$ is an almost complete tilting $\Lambda$-module, and $\Hom_R(M, N_k)$ is a complement. 
Considering a minimal right $(\add_R\bigoplus_{i\in I\setminus\{k\}}N_i)$-approximation of $N_k$: 
\[
0\rightarrow \calK\coloneqq\Ker\varphi\rightarrow\widetilde{N}\rightarrow N_k 
\]
where $\widetilde{N}\in\add_R\bigoplus_{i\in I\setminus\{k\}}N_i$, we have that 
\begin{equation}
\label{rightapprox}
0\rightarrow\Hom_R(M, \calK)\rightarrow\Hom_R(M, \widetilde{N})\rightarrow\Hom_R(M, N_k)\rightarrow0. 
\end{equation}
Then $\Hom_R(M, \mu_k(N))=\Hom_R(M, \bigoplus_{i\in I\setminus\{k\}}N_i)\oplus\Hom_R(M, \calK)$ is a tilting $\Lambda$-module by Theorem~\ref{MM_tilting}, 
and hence $\Hom_R(M, \calK)$ is a complement of $\Hom_R(M, \bigoplus_{i\in I\setminus\{k\}}N_i)$. 
Since $(\ref{rightapprox})$ is a minimal right $\add_\Lambda T$-approximation where $T=\Hom_R(M,\bigoplus_{i\in I\setminus\{k\}}N_i)$, 
we see that $\Hom_R(M, \mu_k(N))$ is just the mutation of tilting $\Lambda$-module $\Hom_R(M, N)$ in the sense of \cite[Section~5]{IR}. 
Thus, we define the \emph{mutation of tilting module} $\Hom_R(M, N)$ at $k\in I$ as 
\[
\mu_k(\Hom_R(M, N))\coloneqq \Hom_R(M, \mu_k(N)). 
\]
For now, we have shown that $\Hom_R(M, N_k)$ and $\Hom_R(M, \calK)$ are complements of $\Hom_R(M,\bigoplus_{i\in I\setminus\{k\}}N_i)$. 
Since the number of such complements is at most two, they are exactly two complements if $N_k\not\cong\calK$. 
Therefore, if $N=\bigoplus_{i\in I}N_i$ and $N^\prime=N_k^\prime\bigoplus_{i\in I\setminus\{k\}}N_i$ are basic MM modules with $N_k\not\cong N_k^\prime$ 
(that is, $N$ and $N^\prime$ have the same direct summands except one component), then we have that $N^\prime\cong\mu_k(N)$. 
%and thus $N$ and $N^\prime$ are connected in $\EG(\MMG(R))$. 

%%%%%%%%%%%%%%%%%%%%%%%%%%%%%%%%%%%%%%%%%%%%%%%%%%%%%%%%%%%%%%%%%%%%%%
\subsection{Preliminaries of toric singularities}
\label{toric_pre}
In this subsection, we review about toric singularities and fix the notation used in later sections. 
For more about toric singularities, see e.g, \cite{BG2,CLS}. 

Let $\sfN\cong\ZZ^d$ be a lattice of rank $d$ and let $\sfM\coloneqq\Hom_\ZZ(\sfN, \ZZ)$ be the dual lattice of $\sfN$. 
We set $\sfN_\RR\coloneqq\sfN\otimes_\ZZ\RR$ and $\sfM_\RR\coloneqq\sfM\otimes_\ZZ\RR$. 
We denote the natural inner product by $\langle\;,\;\rangle:\sfM_\RR\times\sfN_\RR\rightarrow\RR$. 
Let 
\[
\sigma\coloneqq\mathrm{Cone}\{v_1, \cdots, v_n\}=\RR_{\ge 0}v_1+\cdots +\RR_{\ge 0}v_n
\subset\sfN_\RR 
\]
be a strongly convex rational polyhedral cone of rank $d$ 
generated by $v_1, \cdots, v_n\in\ZZ^d$ and suppose that this system of generators is minimal. 
For each generator, we define the linear form $\lambda_i(-)\coloneqq\langle-, v_i\rangle$ and 
the half space $H^+_i\coloneqq \{x\in\sfM_\RR\,|\,\lambda_i(x)\ge0 \}$ for $i=1, \cdots, n$. 
Then, the dual cone $\sigma^\vee$ is determined as the intersection of these half spaces: 
\[
\sigma^\vee\coloneqq\{x\in\sfM_\RR\,|\,\langle x,y\rangle\ge0 \text{ for all } y\in\sigma \}=H^+_1\cap\cdots\cap H^+_n\subset\sfM_\RR. 
\]
Moreover $\sigma^\vee\cap\sfM$ is a positive affine normal semigroup. 
We define the toric singularity (or affine normal semigroup ring) associated with $\sigma$ as follows: 
\[
R\coloneqq K[\sigma^\vee\cap\sfM]=K[t_1^{\alpha_1}\cdots t_d^{\alpha_d}\mid (\alpha_1, \cdots, \alpha_d)\in\sigma^\vee\cap\sfM]
\]
where $K$ is a field. 
In the rest, we always assume that $K$ is an algebraically closed field of characteristic zero. 
In this situation, $R$ is a $d$-dimensional CM normal domain. 
Note that it is known that $R$ is Gorenstein if and only if there is $x\in\sigma^\vee\cap\sfM$ such that $\lambda_i(x)=1$ for all $i=1, \cdots, n$ 
(e.g., \cite[Theorem~6.33]{BG2}). 
In particular, $d$-dimensional Gorenstein toric singularities are $d$-sCY normal domains. 

For each $\mathbf{u}=(u_1, \cdots, u_n)\in\RR^n$, we set 
\[
\mathbb{T}(\mathbf{u})\coloneqq\{x\in\sfM\cong\ZZ^d\,|\,(\lambda_1(x), \cdots, \lambda_n(x))\ge(u_1, \cdots, u_n)\}, 
\]
and define the divisorial ideal $T(\mathbf{u})$ generated by all monomials whose exponent vector is in $\mathbb{T}(\mathbf{u})$. 
From definition, we have that $T(\mathbf{u})=T(\ulcorner \mathbf{u}\urcorner)$ where $\ulcorner \; \urcorner$ implies the round up 
and $\ulcorner \mathbf{u}\urcorner=(\ulcorner u_1\urcorner, \cdots, \ulcorner u_n\urcorner)$. 
Thus, in the rest, we assume that $\mathbf{u}\in\ZZ^n$. 
Clearly, we have that $\mathbb{T}(0)=\sigma^\vee\cap\sfM$ and $T(0)=R$. 
Set $I_i\coloneqq T(\delta_{i1},\cdots,\delta_{in})$ where $\delta_{ij}$ is the Kronecker delta. 
Then, a divisorial ideal $T(u_1,\cdots,u_n)$ corresponds to the element $u_1[I_1]+\cdots+u_n[I_n]$ in $\Cl(R)$. 
In general, a divisorial ideal of $R$ takes this form (see e.g., \cite[Theorem~4.54]{BG2}) and we have the following. 

\begin{lemma}[{see e.g., \cite[Corollary~4.56]{BG2}}] 
\label{div_eq}
For $\mathbf{u}, \mathbf{u}^\prime\in\ZZ^n$, $T(\mathbf{u})\cong T(\mathbf{u}^\prime)$ as an $R$-module 
if and only if there exists $y\in\sfM$ such that $u_i=u_i^\prime+\lambda_i(y)$ for all $i=1, \cdots, n$. 
Especially, we have that $\Cl(R)\cong \ZZ^n/\lambda(\ZZ^d)$. 
\end{lemma}

%%%%%%%%%%%%%%%%%%%%%%%%%%%%%%%%%%%%%%%%%%%%%%%%%%%%%%%%%%%%%%%%%%%%%%
\section{NCCRs arising from consistent dimer models}  
\label{sec2}

In this section, we introduce the notion of dimer models and discuss related topics. 
Especially, we will show that we can obtain NCCRs of a three dimensional Gorenstein toric singularity from dimer models that satisfy the consistency condition. 
For more results concerning dimer models, we refer to a survey article \cite{Boc5} and references quoted in this section. 

\subsection{Dimer models and Perfect mathchings}
\label{subsec_dimer}

\begin{figure}
\[\scalebox{0.65}{
\begin{tikzpicture}
%vertex
\node (P1) at (1,1){$$}; \node (P2) at (3,1){$$}; \node (P3) at
(3,3){$$}; \node (P4) at (1,3){$$};
\draw[thick]  (0,0) rectangle (4,4);
%edge
\draw[line width=0.05cm]  (P1)--(P2)--(P3)--(P4)--(P1);\draw[line width=0.05cm] (0,1)--(P1)--(1,0); \draw[line width=0.05cm]  (4,1)--(P2)--(3,0);
\draw[line width=0.05cm]  (0,3)--(P4)--(1,4);\draw[line width=0.05cm]  (3,4)--(P3)--(4,3);
%black
\filldraw  [ultra thick, fill=black] (1,1) circle [radius=0.18] ;\filldraw  [ultra thick, fill=black] (3,3) circle [radius=0.18] ;
%white
\draw  [ultra thick,fill=white] (3,1) circle [radius=0.18] ;\draw  [ultra thick, fill=white] (1,3) circle [radius=0.18] ;

\end{tikzpicture}
}\]
\caption{Dimer model $\Gamma_{4a}$}
\label{ex_dimer4a}
\end{figure}
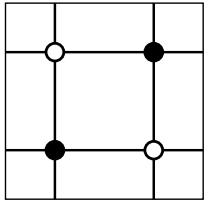

A \emph{dimer model} (or \emph{brane tiling}) is a polygonal cell decomposition of the real two-torus 
$\sfT\coloneqq\RR^2/\ZZ^2$ whose vertices and edges form a finite bipartite graph. 
Therefore, each vertex is colored either black or white so that each edge connects a black vertex to a white vertex. 
For example,  $\Gamma_{4a}$ (see Figure~\ref{ex_dimer4a}) is a dimer model where the outer frame is the fundamental domain of the torus $\sfT$. 
Let $\Gamma$ be a dimer model. We denote the set of vertices (resp. edges, faces) of $\Gamma$ by $\Gamma_0$ (resp. $\Gamma_1$, $\Gamma_2$). 
Since $\Gamma$ is a bipartite graph, we can divide $\Gamma_0$ into the disjoint union of the set of black vertices and the set of white vertices. 

In what follows, we construct a toric singularity from a dimer model $\Gamma$ by paying attention to perfect matchings on $\Gamma$. 

\begin{definition}
A \emph{perfect matching} (or \emph{dimer configuration}) on a dimer model $\Gamma$ is a subset $\sfP$ of $\Gamma_1$ such that each vertex is 
the end point of precisely one edge in $\sfP$. 
\end{definition}

%%% type 4a-1 perfect matching%%%

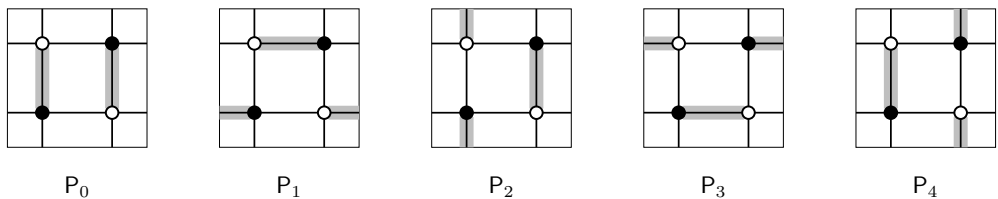
\begin{figure}[h]
\begin{center}
{\scalebox{0.9}{
\begin{tikzpicture} 
\node at (0,-1.6) {$\sfP_0$};\node at (3.1,-1.6) {$\sfP_1$}; \node at (6.2,-1.6) {$\sfP_2$}; \node at (9.3,-1.6) {$\sfP_3$}; \node at (12.4,-1.6) {$\sfP_4$};

\node (PM0) at (0,0) 
{\scalebox{0.51}{
\begin{tikzpicture}
%vertex
\node (P1) at (1,1){$$}; \node (P2) at (3,1){$$}; \node (P3) at
(3,3){$$}; \node (P4) at (1,3){$$};
\draw[thick]  (0,0) rectangle (4,4);

%perfect matching
\draw[line width=0.4cm,color=lightgray] (P1)--(P4);\draw[line width=0.4cm,color=lightgray] (P2)--(P3);

%edge
\draw[line width=0.05cm]  (P1)--(P2)--(P3)--(P4)--(P1);\draw[line width=0.05cm] (0,1)--(P1)--(1,0); \draw[line width=0.05cm]  (4,1)--(P2)--(3,0);
\draw[line width=0.05cm]  (0,3)--(P4)--(1,4);\draw[line width=0.05cm]  (3,4)--(P3)--(4,3);
%black
\filldraw  [ultra thick, fill=black] (1,1) circle [radius=0.18] ;\filldraw  [ultra thick, fill=black] (3,3) circle [radius=0.18] ;
%white
\draw  [ultra thick,fill=white] (3,1) circle [radius=0.18] ;\draw  [ultra thick,fill=white] (1,3) circle [radius=0.18] ;
\end{tikzpicture} }}; 

\node (PM1) at (3.1,0) 
{\scalebox{0.51}{
\begin{tikzpicture}
%vertex
\node (P1) at (1,1){$$}; \node (P2) at (3,1){$$}; \node (P3) at (3,3){$$}; \node (P4) at (1,3){$$};
\draw[thick]  (0,0) rectangle (4,4);

%perfect matching
\draw[line width=0.4cm,color=lightgray] (P3)--(P4);\draw[line width=0.4cm,color=lightgray] (P1)--(0,1);\draw[line width=0.4cm,color=lightgray] (P2)--(4,1);

%edge
\draw[line width=0.05cm]  (P1)--(P2)--(P3)--(P4)--(P1);\draw[line width=0.05cm] (0,1)--(P1)--(1,0); \draw[line width=0.05cm]  (4,1)--(P2)--(3,0);
\draw[line width=0.05cm]  (0,3)--(P4)--(1,4);\draw[line width=0.05cm]  (3,4)--(P3)--(4,3);
%black
\filldraw  [ultra thick, fill=black] (1,1) circle [radius=0.18] ;\filldraw  [ultra thick, fill=black] (3,3) circle [radius=0.18] ;
%white
\draw  [ultra thick,fill=white] (3,1) circle [radius=0.18] ;\draw  [ultra thick,fill=white] (1,3) circle [radius=0.18] ;
\end{tikzpicture} }}; 

\node (PM2) at (6.2,0) 
{\scalebox{0.51}{
\begin{tikzpicture}
%vertex
\node (P1) at (1,1){$$}; \node (P2) at (3,1){$$}; \node (P3) at (3,3){$$}; \node (P4) at (1,3){$$};
\draw[thick]  (0,0) rectangle (4,4);

%perfect matching
\draw[line width=0.4cm,color=lightgray] (P3)--(P2);\draw[line width=0.4cm,color=lightgray] (P4)--(1,4);\draw[line width=0.4cm,color=lightgray] (P1)--(1,0);

%edge
\draw[line width=0.05cm]  (P1)--(P2)--(P3)--(P4)--(P1);\draw[line width=0.05cm] (0,1)--(P1)--(1,0); \draw[line width=0.05cm]  (4,1)--(P2)--(3,0);
\draw[line width=0.05cm]  (0,3)--(P4)--(1,4);\draw[line width=0.05cm]  (3,4)--(P3)--(4,3);
%black
\filldraw  [ultra thick, fill=black] (1,1) circle [radius=0.18] ;\filldraw  [ultra thick, fill=black] (3,3) circle [radius=0.18] ;
%white
\draw  [ultra thick,fill=white] (3,1) circle [radius=0.18] ;\draw  [ultra thick,fill=white] (1,3) circle [radius=0.18] ;
\end{tikzpicture} }} ;  

\node (PM3) at (9.3,0) 
{\scalebox{0.51}{
\begin{tikzpicture}
%vertex
\node (P1) at (1,1){$$}; \node (P2) at (3,1){$$}; \node (P3) at (3,3){$$}; \node (P4) at (1,3){$$};

\draw[thick]  (0,0) rectangle (4,4);

%perfect matching
\draw[line width=0.4cm,color=lightgray] (P2)--(P1);\draw[line width=0.4cm,color=lightgray] (P4)--(0,3);\draw[line width=0.4cm,color=lightgray] (P3)--(4,3);
%edge
\draw[line width=0.05cm]  (P1)--(P2)--(P3)--(P4)--(P1);\draw[line width=0.05cm] (0,1)--(P1)--(1,0); \draw[line width=0.05cm]  (4,1)--(P2)--(3,0);
\draw[line width=0.05cm]  (0,3)--(P4)--(1,4);\draw[line width=0.05cm]  (3,4)--(P3)--(4,3);
%black
\filldraw  [ultra thick, fill=black] (1,1) circle [radius=0.18] ;\filldraw  [ultra thick, fill=black] (3,3) circle [radius=0.18] ;
%white
\draw  [ultra thick,fill=white] (3,1) circle [radius=0.18] ;\draw  [ultra thick,fill=white] (1,3) circle [radius=0.18] ;
\end{tikzpicture} }}; 

\node (PM4) at (12.4,0) 
{\scalebox{0.51}{
\begin{tikzpicture}
%vertex
\node (P1) at (1,1){$$}; \node (P2) at (3,1){$$}; \node (P3) at (3,3){$$}; \node (P4) at (1,3){$$};
\draw[thick]  (0,0) rectangle (4,4);

%perfect matching
\draw[line width=0.4cm,color=lightgray] (P4)--(P1);\draw[line width=0.4cm,color=lightgray] (P3)--(3,4);\draw[line width=0.4cm,color=lightgray] (P2)--(3,0);

%edge
\draw[line width=0.05cm]  (P1)--(P2)--(P3)--(P4)--(P1);\draw[line width=0.05cm] (0,1)--(P1)--(1,0); \draw[line width=0.05cm]  (4,1)--(P2)--(3,0);
\draw[line width=0.05cm]  (0,3)--(P4)--(1,4);\draw[line width=0.05cm]  (3,4)--(P3)--(4,3);
%black
\filldraw  [ultra thick, fill=black] (1,1) circle [radius=0.18] ;\filldraw  [ultra thick, fill=black] (3,3) circle [radius=0.18] ;
%white
\draw  [ultra thick,fill=white] (3,1) circle [radius=0.18] ;\draw  [ultra thick,fill=white] (1,3) circle [radius=0.18] ;
\end{tikzpicture} }} ;

\end{tikzpicture}
}}
\end{center}
\caption{A part of perfect matchings of $\Gamma_{4a}$}
\label{pm_4a}
\end{figure}

For example, Figure~\ref{pm_4a} is a part of perfect matchings on $\Gamma_{4a}$. (In addition, there are three perfect matchings.)  
Every dimer model does not necessarily have a perfect matching. 
A typical case is a dimer model whose number of black vertices is not equal to that of white ones. 
Also, there is a dimer model which has the same number of black and white vertices, but do not have a perfect matching (cf. \cite[(2.8)]{Bro}).  
In this paper, we will not treat such a pathological one. 
Then, we will impose the extra assumption so-called ``non-degeneracy condition", that is, we say that a dimer model is \emph{non-degenerate} if every edge is contained in some perfect matchings. 
For example, $\Gamma_{4a}$ is non-degenerate (see Figure~\ref{ex_dimer4a} and \ref{pm_4a}). 
In the rest of this paper, we assume that a dimer model is non-degenerate, 
and we also assume that a dimer model has no bivalent vertices (i.e., a vertex which connects only two distinct vertices). 
If there are bivalent vertices, we can remove them (see e.g., \cite[Figure~5.1]{IU1}) because this process does not affect our problems. 
To be specific, we mainly investigate the Jacobian algebra associated with a dimer model (see subsection~\ref{sec_QP_from_dimer}) 
to consider our problems, and a removing bivalent vertices does not change the Jacobian algebra up to isomorphism.  

\medskip

For each perfect matching, we give the orientation from a white vertex to a black one. 
Then, the difference of two perfect matchings forms a $1$-cycle. 
We will consider this $1$-cycle as an element in the homology group $\rmH_1(\sfT)\cong\ZZ^2$. 
Precisely, for two perfect matchings $\sfP_i, \sfP_j$ we define 
\[
h(\sfP_i, \sfP_j)\coloneqq [\sfP_i-\sfP_j]\in \rmH_1(\sfT)\cong\ZZ^2. 
\]
We consider this homology class for any pair of perfect matchings, and have finitely many elements in $\ZZ^2$. 
Then, we define the lattice polygon as in Definition~\ref{polygon}. 
Here, note that we have 
\[
h(\sfP_i, \sfP_j)=h(\sfP_i, \sfP_k)-h(\sfP_j, \sfP_k). 
\]
Thus, we fix one perfect matching, and we may consider only the difference between each perfect matching and the fixed one. 

\begin{definition}
\label{polygon}
For a dimer model $\Gamma$, we fix a perfect matching $\sfP^\prime$ of $\Gamma$. 
We define the lattice polygon $\Delta$ as the convex hull of 
\[
\{h(\sfP, \sfP^\prime)\in\ZZ^2 \;|\; \sfP \text{ is a perfect matching of $\Gamma$} \} 
\]
in $\RR^2$. $\Delta$ is called the \emph{perfect matching polygon} (or \emph{characteristic polygon}) of $\Gamma$. 
\end{definition}

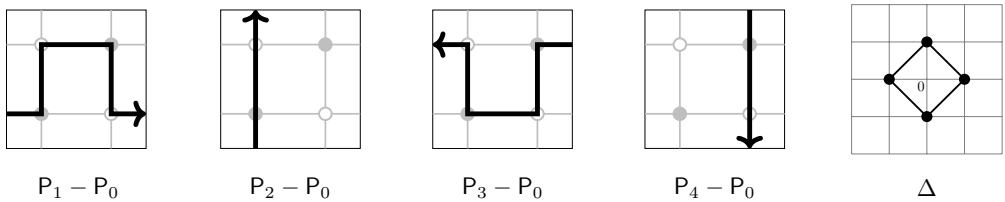
\begin{figure}[h]
\begin{center}
{\scalebox{0.9}{
\begin{tikzpicture} 
\node at (0,-1.6) {$\sfP_1-\sfP_0$};\node at (3.1,-1.6) {$\sfP_2-\sfP_0$}; \node at (6.2,-1.6) {$\sfP_3-\sfP_0$}; 
\node at (9.3,-1.6) {$\sfP_4-\sfP_0$}; \node at (12.4,-1.6) {$\Delta$};

\node (PM1) at (0,0) %pm1
{\scalebox{0.51}{\begin{tikzpicture}
%vertex
\node (P1) at (1,1){$$}; \node (P2) at (3,1){$$}; \node (P3) at (3,3){$$}; \node (P4) at (1,3){$$};

\draw[thick]  (0,0) rectangle (4,4);

\draw[lightgray,line width=0.05cm]  (P1)--(P2);\draw[lightgray,line width=0.05cm] (P1)--(1,0); \draw[lightgray,line width=0.05cm]  (P2)--(3,0);
\draw[lightgray,line width=0.05cm]  (0,3)--(P4)--(1,4);\draw[lightgray,line width=0.05cm]  (3,4)--(P3)--(4,3);
%black
\filldraw  [lightgray, ultra thick, fill=lightgray,] (1,1) circle [radius=0.18] ;\filldraw  [lightgray,ultra thick, fill=lightgray,] (3,3) circle [radius=0.18] ;
%white
\draw  [lightgray,ultra thick, fill=white] (3,1) circle [radius=0.18] ;\draw  [lightgray,ultra thick, fill=white] (1,3) circle [radius=0.18] ;
%edge
\draw[->, line width=0.14cm] (0,1)--(1,1)--(1,3)--(3,3)--(3,1)--(4,1) ;
\end{tikzpicture} }}; 

\node (PM2) at (3.1,0) %pm2
{\scalebox{0.51}{
\begin{tikzpicture}
%vertex
\node (P1) at (1,1){$$}; \node (P2) at (3,1){$$}; \node (P3) at (3,3){$$}; \node (P4) at (1,3){$$};
\draw[thick]  (0,0) rectangle (4,4);

\draw[lightgray,line width=0.05cm]  (P1)--(P2)--(P3)--(P4);\draw[lightgray,line width=0.05cm] (0,1)--(P1); \draw[lightgray,line width=0.05cm]  (4,1)--(P2)--(3,0);
\draw[lightgray,line width=0.05cm]  (0,3)--(P4);\draw[lightgray,line width=0.05cm]  (3,4)--(P3)--(4,3);
%black
\filldraw  [lightgray, ultra thick, fill=lightgray,] (1,1) circle [radius=0.18] ;\filldraw  [lightgray,ultra thick, fill=lightgray,] (3,3) circle [radius=0.18] ;
%white
\draw  [lightgray,ultra thick, fill=white] (3,1) circle [radius=0.18] ;\draw  [lightgray,ultra thick, fill=white] (1,3) circle [radius=0.18] ;
%edge
\draw[->, line width=0.14cm] (1,0)--(1,4) ;
\end{tikzpicture} }}; 

\node (PM3) at (6.2,0) %pm3
{\scalebox{0.51}{
\begin{tikzpicture}
%vertex
\node (P1) at (1,1){$$}; \node (P2) at (3,1){$$}; \node (P3) at (3,3){$$}; \node (P4) at (1,3){$$};

\draw[thick]  (0,0) rectangle (4,4);

\draw[lightgray,line width=0.05cm]  (P3)--(P4);\draw[lightgray,line width=0.05cm] (0,1)--(P1)--(1,0); \draw[lightgray,line width=0.05cm]  (4,1)--(P2)--(3,0);
\draw[lightgray,line width=0.05cm]  (P4)--(1,4);\draw[lightgray,line width=0.05cm]  (3,4)--(P3);
%black
\filldraw  [lightgray, ultra thick, fill=lightgray,] (1,1) circle [radius=0.18] ;\filldraw  [lightgray,ultra thick, fill=lightgray,] (3,3) circle [radius=0.18] ;
%white
\draw  [lightgray,ultra thick, fill=white] (3,1) circle [radius=0.18] ;\draw  [lightgray,ultra thick, fill=white] (1,3) circle [radius=0.18] ;
%edge
\draw[->, line width=0.14cm] (4,3)--(3,3)--(3,1)--(1,1)--(1,3)--(0,3) ;
\end{tikzpicture}
 }} ;  

\node (PM4) at (9.3,0) %pm4
{\scalebox{0.51}{
\begin{tikzpicture}
%vertex
\node (P1) at (1,1){$$}; \node (P2) at (3,1){$$}; \node (P3) at (3,3){$$}; \node (P4) at (1,3){$$};

\draw[thick]  (0,0) rectangle (4,4);

\draw[lightgray,line width=0.05cm]  (P2)--(P1)--(P4)--(P3);\draw[lightgray,line width=0.05cm] (0,1)--(P1)--(1,0); \draw[lightgray,line width=0.05cm]  (4,1)--(P2);
\draw[lightgray,line width=0.05cm]  (0,3)--(P4)--(1,4);\draw[lightgray,line width=0.05cm]  (P3)--(4,3);

%black
\filldraw  [lightgray, ultra thick, fill=lightgray,] (1,1) circle [radius=0.18] ;\filldraw  [lightgray,ultra thick, fill=lightgray,] (3,3) circle [radius=0.18] ;
%white
\draw  [lightgray,ultra thick, fill=white] (3,1) circle [radius=0.18] ;\draw  [lightgray,ultra thick, fill=white] (1,3) circle [radius=0.18] ;
%edge
\draw[->, line width=0.14cm] (3,4)--(3,0) ;
\end{tikzpicture} }}; 

\node (PMP) at (12.4,0) %pm_polygon
{\scalebox{0.55}{
\begin{tikzpicture}
\draw [step=1,thin, gray] (-2,-2) grid (2,2);
\draw [line width=0.05cm] (1,0)--(0,1)--(-1,0)--(0,-1)--cycle ;
\filldraw  [ultra thick, fill=black] (1,0) circle [radius=0.12] ;\filldraw  [ultra thick, fill=black] (0,1) circle [radius=0.12] ;
\filldraw  [ultra thick, fill=black] (-1,0) circle [radius=0.12] ;\filldraw  [ultra thick, fill=black] (0,-1) circle [radius=0.12] ;
%\node  at (0,0) {\Large{4a}} ;
\node at (-0.15,-0.2){{\large$0$}};
\end{tikzpicture} }} ;

\end{tikzpicture}
}}
\end{center}
\caption{Differences of perfect matchings and the perfect matching polygon of $\Gamma_{4a}$}
\label{pm_polygon}
\end{figure}

Although such a convex hull depends on a choice of a fixed perfect matching, it is determined up to isometric transformations, 
and the difference of transformation on $\Delta$ does not affect the toric geometry up to isomorphism. 
For example, for the dimer model $\Gamma_{4a}$, we fix the perfect matching $\sfP_0$. 
Then, we have the perfect matching polygon $\Delta$ as in Figure~\ref{pm_polygon}. 

%%%%%%%%%%%%%%%%%%%%
%%% type 4a-1 p.m. polygon%%%
%%%%%%%%%%%%%%%%%%%%

\begin{definition}
Fix a perfect matching $\sfP^\prime$ for a given dimer model. 
We say that $\sfP$ is an \emph{extremal} (or a \emph{corner}) \emph{perfect matching} if the corresponding point $h(\sfP, \sfP^\prime)$ lies at a vertex of 
the perfect matching polygon $\Delta$. 

We define the \emph{multiplicity} of a lattice point $p\in\Delta$ as the number of perfect matching corresponding to a point $p$. 
\end{definition}

Let $v_1^\prime, \cdots, v_n^\prime\in \ZZ^2$ be all the lattice points lie at a vertex of $\Delta$. 
Since they form $\Delta$ as the convex hull of them, we give a cyclic order to these lattice points along $\Delta$. 
Also, we can take an extremal perfect matching $\sfP_i$ for each lattice point $v_i^\prime$ ($i=1, \cdots, n$).  
We remark that a corresponding perfect matching is not unique for now (see Proposition~\ref{ex_pm}). 
Then, we put $\Delta$ on the hyperplane $z=1$, and construct the cone $\sigma_\Delta$. 
%That is, we add the third coordinate $z=1$ to each primitive vector $v_i^\prime$. 
Namely, we set vectors $v_i=(v_i^\prime, 1)\in \ZZ^3$ for $i=1, \cdots, n$, and construct the cone $\sigma_\Delta=\mathrm{Cone}\{v_1, \cdots, v_n\}$. 
In the same way shown in Section~\ref{toric_pre}, we define the toric singularity $R\coloneqq K[\sigma_\Delta^\vee\cap\sfM]$ where $\sfM\cong\ZZ^3$. 
By the construction, $R$ is a three dimensional Gorenstein normal domain. 

%%%%%%%%%%%%%%%%%%%%%%%%%%%%%%%%%%%%%%%%%%%%%%%%%%%%%%%%%%%%%%%%%%%%%%
\subsection{Quivers associated with dimer models}
\label{sec_QP_from_dimer}

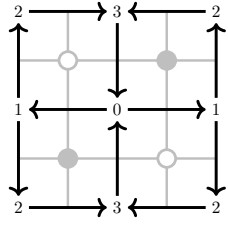
\begin{figure}[h]
\[\scalebox{0.65}{
\begin{tikzpicture}
%vertex
\node (Q1) at (2,2){$0$};\node (Q2a) at (0,2){$1$}; \node(Q2b) at (4,2){$1$};\node (Q3a) at (0,0){$2$};
\node(Q3c) at (4,4){$2$};\node(Q3b) at (4,0){$2$};\node(Q3d) at (0,4){$2$};\node (Q4a) at (2,0){$3$};\node (Q4b) at (2,4){$3$};
%edge
\draw[lightgray, line width=0.05cm]  (P1)--(P2)--(P3)--(P4)--(P1);\draw[lightgray, line width=0.05cm] (0,1)--(P1)--(1,0); 
\draw[lightgray, line width=0.05cm]  (4,1)--(P2)--(3,0);\draw[lightgray, line width=0.05cm]  (0,3)--(P4)--(1,4);\draw[lightgray, line width=0.05cm]  (3,4)--(P3)--(4,3);

%black
\filldraw  [ultra thick, draw=lightgray, fill=lightgray] (1,1) circle [radius=0.18] ;\filldraw  [ultra thick, draw=lightgray, fill=lightgray] (3,3) circle [radius=0.18] ;
%white
\draw  [ultra thick, draw=lightgray,fill=white] (3,1) circle [radius=0.18] ;\draw  [ultra thick, draw=lightgray,fill=white] (1,3) circle [radius=0.18] ;

\draw[->, line width=0.06cm] (Q1)--(Q2a);\draw[->, line width=0.06cm] (Q2a)--(Q3a);\draw[->, line width=0.06cm] (Q3a)--(Q4a);
\draw[->, line width=0.06cm] (Q4a)--(Q1);\draw[->, line width=0.06cm] (Q2a)--(Q3d);\draw[->, line width=0.06cm] (Q3d)--(Q4b);
\draw[->, line width=0.06cm] (Q4b)--(Q1);\draw[->, line width=0.06cm] (Q1)--(Q2b);\draw[->, line width=0.06cm] (Q2b)--(Q3b);
\draw[->, line width=0.06cm] (Q3b)--(Q4a);\draw[->, line width=0.06cm] (Q2b)--(Q3c);\draw[->, line width=0.06cm] (Q3c)--(Q4b);
\end{tikzpicture}

}\]

\caption{Quiver associated with $\Gamma_{4a}$}
\label{ex_quiver4a}
\end{figure}

As the dual notion of a dimer model $\Gamma$, we define the quiver $Q_\Gamma$ associated with $\Gamma$. 
Namely, we assign a vertex dual to each face in $\Gamma_2$, an arrow dual to each edge in $\Gamma_1$.  
The orientation of arrows is determined so that the white vertex is on the right of the arrow (equivalently, the black vertex is on the left of the arrow). 
For example, Figure~\ref{ex_quiver4a} is the quiver associated with $\Gamma_{4a}$ appearing in Figure~\ref{ex_dimer4a}. 
Sometimes we simply denote such a quiver $Q_\Gamma$ by $Q$ if a situation is clear from the context. 
We denote the sets of vertices and arrows of the quiver by $Q_0$ and $Q_1$ respectively. 
In addition, we consider the set of oriented faces $Q_F$ as the dual of vertices in a dimer model $\Gamma$. 
The orientation of faces is determined by its boundary, that is, faces dual to white (resp. black) vertices are oriented clockwise (resp. anti-clockwise). 
Therefore, we decompose the set of faces as $Q_F=Q^+_F\sqcup Q^-_F$ where $Q^+_F$ (resp. $Q^-_F$) stands for the set of faces oriented clockwise (resp. anti-clockwise). 

%%%%%%%%%%%%%%%%%%%%%
%%% type 4a-1 %%%%%%%%%%%%%
%%%%%%%%%%%%%%%%%%%%%

For each arrow, we define the maps $h, t:Q_1\rightarrow Q_0$ sending an arrow to its head and tail respectively. 
A nontrivial path in $Q$ is a finite sequence of arrows $a=a_1\cdots a_r$ satisfying $h(a_\ell)=t(a_{\ell+1})$ for $\ell=1, \cdots r-1$. 
In this situation, we define the length of path $a=a_1\cdots a_r$ as $r \,(\ge 1)$, and denote by $Q_r$ the set of paths of length $r$ in $Q$. 
From this viewpoint, we extend the maps $h, t$ to the maps on paths, that is, we define $t(a)=t(a_1), h(a)=h(a_r)$. 
We say that a path $a$ is a \emph{cycle} if it satisfies $h(a)=t(a)$. We consider each vertex $i\in Q_0$ as a trivial path $e_i$ where $h(e_i)=t(e_i)=i$. 
Also, we denote by $Q^{\rm op}$ the opposite quiver of $Q$ which is obtained from $Q$ by reversing all arrows. 
By replacing white vertices on a dimer model by black vertices and vice versa, we obtain the opposite quiver associated with the original dimer model. 

From this quiver $Q$, we define the path algebra as follows. 
The \emph{path algebra} $KQ$ is a $K$-algebra whose $K$-basis consists of paths in $Q$ and the product is defined as 
$a\cdot b=ab$ (resp. $a\cdot b=0$) if $h(a)=t(b)$ (resp. $h(a)\neq t(b)$) for paths $a$ and $b$, and we extend this product linearly. 
Next, we define a potential, and give relations on paths in $KQ$. 
We denote the $K$-vector space generated by all commutators in $KQ$ by $[KQ, KQ]$ and set $KQ_{\mathrm{cyc}}\coloneqq KQ/[KQ, KQ]$. 
Thus, $KQ_{\mathrm{cyc}}$ has a basis consists of cycles in $Q$. 
Also, we denote $(KQ_{\mathrm{cyc}})_r$ the subspace of $KQ_{\mathrm{cyc}}$ spanned by cycles of length at least $r$. 

\begin{definition}
An element $W\in (KQ_{\mathrm{cyc}})_2$ is called a \emph{potential} (or \emph{superpotential}), 
and a pair $(Q,W)$ consisting of a quiver $Q$ and a potential $W$ is called a \emph{quiver with potential} (= \emph{QP}). 
\end{definition}

For each face $f\in Q_F$, we associate the \emph{small cycle} $\omega_f\in (KQ_{\mathrm{cyc}})_3$ obtained by the product of all arrows 
around the boundary of $f$. 
(Since we assume that a dimer model has no bivalent vertices in this paper, the length of each small cycle is at least $3$.) 
For the quiver $Q$ associated with a dimer model, we define the potential $W_Q$ as 
\[
W_Q\coloneqq \sum_{f\in Q^+_F}\omega_f-\sum_{f\in Q^-_F}\omega_f. 
\]
In the rest of this paper, we consider the QP $(Q,W_Q)$ associated with a dimer model. 

\begin{figure}[h]
\[\scalebox{0.75}{
\begin{tikzpicture}
%vertex
\node (P1) at (0,0){$$}; \node (P2) at (2,0){$$}; 

%edge
\draw[line width=0.03cm]  (P1)--(P2); \draw[line width=0.03cm]  (P1)--(-1.2,1); \draw[line width=0.03cm]  (P1)--(-1.2,-1);
\draw[line width=0.03cm]  (P2)--(3.2,1.2);  \draw[line width=0.03cm]  (P2)--(3.8,0.7); \draw[line width=0.03cm]  (P2)--(3.2,-1); 
%black
\filldraw  [thick, fill=black] (0,0) circle [radius=0.14] ;
%white
\draw  [ thick, fill=white] (2,0) circle [radius=0.14] ;
%arrow
\draw[->, line width=0.05cm]  (1,-1)--(1,1); 
\draw[->, line width=0.05cm]  (0.8,1)--(-1.2,0.1); 
%\draw[dotted, line width=0.05cm]  (-1.2,0.4)--(-1.2,-0.4); 
\draw[->, line width=0.05cm]  (-1.2,-0.1)--(0.8,-1); 
\draw[->, line width=0.05cm]  (1.2,1)--(3.2,0.7); 
\draw[->, line width=0.05cm]  (3.3,0.7)--(3.4, -0.2); 
%\draw[dotted, line width=0.05cm]  (3.6,-0.1)--(3.3,-0.5); 
\draw[->, line width=0.05cm]  (3.4,-0.3)--(1.2,-1); 

\node  at (1.2,0.3) {$a$} ;\node  at (-1.7,0.2) {$p^-_a$} ;\node  at (4,0.2) {$p^+_a$} ;
\end{tikzpicture}
}\]
\caption{An example of paths $p_a^{\pm}$}
\label{relation}
\end{figure}
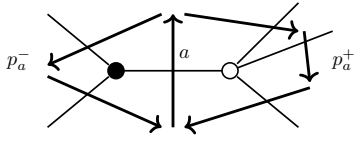

For every face $f\in Q_F$ and arrow $a\in \omega_f$, we choose $h(a)$ as the starting point of $\omega_f$ 
and set $e_{h(a)}\omega_fe_{h(a)}\coloneqq a_1\cdots a_ra$ for some path $a_1\cdots a_r$. 
Then, the partial derivative of $\omega_f$ with respect to $a$ is defined by the path $\partial \omega_f/\partial a\coloneqq a_1\cdots a_r$. 
We note that $\partial \omega_f/\partial a=0$ if $\omega_f$ does not contain the arrow $a$. 
Extending this derivative linearly, we have $\partial W_Q/\partial a\in KQ$ for any $a\in Q_1$ and define the two-sided ideal 
$J(W_Q)\coloneqq \langle\partial W_Q/\partial a\,|\,a\in Q_1\rangle\subset KQ$. 
Then, we define the \emph{Jacobian algebra} (or \emph{superpotential algebra}) of a dimer model as 
\[
\calP(Q, W_Q)\coloneqq KQ/J(W_Q) .
%=kQ/\langle p^+_a-p^-_a \; | \; a\in Q_1 \rangle. 
\]
From the construction, $\partial W_Q/\partial a$ gives the relation on paths in $\calP(Q, W_Q)$ for each arrow $a\in Q_1$. 
Namely, for each arrow $a\in Q_1$, there are precisely two oppositely oriented cycles which contain the arrow $a$ as a boundary. 
We denote them by $f_a^+, f_a^-\in Q_F$ respectively. Let $p_a^{\pm}$ be the path from $h(a)$ around the boundary of $f^{\pm}_a$ to $t(a)$ 
(see Figure~\ref{relation}). 
Then, we can describe $\partial W_Q/\partial a$ as a difference of $p_a^{\pm}$, that is, $\partial W_Q/\partial a=p_a^+-p_a^-$. 
Thus, we have $p_a^+=p_a^-$ in $\calP(Q, W_Q)$ for each arrow $a\in Q_1$.  

We say that two paths $a$ and $b$ in $Q$ are \emph{equivalent} (denoted by $a\sim b$) if they give the same element in the Jacobian algebra $\calP(Q, W_Q)$. 
For any face $f\in Q_F$ whose boundary factor through a vertex $i\in Q_0$, we define a small cycle starting from $i$ as $\omega_{i, f}\coloneqq e_i\omega_fe_i$. 
It is easy to see that small cycles starting from the same vertex are equivalent to each other.

%%%%%%%%%%%%%%%%%%%%%%%%%%%%%%%%%%%%%%%%%%%%%%%%%%%%%%%%%%%%%%%%%%%%%%
\subsection{Consistency conditions and NCCRs}
\label{sec_consist}

In this subsection, we impose the extra condition called consistency condition. 
Under this assumption, a dimer model will give an NCCR of a three dimensional Gorenstein toric singularity associated with a dimer model (see Theorem~\ref{NCCR1}). 
In the literature, there are several consistency conditions, for example, 
\begin{enumerate} [\rm(a)]
\item the consistency condition in the sense of \cite[Definition~3.5]{IU1}, 
\item a properly-ordered dimer model in the sense of \cite{Gul}, 
\item the first consistency condition (or cancellation) in the sense of \cite{MR} (see also \cite{Dav}), 
\item an algebraically consistent dimer model in the sense of \cite{Bro}, 
\item the existence of a consistent $\sfR$-charge in the sense of \cite{Kenn},
\item a toric order in the sense of \cite{Boc3},  
\item an isoradial (or a geometrically consistent) dimer model (e.g., \cite{KS,Bro}).  
\end{enumerate}
It is known that {\rm(a)} and {\rm(b)} are equivalent \cite{IU1}, {\rm(a)} and {\rm(c)} are equivalent under the non-degeneracy condition \cite{IU1}, 
${\rm(c)}, {\rm(d)}, {\rm(e)}$ and {\rm(f)} are equivalent \cite{Boc1}. 
Since we assume that the non-degeneracy, the conditions {\rm(a)}--{\rm(f)} are equivalent. 
Furthermore, the condition {\rm(g)} implies these equivalent conditions (see e.g., \cite{IU1,Bro}). 
We remark that a consistent dimer model is always non-degenerate (see e.g., \cite[Proposition~8.1]{IU2}). 
Here, we note the condition {\rm(e)}.  

\begin{definition}
Let $Q$ be a quiver associated with a dimer model $\Gamma$. A positively grading $\sfR:Q_1\rightarrow \RR_{>0}$ that 
satisfies the following conditions is called a \emph{consistent $\sfR$-charge} (cf. \cite{Kenn}): 
\begin{itemize}
\item [(1)] $\displaystyle\sum_{a\in\partial f}\sfR(a)=2$ \;for any $f\in Q_F$,  
\item [(2)] $\displaystyle\sum_{h(a)=i}(1-\sfR(a))+\displaystyle\sum_{t(a)=i}(1-\sfR(a))=2$ \;for any $i\in Q_0$. 
\end{itemize}
We say that a dimer model $\Gamma$ is \emph{consistent} if it admits a consistent $\sfR$-charge.  
\end{definition} 

As we showed in subsection~\ref{subsec_dimer}, for a given dimer model $\Gamma$, 
we can construct the three dimensional Gorenstein toric singularity $R=K[\sigma_\Delta^\vee\cap\sfM]$. 
In the rest of this article, we only consider such a toric singularity $R$ obtained from a consistent dimer model. 

\medskip

By the dual point of view, we can consider a perfect matching to be a function on $Q_1$. 
That is, for each arrow $a\in Q_1$, we define 
\begin{eqnarray*}
\sfP(a)=\left \{\begin{array}{ll}
1&\text{if the corresponding edge to $a$ is in $\sfP$}\\
0&\text{otherwise}. \\
\end{array} \right. 
\end{eqnarray*}
We note that by the definition of perfect matching, we have that $\sfP(\omega)=1$ for each small cycle $\omega$. 
Also, we extend this function on the double quiver of $Q$. 
The double quiver $\widetilde{Q}$ is obtained from $Q$ by adding the oppositely directed arrow $a^*$ for every $a\in Q_1$. 
%and we introduce the new relations on the path algebra $K\widetilde{Q}$ as $aa^*=e_{t(a)}$, $a^*a=e_{h(a)}$ for any $a\in Q_1$ together with the original ones. 
Then, we define $\sfP(a^*)=-\sfP(a)$. 

\medskip

Here, we note the important property of consistent dimer models. 

\begin{proposition}[{\cite[Corollary~4.27 ]{Bro}, \cite[Proposition~9.2]{IU2}}] 
\label{ex_pm}
If a dimer model is consistent, then the multiplicity of a vertex of the perfect matching polygon $\Delta$ is equal to one. 
That is, for each vertex of $\Delta$, there is the extremal perfect matching which corresponds to it. 
Furthermore, every edge is contained in some extremal perfect matchings. 
\end{proposition}

Let $\sfP_1, \cdots, \sfP_n$ be the extremal perfect matchings of a consistent dimer model corresponding to $v_1, \cdots, v_n$ respectively. 
For $i, j\in Q_0$, let $a_{ij}$ be a path from $i$ to $j$ (i.e., $h(a_{ij})=j$ and $t(a_{ij})=i$). 
For a three dimensional Gorenstein toric singularity $R=K[\sigma_\Delta^\vee\cap\ZZ^3]$ which is constructed from a dimer model (see subsection~\ref{subsec_dimer}), we define the divisorial ideal of $R$ associated with $a_{ij}$ as 
\[
T_{a_{ij}}\coloneqq T(\sfP_1(a_{ij}), \cdots, \sfP_n(a_{ij})). 
\]
As the following lemma shows, this ideal only depends on the starting point $i$ and the ending point $j$, whereas a path from $i$ to $j$ is not unique.  
Thus, we simply denote it by $T_{ij}$. When we want to specify the considering path, we use the former notation. 

\begin{lemma}
\label{equiv_div_ideal}
Let $Q$ be the quiver associated with a consistent dimer model. 
\begin{itemize}
\item [(1)] If two paths $a, b\in Q$ are equivalent, we have that $T_a\cong T_b$. 
\item [(2)] Let $\omega$ be a small cycle. Then, we have that $T_{\omega^m}\cong R$ for any $m\in\ZZ_{\ge0}$. 
\item [(3)] Let $a_{ij}, b_{ij}$ be paths from $i$ to $j$. Then, we have that $T_{a_{ij}}\cong T_{b_{ij}}$. 
\end{itemize}
\end{lemma}

\begin{proof}
$(1)$ We may assume that $a=p^+_c, b=p^-_c$ for some path $c\in Q_1$. 
For each perfect matching $\sfP$, we have that 
\[
\sfP(p^+_c)-\sfP(p^-_c)=\sfP(p^+_c)+\sfP(c)-\sfP(c)-\sfP(p^-_c)=\sfP(\omega)-\sfP(\omega^\prime)=1-1=0, 
\] 
where $\omega, \omega^\prime$ are small cycles starting from the vertex $t(c)$. 
Therefore, we have that $(\sfP_1(a), \cdots, \sfP_n(a))=(\sfP_1(b), \cdots, \sfP_n(b))$, and hence $T_a\cong T_b$ by Lemma~\ref{div_eq}. 
%This implies $\Hom_R(T_b, T_a)\cong R$. 

$(2)$ We have that $(\sfP_1(\omega^m), \cdots, \sfP_n(\omega^m))=(m, \cdots, m)=(\lambda_1(0,0,m), \cdots, \lambda_n(0,0,m))$. 
By Lemma~\ref{div_eq}, we have the conclusion.  

$(3)$ Let $(Q,W_Q)$ be the QP associated with a consistent dimer model $\Gamma$. 
Let $R$ be the three dimensional Gorenstein toric singularity arising from $\Gamma$. 
By \cite[Chapter~5 and 8]{Bro}, we have a basic splitting reflexive $R$-module $M=\bigoplus_{i\in Q_0}M_i$ such that $\calP(Q,W_Q)\cong\End_R(M)$, 
and this isomorphism induces $e_i\calP(Q,W_Q)e_j\cong\Hom_R(M_i,M_j)$. 
In particular, we have $e_i\calP(Q,W_Q)e_i\cong\Hom_R(M_i,M_i)\cong R$ for any $i\in Q_0$. 
By the same argument as \cite[The proof of Theorem~6.6]{Boc2}, we have that $T_{a_{ij}}\cong e_i\calP(Q,W_Q)e_j$ for any path $a_{ij}$ from $i$ to $j$. 
\end{proof}

Using these divisorial ideals, we obtain an NCCR of $R$ as follows. 

\begin{theorem}[{\cite[Theorem~6.6]{Boc2}, see also \cite{Bro,IU2}}] 
\label{NCCR1}
Suppose that $(Q, W_Q)$ is the QP associated with a consistent dimer model and $\calP(Q, W_Q)$ is the Jacobian algebra.  
Fix a vertex $i\in Q_0$, then we have that 
\[
\calP(Q, W_Q)\cong\End_R(\bigoplus_{j\in Q_0}T_{ij}), 
\]
and the center of $\calP(Q, W_Q)$ is just a three dimensional Gorenstein toric singularity $R$ constructed from a given dimer model. 
Moreover, this endomorphism ring is a splitting NCCR of $R$. 
\end{theorem}
 
Note that $T^i\coloneqq\oplus_{j\in Q_0}T_{ij}$ is an MM module by Propostion~\ref{maximal_NCCR}. 
Since $T^i$ contains $R\cong T_{ii}$ as a direct summand for any fixed vertex $i\in Q_0$, 
$T^i$ is a generator and we have that $T_{ij}\in\CM R$ for any $i, j \in Q_0$. 
Since $e_i$ is the idempotent corresponding to $i\in Q_0$, we have that $T^i\cong e_i\calP(Q, W_Q)$. 
Also, since $T_{ij}^* \cong T_{ji}$, we have that 
\[
\calP(Q, W_Q)\cong\End_R(T^i)\cong\End_R((T^i)^*)\cong\calP(Q^{\rm op}, W_{Q^{\rm op}}). 
\]

In this way, we obtain three dimensional Gorenstein toric singularities and their NCCRs from consistent dimer models. 
Conversely, for every lattice polygon $\Delta$ in $\RR^2$ (equivalently, for every three dimensional Gorenstein toric singularity), 
there exists a consistent dimer model whose perfect matching polygon coincides with $\Delta$ \cite{Gul,IU2}. 
(Note that such a consistent dimer model is not unique in general.) 
Therefore, we have the following. 

\begin{corollary}
Every three dimensional Gorenstein toric singularity has NCCRs arising from consistent dimer models. 
\end{corollary}

%%%%%%%%%%%%%%%%%%%%%%%%%%%%%%%%%%%%%%%%%%%%%%%%%%%%%%%%%%%%%%%%%%%%%%
\section{Mutations of dimer models and splitting MM generators} 
\label{sec_toric_mutation}

In this section, we will apply the mutations of MM modules introduced in subsection~\ref{subsec_MMmutation} to modules giving NCCRs of 
three dimensional Gorenstein toric singularities. 
Although the basic idea is the same as \cite{Boc2}, we will discuss details for the sake of completeness, 
and for the convenience of the readers. 

As we saw in subsection~\ref{sec_consist}, there is a consistent dimer model such that $\calP(Q, W_Q)\cong\End_R(\bigoplus_{j\in Q_0}T_{ij})$ 
is an NCCR of a given three dimensional Gorenstein toric singularity $R$. 
Especially, $\bigoplus_{j\in Q_0}T_{ij}$ is a splitting MM generator which gives a splitting NCCR. 
On the other hand, Bocklandt showed that every splitting NCCR is coming from a consistent dimer model \cite{Boc3}. 
%Thus, using these facts and Proposition~\ref{maximal_NCCR}, we see that splitting MM modules are precisely modules giving NCCRs arising from consistent dimer models. 
In this section, we consider such a splitting MM generator and a splitting NCCR (which is also called a \emph{toric NCCR} in \cite{Boc2}) 
arising from a consistent dimer model. 
Also, we remark that we always assume that a splitting MM generator is basic. 
Meanwhile, Bruns and Gubeladze proved that 
every toric singularity has only finitely many rank one MCM modules up to isomorphism \cite[Corollary~5.2]{BG1}. 
Collectively, we have the following corollary. 

\begin{corollary}
For any three dimensional Gorenstein toric singularity, there are only finitely many splitting NCCRs, and they arise from consistent dimer models. 
\end{corollary}

Therefore, it is possible to classify all splitting NCCRs of $R$. 
In fact, Bocklandt discussed an algorithm for finding all of them \cite{Boc2}. 
The point is to find maximal pairs of rank one MCM modules $\{M_0, \cdots, M_r\}$ which satisfy $\End_R(M_0\oplus\cdots\oplus M_r)\in\CM R$ 
because of the maximality of an NCCR (see Proposition~\ref{maximal_NCCR}). 
Then, we can reconstruct a consistent dimer model from such a pair $\{M_0, \cdots, M_r\}$ \cite[Theorem~5.9]{CV} (see also \cite[Section~6]{Boc2}). 
Precisely, suppose that $M\coloneqq M_0\oplus\cdots\oplus M_r$ gives a splitting NCCR, 
by using the methods as in \cite{CV,Boc2}, we obtain a consistent dimer model $\Gamma$ such that 
the associated Jacobian algebra $\calP(Q_\Gamma, W_{Q_\Gamma})$ is isomorphic to $\End_R(M)$ as an $R$-algebra. 
%Also, there is a one to one correspondence between direct summands $M_i$ of $M$ and vertices of $Q_\Gamma$ by a construction. 
We remark that different pairs of rank one MCM modules often give distinct consistent dimer models associated with $R$. 
Under these backgrounds, our first question is the following. 

\begin{question}
\label{que_mutation_dimer}
Let $\{\Gamma_1, \cdots, \Gamma_s\}$ be the complete set of consistent dimer models giving splitting NCCRs of $R$. 
Is there a relationship between these dimer models ?
\end{question}

In order to investigate this question, we will try to define the mutation of dimer models. 
As the first step, we consider the mutation of QPs introduced in \cite{DWZ}. To define this mutation, we prepare some notations. 
In what follows, we consider the complete path algebra of $Q$. 
Namely, let $Q_r$ be the set of paths of length $r$ in $Q$. Then, the \emph{complete path algebra} is defined as 
\[
\widehat{KQ}\coloneqq\prod_{r\ge 0}KQ_r 
\]
where $KQ_r$ is the vector space with a basis $Q_r$. 
Note that the multiplication is induced by the same way as $KQ$. Also,we set 
$\fkm_Q\coloneqq \prod_{r\ge 1}KQ_r$. 
For any subset $U\subseteq\widehat{KQ}$, we define the \emph{$\fkm_Q$-adic closure} of $U$ by 
$\overline{U}\coloneqq \bigcap_{n\ge 0}(U+\fkm_Q^n)$. 
Then, we define the notion of right equivalence introduced in \cite{DWZ}. 

\begin{definition} 
We say that two potentials $W_1, W_2$ associated with a quiver $Q$ are \emph{cyclically equivalent} if $W_1-W_2\in\overline{[KQ, KQ]}$. 

Let $(Q, W), (Q^\prime, W^\prime)$ be QPs. 
We say that $(Q, W)$ and $(Q^\prime, W^\prime)$ are \emph{right equivalent} if $Q_0=Q_0^\prime$ and 
there exists a isomorphism $\varphi:\widehat{KQ}\rightarrow\widehat{KQ^\prime}$ such that $\varphi\mid_{Q_0}=id$, 
and $\varphi(W)$ and $W^\prime$ are cyclically equivalent. 
\end{definition}

In the rest of this paper, by the abuse of the notation, we denote the complete Jacobian algebra $\widehat{KQ}/\overline{J(W)}$ by $\calP(Q, W)$. 
We remark that if QPs $(Q, W), (Q^\prime, W^\prime)$ are right equivalent, 
then we have that $\calP(Q, W)\cong\calP(Q^\prime, W^\prime)$ \cite[Proposition~3.3 and 3.7]{DWZ}. 
Therefore, we may only consider right equivalence classes of QPs. 

\begin{definition}
We say that a QP $(Q, W)$ is \emph{trivial} if $W$ is consisting of cycles of length $2$ and $\calP(Q,W)$ is isomorphic to $\widehat{KQ_0}$, 
and a QP $(Q, W)$ is \emph{reduced} if $W$ is a linear combination of cycles of length $3$ or more. 
\end{definition}

By the splitting theorem \cite[Theorem~4.6]{DWZ}, every QP $(Q,W)$ is decomposed as a direct sum of a trivial QP $(Q_{\rm triv}, W_{\rm triv})$ and 
a reduced QP $(Q_{\rm red}, W_{\rm red})$: 
\[
(Q, W)\cong (Q_{\rm triv}, W_{\rm triv})\oplus(Q_{\rm red}, W_{\rm red}), 
\]
up to right equivalence where $Q_0=(Q_{\rm triv})_0=(Q_{\rm red})_0$, $Q_1=(Q_{\rm triv})_1\sqcup(Q_{\rm red})_1$, 
and $W=W_{\rm triv}+W_{\rm red}$. 

\medskip

Now we start to introduce the mutation of QPs (see also \cite[Section~5]{DWZ}, \cite[subsection~1.2]{BIRS} for more details). 
Let $(Q, W)$ be a QP. 
First, we define $(Q^\prime, W^\prime)=\widetilde{\mu_k}(Q, W)$ by the procedure below for each vertex $k$ not lying on a $2$-cycle. 

Since a potential is defined as an element in $KQ_{\mathrm{cyc}}$, we may assume that no cycles in $W$ start at $k$. 
% (equivalently, end at $k$). 

\begin{itemize}
\item [(A)] $Q^\prime$ is a quiver obtained from $Q$ as follows. Fix a vertex $k\in Q_0$ not lying on a $2$-cycle. 
 \begin{itemize}
 \item [(A-1)] For each pair of arrows $a: i\rightarrow k$ and $b: k\rightarrow j$ in $Q_1$, we add a new arrow $ab: i\rightarrow j$. 
 \item [(A-2)] Replace each arrow $a: i\rightarrow k \in Q_1$ by a new arrow $a^*: k\rightarrow i$.  
 \item [(A-3)] Replace each arrow $b: k\rightarrow j \in Q_1$ by a new arrow $b^*: j\rightarrow k$. 
 \end{itemize}

\item [(B)] $W^\prime=[W]+\Omega$ where $[W]$ and $\Omega$ are obtained by the following fashion. 
 \begin{itemize}
 \item [(B-1)] $[W]$ is obtained by substituting $[ab]$ for each part $ab$ in $W$ which satisfies $a:i\rightarrow k$ and $b:k\rightarrow j$. 
 \item [(B-2)] $\Omega=\displaystyle\sum_{\genfrac{}{}{0pt}{}{a, b\in Q_1,}{h(a)=k=t(b)}} a^*[ab]b^*$. 
 \end{itemize}
\end{itemize}

\begin{definition}
Let $k\in Q_0$ be a vertex not lying on a $2$-cycle. 
We define the \emph{mutation $\mu_k(Q, W)$ of a QP $(Q, W)$ at $k$} as a reduced part of $(Q^\prime, W^\prime)$, 
that is, $\mu_k(Q, W)\coloneqq (Q^\prime_{\rm red}, W^\prime_{\rm red})$. 
\end{definition}

In this situation, $k$ is not lying on a $2$-cycle in $\mu_k(Q, W)$, and $\mu_{k}\left(\mu_k(Q, W)\right)$ is right equivalent to $(Q, W)$ \cite{DWZ}. 
If both $Q$ and the quiver part of $\mu_k(Q, W)$ do not contain loops and 2-cycles, then the above mutation coincides with Fomin-Zelevinsky's mutation \cite{FZ1}. 

Until now, we introduced the mutation of a QP. 
Then, we apply these procedure to the QP associated with a dimer model. 
However, the mutated QP at some vertex is not the dual of a dimer model in general. 
To make the mutated QP the dual of a dimer model, we have to impose some restrictions. 

\begin{lemma}
Let $\Gamma$ be a dimer model, and $(Q, W_Q)$ be the QP associated with $\Gamma$.
We assume that $k\in Q_0$ is not lying on $2$-cycles, and has no loops. 
Then, the mutation of $(Q, W_Q)$ at $k\in Q_0$ is the dual of a dimer model if and only if the vertex $k$ satisfies the following conditions: 
\[
(\clubsuit)\, \text{a vertex $k\in Q_0$ has exactly two incoming arrows, and exactly two outgoing arrows.}
\] 
\end{lemma}

\begin{proof}
First, we have the following by the definition of a dimer model. 
\begin{itemize}
\item[$\cdot$] For each vertex, the number of incoming arrows coincides with that of outgoing arrows. 
%(Thus, in the condition $(\clubsuit)$ we may only assume the one side.) 
\item[$\cdot$] For each vertex, there are at least two incoming arrows (thus, there are at least two outgoing arrows). 
\item[$\cdot$] Incoming arrows and outgoing arrows appear alternately around a vertex. 
\end{itemize}
Therefore, if the number of incoming arrows (= that of outgoing arrows) around $k\in Q_0$ is greater than or equal to $3$, 
then some arrows added in the process of the mutation intersect with other arrows. Thus, it is no longer the dual of a dimer model. 
Conversely, we suppose that $k\in Q_0$ satisfies $(\clubsuit)$. 
Then, we see that the mutated QP at $k$ is tiling the two-torus $\sfT$ by combining clockwise faces and anti-clockwise faces, 
so it is the dual of a dimer model. 
\end{proof}

We denote by $Q_0^\mu$ the set of vertices in a quiver $Q$ that are not lying on $2$-cycles, have no loops, and satisfy the condition $(\clubsuit)$. 
We will consider the mutation at a vertex $k\in Q_0^\mu$. 
Meanwhile, each cycle in the potential obtained from a dimer model has the sign corresponding to its orientation. 
Thus, we also have to modify the sign of the mutated QP. 
Since a vertex $k\in Q_0^\mu$ has two incoming arrows and two outgoing arrows, we can divide them into two pairs consisting of 
an outgoing arrow and an incoming arrow.  
We denote them $\{a,b\}, \{a^\prime,b^\prime\}$ respectively where $a, a^\prime$ (resp. $b, b^\prime$) are incoming (resp. outgoing) arrows. 
Notice that there are two possibilities of choice of such pairs. (They belong to clockwise cycles or anti-clockwise cycles.)  
In the rest, we assume that $\{a,b\}, \{a^\prime,b^\prime\}$ belong to anti-clockwise cycles. 
Then, we replace the procedure (B-2) by 
\[
\text{(B-2$^\prime$)}\,\,\, \Omega=a^*[ab]b^*-a^*[ab^\prime](b^\prime)^*+(a^\prime)^*[a^\prime b^\prime](b^\prime)^*-(a^\prime)^*[a^\prime b]b^*. 
\]
Note that the potential $\Omega$ in (B-2$^\prime$) is obtained from the potential $\Omega$ in (B-2) by sending $\{a^\prime,b^\prime\}$ to $\{-a^\prime,-b^\prime\}$, 
and it induces an automorphism on the Jacobian algebra. 
Therefore, we define the mutation of a dimer model as follows. 

\begin{definition}
Let $\Gamma$ be a dimer model, and $(Q, W_Q)$ be the QP associated with $\Gamma$. 
For a vertex $k\in Q_0^\mu$, we define the \emph{mutation of a dimer model} $\Gamma$ at $k$ as the mutation of the associated QP $(Q, W_Q)$ at $k$, 
that is, it consists of the procedures (A-1),(A-2),(A-3),(B-1) and (B-2$^\prime$).  
\end{definition}

Next, we show examples of the mutation of a dimer model. 

\begin{example}
\label{mutation_QP1}
Let $Q$ be the quiver as in the left hand side of the figure below, and 
\[
W\coloneqq d_2c_1a_2b_1-d_1c_1a_1b_1+d_1c_2a_1b_2-d_2c_2a_2b_2
\]
be a potential of $Q$. 
Note that this QP is the one associated with $\Gamma_{4a}$ (see Figure~\ref{ex_dimer4a} and \ref{ex_quiver4a}), and all vertices are in $Q_0^\mu$.  
For example, we consider the mutation of $(Q, W)$ at $0\in Q_0^\mu$.  

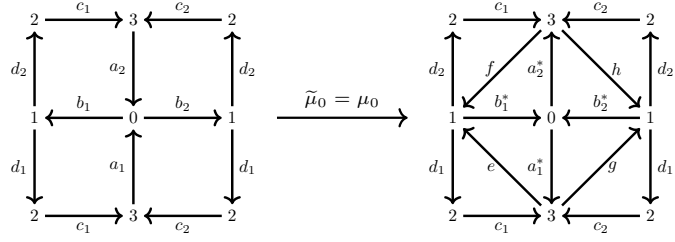
\begin{center} 
\begin{tikzpicture} 
\node (QPa) at (0,0) 
{\scalebox{0.65}{
\begin{tikzpicture}
%vertex
\node (Q0) at (2,2){$0$};\node (Q1a) at (0,2){$1$}; \node(Q1b) at (4,2){$1$};\node (Q2a) at (0,0){$2$};
\node(Q2c) at (4,4){$2$};\node(Q2b) at (4,0){$2$};\node(Q2d) at (0,4){$2$};\node (Q3a) at (2,0){$3$};\node (Q3b) at (2,4){$3$};

\draw[->, line width=0.05cm] (Q0)--node[above] {$b_1$}(Q1a);\draw[->, line width=0.05cm] (Q1a)--node[midway,left] {$d_1$}(Q2a);
\draw[->, line width=0.05cm] (Q2a)--node[below] {$c_1$}(Q3a);\draw[->, line width=0.05cm] (Q3a)--node[midway,left] {$a_1$}(Q0);
\draw[->, line width=0.05cm] (Q1a)--node[midway,left] {$d_2$}(Q2d);\draw[->, line width=0.05cm] (Q2d)--node[above] {$c_1$}(Q3b);
\draw[->, line width=0.05cm] (Q3b)--node[midway,left] {$a_2$}(Q0);\draw[->, line width=0.05cm] (Q0)--node[above] {$b_2$}(Q1b);
\draw[->, line width=0.05cm] (Q1b)--node[midway,right] {$d_1$}(Q2b);\draw[->, line width=0.05cm] (Q2b)--node[below] {$c_2$}(Q3a);
\draw[->, line width=0.05cm] (Q1b)--node[midway,right] {$d_2$}(Q2c);\draw[->, line width=0.05cm] (Q2c)--node[above] {$c_2$}(Q3b);
\end{tikzpicture} }}; 

\node (QPb) at (5.5,0) 
{\scalebox{0.65}{
\begin{tikzpicture}
%vertex
\node (Q0) at (2,2){$0$};\node (Q1a) at (0,2){$1$}; \node(Q1b) at (4,2){$1$};\node (Q2a) at (0,0){$2$};
\node(Q2c) at (4,4){$2$};\node(Q2b) at (4,0){$2$};\node(Q2d) at (0,4){$2$};\node (Q3a) at (2,0){$3$};\node (Q3b) at (2,4){$3$};

\draw[->, line width=0.05cm] (Q1a)--node[above] {$b_1^*$}(Q0);
\draw[->, line width=0.05cm] (Q1a)--node[midway,left] {$d_1$}(Q2a);
\draw[->, line width=0.05cm] (Q2a)--node[below] {$c_1$}(Q3a);
\draw[->, line width=0.05cm] (Q0)--node[midway,left] {$a_1^*$}(Q3a);
\draw[->, line width=0.05cm] (Q1a)--node[midway,left] {$d_2$}(Q2d);
\draw[->, line width=0.05cm] (Q2d)--node[above] {$c_1$}(Q3b);
\draw[->, line width=0.05cm] (Q0)--node[midway,left] {$a_2^*$}(Q3b);
\draw[->, line width=0.05cm] (Q1b)--node[above] {$b_2^*$}(Q0);
\draw[->, line width=0.05cm] (Q1b)--node[midway,right] {$d_1$}(Q2b);
\draw[->, line width=0.05cm] (Q2b)--node[below] {$c_2$}(Q3a);
\draw[->, line width=0.05cm] (Q1b)--node[midway,right] {$d_2$}(Q2c);
\draw[->, line width=0.05cm] (Q2c)--node[above] {$c_2$}(Q3b);

\draw[->, line width=0.05cm] (Q3a)--node[midway,left] {$e$}(Q1a);
\draw[->, line width=0.05cm] (Q3b)--node[midway,left] {$f$}(Q1a);
\draw[->, line width=0.05cm] (Q3a)--node[midway,right] {$g$}(Q1b);
\draw[->, line width=0.05cm] (Q3b)--node[midway,right,xshift=2pt] {$h$}(Q1b);

\end{tikzpicture} }}; 

\draw[->, line width=0.03cm] (QPa)--node[above] {{\scriptsize$\widetilde{\mu}_0=\mu_0$}}(QPb);
\end{tikzpicture}
\end{center}

Then, we obtain the quiver $Q^\prime$ as in the right hand side and the potential $W^\prime=[W]+\Omega$ where 
\[ [W]=d_2c_1f-d_1c_1e+d_1c_2g-d_2c_2h, \]
\[\Omega=a_1^*eb_1^*-a_2^*fb_1^*+a_2^*hb_2^*-a_1^*gb_2^*. \]
Since $\widetilde{\mu}_0(Q, W)=(Q^\prime, W^\prime)$ is reduced, we have that $\widetilde{\mu}_0(Q, W)=\mu_0(Q, W)$. 
%Remark that the mutated QP is just the QP appearing in subsection~\ref{type4a}. 
\end{example}

\begin{example}
\label{mutation_QP2}
Let $Q$ be the quiver as in the left hand side of the figure below, and 
\[W\coloneqq a_1b_1c_1-a_1b_2f_1+a_3f_3d_1-a_3f_2c_1+e_1e_2f_1-e_1a_2b_1c_2d_1+d_2a_2b_2f_2c_2-d_2e_2f_3 \]
be a potential (see subsection~\ref{type6b}). 
We consider the mutation of $(Q, W)$ at $1\in Q_0$.

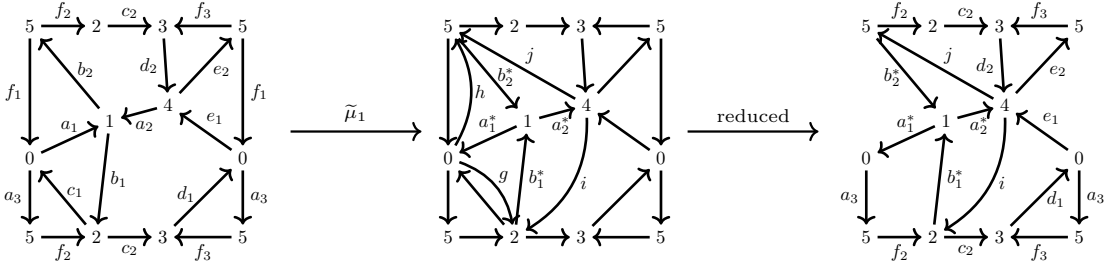
\begin{center} 
\begin{tikzpicture} 
\node (QPa) at (0,0) 
{\scalebox{0.7}{
\begin{tikzpicture}
\node (Q0a) at (0,1.5){$0$}; \node (Q0b) at (4,1.5){$0$}; 
\node (Q1) at (1.5,2.2){$1$}; \node (Q2a) at (1.25,0){$2$}; \node (Q2b) at (1.25,4){$2$}; 
\node (Q3a) at (2.5,0){$3$}; \node (Q3b) at (2.5,4){$3$}; \node (Q4) at (2.6,2.5){$4$}; 
\node (Q5a) at (0,0){$5$}; \node (Q5b) at (4,0){$5$}; \node (Q5c) at (4,4){$5$}; \node (Q5d) at (0,4){$5$}; 

\draw[->, line width=0.05cm] (Q0a)--node[above] {$a_1$}(Q1); \draw[->, line width=0.05cm] (Q0a)--node[midway,left] {$a_3$}(Q5a); 
\draw[->, line width=0.05cm] (Q0b)--node[midway,right] {$a_3$}(Q5b); \draw[->, line width=0.05cm] (Q0b)--node[above,xshift=5pt] {$e_1$}(Q4);  
\draw[->, line width=0.05cm] (Q1)--node[midway,right,yshift=2pt] {$b_2$}(Q5d); \draw[->, line width=0.05cm] (Q1)--node[midway,right] {$b_1$}(Q2a);
\draw[->, line width=0.05cm] (Q2a)--node[below] {$c_2$}(Q3a); \draw[->, line width=0.05cm] (Q2a)--node[midway,right,xshift=-2pt,yshift=3pt] {$c_1$}(Q0a); 
\draw[->, line width=0.05cm] (Q2b)--node[above] {$c_2$}(Q3b); \draw[->, line width=0.05cm] (Q3a)--node[midway,left,yshift=2pt] {$d_1$}(Q0b); 
\draw[->, line width=0.05cm] (Q3b)--node[midway,left] {$d_2$}(Q4); \draw[->, line width=0.05cm] (Q4)--node[below,xshift=3pt] {$a_2$}(Q1); 
\draw[->, line width=0.05cm] (Q4)--node[midway,right,yshift=-3pt] {$e_2$}(Q5c);  
\draw[->, line width=0.05cm] (Q5c)--node[midway,right] {$f_1$}(Q0b); \draw[->, line width=0.05cm] (Q5d)--node[midway,left] {$f_1$}(Q0a);  
\draw[->, line width=0.05cm] (Q5a)--node[below] {$f_2$}(Q2a); \draw[->, line width=0.05cm] (Q5d)--node[above] {$f_2$}(Q2b); 
\draw[->, line width=0.05cm] (Q5c)--node[above] {$f_3$}(Q3b); \draw[->, line width=0.05cm] (Q5b)--node[below] {$f_3$}(Q3a); 
 
\end{tikzpicture}
 }}; 

\node (QPb) at (5.5,0) 
{\scalebox{0.7}{
\begin{tikzpicture}
\node (Q0a) at (0,1.5){$0$}; \node (Q0b) at (4,1.5){$0$}; 
\node (Q1) at (1.5,2.2){$1$}; \node (Q2a) at (1.25,0){$2$}; \node (Q2b) at (1.25,4){$2$}; 
\node (Q3a) at (2.5,0){$3$}; \node (Q3b) at (2.5,4){$3$}; \node (Q4) at (2.6,2.5){$4$}; 
\node (Q5a) at (0,0){$5$}; \node (Q5b) at (4,0){$5$}; \node (Q5c) at (4,4){$5$}; \node (Q5d) at (0,4){$5$}; 

\draw[->, line width=0.05cm] (Q5d)--node[midway,right,xshift=1pt,yshift=-2pt] {$b_2^*$}(Q1); \draw[->, line width=0.05cm] (Q2a)--node[midway,right] {$b_1^*$}(Q1);
\draw[->, line width=0.05cm] (Q1)--node[above] {$a_1^*$}(Q0a); \draw[->, line width=0.05cm] (Q1)--node[below,xshift=2pt,yshift=2pt] {$a_2^*$}(Q4); 

\draw[->, line width=0.05cm] (Q0a) to [bend left]  node [midway, right,xshift=-2pt, yshift=5pt] {$g$}(Q2a);  
\draw[->, line width=0.05cm] (Q0a) to [bend right]  node [midway, right,xshift=-2pt] {$h$}(Q5d);  
\draw[->, line width=0.05cm] (Q4) to [bend left]  node [midway, right] {$i$}(Q2a);  
\draw[->, line width=0.05cm] (Q4) to node [above,xshift=7pt, yshift=-3pt] {$j$}(Q5d);  

\draw[->, line width=0.05cm] (Q0a)--(Q5a); \draw[->, line width=0.05cm] (Q0b)--(Q5b); \draw[->, line width=0.05cm] (Q0b)--(Q4);  
\draw[->, line width=0.05cm] (Q2a)--(Q3a); \draw[->, line width=0.05cm] (Q2a)--(Q0a); \draw[->, line width=0.05cm] (Q2b)--(Q3b); 
\draw[->, line width=0.05cm] (Q3a)--(Q0b); \draw[->, line width=0.05cm] (Q3b)--(Q4); 
\draw[->, line width=0.05cm] (Q4)--(Q5c);  \draw[->, line width=0.05cm] (Q5c)--(Q0b); \draw[->, line width=0.05cm] (Q5d)--(Q0a);  
\draw[->, line width=0.05cm] (Q5a)--(Q2a); \draw[->, line width=0.05cm] (Q5d)--(Q2b); 
\draw[->, line width=0.05cm] (Q5c)--(Q3b); \draw[->, line width=0.05cm] (Q5b)--(Q3a);

\end{tikzpicture}
 }}; 

\node (QPc) at (11,0) 
{\scalebox{0.7}{
\begin{tikzpicture}
\node (Q0a) at (0,1.5){$0$}; \node (Q0b) at (4,1.5){$0$}; 
\node (Q1) at (1.5,2.2){$1$}; \node (Q2a) at (1.25,0){$2$}; \node (Q2b) at (1.25,4){$2$}; 
\node (Q3a) at (2.5,0){$3$}; \node (Q3b) at (2.5,4){$3$}; \node (Q4) at (2.6,2.5){$4$}; 
\node (Q5a) at (0,0){$5$}; \node (Q5b) at (4,0){$5$}; \node (Q5c) at (4,4){$5$}; \node (Q5d) at (0,4){$5$}; 

\draw[->, line width=0.05cm] (Q5d)--node[midway,left,xshift=1pt,yshift=-2pt] {$b_2^*$}(Q1); \draw[->, line width=0.05cm] (Q2a)--node[midway,right] {$b_1^*$}(Q1);
\draw[->, line width=0.05cm] (Q1)--node[above] {$a_1^*$}(Q0a); \draw[->, line width=0.05cm] (Q1)--node[below,xshift=2pt,yshift=2pt] {$a_2^*$}(Q4); 

\draw[->, line width=0.05cm] (Q4) to [bend left]  node [midway, right] {$i$}(Q2a);  
\draw[->, line width=0.05cm] (Q4) to node [above,xshift=7pt, yshift=-3pt] {$j$}(Q5d);

\draw[->, line width=0.05cm] (Q0a)--node[midway,left] {$a_3$}(Q5a); 
\draw[->, line width=0.05cm] (Q0b)--node[midway,right] {$a_3$}(Q5b); \draw[->, line width=0.05cm] (Q0b)--node[above,xshift=5pt] {$e_1$}(Q4);  
\draw[->, line width=0.05cm] (Q2a)--node[below] {$c_2$}(Q3a);
\draw[->, line width=0.05cm] (Q2b)--node[above] {$c_2$}(Q3b); \draw[->, line width=0.05cm] (Q3a)--node[midway,right,yshift=-1pt] {$d_1$}(Q0b); 
\draw[->, line width=0.05cm] (Q3b)--node[midway,left] {$d_2$}(Q4); 
\draw[->, line width=0.05cm] (Q4)--node[midway,right,yshift=-3pt] {$e_2$}(Q5c);  
\draw[->, line width=0.05cm] (Q5a)--node[below] {$f_2$}(Q2a); \draw[->, line width=0.05cm] (Q5d)--node[above] {$f_2$}(Q2b); 
\draw[->, line width=0.05cm] (Q5c)--node[above] {$f_3$}(Q3b); \draw[->, line width=0.05cm] (Q5b)--node[below] {$f_3$}(Q3a); 
\end{tikzpicture}
 }}; 

\draw[->, line width=0.03cm] (QPa)--node[above] {{\scriptsize$\widetilde{\mu}_1$}} (QPb);
\draw[->, line width=0.03cm] (QPb)--node[above] {{\scriptsize\text{reduced}}} (QPc);

\end{tikzpicture}
\end{center}

Then, we obtain the QP $\widetilde{\mu}_1(Q, W)=(Q^\prime, W^\prime)$ where $Q^\prime$ is a quiver as in the middle of the above figure 
and $W^\prime=[W]+\Omega :$ 
\[ [W]=gc_1-hf_1+a_3f_3d_1-a_3f_2c_1+e_1e_2f_1-e_1ic_2d_1+d_2jf_2c_2-d_2e_2f_3, \]
\[\Omega=a_1^*hb_2^*-a_1^*gb_1^*+a_2^*ib_1^*-a_2^*jb_2^*. \]

Then, we remove arrows constructing $2$-cycles in $W^\prime$. Finally, we have the QP 
$\mu_1(Q, W)=(Q^\prime_{\rm red}, W^\prime_{\rm red})$ where $Q^\prime_{\rm red}$ is the quiver as in the right hand side, and 
\[
W^\prime_{\rm red}=a_3f_3d_1-a_3f_2b_1^*a_1^*+e_1e_2b_2^*a_1^*-e_1ic_2d_1+d_2jf_2c_2-d_2e_2f_3+a_2^*ib_1^*-a_2^*jb_2^*.
\]
\end{example}

Later, we will show that the mutated consistent dimer model is again consistent (see Corollary~\ref{mutated_consistent}). 
Now we are ready to consider Question~\ref{que_mutation_dimer}. 
For this question, Bocklandt pointed out the following. 

\begin{theorem}[{\cite[Theorem~7.6]{Boc2}}]
\label{mutation_Boc}
Let $R$ be a three dimensional Gorenstein toric singularity associated with a reflexive polygon (see Section~\ref{mutation_reflexive}). 
Suppose that $\fkD=\{\Gamma_1, \cdots, \Gamma_m\}$ is the complete set of consistent dimer models giving splitting NCCRs of $R$. 
Then, any two dimer models in $\fkD$ are transformed into each other by repeating mutations of dimer models. 
\end{theorem}

Note that a similar statement also holds for other special situations (see \cite[subsection~7.3, 7.4]{Boc2}). 

In this way, we could have the answer to Question~\ref{que_mutation_dimer} for the case of reflexive polygons. 
However, a splitting MM generator giving an NCCR obtained from a consistent dimer model is not unique. 
More precisely, we fix a consistent dimer model $\Gamma$, then there exists a splitting MM generator $M$ such that $\End_R(M)\cong\calP(Q_\Gamma, W_{Q_\Gamma})$, but such a module $M$ is not unique. Therefore, we should also ask the following. 

\begin{question}
\label{que_mutation_MM}
Let $\{\Gamma_1, \cdots, \Gamma_m\}$ be the complete set of consistent dimer models giving splitting NCCRs of $R$.  
We suppose that $\{M_{t1}, \cdots, M_{ts_t} \}$ is the set of non-isomorphic splitting MM generators such that $\End_R(M_{ts})\cong\calP(Q_{\Gamma_t}, W_{Q_{\Gamma_t}})$ 
for any $t=1, \cdots, m$ and $s=1, \cdots, s_t$. 
Is there a relationship between these splitting MM generators $M_{ts}$? 
\end{question}

In the rest, we suppose that $R=K[[ t_1^{\alpha_1}t_2^{\alpha_2}t_3^{\alpha_3}\mid (\alpha_1,\alpha_2,\alpha_3)\in\sigma^\vee\cap\ZZ^3 ]]$ 
is the $\fkm$-adic completion of a three dimensional Gorenstein toric singularity where $\fkm$ is the irrelevant maximal ideal. 
Note that we can obtain this singularity as the center of the complete Jacobian algebra of a certain dimer model. 
Let $M$ be a splitting MM generator. As we mentioned, we can construct a consistent dimer model $\Gamma$ and the associated QP $(Q,W_Q)$ 
that satisfies $\calP(Q,W_Q)\cong\End_R(M)$, and we may write $M=\bigoplus_{i\in Q_0}M_i$. Especially, we set $M_0=R$. 
In order to consider Question~\ref{que_mutation_MM}, we also consider the mutation of a splitting MM generator $M$ at $0\neq k\in Q_0$. 
By Proposition~\ref{prop_mutation_MM}(5), $\mu_k(M)$ is also an MM module, and it gives an NCCR. 
However, we again remark that even if $M$ is a splitting MM generator, $\mu_k(M)$ is not a splitting MM generator in general. 
Thus, we will consider the mutation of $M$ at $0\neq k\in Q_0^\mu$ to achieve our purpose. 
%Let $\Gamma$ be a consistent dimer model giving a splitting NCCR $\End_R(M)$, and $(Q, W_Q)$ be the associated QP. 
%We denote the vertex of $Q$ corresponding to an indecomposable direct summand $M_i\in\add_RM$ by $i\in Q_0$. 
First, for a vertex $0\neq k\in Q_0^\mu$, we suppose that $a_1, a_2\in Q_1$ (resp. $b_1, b_2\in Q_1$) are two incoming (resp. outgoing) arrows. 
Since an isomorphism $\calP(Q,W_Q)\cong\End_R(M)$ is established by sending a path $i\rightarrow j$ in $Q$ to an $R$-linear map $M_i\rightarrow M_j$ 
(see \cite[Section~5]{Bro}), we have a morphism 
\begin{equation}
\label{Mk_approx}
\varphi\coloneqq\varphi_{a_1}\oplus\varphi_{a_2}:M_{t(a_1)}\oplus M_{t(a_2)}\rightarrow M_k
\end{equation}
corresponding to arrows $a_1,a_2$. 
Since $R$-linear maps corresponding to arrows in $Q$ generate $\End_R(M)$, they contain generators of $\Hom_R(\oplus_{i\in Q_0\setminus\{k\}}M_i,M_k)$. 
Also, since arrows whose target is $k$ are only $a_1,a_2$, elements in $\Hom_R(\oplus_{i\in Q_0\setminus\{k\}}M_i,M_k)$ certainly factor through $\varphi_{a_1}$ or $\varphi_{a_2}$. 
Thus, we see that the morphism (\ref{Mk_approx}) is a right $(\add_R \oplus_{i\in Q_0\setminus\{k\}}M_i)$-approximation of $M_k$, 
and clearly this is minimal. 
Let $\calK_k\coloneqq\Ker\varphi$ be  the kernel of this morphism. 
Then, we have that $\mu_k(M)=\oplus_{i\in Q_0\setminus\{k\}}M_i\oplus \calK_k$.  
By counting the rank, we see that $\rank_R\calK_k=1$, and hence $\mu_k(M)$ is also a splitting MM generator. 
Conversely, if $\mu_k(M)$ is a splitting MM generator, then $k\in Q_0^\mu$ because $\mu_k(\mu_k(M))=M$. 
Thus, we have the following lemma. 

\begin{lemma}
\label{mutatable_MM} 
Let the notation be the same as above. 
Then, $\mu_k(M)$ is a splitting MM generator if and only if a vertex $0\neq k$ is in $Q_0^\mu$. 
\end{lemma}

For each MM module $M^\prime$, $T\coloneqq\Hom_R(M^\prime, M)$ is a tilting $\Lambda$-module where $\Lambda=\End_R(M^\prime)$ 
(see Theorem~\ref{MM_tilting}), and we have that 
\[
\End_\Lambda(T)\cong\End_R(M)\cong\calP(Q, W_Q). 
\]
Therefore, by \cite[Theorem~5.1]{BIRS} (see also \cite{KY,Vit}), we have the following. 

\begin{theorem}
Let the notations be the same as above. For a vertex $k\in Q_0^\mu$, we have that 
\[
\End_\Lambda(\mu_k(T))\cong\End_R(\mu_k(M))\cong\calP(\mu_k(Q, W_Q)), 
\]
where the above mutations stand for the mutation of tilting modules, that of splitting MM modules and that of dimer models (i.e., the mutation of QPs) respectively. 
\end{theorem}

Since the mutation of a splitting MM generator at $0\neq k\in Q_0^\mu$ is a splitting MM generator, we also have the following corollary. 

\begin{corollary}
\label{mutated_consistent}
The mutation of a consistent dimer model is also a consistent dimer model. 
\end{corollary}

Now, we consider Question~\ref{que_mutation_MM}.   
Especially, we observe the exchange graph $\EG(\TMMG(R))$, that is, 
the vertices of this graph are elements in $\TMMG(R)$, and draw an edge between $M$ and $\mu_k(M)$ for each $M\in\TMMG(R)$ and $0\neq k\in Q_0^\mu$. 
In the rest, we discuss whether $\EG(\TMMG(R))$ is connected or not. 
%With this terminology, our question is rephrased as ``\,is the exchange graph $\EG(\TMMG(R))$ connected?\," 

For the case of simplicial cones and the $A_1$-singularity (or conifold), we easily have the affirmative answer as in examples below. 
We note that the same statement as in Theorem~\ref{mutation_Boc} also holds for these cases, 
because they have a unique consistent dimer model. 

\begin{example}
\label{mutation_triangle}
Let $\Delta$ be a triangle polygon. Then, the associated cone $\sigma_\Delta$ is simplicial, 
and $R$ is a quotient singularity by a finite abelian group $G$ (see e.g., \cite[Example~1.3.20]{CLS}). 
In this case, $R$ has the unique basic splitting MM generator (see \cite[Theorem~3.1]{IN}). 
Especially, $\EG(\TMMG(R))$ is the single point, and there is the unique consistent dimer model giving such a splitting MM generator. 
Note that the associated quiver is the McKay quiver of $G$. (For more details, see \cite{UY}, \cite[Corollary~1.7]{IN}.) 
\end{example}

\begin{example}
\label{mutation_A1}
Let $R$ be the three dimensional $A_1$-singularity (i.e., $R\cong K[[x, y, z, w]]/(xy-zw)$). 
Remark that $R$ is of finite CM representation type, that is, it has only finitely many non-isomorphic indecomposable MCM modules, 
and finitely many MCM $R$-modules are $R$, $I=(x, z)$, $I^*=(x, w)$ (see e.g., \cite{Yos}). 
Then, modules giving an NCCR are only $R\oplus I$ and $R\oplus I^*$ \cite{VdB2}, and they are splitting MM generators. 
By taking a minimal right $\add R$-approximation of $I$: 
$0\rightarrow I^*\rightarrow R^2\rightarrow I\rightarrow 0$, we can connect $R\oplus I$ to $R\oplus I^*$ in $\EG(\TMMG(R))$, 
hence $\EG(\TMMG(R))$ is connected. 

Also, for these splitting MM generators, there is the unique consistent dimer model $\Gamma$ such that 
$\calP(Q_\Gamma, W_{Q_\Gamma})\cong\End_R(R\oplus I)\cong\End_R(R\oplus I^*)$. 
Here, $\Gamma$ takes the following form. 

\medskip

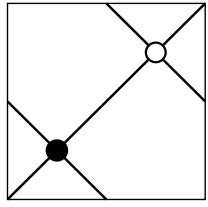
\begin{center}
\scalebox{0.65}{
\begin{tikzpicture}
%vertex
\node (B1) at (1,1){$$}; \node (W1) at (3,3){$$}; 
\draw[thick]  (0,0) rectangle (4,4);

%edge
\draw[line width=0.05cm]  (B1)--(W1); \draw[line width=0.05cm]  (B1)--(0,0);  \draw[line width=0.05cm]  (B1)--(2,0);  \draw[line width=0.05cm]  (B1)--(0,2); 
\draw[line width=0.05cm]  (W1)--(4,4);  \draw[line width=0.05cm]  (W1)--(2,4);  \draw[line width=0.05cm]  (W1)--(4,2); 
%black
\filldraw  [ultra thick, fill=black] (1,1) circle [radius=0.2] ;
%white
\draw  [line width=1.3pt,fill=white] (3,3) circle [radius=0.2] ;
\end{tikzpicture} }
\end{center}
\end{example}

%%%%%%%%%%%%%%%%%%%%%%%%%%%%%%%%%%%%%%%%%%%%%%%%%%%%%%%%%%%%%%%%%
\section{Mutations of splitting MM generators associated with reflexive polygons} 
\label{mutation_reflexive}

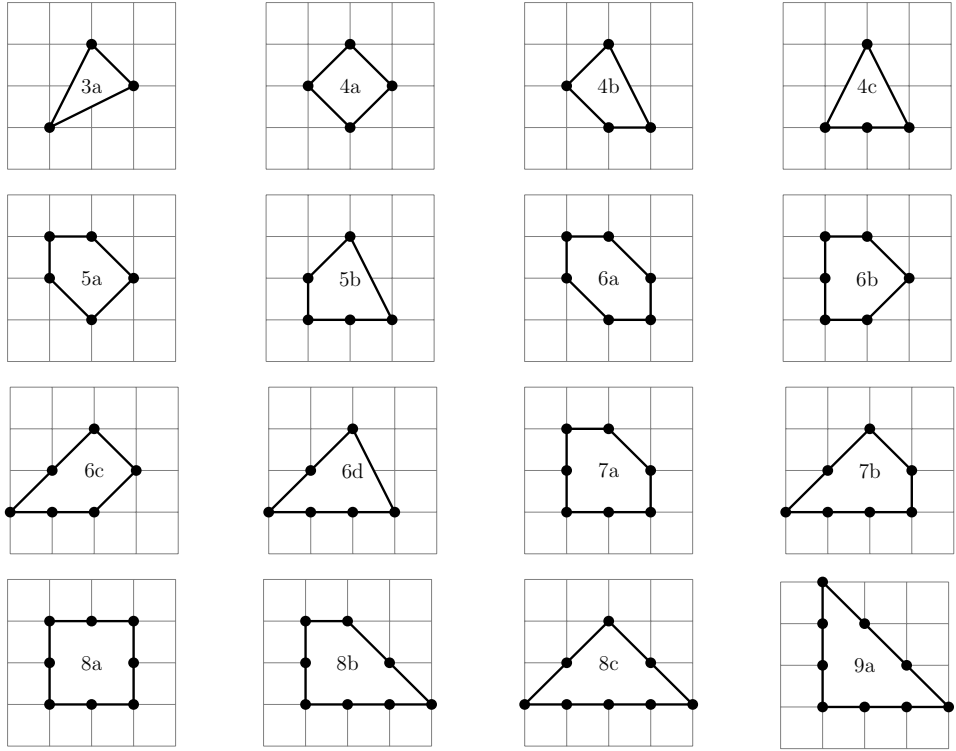
\begin{figure}[h]
\begin{center}
\scalebox{0.85}{
\begin{tikzpicture}

\node (3a) at (0,0){
%%% 3a %%%
\scalebox{0.65}{
\begin{tikzpicture}
\draw [step=1,thin, gray] (-2,-2) grid (2,2);
\filldraw [ultra thick, fill=white] (1,0)--(0,1)--(-1,-1)--cycle ;
\filldraw  [ultra thick, fill=black] (1,0) circle [radius=0.1] ;
\filldraw  [ultra thick, fill=black] (0,1) circle [radius=0.1] ;
\filldraw  [ultra thick, fill=black] (-1,-1) circle [radius=0.1] ;

\node  at (0,0) {\LARGE{3a}} ;
\end{tikzpicture} }} ;

%%% 4a %%%
\node (4a) at (4,0){
\scalebox{0.65}{
\begin{tikzpicture}
\draw [step=1,thin, gray] (-2,-2) grid (2,2);
\filldraw [ultra thick, fill=white] (1,0)--(0,1)--(-1,0)--(0,-1)--cycle ;
\filldraw  [ultra thick, fill=black] (1,0) circle [radius=0.1] ;
\filldraw  [ultra thick, fill=black] (0,1) circle [radius=0.1] ;
\filldraw  [ultra thick, fill=black] (-1,0) circle [radius=0.1] ;
\filldraw  [ultra thick, fill=black] (0,-1) circle [radius=0.1] ;

\node  at (0,0) {\LARGE{4a}} ;
\end{tikzpicture} }};

%%% 4b %%%
\node (4b) at (8,0){
\scalebox{0.65}{
\begin{tikzpicture}
\draw [step=1,thin, gray] (-2,-2) grid (2,2);
\filldraw [ultra thick, fill=white] (0,1)--(-1,0)--(0,-1)--(1,-1)--cycle ;
\filldraw  [ultra thick, fill=black] (0,1) circle [radius=0.1] ;
\filldraw  [ultra thick, fill=black] (-1,0) circle [radius=0.1] ;
\filldraw  [ultra thick, fill=black] (0,-1) circle [radius=0.1] ;
\filldraw  [ultra thick, fill=black] (1,-1) circle [radius=0.1] ;

\node  at (0,0) {\LARGE{4b}} ;
\end{tikzpicture} }}; 

%%% 4c %%%
\node (4c) at (12,0){
\scalebox{0.65}{
\begin{tikzpicture}
\draw [step=1,thin, gray] (-2,-2) grid (2,2);
\filldraw [ultra thick, fill=white] (0,1)--(-1,-1)--(1,-1)--cycle ;
\filldraw  [ultra thick, fill=black] (0,1) circle [radius=0.1] ;
\filldraw  [ultra thick, fill=black] (-1,-1) circle [radius=0.1] ;
\filldraw  [ultra thick, fill=black] (1,-1) circle [radius=0.1] ;
\filldraw  [ultra thick, fill=black] (0,-1) circle [radius=0.1] ;

\node  at (0,0) {\LARGE{4c}} ;
\end{tikzpicture} }};

%%% 5a %%%
\node (5a) at (0,-3){
\scalebox{0.65}{
\begin{tikzpicture}
\draw [step=1, thin, gray] (-2,-2) grid (2,2);
\filldraw [ultra thick, fill=white] (1,0)--(0,1)--(-1,1)--(-1,0)--(0,-1)--cycle ;
\filldraw  [ultra thick, fill=black] (1,0) circle [radius=0.1] ;
\filldraw  [ultra thick, fill=black] (0,1) circle [radius=0.1] ;
\filldraw  [ultra thick, fill=black] (-1,1) circle [radius=0.1] ;
\filldraw  [ultra thick, fill=black] (-1,0) circle [radius=0.1] ;
\filldraw  [ultra thick, fill=black] (0,-1) circle [radius=0.1] ;

\node  at (0,0) {\LARGE{5a}} ;
\end{tikzpicture} }}; 

%%% 5b %%%
\node (5b) at (4,-3){
\scalebox{0.65}{
\begin{tikzpicture}
\draw [step=1,thin, gray] (-2,-2) grid (2,2);
\filldraw [ultra thick, fill=white] (0,1)--(-1,0)--(-1,-1)--(1,-1)--cycle ;
\filldraw  [ultra thick, fill=black] (0,1) circle [radius=0.1] ;
\filldraw  [ultra thick, fill=black] (-1,0) circle [radius=0.1] ;
\filldraw  [ultra thick, fill=black] (-1,-1) circle [radius=0.1] ;
\filldraw  [ultra thick, fill=black] (0,-1) circle [radius=0.1] ;
\filldraw  [ultra thick, fill=black] (1,-1) circle [radius=0.1] ;

\node  at (0,0) {\LARGE{5b}} ;
\end{tikzpicture} }};

%%% 6a %%%
\node (6a) at (8,-3){
\scalebox{0.65}{
\begin{tikzpicture}
\draw [step=1,thin, gray] (-2,-2) grid (2,2);
\filldraw [ultra thick, fill=white] (1,0)--(0,1)--(-1,1)--(-1,0)--(0,-1)--(1,-1)--cycle ;
\filldraw  [ultra thick, fill=black] (1,0) circle [radius=0.1] ;
\filldraw  [ultra thick, fill=black] (0,1) circle [radius=0.1] ;
\filldraw  [ultra thick, fill=black] (-1,1) circle [radius=0.1] ;
\filldraw  [ultra thick, fill=black] (-1,0) circle [radius=0.1] ;
\filldraw  [ultra thick, fill=black] (0,-1) circle [radius=0.1] ;
\filldraw  [ultra thick, fill=black] (1,-1) circle [radius=0.1] ;

\node  at (0,0) {\LARGE{6a}} ;
\end{tikzpicture} }};

%%% 6b %%%
\node (6b) at (12,-3){
\scalebox{0.65}{
\begin{tikzpicture}
\draw [step=1,thin, gray] (-2,-2) grid (2,2);
\filldraw [ultra thick, fill=white] (1,0)--(0,1)--(-1,1)--(-1,-1)--(0,-1)--cycle ;
\filldraw  [ultra thick, fill=black] (1,0) circle [radius=0.1] ;
\filldraw  [ultra thick, fill=black] (0,1) circle [radius=0.1] ;
\filldraw  [ultra thick, fill=black] (-1,1) circle [radius=0.1] ;
\filldraw  [ultra thick, fill=black] (-1,0) circle [radius=0.1] ;
\filldraw  [ultra thick, fill=black] (-1,-1) circle [radius=0.1] ;
\filldraw  [ultra thick, fill=black] (0,-1) circle [radius=0.1] ;

\node  at (0,0) {\LARGE{6b}} ;
\end{tikzpicture} }};

%%% 6c %%%
\node (6c) at (0,-6){
\scalebox{0.65}{
\begin{tikzpicture}
\draw [step=1,thin, gray] (-2,-2) grid (2,2);
\filldraw [ultra thick, fill=white] (1,0)--(0,1)--(-2,-1)--(0,-1)--cycle ;
\filldraw  [ultra thick, fill=black] (1,0) circle [radius=0.1] ;
\filldraw  [ultra thick, fill=black] (0,1) circle [radius=0.1] ;
\filldraw  [ultra thick, fill=black] (-1,0) circle [radius=0.1] ;
\filldraw  [ultra thick, fill=black] (-2,-1) circle [radius=0.1] ;
\filldraw  [ultra thick, fill=black] (-1,-1) circle [radius=0.1] ;
\filldraw  [ultra thick, fill=black] (0,-1) circle [radius=0.1] ;

\node  at (0,0) {\LARGE{6c}} ;
\end{tikzpicture} }};

%%% 6d %%%
\node (6d) at (4,-6){
\scalebox{0.65}{
\begin{tikzpicture}
\draw [step=1,thin, gray] (-2,-2) grid (2,2);
\filldraw [ultra thick, fill=white] (0,1)--(-2,-1)--(1,-1)--cycle ;
\filldraw  [ultra thick, fill=black] (0,1) circle [radius=0.1] ;
\filldraw  [ultra thick, fill=black] (-1,0) circle [radius=0.1] ;
\filldraw  [ultra thick, fill=black] (-2,-1) circle [radius=0.1] ;
\filldraw  [ultra thick, fill=black] (-1,-1) circle [radius=0.1] ;
\filldraw  [ultra thick, fill=black] (0,-1) circle [radius=0.1] ;
\filldraw  [ultra thick, fill=black] (1,-1) circle [radius=0.1] ;

\node  at (0,0) {\LARGE{6d}} ;
\end{tikzpicture} }};

%%% 7a %%%
\node (7a) at (8,-6){
\scalebox{0.65}{
\begin{tikzpicture}
\draw [step=1,thin, gray] (-2,-2) grid (2,2);
\filldraw [ultra thick, fill=white] (1,0)--(0,1)--(-1,1)--(-1,-1)--(1,-1)--cycle ;
\filldraw  [ultra thick, fill=black] (1,0) circle [radius=0.1] ;
\filldraw  [ultra thick, fill=black] (0,1) circle [radius=0.1] ;
\filldraw  [ultra thick, fill=black] (-1,1) circle [radius=0.1] ;
\filldraw  [ultra thick, fill=black] (-1,0) circle [radius=0.1] ;
\filldraw  [ultra thick, fill=black] (-1,-1) circle [radius=0.1] ;
\filldraw  [ultra thick, fill=black] (0,-1) circle [radius=0.1] ;
\filldraw  [ultra thick, fill=black] (1,-1) circle [radius=0.1] ;

\node  at (0,0) {\LARGE{7a}} ;
\end{tikzpicture} }};

%%% 7b %%%
\node (7b) at (12,-6){
\scalebox{0.65}{
\begin{tikzpicture}
\draw [step=1,thin, gray] (-2,-2) grid (2,2);
\filldraw [ultra thick, fill=white] (1,0)--(0,1)--(-2,-1)--(1,-1)--cycle ;
\filldraw  [ultra thick, fill=black] (1,0) circle [radius=0.1] ;
\filldraw  [ultra thick, fill=black] (0,1) circle [radius=0.1] ;
\filldraw  [ultra thick, fill=black] (-1,0) circle [radius=0.1] ;
\filldraw  [ultra thick, fill=black] (-2,-1) circle [radius=0.1] ;
\filldraw  [ultra thick, fill=black] (-1,-1) circle [radius=0.1] ;
\filldraw  [ultra thick, fill=black] (0,-1) circle [radius=0.1] ;
\filldraw  [ultra thick, fill=black] (1,-1) circle [radius=0.1] ;

\node  at (0,0) {\LARGE{7b}} ;
\end{tikzpicture} }};

%%% 8a %%%
\node (8a) at (0,-9){
\scalebox{0.65}{
\begin{tikzpicture}
\draw [step=1,thin, gray] (-2,-2) grid (2,2);
\filldraw [ultra thick, fill=white] (1,1)--(-1,1)--(-1,-1)--(1,-1)--cycle ;
\filldraw  [ultra thick, fill=black] (1,0) circle [radius=0.1] ;
\filldraw  [ultra thick, fill=black] (1,1) circle [radius=0.1] ;
\filldraw  [ultra thick, fill=black] (0,1) circle [radius=0.1] ;
\filldraw  [ultra thick, fill=black] (-1,1) circle [radius=0.1] ;
\filldraw  [ultra thick, fill=black] (-1,0) circle [radius=0.1] ;
\filldraw  [ultra thick, fill=black] (-1,-1) circle [radius=0.1] ;
\filldraw  [ultra thick, fill=black] (0,-1) circle [radius=0.1] ;
\filldraw  [ultra thick, fill=black] (1,-1) circle [radius=0.1] ;

\node  at (0,0) {\LARGE{8a}} ;
\end{tikzpicture} }};

%%% 8b %%%
\node (8b) at (4,-9){
\scalebox{0.65}{
\begin{tikzpicture}
\draw [step=1,thin, gray] (-2,-2) grid (2,2);
\filldraw [ultra thick, fill=white] (0,1)--(-1,1)--(-1,-1)--(2,-1)--cycle ;
\filldraw  [ultra thick, fill=black] (1,0) circle [radius=0.1] ;
\filldraw  [ultra thick, fill=black] (0,1) circle [radius=0.1] ;
\filldraw  [ultra thick, fill=black] (-1,1) circle [radius=0.1] ;
\filldraw  [ultra thick, fill=black] (-1,0) circle [radius=0.1] ;
\filldraw  [ultra thick, fill=black] (-1,-1) circle [radius=0.1] ;
\filldraw  [ultra thick, fill=black] (0,-1) circle [radius=0.1] ;
\filldraw  [ultra thick, fill=black] (1,-1) circle [radius=0.1] ;
\filldraw  [ultra thick, fill=black] (2,-1) circle [radius=0.1] ;

\node  at (0,0) {\LARGE{8b}} ;
\end{tikzpicture} }};

%%% 8c %%%
\node (8c) at (8,-9){
\scalebox{0.65}{
\begin{tikzpicture}
\draw [step=1,thin, gray] (-2,-2) grid (2,2);
\filldraw [ultra thick, fill=white] (0,1)--(-2,-1)--(2,-1)--cycle ;
\filldraw  [ultra thick, fill=black] (1,0) circle [radius=0.1] ;
\filldraw  [ultra thick, fill=black] (0,1) circle [radius=0.1] ;
\filldraw  [ultra thick, fill=black] (-1,0) circle [radius=0.1] ;
\filldraw  [ultra thick, fill=black] (-2,-1) circle [radius=0.1] ;
\filldraw  [ultra thick, fill=black] (-1,-1) circle [radius=0.1] ;
\filldraw  [ultra thick, fill=black] (0,-1) circle [radius=0.1] ;
\filldraw  [ultra thick, fill=black] (1,-1) circle [radius=0.1] ;
\filldraw  [ultra thick, fill=black] (2,-1) circle [radius=0.1] ;

\node  at (0,0) {\LARGE{8c}} ;
\end{tikzpicture} }};

%%% 9a %%%
\node (9a) at (12,-9){
\scalebox{0.65}{
\begin{tikzpicture}
\draw [step=1,thin, gray] (-2,-2) grid (2,2);
\filldraw [ultra thick, fill=white] (-1,2)--(-1,-1)--(2,-1)--cycle ;
\filldraw  [ultra thick, fill=black] (1,0) circle [radius=0.1] ;
\filldraw  [ultra thick, fill=black] (0,1) circle [radius=0.1] ;
\filldraw  [ultra thick, fill=black] (-1,2) circle [radius=0.1] ;
\filldraw  [ultra thick, fill=black] (-1,1) circle [radius=0.1] ;
\filldraw  [ultra thick, fill=black] (-1,0) circle [radius=0.1] ;
\filldraw  [ultra thick, fill=black] (-1,-1) circle [radius=0.1] ;
\filldraw  [ultra thick, fill=black] (0,-1) circle [radius=0.1] ;
\filldraw  [ultra thick, fill=black] (1,-1) circle [radius=0.1] ;
\filldraw  [ultra thick, fill=black] (2,-1) circle [radius=0.1] ;

\node  at (0,0) {\LARGE{9a}} ;
\end{tikzpicture} }};

\end{tikzpicture} }
\end{center}
\caption{The classification of reflexive polygons}
\label{class_ref}
\end{figure}

In this section, we consider Question~\ref{que_mutation_MM} for the case of three dimensional Gorenstein toric singularities 
associated with reflexive polygons. 
We recall that $\Delta$ is called a reflexive polygon (or Fano polygon) if the origin is the unique interior point of $\Delta$. 
Reflexive polygons are classified in $16$ types up to integral unimodular transformations as in Figure~\ref{class_ref} 
(see e.g., \cite[Theorem~8.3.7]{CLS}, \cite[Appendix]{Boc2}). 
Consistent dimer models giving a reflexive polygon as the perfect matching polygon 
are well-studied in several papers (see e.g., \cite{Boc2,Boc4,HS}), 
and such dimer models are classified up to right equivalence of the associated QPs. 
Thus, we can obtain all splitting MM generators from those consistent dimer models. 
For further studies on rank one MCM modules appearing in subsections below, see \cite{Nak}. 

%%%%%%%%%%%%%%%%%%%%%%%%%%%%%%%%%%%%%%%%%%%%%%%%%%%%%%%
In what follows, we show the following theorem. 

\begin{theorem}
\label{main_reflexive}
Let $\Delta$ be a reflexive polygon, and $R$ be a three dimensional complete local Gorenstein toric singularity associated with $\Delta$. 
Then $\EG(\TMMG(R))$ is connected. 
\end{theorem}

The proof is a case-by-case check for all classified dimer models, and the strategy is similar for each type. 
Therefore, we mainly explain the case of type 4a (see subsection~\ref{type4a}). For other types, we only mention an outline. 
We remark that Theorem~\ref{mutation_Boc} is recovered from those arguments. 

\begin{remark}
\label{rem_auto}
Any splitting MM generator takes the form $e_i\calP(Q, W_Q)$ for some QP $(Q,W_Q)$ associated with a consistent dimer model and $i\in Q_0$. 
When we check each splitting MM generator, we have to take care of the following points. 
\begin{enumerate}[\rm (a)] 
\item Since there is an isomorphism  
\[
\calP(Q, W_Q)\cong\End_R(e_i\calP(Q,W_Q))\cong\End_R((e_j\calP(Q,W_Q))^*)\cong\calP(Q^{\rm op}, W_{Q^{\rm op}})
\]
for any $i,j\in Q_0$, the $R$-dual of a splitting MM generator also gives a splitting NCCR. 
Thus, we should also consider $(e_i\calP(Q,W_Q))^*$ if there is no vertex $j\in Q_0$ such that $(e_i\calP(Q,W_Q))^*\cong e_j\calP(Q,W_Q)$. 

\item Even if two Jacobian algebras are isomorphic as an $R$-algebra, they are not isomorphic as an $R$-module in general. 
More precisely, a different choice of generators of $\rmH_1(\sfT)$ induces an automorphism on $R$, 
and it sometimes gives another $R$-module structure. 
Thus, we need to take care about such automorphisms in discussions below. 
\end{enumerate}
\end{remark}

%%%%%%%%%%%%%%%%%%%%%%%%%%%%%%%%%%%%%%%%%%%%%%%%%%%%%%%
\subsection{Type 3a, 4c, 6d, 8c and 9a}
\label{ref_3a}
For these types, we have the assertion by Example~\ref{mutation_triangle}, because the associated cones are simplicial. 
For example, the following figure is the unique consistent dimer model associated with type 3a. 
For this case, the associated toric singularity is the quotient singularity by the cyclic group $G=\left<\mathrm{diag}(\omega, \omega, \omega)\right>$ 
where $\omega$ is a primitive cubic root of unity, and the associated quiver is the McKay quiver of $G$. 

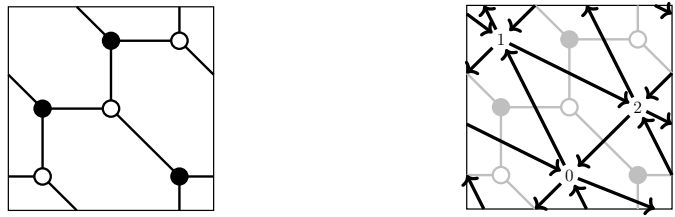
\medskip

%%%%%%%%%%%%%%%%%%%%%
%%% type 3a %%%%%%%%%%%%%
%%%%%%%%%%%%%%%%%%%%%

\begin{center}
\begin{tikzpicture}
\node (DM) at (0,0) 
{\scalebox{0.45}{
\begin{tikzpicture}
%vertex
\node (W1) at (1,1){$$}; \node (B1) at (5,1){$$}; \node (B2) at (1,3){$$}; \node (W2) at (3,3){$$};
\node (B3) at (3,5){$$}; \node (W3) at (5,5){$$};
\draw[line width=0.036cm]  (0,0) rectangle (6,6);

%edge
\draw[line width=0.07cm]  (0,1)--(W1)--(B2)--(W2)--(B3)--(W3)--(5,6);
\draw[line width=0.07cm]  (5,0)--(B1)--(6,1); \draw[line width=0.07cm]  (B1)--(W2); 
\draw[line width=0.07cm]  (W1)--(2,0);\draw[line width=0.07cm]  (B2)--(0,4);
\draw[line width=0.07cm]  (B3)--(2,6);\draw[line width=0.07cm]  (W3)--(6,4);
%black
\filldraw  [ultra thick, fill=black] (5,1) circle [radius=0.24] ;\filldraw  [ultra thick, fill=black] (1,3) circle [radius=0.24] ;
\filldraw  [ultra thick, fill=black] (3,5) circle [radius=0.24] ;
%white
\draw  [line width=1.92pt,fill=white] (1,1) circle [radius=0.24] ;\draw  [line width=1.92pt, fill=white] (3,3) circle [radius=0.24] ;
\draw  [line width=1.92pt, fill=white] (5,5) circle [radius=0.24] ;
\end{tikzpicture}
} }; 

\node (QV) at (6,0) 
{\scalebox{0.45}{
\begin{tikzpicture}
%vertex
\node (W1) at (1,1){$$}; \node (B1) at (5,1){$$}; \node (B2) at (1,3){$$}; \node (W2) at (3,3){$$}; \node (B3) at (3,5){$$}; \node (W3) at (5,5){$$};
\node (Q0) at (3,1){{\LARGE$0$}}; \node (Q1) at (1,5){{\LARGE$1$}}; \node (Q2) at (5,3){{\LARGE$2$}}; 
\draw[line width=0.036cm]  (0,0) rectangle (6,6);
%edge
\draw[lightgray,line width=0.07cm]  (0,1)--(W1)--(B2)--(W2)--(B3)--(W3)--(5,6);
\draw[lightgray,line width=0.07cm]  (5,0)--(B1)--(6,1); \draw[lightgray,line width=0.07cm]  (B1)--(W2); 
\draw[lightgray,line width=0.07cm]  (W1)--(2,0);\draw[lightgray,line width=0.07cm]  (B2)--(0,4);
\draw[lightgray,line width=0.07cm]  (B3)--(2,6);\draw[lightgray,line width=0.07cm]  (W3)--(6,4);
%black
\filldraw  [ultra thick, draw=lightgray, fill=lightgray] (5,1) circle [radius=0.24] ;\filldraw  [ultra thick, draw=lightgray, fill=lightgray] (1,3) circle [radius=0.24] ;
\filldraw  [ultra thick, draw=lightgray, fill=lightgray] (3,5) circle [radius=0.24] ;
%white
\draw  [line width=1.92pt, draw=lightgray,fill=white] (1,1) circle [radius=0.24] ;\draw  [line width=1.92pt, draw=lightgray,fill=white] (3,3) circle [radius=0.24] ;
\draw  [line width=1.92pt, draw=lightgray,fill=white] (5,5) circle [radius=0.24] ;

\draw[->, line width=0.096cm] (Q0)--(Q1); \draw[->, line width=0.096cm] (Q1)--(Q2); \draw[->, line width=0.096cm] (Q2)--(Q0);
\draw[->, line width=0.096cm] (Q0)--(2,0); \draw[->, line width=0.096cm] (2,6)--(Q1); \draw[->, line width=0.096cm] (Q2)--(3.5,6);
\draw[->, line width=0.096cm] (3.5,0)--(Q0); \draw[->, line width=0.096cm] (Q2)--(6,2.5); \draw[->, line width=0.096cm] (0,2.5)--(Q0);
\draw[->, line width=0.096cm] (Q1)--(0,4); \draw[->, line width=0.096cm] (6,4)--(Q2);
\draw[->, line width=0.096cm] (Q0)--(5.5,0); \draw[->, line width=0.096cm] (5.5,6)--(6,5.75);  \draw[->, line width=0.096cm] (0,5.75)--(Q1);
\draw[<-, line width=0.096cm] (Q2)--(6,1); \draw[<-, line width=0.096cm] (0,1)--(0.5,0); \draw[<-, line width=0.096cm] (0.5,6)--(Q1);
\end{tikzpicture}
} }; 
\end{tikzpicture}
\end{center}

\medskip

%%%%%%%%%%%%%%%%%%%%%%%%%%%%%%%%%%%%%%%%%%%%%%%%%%%%%%%
\subsection{Type 4a}
\label{type4a}

In this subsection, we consider the reflexive polygon of type 4a. 
Thus, let $R$ be the three dimensional complete local Gorenstein toric singularity defined by the cone $\sigma$: 
\[
\sigma=\mathrm{Cone}\{v_1=(1,0,1), v_2=(0,1,1), v_3=(-1,0,1), v_4=(0,-1,1) \}. 
\]
Every divisorial ideal (i.e., rank one reflexive $R$-module) takes the form $T(u_1,u_2,u_3,u_4)$, and it corresponds to  
the element $u_1[I_1]+\cdots+u_4[I_4]$ in $\Cl(R)$ (see subsection~\ref{toric_pre}). 
As an element in $\Cl(R)$, we have that $[I_1]=[I_3]$, $[I_2]=[I_4]$, $2[I_1]+2[I_2]=0$ (see Lemma~\ref{div_eq}). 
Therefore, we have that $\Cl(R)\cong\ZZ\times\ZZ/2\ZZ$, and each divisorial ideal is represented by $T(a,b,0,0)$ where $a\in\ZZ, b\in\ZZ/2\ZZ$. 
Also, there are two consistent dimer models (up to right equivalence) written below which give the reflexive polygon of type 4a as the perfect matching polygon. 

\subsubsection{Type 4a-1}

First, we consider the following consistent dimer model and the associated QP $(Q, W_Q)$. 
This is just $\Gamma_{4a}$ used in Section~\ref{sec2} and Example~\ref{mutation_QP1}. 
For simplicity, we denote by $\sfA$ the complete Jacobian algebra $\calP(Q, W_Q)$. 
Then, we have that $\sfA\cong\End_R(\bigoplus_{j\in Q_0}T_{ij})$ and $e_i\sfA\cong\bigoplus_{j\in Q_0}T_{ij}$ for any $i\in Q_0$ by Theorem~\ref{NCCR1}. 
Here, $e_i$ is a primitive idempotent (i.e., it is a trivial path satisfying $h(e_i)=t(e_i)=i\in Q_0$). 
In what follows, we give the concrete description of these splitting MM generators $\bigoplus_{j\in Q_0}T_{ij}$. 

\medskip

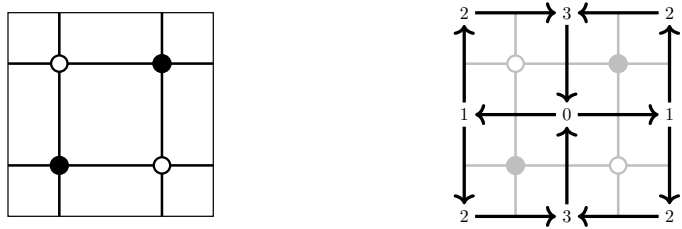
\begin{center}
\begin{tikzpicture}
\node (DM) at (0,0) 
{\scalebox{0.675}{
\begin{tikzpicture}
%vertex
\node (P1) at (1,1){$$}; \node (P2) at (3,1){$$}; \node (P3) at
(3,3){$$}; \node (P4) at (1,3){$$};
\draw[line width=0.024cm]  (0,0) rectangle (4,4);

%edge
\draw[line width=0.05cm]  (P1)--(P2)--(P3)--(P4)--(P1);\draw[line width=0.05cm] (0,1)--(P1)--(1,0); \draw[line width=0.05cm]  (4,1)--(P2)--(3,0);
\draw[line width=0.05cm]  (0,3)--(P4)--(1,4);\draw[line width=0.05cm]  (3,4)--(P3)--(4,3);
%black
\filldraw  [ultra thick, fill=black] (1,1) circle [radius=0.16] ;\filldraw  [ultra thick, fill=black] (3,3) circle [radius=0.16] ;
%white
\draw  [line width=1.28pt,fill=white] (3,1) circle [radius=0.16] ;\draw  [line width=1.28pt, fill=white] (1,3) circle [radius=0.16] ;
\end{tikzpicture}
} }; 

\node (QV) at (6,0) 
{\scalebox{0.675}{
\begin{tikzpicture}
%vertex
\node (Q1) at (2,2){$0$};\node (Q2a) at (0,2){$1$}; \node(Q2b) at (4,2){$1$};\node (Q3a) at (0,0){$2$};
\node(Q3c) at (4,4){$2$};\node(Q3b) at (4,0){$2$};\node(Q3d) at (0,4){$2$};\node (Q4a) at (2,0){$3$};
\node (Q4b) at (2,4){$3$};

%edge
\draw[lightgray, line width=0.05cm]  (P1)--(P2)--(P3)--(P4)--(P1);\draw[lightgray, line width=0.05cm] (0,1)--(P1)--(1,0); 
\draw[lightgray, line width=0.05cm]  (4,1)--(P2)--(3,0);\draw[lightgray, line width=0.05cm]  (0,3)--(P4)--(1,4);\draw[lightgray, line width=0.05cm]  (3,4)--(P3)--(4,3);
%black
\filldraw  [ultra thick, draw=lightgray, fill=lightgray] (1,1) circle [radius=0.16] ;\filldraw  [ultra thick, draw=lightgray, fill=lightgray] (3,3) circle [radius=0.16] ;
%white
\draw  [line width=1.28pt, draw=lightgray,fill=white] (3,1) circle [radius=0.16] ;\draw  [line width=1.28pt, draw=lightgray,fill=white] (1,3) circle [radius=0.16] ;

\draw[->, line width=0.064cm] (Q1)--(Q2a);\draw[->, line width=0.064cm] (Q2a)--(Q3a);\draw[->, line width=0.064cm] (Q3a)--(Q4a);
\draw[->, line width=0.064cm] (Q4a)--(Q1);\draw[->, line width=0.064cm] (Q2a)--(Q3d);\draw[->, line width=0.064cm] (Q3d)--(Q4b);
\draw[->, line width=0.064cm] (Q4b)--(Q1);\draw[->, line width=0.064cm] (Q1)--(Q2b);\draw[->, line width=0.064cm] (Q2b)--(Q3b);
\draw[->, line width=0.064cm] (Q3b)--(Q4a);\draw[->, line width=0.064cm] (Q2b)--(Q3c);\draw[->, line width=0.064cm] (Q3c)--(Q4b);
\end{tikzpicture}
} }; 
\end{tikzpicture}
\end{center}

\medskip

Since this dimer model is consistent, there are extremal perfect matchings corresponding to $v_1, \cdots, v_4$ (see Proposition~\ref{ex_pm}), 
and the following perfect matchings $\sfP_1, \cdots, \sfP_4$ are desired ones. 
Namely, we fix the perfect matching $\sfP_0$, then we have that $v_i=(h(\sfP_i, \sfP_0), 1)\in\ZZ^3$. 
(Note that a choice of perfect matchings is not unique, but we may only check one of them by arguments in subsection~\ref{sec_consist}.) 

\medskip

%%%%%perfect_matching%%%%%
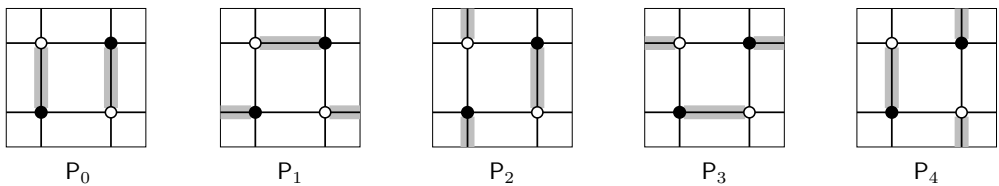
\begin{center}
{\scalebox{0.9}{
\begin{tikzpicture} 
\node at (0,-1.4) {$\sfP_0$};\node at (3.1,-1.4) {$\sfP_1$}; \node at (6.2,-1.4) {$\sfP_2$}; 
\node at (9.3,-1.4) {$\sfP_3$}; \node at (12.4,-1.4) {$\sfP_4$};

\node (PM0) at (0,0) 
{\scalebox{0.51}{\begin{tikzpicture}
%vertex
\node (P1) at (1,1){$$}; \node (P2) at (3,1){$$}; \node (P3) at (3,3){$$}; \node (P4) at (1,3){$$};
\draw[line width=0.024cm]  (0,0) rectangle (4,4);
%perfect matching
\draw[line width=0.4cm,color=lightgray] (P1)--(P4);\draw[line width=0.4cm,color=lightgray] (P2)--(P3);
%edge
\draw[line width=0.05cm]  (P1)--(P2)--(P3)--(P4)--(P1);\draw[line width=0.05cm] (0,1)--(P1)--(1,0); \draw[line width=0.05cm]  (4,1)--(P2)--(3,0);
\draw[line width=0.05cm]  (0,3)--(P4)--(1,4);\draw[line width=0.05cm]  (3,4)--(P3)--(4,3);
%black
\filldraw  [ultra thick, fill=black] (1,1) circle [radius=0.16] ;\filldraw  [ultra thick, fill=black] (3,3) circle [radius=0.16] ;
%white
\draw  [line width=1.28pt,fill=white] (3,1) circle [radius=0.16] ;\draw  [line width=1.28pt, fill=white] (1,3) circle [radius=0.16] ;
\end{tikzpicture} }}; 

\node (PM1) at (3.1,0) 
{\scalebox{0.51}{
\begin{tikzpicture}
%vertex
\node (P1) at (1,1){$$}; \node (P2) at (3,1){$$}; \node (P3) at (3,3){$$}; \node (P4) at (1,3){$$};
\draw[line width=0.024cm]  (0,0) rectangle (4,4);

%perfect matching
\draw[line width=0.4cm,color=lightgray] (P3)--(P4);\draw[line width=0.4cm,color=lightgray] (P1)--(0,1);\draw[line width=0.4cm,color=lightgray] (P2)--(4,1);

%edge
\draw[line width=0.05cm]  (P1)--(P2)--(P3)--(P4)--(P1);\draw[line width=0.05cm] (0,1)--(P1)--(1,0); \draw[line width=0.05cm]  (4,1)--(P2)--(3,0);
\draw[line width=0.05cm]  (0,3)--(P4)--(1,4);\draw[line width=0.05cm]  (3,4)--(P3)--(4,3);
%black
\filldraw  [ultra thick, fill=black] (1,1) circle [radius=0.16] ;\filldraw  [ultra thick, fill=black] (3,3) circle [radius=0.16] ;
%white
\draw  [line width=1.28pt,fill=white] (3,1) circle [radius=0.16] ;\draw  [line width=1.28pt, fill=white] (1,3) circle [radius=0.16] ;
\end{tikzpicture} }}; 

\node (PM2) at (6.2,0) 
{\scalebox{0.51}{
\begin{tikzpicture}
%vertex
\node (P1) at (1,1){$$}; \node (P2) at (3,1){$$}; \node (P3) at (3,3){$$}; \node (P4) at (1,3){$$};
\draw[line width=0.024cm]  (0,0) rectangle (4,4);

%perfect matching
\draw[line width=0.4cm,color=lightgray] (P3)--(P2);\draw[line width=0.4cm,color=lightgray] (P4)--(1,4);\draw[line width=0.4cm,color=lightgray] (P1)--(1,0);

%edge
\draw[line width=0.05cm]  (P1)--(P2)--(P3)--(P4)--(P1);\draw[line width=0.05cm] (0,1)--(P1)--(1,0); \draw[line width=0.05cm]  (4,1)--(P2)--(3,0);
\draw[line width=0.05cm]  (0,3)--(P4)--(1,4);\draw[line width=0.05cm]  (3,4)--(P3)--(4,3);

%black
\filldraw  [ultra thick, fill=black] (1,1) circle [radius=0.16] ;\filldraw  [ultra thick, fill=black] (3,3) circle [radius=0.16] ;
%white
\draw  [line width=1.28pt,fill=white] (3,1) circle [radius=0.16] ;\draw  [line width=1.28pt, fill=white] (1,3) circle [radius=0.16] ;
\end{tikzpicture} }} ;  

\node (PM3) at (9.3,0) 
{\scalebox{0.51}{
\begin{tikzpicture}
%vertex
\node (P1) at (1,1){$$}; \node (P2) at (3,1){$$}; \node (P3) at (3,3){$$}; \node (P4) at (1,3){$$};
\draw[line width=0.024cm]  (0,0) rectangle (4,4);

%perfect matching
\draw[line width=0.4cm,color=lightgray] (P2)--(P1);\draw[line width=0.4cm,color=lightgray] (P4)--(0,3);\draw[line width=0.4cm,color=lightgray] (P3)--(4,3);

%edge
\draw[line width=0.05cm]  (P1)--(P2)--(P3)--(P4)--(P1);\draw[line width=0.05cm] (0,1)--(P1)--(1,0); \draw[line width=0.05cm]  (4,1)--(P2)--(3,0);
\draw[line width=0.05cm]  (0,3)--(P4)--(1,4);\draw[line width=0.05cm]  (3,4)--(P3)--(4,3);
%black
\filldraw  [ultra thick, fill=black] (1,1) circle [radius=0.16] ;\filldraw  [ultra thick, fill=black] (3,3) circle [radius=0.16] ;
%white
\draw  [line width=1.28pt,fill=white] (3,1) circle [radius=0.16] ;\draw  [line width=1.28pt, fill=white] (1,3) circle [radius=0.16] ;
\end{tikzpicture} }}; 

\node (PM4) at (12.4,0) 
{\scalebox{0.51}{
\begin{tikzpicture}
%vertex
\node (P1) at (1,1){$$}; \node (P2) at (3,1){$$}; \node (P3) at (3,3){$$}; \node (P4) at (1,3){$$};
\draw[line width=0.024cm]  (0,0) rectangle (4,4);

%perfect matching
\draw[line width=0.4cm,color=lightgray] (P4)--(P1);\draw[line width=0.4cm,color=lightgray] (P3)--(3,4);\draw[line width=0.4cm,color=lightgray] (P2)--(3,0);
%edge
\draw[line width=0.05cm]  (P1)--(P2)--(P3)--(P4)--(P1);\draw[line width=0.05cm] (0,1)--(P1)--(1,0); \draw[line width=0.05cm]  (4,1)--(P2)--(3,0);
\draw[line width=0.05cm]  (0,3)--(P4)--(1,4);\draw[line width=0.05cm]  (3,4)--(P3)--(4,3);
%black
\filldraw  [ultra thick, fill=black] (1,1) circle [radius=0.16] ;\filldraw  [ultra thick, fill=black] (3,3) circle [radius=0.16] ;
%white
\draw  [line width=1.28pt,fill=white] (3,1) circle [radius=0.16] ;\draw  [line width=1.28pt, fill=white] (1,3) circle [radius=0.16] ;
\end{tikzpicture} }} ;

\end{tikzpicture}
}}
\end{center}

%\medskip

Then, for any $i, j\in Q_0$, we consider the divisorial ideal $T_{ij}= T(\sfP_1(a_{ij}), \cdots, \sfP_4(a_{ij}))$ as in subsection~\ref{sec_consist}. 
Note that $T_{ij}$ depends only on the starting point $i$ and the ending point $j$. 
For example, taking a path from $0\in Q_0$ to each vertex: 
\[0\xrightarrow{b_2} 1,\quad 0\xrightarrow{b_2d_2} 2,\quad 0\xrightarrow{b_2d_2c_2} 3 \quad\text{(see the notation in Example~\ref{mutation_QP1})}, \]
we have the following divisorial ideals and the direct sum of these ideals is just $e_0\sfA$. 

\begin{itemize}
\item [] $T_{00}\cong T(\sfP_1(e_0),\cdots,\sfP_4(e_0))= T(0,0,0,0)\cong R$, 
\item [] $T_{01}\cong T(\sfP_1(b_2),\cdots,\sfP_4(b_2))= T(0,1,0,0)$, 
\item [] $T_{02}\cong T(\sfP_1(b_2d_2),\cdots,\sfP_4(b_2d_2))= T(0,1,1,0)\cong T(1,1,0,0)$, 
\item [] $T_{03}\cong T(\sfP_1(b_2d_2c_2),\cdots,\sfP_4(b_2d_2c_2))= T(0,1,1,1)\cong T(-1,0,0,0)$.
\end{itemize}

%Note that the last isomorphism of each ideal follows from Lemma~\ref{div_eq}. 
Similarly, we have the following table, and easy to see that $e_0\sfA\cong e_2\sfA$, $e_1\sfA\cong e_3\sfA$, $(e_0\sfA)^*\cong e_1\sfA$. 

\medskip

\begin{center}
{\small{
\begin{tabular}{|l|l|l|l|} \hline 
$e_0\sfA$&$e_1\sfA$&$e_2\sfA$&$e_3\sfA$ \\ \hline
$T_{00}\cong R$&$T_{10}\cong T(2,1,0,0)$&$T_{20}\cong T(1,1,0,0)$&$T_{30}\cong T(1,0,0,0)$ \\
$T_{01}\cong T(0,1,0,0)$&$T_{11}\cong R$&$T_{21}\cong T(-1,0,0,0)$&$T_{31}\cong T(1,1,0,0)$ \\
$T_{02}\cong T(1,1,0,0)$&$T_{12}\cong T(1,0,0,0)$&$T_{22}\cong R$&$T_{32}\cong T(2,1,0,0)$ \\
$T_{03}\cong T(-1,0,0,0)$&$T_{13}\cong T(1,1,0,0)$&$T_{23}\cong T(0,1,0,0)$&$T_{33}\cong R$\\ \hline 
\end{tabular}
}}
\end{center}

\subsubsection{Type 4a-2}
Next, we consider the following consistent dimer model and the associated QP $(Q, W_Q)$. 
As we showed in Example~\ref{mutation_QP1}, we can obtain this dimer model 
by mutating $\Gamma_{4a}$ at $0\in Q_0$. We denote the complete Jacobian algebra $\calP(Q, W_Q)$ by $\sfB$. 
In the same way as the above case, we give the concrete description of splitting MM generators $e_i\sfB\cong\bigoplus_{j\in Q_0}T_{ij}$. 

\medskip

%%%%%%%%%%%%%%%%%%%%%
%%% type 4a-2 %%%%%%%%%%%%%
%%%%%%%%%%%%%%%%%%%%%

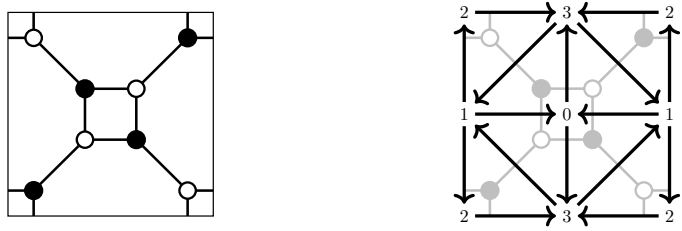
\begin{center}
\begin{tikzpicture}
\node (DM) at (0,0) 
{\scalebox{0.675}{
\begin{tikzpicture}
%vertex
\node (B1) at (0.5,0.5){$$}; \node (W1) at (3.5,0.5){$$}; \node (B2) at (2.5,1.5){$$}; 
\node (W2) at (1.5,1.5){$$}; \node (B3) at (1.5,2.5){$$}; \node (W3) at (2.5,2.5){$$}; \node (B4) at (3.5,3.5){$$}; \node (W4) at (0.5,3.5){$$}; 
\draw[line width=0.024cm]  (0,0) rectangle (4,4);
%edge
\draw[line width=0.05cm]  (W2)--(B2)--(W3)--(B3)--(W2);\draw[line width=0.05cm] (B1)--(W2); \draw[line width=0.05cm] (B2)--(W1); 
\draw[line width=0.05cm] (B3)--(W4); \draw[line width=0.05cm] (B4)--(W3);  
\draw[line width=0.05cm] (B1)--(0.5,0); \draw[line width=0.05cm] (B1)--(0,0.5); \draw[line width=0.05cm] (W1)--(3.5,0);   
\draw[line width=0.05cm] (W1)--(4,0.5); \draw[line width=0.05cm] (W4)--(0,3.5); \draw[line width=0.05cm] (W4)--(0.5,4);  
\draw[line width=0.05cm] (B4)--(4,3.5); \draw[line width=0.05cm] (B4)--(3.5,4);   
%black
\filldraw  [ultra thick, fill=black] (0.5,0.5) circle [radius=0.16] ;\filldraw  [ultra thick, fill=black] (2.5,1.5) circle [radius=0.16] ;
\filldraw  [ultra thick, fill=black] (1.5,2.5) circle [radius=0.16] ;\filldraw  [ultra thick, fill=black] (3.5,3.5) circle [radius=0.16] ;
%white
\draw  [line width=1.28pt,fill=white] (3.5,0.5) circle [radius=0.16] ;\draw  [line width=1.28pt, fill=white] (1.5,1.5) circle [radius=0.16] ;
\draw  [line width=1.28pt, fill=white] (2.5,2.5) circle [radius=0.16] ;\draw  [line width=1.28pt, fill=white] (0.5,3.5) circle [radius=0.16] ;
\end{tikzpicture}
} }; 

\node (QV) at (6,0) 
{\scalebox{0.675}{
\begin{tikzpicture}
%vertex
\node (Q0) at (2,2){$0$};\node (Q1a) at (0,2){$1$}; \node(Q1b) at (4,2){$1$};\node (Q2a) at (0,0){$2$};\node(Q2c) at (4,4){$2$};
\node(Q2b) at (4,0){$2$};\node(Q2d) at (0,4){$2$};\node (Q3a) at (2,0){$3$};\node (Q3b) at (2,4){$3$};

%edge
\draw[lightgray,line width=0.05cm]  (W2)--(B2)--(W3)--(B3)--(W2);
\draw[lightgray,line width=0.05cm] (B1)--(W2); \draw[lightgray,line width=0.05cm] (B2)--(W1); 
\draw[lightgray,line width=0.05cm] (B3)--(W4); \draw[lightgray,line width=0.05cm] (B4)--(W3);  
\draw[lightgray,line width=0.05cm] (B1)--(0.5,0); \draw[lightgray,line width=0.05cm] (B1)--(0,0.5); \draw[lightgray,line width=0.05cm] (W1)--(3.5,0);   
\draw[lightgray,line width=0.05cm] (W1)--(4,0.5); \draw[lightgray,line width=0.05cm] (W4)--(0,3.5); \draw[lightgray,line width=0.05cm] (W4)--(0.5,4);  
\draw[lightgray,line width=0.05cm] (B4)--(4,3.5); \draw[lightgray,line width=0.05cm] (B4)--(3.5,4);   

%black
\filldraw  [ultra thick, draw=lightgray, fill=lightgray] (0.5,0.5) circle [radius=0.16] ;\filldraw  [ultra thick, draw=lightgray, fill=lightgray] (2.5,1.5) circle [radius=0.16] ;
\filldraw  [ultra thick, draw=lightgray, fill=lightgray] (1.5,2.5) circle [radius=0.16] ;\filldraw  [ultra thick, draw=lightgray, fill=lightgray] (3.5,3.5) circle [radius=0.16] ;

%white
\draw  [line width=1.28pt,draw=lightgray,fill=white] (3.5,0.5) circle [radius=0.16] ;\draw  [line width=1.28pt, draw=lightgray,fill=white] (1.5,1.5) circle [radius=0.16] ;
\draw  [line width=1.28pt, draw=lightgray,fill=white] (2.5,2.5) circle [radius=0.16] ;\draw  [line width=1.28pt, draw=lightgray,fill=white] (0.5,3.5) circle [radius=0.16] ;

\draw[->, line width=0.064cm] (Q0)--(Q3a); \draw[->, line width=0.064cm] (Q0)--(Q3b); \draw[->, line width=0.064cm] (Q1a)--(Q0); 
\draw[->, line width=0.064cm] (Q1b)--(Q0); \draw[->, line width=0.064cm] (Q1a)--(Q2a); \draw[->, line width=0.064cm] (Q1a)--(Q2d); 
\draw[->, line width=0.064cm] (Q1b)--(Q2b); \draw[->, line width=0.064cm] (Q1b)--(Q2c);  
\draw[->, line width=0.064cm] (Q2a)--(Q3a); \draw[->, line width=0.064cm] (Q2b)--(Q3a); \draw[->, line width=0.064cm] (Q2c)--(Q3b);
\draw[->, line width=0.064cm] (Q2d)--(Q3b); \draw[->, line width=0.064cm] (Q3a)--(Q1a); \draw[->, line width=0.064cm] (Q3a)--(Q1b); 
\draw[->, line width=0.064cm] (Q3b)--(Q1a);  \draw[->, line width=0.064cm] (Q3b)--(Q1b); 
\end{tikzpicture}
} }; 
\end{tikzpicture}
\end{center}

\medskip

We fix the perfect matching $\sfP_0$, 
then the following $\sfP_1, \cdots, \sfP_4$ are extremal perfect matchings corresponding to $v_1, \cdots, v_4$. 
Namely, $v_i=(h(\sfP_i, \sfP_0), 1)\in\ZZ^3$. 

\medskip

%%%%%perfect_matching%%%%%
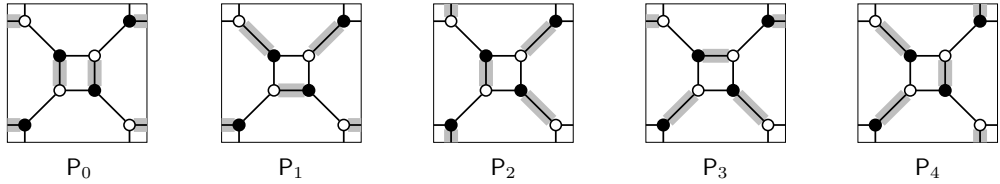
\begin{center}
{\scalebox{0.9}{
\begin{tikzpicture} 
\node at (0,-1.4) {$\sfP_0$};\node at (3.1,-1.4) {$\sfP_1$}; \node at (6.2,-1.4) {$\sfP_2$}; 
\node at (9.3,-1.4) {$\sfP_3$}; \node at (12.4,-1.4) {$\sfP_4$};

\node (PM0) at (0,0) 
{\scalebox{0.51}{\begin{tikzpicture}
%vertex
\node (B1) at (0.5,0.5){$$}; \node (W1) at (3.5,0.5){$$}; \node (B2) at (2.5,1.5){$$}; \node (W2) at (1.5,1.5){$$}; \node (B3) at (1.5,2.5){$$}; \node (W3) at (2.5,2.5){$$}; 
\node (B4) at (3.5,3.5){$$}; \node (W4) at (0.5,3.5){$$}; 
\draw[line width=0.024cm]  (0,0) rectangle (4,4);

%perfect matching
\draw[line width=0.4cm,color=lightgray] (W2)--(B3); \draw[line width=0.4cm,color=lightgray] (W3)--(B2);
\draw[line width=0.4cm,color=lightgray] (B1)--(0,0.5); \draw[line width=0.4cm,color=lightgray] (W1)--(4,0.5); 
\draw[line width=0.4cm,color=lightgray] (W4)--(0,3.5); \draw[line width=0.4cm,color=lightgray] (B4)--(4,3.5);
%edge
\draw[line width=0.05cm]  (W2)--(B2)--(W3)--(B3)--(W2);\draw[line width=0.05cm] (B1)--(W2); \draw[line width=0.05cm] (B2)--(W1); 
\draw[line width=0.05cm] (B3)--(W4); \draw[line width=0.05cm] (B4)--(W3);  
\draw[line width=0.05cm] (B1)--(0.5,0); \draw[line width=0.05cm] (B1)--(0,0.5); \draw[line width=0.05cm] (W1)--(3.5,0);   
\draw[line width=0.05cm] (W1)--(4,0.5); \draw[line width=0.05cm] (W4)--(0,3.5); \draw[line width=0.05cm] (W4)--(0.5,4);  
\draw[line width=0.05cm] (B4)--(4,3.5); \draw[line width=0.05cm] (B4)--(3.5,4);   
%black
\filldraw  [ultra thick, fill=black] (0.5,0.5) circle [radius=0.16] ;\filldraw  [ultra thick, fill=black] (2.5,1.5) circle [radius=0.16] ;
\filldraw  [ultra thick, fill=black] (1.5,2.5) circle [radius=0.16] ;\filldraw  [ultra thick, fill=black] (3.5,3.5) circle [radius=0.16] ;
%white
\draw  [line width=1.28pt,fill=white] (3.5,0.5) circle [radius=0.16] ;\draw  [line width=1.28pt, fill=white] (1.5,1.5) circle [radius=0.16] ;
\draw  [line width=1.28pt, fill=white] (2.5,2.5) circle [radius=0.16] ;\draw  [line width=1.28pt, fill=white] (0.5,3.5) circle [radius=0.16] ;

\end{tikzpicture} }}; 

\node (PM1) at (3.1,0) 
{\scalebox{0.51}{
\begin{tikzpicture}
%vertex
\node (B1) at (0.5,0.5){$$}; \node (W1) at (3.5,0.5){$$}; \node (B2) at (2.5,1.5){$$}; 
\node (W2) at (1.5,1.5){$$}; \node (B3) at (1.5,2.5){$$}; \node (W3) at (2.5,2.5){$$}; \node (B4) at (3.5,3.5){$$}; \node (W4) at (0.5,3.5){$$}; 
\draw[line width=0.024cm]  (0,0) rectangle (4,4);

%perfect matching
\draw[line width=0.4cm,color=lightgray] (W2)--(B2); \draw[line width=0.4cm,color=lightgray] (B1)--(0,0.5); \draw[line width=0.4cm,color=lightgray] (W1)--(4,0.5); 
\draw[line width=0.4cm,color=lightgray] (W4)--(B3); \draw[line width=0.4cm,color=lightgray] (W3)--(B4); 

%edge
\draw[line width=0.05cm]  (W2)--(B2)--(W3)--(B3)--(W2);
\draw[line width=0.05cm] (B1)--(W2); \draw[line width=0.05cm] (B2)--(W1); 
\draw[line width=0.05cm] (B3)--(W4); \draw[line width=0.05cm] (B4)--(W3);  
\draw[line width=0.05cm] (B1)--(0.5,0); \draw[line width=0.05cm] (B1)--(0,0.5); \draw[line width=0.05cm] (W1)--(3.5,0);   
\draw[line width=0.05cm] (W1)--(4,0.5); \draw[line width=0.05cm] (W4)--(0,3.5); \draw[line width=0.05cm] (W4)--(0.5,4);  
\draw[line width=0.05cm] (B4)--(4,3.5); \draw[line width=0.05cm] (B4)--(3.5,4);   
%black
\filldraw  [ultra thick, fill=black] (0.5,0.5) circle [radius=0.16] ;\filldraw  [ultra thick, fill=black] (2.5,1.5) circle [radius=0.16] ;
\filldraw  [ultra thick, fill=black] (1.5,2.5) circle [radius=0.16] ;\filldraw  [ultra thick, fill=black] (3.5,3.5) circle [radius=0.16] ;
%white
\draw  [line width=1.28pt,fill=white] (3.5,0.5) circle [radius=0.16] ;\draw  [line width=1.28pt, fill=white] (1.5,1.5) circle [radius=0.16] ;
\draw  [line width=1.28pt, fill=white] (2.5,2.5) circle [radius=0.16] ;\draw  [line width=1.28pt, fill=white] (0.5,3.5) circle [radius=0.16] ;

\end{tikzpicture} }}; 

\node (PM2) at (6.2,0) 
{\scalebox{0.51}{
\begin{tikzpicture}
%vertex
\node (B1) at (0.5,0.5){$$}; \node (W1) at (3.5,0.5){$$}; \node (B2) at (2.5,1.5){$$}; 
\node (W2) at (1.5,1.5){$$}; \node (B3) at (1.5,2.5){$$}; \node (W3) at (2.5,2.5){$$}; \node (B4) at (3.5,3.5){$$}; \node (W4) at (0.5,3.5){$$}; 
\draw[line width=0.024cm]  (0,0) rectangle (4,4);

%perfect matching
\draw[line width=0.4cm,color=lightgray] (W3)--(B4); \draw[line width=0.4cm,color=lightgray] (W1)--(B2); 
\draw[line width=0.4cm,color=lightgray] (W2)--(B3); \draw[line width=0.4cm,color=lightgray] (B1)--(0.5,0); \draw[line width=0.4cm,color=lightgray] (W4)--(0.5,4);  

%edge
\draw[line width=0.05cm]  (W2)--(B2)--(W3)--(B3)--(W2);
\draw[line width=0.05cm] (B1)--(W2); \draw[line width=0.05cm] (B2)--(W1); 
\draw[line width=0.05cm] (B3)--(W4); \draw[line width=0.05cm] (B4)--(W3);  
\draw[line width=0.05cm] (B1)--(0.5,0); \draw[line width=0.05cm] (B1)--(0,0.5); \draw[line width=0.05cm] (W1)--(3.5,0);   
\draw[line width=0.05cm] (W1)--(4,0.5); \draw[line width=0.05cm] (W4)--(0,3.5); \draw[line width=0.05cm] (W4)--(0.5,4);  
\draw[line width=0.05cm] (B4)--(4,3.5); \draw[line width=0.05cm] (B4)--(3.5,4);   
%black
\filldraw  [ultra thick, fill=black] (0.5,0.5) circle [radius=0.16] ;\filldraw  [ultra thick, fill=black] (2.5,1.5) circle [radius=0.16] ;
\filldraw  [ultra thick, fill=black] (1.5,2.5) circle [radius=0.16] ;\filldraw  [ultra thick, fill=black] (3.5,3.5) circle [radius=0.16] ;
%white
\draw  [line width=1.28pt,fill=white] (3.5,0.5) circle [radius=0.16] ;\draw  [line width=1.28pt, fill=white] (1.5,1.5) circle [radius=0.16] ;
\draw  [line width=1.28pt, fill=white] (2.5,2.5) circle [radius=0.16] ;\draw  [line width=1.28pt, fill=white] (0.5,3.5) circle [radius=0.16] ;

\end{tikzpicture} }} ;  

\node (PM3) at (9.3,0) 
{\scalebox{0.51}{
\begin{tikzpicture}
%vertex
\node (B1) at (0.5,0.5){$$}; \node (W1) at (3.5,0.5){$$}; \node (B2) at (2.5,1.5){$$}; 
\node (W2) at (1.5,1.5){$$}; \node (B3) at (1.5,2.5){$$}; \node (W3) at (2.5,2.5){$$}; \node (B4) at (3.5,3.5){$$}; \node (W4) at (0.5,3.5){$$}; 
\draw[line width=0.024cm]  (0,0) rectangle (4,4);

%perfect matching
\draw[line width=0.4cm,color=lightgray] (W2)--(B1); \draw[line width=0.4cm,color=lightgray] (W1)--(B2);
\draw[line width=0.4cm,color=lightgray] (W3)--(B3); \draw[line width=0.4cm,color=lightgray] (W4)--(0,3.5);\draw[line width=0.4cm,color=lightgray] (B4)--(4,3.5);

%edge
\draw[line width=0.05cm]  (W2)--(B2)--(W3)--(B3)--(W2);
\draw[line width=0.05cm] (B1)--(W2); \draw[line width=0.05cm] (B2)--(W1); 
\draw[line width=0.05cm] (B3)--(W4); \draw[line width=0.05cm] (B4)--(W3);  
\draw[line width=0.05cm] (B1)--(0.5,0); \draw[line width=0.05cm] (B1)--(0,0.5); \draw[line width=0.05cm] (W1)--(3.5,0);   
\draw[line width=0.05cm] (W1)--(4,0.5); \draw[line width=0.05cm] (W4)--(0,3.5); \draw[line width=0.05cm] (W4)--(0.5,4);  
\draw[line width=0.05cm] (B4)--(4,3.5); \draw[line width=0.05cm] (B4)--(3.5,4);   
%black
\filldraw  [ultra thick, fill=black] (0.5,0.5) circle [radius=0.16] ;\filldraw  [ultra thick, fill=black] (2.5,1.5) circle [radius=0.16] ;
\filldraw  [ultra thick, fill=black] (1.5,2.5) circle [radius=0.16] ;\filldraw  [ultra thick, fill=black] (3.5,3.5) circle [radius=0.16] ;
%white
\draw  [line width=1.28pt,fill=white] (3.5,0.5) circle [radius=0.16] ;\draw  [line width=1.28pt, fill=white] (1.5,1.5) circle [radius=0.16] ;
\draw  [line width=1.28pt, fill=white] (2.5,2.5) circle [radius=0.16] ;\draw  [line width=1.28pt, fill=white] (0.5,3.5) circle [radius=0.16] ;

\end{tikzpicture} }};

\node (PM4) at (12.4,0) 
{\scalebox{0.51}{
\begin{tikzpicture}
%vertex
\node (B1) at (0.5,0.5){$$}; \node (W1) at (3.5,0.5){$$}; \node (B2) at (2.5,1.5){$$}; 
\node (W2) at (1.5,1.5){$$}; \node (B3) at (1.5,2.5){$$}; \node (W3) at (2.5,2.5){$$}; \node (B4) at (3.5,3.5){$$}; \node (W4) at (0.5,3.5){$$}; 
\draw[line width=0.024cm]  (0,0) rectangle (4,4);

%perfect matching
\draw[line width=0.4cm,color=lightgray] (W2)--(B1); \draw[line width=0.4cm,color=lightgray] (W4)--(B3); 
\draw[line width=0.4cm,color=lightgray] (W3)--(B2); \draw[line width=0.4cm,color=lightgray] (W1)--(3.5,0);\draw[line width=0.4cm,color=lightgray] (B4)--(3.5,4);  

%edge
\draw[line width=0.05cm]  (W2)--(B2)--(W3)--(B3)--(W2);\draw[line width=0.05cm] (B1)--(W2); \draw[line width=0.05cm] (B2)--(W1); 
\draw[line width=0.05cm] (B3)--(W4); \draw[line width=0.05cm] (B4)--(W3);  
\draw[line width=0.05cm] (B1)--(0.5,0); \draw[line width=0.05cm] (B1)--(0,0.5); \draw[line width=0.05cm] (W1)--(3.5,0);   
\draw[line width=0.05cm] (W1)--(4,0.5); \draw[line width=0.05cm] (W4)--(0,3.5); \draw[line width=0.05cm] (W4)--(0.5,4);  
\draw[line width=0.05cm] (B4)--(4,3.5); \draw[line width=0.05cm] (B4)--(3.5,4);   
%black
\filldraw  [ultra thick, fill=black] (0.5,0.5) circle [radius=0.16] ;\filldraw  [ultra thick, fill=black] (2.5,1.5) circle [radius=0.16] ;
\filldraw  [ultra thick, fill=black] (1.5,2.5) circle [radius=0.16] ;\filldraw  [ultra thick, fill=black] (3.5,3.5) circle [radius=0.16] ;
%white
\draw  [line width=1.28pt,fill=white] (3.5,0.5) circle [radius=0.16] ;\draw  [line width=1.28pt, fill=white] (1.5,1.5) circle [radius=0.16] ;
\draw  [line width=1.28pt, fill=white] (2.5,2.5) circle [radius=0.16] ;\draw  [line width=1.28pt, fill=white] (0.5,3.5) circle [radius=0.16] ;

\end{tikzpicture} }} ;
\end{tikzpicture}
}}
\end{center}

Then, we have the following table, and $(e_0\sfB)^*\cong e_2\sfB$, $(e_1\sfB)^*\cong e_3\sfB$. 

\medskip

\begin{center}
{\small{
\begin{tabular}{|l|l|l|l|} \hline 
$e_0\sfB$&$e_1\sfB$&$e_2\sfB$&$e_3\sfB$ \\ \hline
$T_{00}\cong R$&$T_{10}\cong T(0,1,0,0)$&$T_{20}\cong T(-1,1,0,0)$&$T_{30}\cong T(-1,0,0,0)$ \\
$T_{01}\cong T(2,1,0,0)$&$T_{11}\cong R$&$T_{21}\cong T(-1,0,0,0)$&$T_{31}\cong T(1,1,0,0)$ \\
$T_{02}\cong T(3,1,0,0)$&$T_{12}\cong T(1,0,0,0)$&$T_{22}\cong R$&$T_{32}\cong T(2,1,0,0)$ \\
$T_{03}\cong T(1,0,0,0)$&$T_{13}\cong T(1,1,0,0)$&$T_{23}\cong T(0,1,0,0)$&$T_{33}\cong R$\\ \hline 
\end{tabular}
}}
\end{center}

\subsubsection{Exchange graph of type 4a}
By combining the above results, we describe the exchange graph $\EG(\TMMG(R))$ for the case of type 4a as follows, 
and it is actually connected. 
In the following picture, a double circle stands for the origin $(0,0)\in\ZZ^2$ and
each point $(a,b)\in\ZZ^2$ corresponds to a divisorial ideal $T(a,b,0,0)$. 
We note that if $M, M^\prime\in\TMMG(R)$ have the same direct summands except one component, then we connect $M$ to $M^\prime$ in $\EG(\TMMG(R))$ (see Section~\ref{subsec_MMmutation}). 
Also, mutations are determined by a minimal right approximation given in (\ref{Mk_approx}), 
and we have that $[\calK_k]=[M_{t(a_1)}]+[M_{t(a_2)}]-[M_k]$. Thus, we also compute the exchange graph using this observation. 

%%%%%%%%%%%%%%%exchange_graph_4a%%%%%%%%%%%%%%%%%%%%%%%%
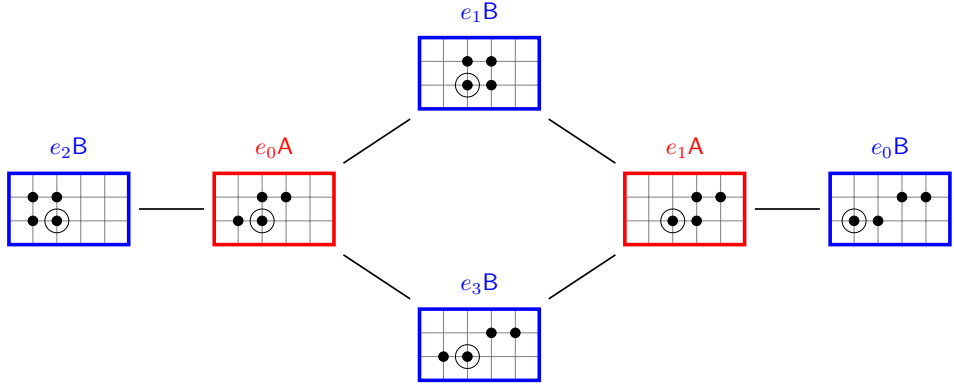
\begin{center}
{\scalebox{0.9}{
\begin{tikzpicture}

\node [blue]  at (0,0.9) {$e_2\sfB$}; 
\node (B2) at (0,0)
{\scalebox{0.35}{\begin{tikzpicture}
\draw [step=1,thin, gray] (-2,-1) grid (3,2); 
\draw [blue, line width=0.15cm] (-2,-1) rectangle (3,2);
\filldraw  [thick, fill=black] (-1,0) circle [radius=0.2] ;
\draw (0,0) circle [ultra thick, radius=0.5]; 
\filldraw [thick, fill=black] (0,0) circle [radius=0.2]; 
\filldraw  [thick, fill=black] (0,1) circle [radius=0.2] ;
\filldraw  [thick, fill=black] (-1,1) circle [radius=0.2] ;
\end{tikzpicture}}};

\node [red]  at (3,0.9) {$e_0\sfA$}; 
\node (A0) at (3,0)
{\scalebox{0.35}{\begin{tikzpicture}
\draw [step=1,thin, gray] (-2,-1) grid (3,2); 
\draw [red, line width=0.15cm] (-2,-1) rectangle (3,2);
\draw (0,0) circle [ultra thick, radius=0.5]; 
\filldraw [thick, fill=black] (0,0) circle [radius=0.2]; 
\filldraw  [thick, fill=black] (0,1) circle [radius=0.2] ;
\filldraw  [thick, fill=black] (-1,0) circle [radius=0.2] ;
\filldraw  [thick, fill=black] (1,1) circle [radius=0.2] ;
\end{tikzpicture}}};

\node [blue]  at (6,2.9) {$e_1\sfB$}; 
\node (B1) at (6,2)
{\scalebox{0.35}{\begin{tikzpicture}
\draw [step=1,thin, gray] (-2,-1) grid (3,2); 
\draw [blue, line width=0.15cm] (-2,-1) rectangle (3,2);
\draw (0,0) circle [ultra thick, radius=0.5]; 
\filldraw [thick, fill=black] (0,0) circle [radius=0.2]; 
\filldraw  [thick, fill=black] (0,1) circle [radius=0.2] ;
\filldraw  [thick, fill=black] (1,0) circle [radius=0.2] ;
\filldraw  [thick, fill=black] (1,1) circle [radius=0.2] ;
\end{tikzpicture}}};

\node [blue]  at (6,-1.1) {$e_3\sfB$}; 
\node (B3) at (6,-2)
{\scalebox{0.35}{\begin{tikzpicture}
\draw [step=1,thin, gray] (-2,-1) grid (3,2); 
\draw [blue, line width=0.15cm] (-2,-1) rectangle (3,2);
\draw (0,0) circle [ultra thick, radius=0.5]; 
\filldraw [thick, fill=black] (0,0) circle [radius=0.2]; 
\filldraw  [thick, fill=black] (-1,0) circle [radius=0.2] ;
\filldraw  [thick, fill=black] (2,1) circle [radius=0.2] ;
\filldraw  [thick, fill=black] (1,1) circle [radius=0.2] ;
\end{tikzpicture}}};

\node [red]  at (9,0.9) {$e_1\sfA$}; 
\node (A1) at (9,0)
{\scalebox{0.35}{\begin{tikzpicture}
\draw [step=1,thin, gray] (-2,-1) grid (3,2); 
\draw [red, line width=0.15cm] (-2,-1) rectangle (3,2);
\draw (0,0) circle [ultra thick, radius=0.5]; 
\filldraw [thick, fill=black] (0,0) circle [radius=0.2]; 
\filldraw  [thick, fill=black] (1,0) circle [radius=0.2] ;
\filldraw  [thick, fill=black] (2,1) circle [radius=0.2] ;
\filldraw  [thick, fill=black] (1,1) circle [radius=0.2] ;
\end{tikzpicture}}};

\node [blue]  at (12,0.9) {$e_0\sfB$}; 
\node (B0) at (12,0)
{\scalebox{0.35}{\begin{tikzpicture}
\draw [step=1,thin, gray] (-1,-1) grid (4,2); 
\draw [blue, line width=0.15cm] (-1,-1) rectangle (4,2);
\draw (0,0) circle [ultra thick, radius=0.5]; 
\filldraw [thick, fill=black] (0,0) circle [radius=0.2]; 
\filldraw  [thick, fill=black] (1,0) circle [radius=0.2] ;
\filldraw  [thick, fill=black] (2,1) circle [radius=0.2] ;
\filldraw  [thick, fill=black] (3,1) circle [radius=0.2] ;
\end{tikzpicture}}};

\draw[line width=0.025cm](B2)--(A0) ;\draw[line width=0.025cm] (A0)--(B1) ; \draw[line width=0.025cm] (A0)--(B3) ; 
\draw[line width=0.025cm] (B1)--(A1) ; \draw[line width=0.025cm] (B3)--(A1) ;\draw[line width=0.025cm] (A1)--(B0) ;

\end{tikzpicture} }} 
\end{center}

%%%%%%%%%%%%%%%%%%%%%%%%%%%%%%%%%%%%%%%%%%%%%%%%%%%%%%%
%\iffalse

\subsection{Type 4b} 
\label{type4b}

In this subsection, we consider the reflexive polygon of type 4b. 
Thus, let $R$ be the three dimensional complete local Gorenstein toric singularity defined by the cone $\sigma$: 
\[
\sigma=\mathrm{Cone}\{v_1=(0,1,1), v_2=(-1,0,1), v_3=(0,-1,1), v_4=(1,-1,1) \}. 
\]
Then, we consider a divisorial ideal $T(u_1,\cdots,u_4)$. 
As an element in $\Cl(R)$, we have that $[I_3]=3[I_1]$, $2[I_1]+[I_2]=0$, $[I_2]=[I_4]$.  
Therefore, we have that $\Cl(R)\cong\ZZ$, and each divisorial ideal is represented by $T(a,0,0,0)$ where $a\in\ZZ$. 
Also, there is the unique consistent dimer model (up to right equivalence) written below which give the reflexive polygon of type 4b as the perfect matching polygon.

\subsubsection{Type 4b-1}
We denote the complete Jacobian algebra associated with the following QP by $\sfA$. 
In the same way as the above case, we give splitting MM generators $e_i\sfA\cong\bigoplus_{j\in Q_0}T_{ij}$. 

\medskip

%%%%%%%%%%%%%%%%%%%%%
%%% type 4b %%%%%%%%%%%%%
%%%%%%%%%%%%%%%%%%%%%

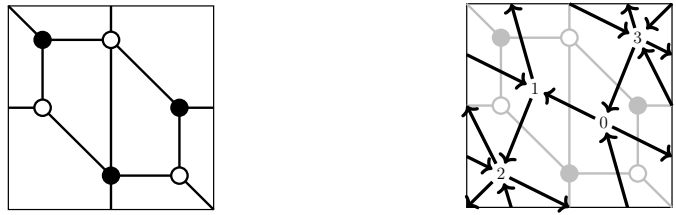
\begin{center}
\begin{tikzpicture}
\node (DM) at (0,0) 
{\scalebox{0.45}{
\begin{tikzpicture}
%vertex
\node (P1) at (1,3){$$}; \node (P2) at (3,1){$$}; \node (P3) at (5,1){$$}; \node (P4) at (5,3){$$}; \node (P5) at (3,5){$$};\node (P6) at (1,5){$$};
\draw[line width=0.036cm]  (0,0) rectangle (6,6);
%edge
\draw[line width=0.07cm]  (P1)--(P2)--(P3)--(P4)--(P5)--(P6)--(P1);\draw[line width=0.07cm]  (P2)--(P5); 
\draw[line width=0.07cm] (0,3)--(P1); \draw[line width=0.07cm]  (3,0)--(P2);\draw[line width=0.07cm]  (6,0)--(P3);
\draw[line width=0.07cm]  (6,3)--(P4);\draw[line width=0.07cm]  (3,6)--(P5);\draw[line width=0.07cm]  (0,6)--(P6);
%black
\filldraw  [ultra thick, fill=black] (3,1) circle [radius=0.24] ;\filldraw  [ultra thick, fill=black] (5,3) circle [radius=0.24] ;
\filldraw  [ultra thick, fill=black] (1,5) circle [radius=0.24] ;
%white
\draw  [line width=1.92pt,fill=white] (1,3) circle [radius=0.24] ;\draw  [line width=1.92pt,fill=white] (5,1) circle [radius=0.24] ;
\draw  [line width=1.92pt,fill=white] (3,5) circle [radius=0.24] ;
\end{tikzpicture}
} }; 

\node (QV) at (6,0) 
{\scalebox{0.45}{
\begin{tikzpicture}
%vertex
\node (P1) at (1,3){$$}; \node (P2) at (3,1){$$}; \node (P3) at
(5,1){$$}; \node (P4) at (5,3){$$}; \node (P5) at (3,5){$$};
\node (P6) at (1,5){$$};

\node (Q1) at (2,3.5){{\LARGE$1$}}; \node (Q0) at (4,2.5){{\LARGE$0$}}; 
\node (Q3) at (5,5){{\LARGE$3$}}; \node (Q2) at (1,1){{\LARGE$2$}}; 
\draw[line width=0.036cm]  (0,0) rectangle (6,6);

%edge
\draw[lightgray,line width=0.07cm]  (P1)--(P2)--(P3)--(P4)--(P5)--(P6)--(P1);\draw[lightgray,line width=0.07cm]  (P2)--(P5); 
\draw[lightgray,line width=0.07cm] (0,3)--(P1); \draw[lightgray,line width=0.07cm]  (3,0)--(P2);
\draw[lightgray,line width=0.07cm]  (6,0)--(P3);\draw[lightgray,line width=0.07cm]  (6,3)--(P4);
\draw[lightgray,line width=0.07cm]  (3,6)--(P5);\draw[lightgray,line width=0.07cm]  (0,6)--(P6);

%black
\filldraw  [ultra thick, draw=lightgray, fill=lightgray] (3,1) circle [radius=0.24] ;\filldraw  [ultra thick, draw=lightgray, fill=lightgray] (5,3) circle [radius=0.24] ;
\filldraw  [ultra thick, draw=lightgray, fill=lightgray] (1,5) circle [radius=0.24] ;

%white
\draw  [line width=1.92pt,draw=lightgray, fill=white] (1,3) circle [radius=0.24] ;\draw  [line width=1.92pt,draw=lightgray, fill=white] (5,1) circle [radius=0.24] ;
\draw  [line width=1.92pt,draw=lightgray, fill=white] (3,5) circle [radius=0.24] ;

\draw[->, line width=0.096cm] (Q0)--(Q1);\draw[->, line width=0.096cm] (Q3)--(Q0);\draw[->, line width=0.096cm] (Q1)--(Q2);
\draw[->, line width=0.096cm] (Q2)--(3,0);\draw[->, line width=0.096cm] (3,6)--(Q3);\draw[->, line width=0.096cm] (Q3)--(4.7,6);
\draw[->, line width=0.096cm] (4.7,0)--(Q0);\draw[->, line width=0.096cm] (Q0)--(6,1.5);\draw[->, line width=0.096cm] (0,1.5)--(Q2);
\draw[->, line width=0.096cm] (Q2)--(0,3);\draw[->, line width=0.096cm] (6,3)--(Q3);\draw[->, line width=0.096cm] (Q3)--(6,4.5);
\draw[->, line width=0.096cm] (0,4.5)--(Q1);
\draw[->, line width=0.096cm] (1.3,0)--(Q2);\draw[->, line width=0.096cm] (Q1)--(1.3,6);\draw[->, line width=0.096cm] (Q2)--(0,0);
\draw[->, line width=0.096cm] (6,6)--(Q3);
\end{tikzpicture}
} }; 
\end{tikzpicture}
\end{center}

\medskip

In the following figure, we fix the perfect matching $\sfP_0$.  
Then, $\sfP_1, \cdots, \sfP_4$ are extremal perfect matchings corresponding to $v_1, \cdots, v_4$. 
Namely, $v_i=(h(\sfP_i, \sfP_0), 1)\in\ZZ^3$. 

\medskip

%%%%%perfect_matching%%%%%
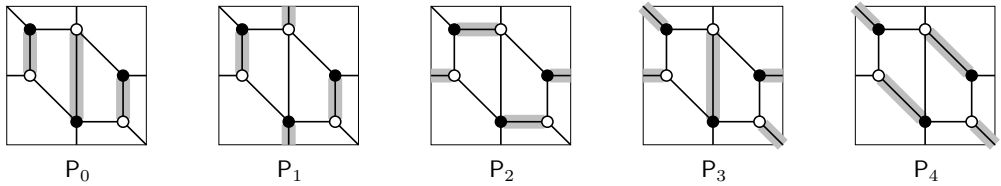
\begin{center}
{\scalebox{0.9}{
\begin{tikzpicture} 
\node at (0,-1.4) {$\sfP_0$};\node at (3.1,-1.4) {$\sfP_1$}; \node at (6.2,-1.4) {$\sfP_2$}; 
\node at (9.3,-1.4) {$\sfP_3$}; \node at (12.4,-1.4) {$\sfP_4$};

\node (PM0) at (0,0) 
{\scalebox{0.34}{\begin{tikzpicture}
%vertex
\node (P1) at (1,3){$$}; \node (P2) at (3,1){$$}; \node (P3) at (5,1){$$}; \node (P4) at (5,3){$$}; \node (P5) at (3,5){$$};\node (P6) at (1,5){$$};
\draw[line width=0.036cm]  (0,0) rectangle (6,6);
%perfect matching
\draw[line width=0.6cm,color=lightgray] (P1)--(P6);\draw[line width=0.6cm,color=lightgray] (P2)--(P5);\draw[line width=0.6cm,color=lightgray] (P3)--(P4);
%edge
\draw[line width=0.07cm]  (P1)--(P2)--(P3)--(P4)--(P5)--(P6)--(P1);\draw[line width=0.07cm]  (P2)--(P5); \draw[line width=0.07cm] (0,3)--(P1); 
\draw[line width=0.07cm]  (3,0)--(P2);\draw[line width=0.07cm]  (6,0)--(P3);\draw[line width=0.07cm]  (6,3)--(P4);
\draw[line width=0.07cm]  (3,6)--(P5);\draw[line width=0.07cm]  (0,6)--(P6);
%black
\filldraw  [ultra thick, fill=black] (3,1) circle [radius=0.24] ;\filldraw  [ultra thick, fill=black] (5,3) circle [radius=0.24] ;
\filldraw  [ultra thick, fill=black] (1,5) circle [radius=0.24] ;
%white
\draw  [line width=1.92pt,fill=white] (1,3) circle [radius=0.24] ;\draw  [line width=1.92pt,fill=white] (5,1) circle [radius=0.24] ;
\draw  [line width=1.92pt,fill=white] (3,5) circle [radius=0.24] ;
\end{tikzpicture} }}; 

\node (PM1) at (3.1,0) 
{\scalebox{0.34}{\begin{tikzpicture}
%vertex
\node (P1) at (1,3){$$}; \node (P2) at (3,1){$$}; \node (P3) at (5,1){$$}; \node (P4) at (5,3){$$}; \node (P5) at (3,5){$$};\node (P6) at (1,5){$$};
\draw[line width=0.036cm]  (0,0) rectangle (6,6);
%perfect matching
\draw[line width=0.6cm,color=lightgray] (P1)--(P6);\draw[line width=0.6cm,color=lightgray] (P2)--(3,0);
\draw[line width=0.6cm,color=lightgray] (P5)--(3,6);\draw[line width=0.6cm,color=lightgray] (P3)--(P4);
%edge
\draw[line width=0.07cm]  (P1)--(P2)--(P3)--(P4)--(P5)--(P6)--(P1);\draw[line width=0.07cm]  (P2)--(P5); \draw[line width=0.07cm] (0,3)--(P1); 
\draw[line width=0.07cm]  (3,0)--(P2);\draw[line width=0.07cm]  (6,0)--(P3);
\draw[line width=0.07cm]  (6,3)--(P4);\draw[line width=0.07cm]  (3,6)--(P5);\draw[line width=0.07cm]  (0,6)--(P6);
%black
\filldraw  [ultra thick, fill=black] (3,1) circle [radius=0.24] ;\filldraw  [ultra thick, fill=black] (5,3) circle [radius=0.24] ;
\filldraw  [ultra thick, fill=black] (1,5) circle [radius=0.24] ;
%white
\draw  [line width=1.92pt,fill=white] (1,3) circle [radius=0.24] ;\draw  [line width=1.92pt,fill=white] (5,1) circle [radius=0.24] ;
\draw  [line width=1.92pt,fill=white] (3,5) circle [radius=0.24] ;
\end{tikzpicture} }}; 

\node (PM2) at (6.2,0) 
{\scalebox{0.34}{\begin{tikzpicture}
%vertex
\node (P1) at (1,3){$$}; \node (P2) at (3,1){$$}; \node (P3) at (5,1){$$}; \node (P4) at (5,3){$$}; \node (P5) at (3,5){$$}; \node (P6) at (1,5){$$};
\draw[line width=0.036cm]  (0,0) rectangle (6,6);
%perfect matching
\draw[line width=0.6cm,color=lightgray] (P1)--(0,3);\draw[line width=0.6cm,color=lightgray] (P2)--(P3);
\draw[line width=0.6cm,color=lightgray] (P5)--(P6);\draw[line width=0.6cm,color=lightgray] (P4)--(6,3);
%edge
\draw[line width=0.07cm]  (P1)--(P2)--(P3)--(P4)--(P5)--(P6)--(P1);\draw[line width=0.07cm]  (P2)--(P5); \draw[line width=0.07cm] (0,3)--(P1); 
\draw[line width=0.07cm]  (3,0)--(P2);\draw[line width=0.07cm]  (6,0)--(P3);\draw[line width=0.07cm]  (6,3)--(P4);
\draw[line width=0.07cm]  (3,6)--(P5);\draw[line width=0.07cm]  (0,6)--(P6);
%black
\filldraw  [ultra thick, fill=black] (3,1) circle [radius=0.24] ;\filldraw  [ultra thick, fill=black] (5,3) circle [radius=0.24] ;
\filldraw  [ultra thick, fill=black] (1,5) circle [radius=0.24] ;
%white
\draw  [line width=1.92pt,fill=white] (1,3) circle [radius=0.24] ;\draw  [line width=1.92pt,fill=white] (5,1) circle [radius=0.24] ;
\draw  [line width=1.92pt,fill=white] (3,5) circle [radius=0.24] ;

\end{tikzpicture} }}; 

\node (PM3) at (9.3,0) 
{\scalebox{0.34}{\begin{tikzpicture}
%vertex
\node (P1) at (1,3){$$}; \node (P2) at (3,1){$$}; \node (P3) at (5,1){$$}; \node (P4) at (5,3){$$}; \node (P5) at (3,5){$$}; \node (P6) at (1,5){$$};
\draw[line width=0.036cm]  (0,0) rectangle (6,6);
%perfect matching
\draw[line width=0.6cm,color=lightgray] (P1)--(0,3);\draw[line width=0.6cm,color=lightgray] (P4)--(6,3);\draw[line width=0.6cm,color=lightgray] (P2)--(P5);
\draw[line width=0.6cm,color=lightgray] (P3)--(6,0);\draw[line width=0.6cm,color=lightgray] (P6)--(0,6);
%edge
\draw[line width=0.07cm]  (P1)--(P2)--(P3)--(P4)--(P5)--(P6)--(P1);\draw[line width=0.07cm]  (P2)--(P5); \draw[line width=0.07cm] (0,3)--(P1); 
\draw[line width=0.07cm]  (3,0)--(P2);\draw[line width=0.07cm]  (6,0)--(P3);\draw[line width=0.07cm]  (6,3)--(P4);
\draw[line width=0.07cm]  (3,6)--(P5);\draw[line width=0.07cm]  (0,6)--(P6);
%black
\filldraw  [ultra thick, fill=black] (3,1) circle [radius=0.24] ;\filldraw  [ultra thick, fill=black] (5,3) circle [radius=0.24] ;
\filldraw  [ultra thick, fill=black] (1,5) circle [radius=0.24] ;
%white
\draw  [line width=1.92pt,fill=white] (1,3) circle [radius=0.24] ;\draw  [line width=1.92pt,fill=white] (5,1) circle [radius=0.24] ;
\draw  [line width=1.92pt,fill=white] (3,5) circle [radius=0.24] ;
\end{tikzpicture} }}; 

\node (PM4) at (12.4,0) 
{\scalebox{0.34}{\begin{tikzpicture}
%vertex
\node (P1) at (1,3){$$}; \node (P2) at (3,1){$$}; \node (P3) at (5,1){$$}; \node (P4) at (5,3){$$}; \node (P5) at (3,5){$$}; \node (P6) at (1,5){$$};
\draw[line width=0.036cm]  (0,0) rectangle (6,6);
%perfect matching
\draw[line width=0.6cm,color=lightgray] (P1)--(P2);\draw[line width=0.6cm,color=lightgray] (P4)--(P5);\draw[line width=0.6cm,color=lightgray] (P3)--(6,0);
\draw[line width=0.6cm,color=lightgray] (P6)--(0,6);
%edge
\draw[line width=0.07cm]  (P1)--(P2)--(P3)--(P4)--(P5)--(P6)--(P1);\draw[line width=0.07cm]  (P2)--(P5); \draw[line width=0.07cm] (0,3)--(P1); 
\draw[line width=0.07cm]  (3,0)--(P2);\draw[line width=0.07cm]  (6,0)--(P3);\draw[line width=0.07cm]  (6,3)--(P4);
\draw[line width=0.07cm]  (3,6)--(P5);\draw[line width=0.07cm]  (0,6)--(P6);
%black
\filldraw  [ultra thick, fill=black] (3,1) circle [radius=0.24] ;\filldraw  [ultra thick, fill=black] (5,3) circle [radius=0.24] ;
\filldraw  [ultra thick, fill=black] (1,5) circle [radius=0.24] ;
%white
\draw  [line width=1.92pt,fill=white] (1,3) circle [radius=0.24] ;\draw  [line width=1.92pt,fill=white] (5,1) circle [radius=0.24] ;
\draw  [line width=1.92pt,fill=white] (3,5) circle [radius=0.24] ;
\end{tikzpicture} }}; 

\end{tikzpicture}
}}
\end{center}

\medskip

Then, we have the following table, and we see that $(e_0\sfA)^*\cong e_1\sfA$, $(e_2\sfA)^*\cong e_3\sfA$. 

\medskip

\begin{center}
{\small{
\begin{tabular}{|l|l|l|l|} \hline 
$e_0\sfA$&$e_1\sfA$&$e_2\sfA$&$e_3\sfA$ \\ \hline
$T_{00}\cong R$&$T_{10}\cong T(-3,0,0,0)$&$T_{20}\cong T(-1,0,0,0)$&$T_{30}\cong T(-2,0,0,0)$ \\
$T_{01}\cong T(3,0,0,0)$&$T_{11}\cong R$&$T_{21}\cong T(2,0,0,0)$&$T_{31}\cong T(1,0,0,0)$ \\
$T_{02}\cong T(1,0,0,0)$&$T_{12}\cong T(-2,0,0,0)$&$T_{22}\cong R$&$T_{32}\cong T(-1,0,0,0)$ \\
$T_{03}\cong T(2,0,0,0)$&$T_{13}\cong T(-1,0,0,0)$&$T_{23}\cong T(1,0,0,0)$&$T_{33}\cong R$\\ \hline 
\end{tabular}
}}
\end{center}

\subsubsection{Exchange graph of type 4b}
By the above results, we describe the exchange graph $\EG(\TMMG(R))$ for the case of type 4b as follows, 
and it is connected. 

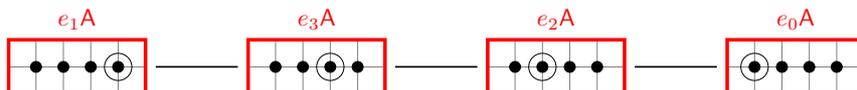
\begin{center}
{\scalebox{0.9}{
\begin{tikzpicture}

\node [red]  at (10.5,0.7) {$e_0\sfA$}; 
\node (A0) at (10.5,0)
{\scalebox{0.4}{\begin{tikzpicture}
\draw [step=1,thin, gray] (-1,-1) grid (4,1); 
\draw [red, line width=0.13cm] (-1,-1) rectangle (4,1);
\draw (0,0) circle [ultra thick, radius=0.5]; 
\filldraw [thick, fill=black] (0,0) circle [radius=0.2]; 
\filldraw  [thick, fill=black] (1,0) circle [radius=0.2] ;
\filldraw  [thick, fill=black] (2,0) circle [radius=0.2] ;
\filldraw  [thick, fill=black] (3,0) circle [radius=0.2] ;
\end{tikzpicture}}};

\node [red]  at (7,0.7) {$e_2\sfA$}; 
\node (A2) at (7,0)
{\scalebox{0.4}{\begin{tikzpicture}
\draw [step=1,thin, gray] (-1,-1) grid (4,1); 
\draw [red, line width=0.13cm] (-1,-1) rectangle (4,1);
\filldraw [thick, fill=black] (0,0) circle [radius=0.2]; 
\draw (1,0) circle [ultra thick, radius=0.5]; 
\filldraw  [thick, fill=black] (1,0) circle [radius=0.2] ;
\filldraw  [thick, fill=black] (2,0) circle [radius=0.2] ;
\filldraw  [thick, fill=black] (3,0) circle [radius=0.2] ;
\end{tikzpicture}}};

\node [red]  at (3.5,0.7) {$e_3\sfA$}; 
\node (A3) at (3.5,0)
{\scalebox{0.4}{\begin{tikzpicture}
\draw [step=1,thin, gray] (-1,-1) grid (4,1); 
\draw [red, line width=0.13cm] (-1,-1) rectangle (4,1);
\filldraw [thick, fill=black] (0,0) circle [radius=0.2]; 
\draw (2,0) circle [ultra thick, radius=0.5]; 
\filldraw  [thick, fill=black] (1,0) circle [radius=0.2] ;
\filldraw  [thick, fill=black] (2,0) circle [radius=0.2] ;
\filldraw  [thick, fill=black] (3,0) circle [radius=0.2] ;
\end{tikzpicture}}};

\node [red]  at (0,0.7) {$e_1\sfA$}; 
\node (A1) at (0,0)
{\scalebox{0.4}{\begin{tikzpicture}
\draw [step=1,thin, gray] (-1,-1) grid (4,1); 
\draw [red, line width=0.13cm] (-1,-1) rectangle (4,1);
\filldraw [thick, fill=black] (0,0) circle [radius=0.2]; 
\draw (3,0) circle [ultra thick, radius=0.5]; 
\filldraw  [thick, fill=black] (1,0) circle [radius=0.2] ;
\filldraw  [thick, fill=black] (2,0) circle [radius=0.2] ;
\filldraw  [thick, fill=black] (3,0) circle [radius=0.2] ;
\end{tikzpicture}}};

\draw[line width=0.025cm] (A0)--(A2)--(A3)--(A1) ;

\end{tikzpicture} }}
\end{center}

%%%%%%%%%%%%%%%%%%%%%%%%%%%%%%%%%%%%%%%%%%%%%%%%%%%%%%%
%%%%%%%%%%%%%%%%%%%%%%%%%%%%%%%%%%%%%%%%%%%%%%%%%%%%%%%
%%%%%%%%%%%%%%%%%%%%%%%%%%%%%%%%%%%%%%%%%%%%%%%%%%%%%%%
\subsection{Type 5a}
\label{type5a}

In this subsection, we consider the reflexive polygon of type 5a. 
Thus, let $R$ be the three dimensional complete local Gorenstein toric singularity defined by the cone $\sigma$: 
\[
\sigma=\mathrm{Cone}\{v_1=(1,0,1), v_2=(0,1,1), v_3=(-1,1,1), v_4=(-1,0,1), v_5=(0,-1,1) \}. 
\]
As an element in $\Cl(R)$, we obtain $2[I_1]+2[I_2]+[I_3]=0$, $3[I_1]+2[I_2]-[I_4]=0$, $2[I_1]+[I_2]+[I_5]=0$.  
Therefore, we have that $\Cl(R)\cong\ZZ^2$, and each divisorial ideal is represented by $T(a,b,0,0)$ where $a, b\in\ZZ$. 
Also, there are two consistent dimer models (up to right equivalence) written below which give the reflexive polygon of type 5a as the perfect matching polygon. 

\subsubsection{Type 5a-1}

First, we consider the following consistent dimer model and the associated QP $(Q, W_Q)$. 
For simplicity, we denote by $\sfA$ the complete Jacobian algebra $\calP(Q, W_Q)$. 
Then, we have that $\sfA\cong\End_R(\bigoplus_{j\in Q_0}T_{ij})$ and $e_i\sfA\cong\bigoplus_{j\in Q_0}T_{ij}$ for any $i\in Q_0$ by Theorem~\ref{NCCR1}. 
We will give these splitting MM generators $\bigoplus_{j\in Q_0}T_{ij}$. 

\medskip

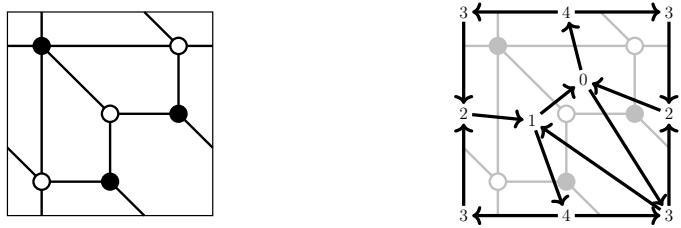
\begin{center}

\end{center}

\medskip

In the following figure, we fix the perfect matching $\sfP_0$.  
Then, $\sfP_1, \cdots, \sfP_5$ are extremal perfect matchings corresponding to $v_1, \cdots, v_5$. 
Namely, $v_i=(h(\sfP_i, \sfP_0), 1)\in\ZZ^3$.

\medskip
%%%%%perfect_matching%%%%%

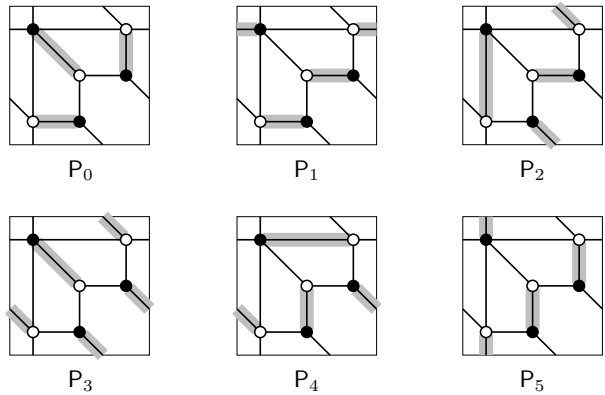
\begin{center}
{\scalebox{0.9}{
%
}}
\end{center}

\medskip

Then, we have the following table. For simplicity, we denote by $T(a,b)$ a divisorial ideal $T(a,b,0,0,0)$. 

\medskip

\begin{center}
{\small{
\begin{tabular}{|l|l|l|l|l|} \hline 
$e_0\sfA$&$e_1\sfA$&$e_2\sfA$&$e_3\sfA$&$e_4\sfA$ \\ \hline
$T_{00}\cong R$&$T_{10}\cong T(-2,-2)$&$T_{20}\cong T(-2,-1)$&$T_{30}\cong T(-1,-1)$&$T_{40}\cong T(-3,-2)$ \\
$T_{01}\cong T(2,2)$&$T_{11}\cong R$&$T_{21}\cong T(0,1)$&$T_{31}\cong T(1,1)$&$T_{41}\cong T(-1,0)$ \\
$T_{02}\cong T(2,1)$&$T_{12}\cong T(0,-1)$&$T_{22}\cong R$&$T_{32}\cong T(1,0)$&$T_{42}\cong T(-1,-1)$ \\
$T_{03}\cong T(1,1)$&$T_{13}\cong T(-1,-1)$&$T_{23}\cong T(-1,0)$&$T_{33}\cong R$&$T_{43}\cong T(-2,-1)$\\ 
$T_{04}\cong T(3,2)$&$T_{14}\cong T(1,0)$&$T_{24}\cong T(1,1)$&$T_{34}\cong T(2,1)$&$T_{44}\cong R$\\ \hline 
\end{tabular} 
}}
\end{center}

\subsubsection{Type 5a-2}

%%%%%%%%%%%%%%%%%%%%%
%%% type 5a-2 %%%%%%%%%%%%%
%%%%%%%%%%%%%%%%%%%%%

Next, we consider the following consistent dimer model and the associated QP $(Q, W_Q)$. 
We denote by $\sfB$ the complete Jacobian algebra $\calP(Q, W_Q)$. 
We will give splitting MM generators $e_i\sfB$. 

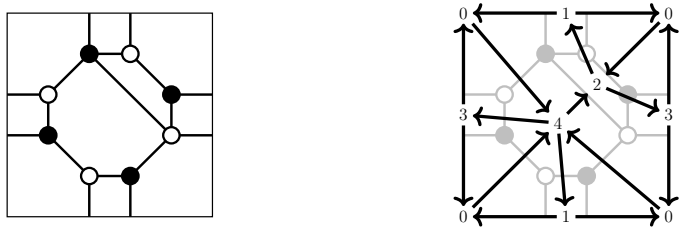
\medskip

\begin{center}
%
\end{center}

\medskip

In the following figure, we fix the perfect matching $\sfP_0$.  
Then, $\sfP_1, \cdots, \sfP_5$ are extremal perfect matchings corresponding to $v_1, \cdots, v_5$. 

\medskip

%%%%%perfect_matching%%%%%
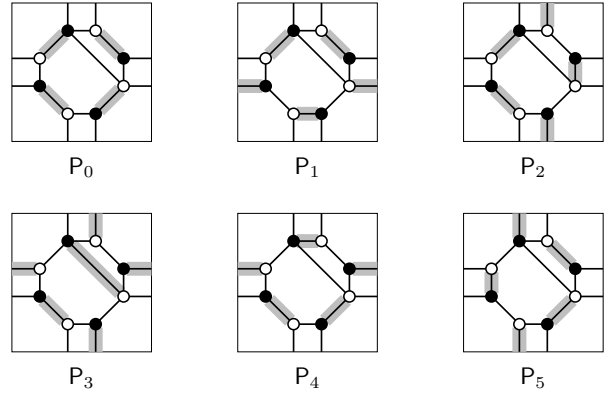
\begin{center}
{\scalebox{0.9}{
%
}}
\end{center}

\medskip

Then, we have the following tables. For simplicity, we denote by $T(a,b)$ a divisorial ideal $T(a,b,0,0,0)$. 

\medskip

\begin{center}
{\small{
\begin{tabular}{|l|l|l|l|l|} \hline 
$e_0\sfB$&$e_1\sfB$&$e_2\sfB$&$e_3\sfB$&$e_4\sfB$ \\ \hline
$T_{00}\cong R$&$T_{10}\cong T(-2,-1)$&$T_{20}\cong T(1,1)$&$T_{30}\cong T(1,0)$&$T_{40}\cong T(-1,-1)$ \\
$T_{01}\cong T(2,1)$&$T_{11}\cong R$&$T_{21}\cong T(3,2)$&$T_{31}\cong T(3,1)$&$T_{41}\cong T(1,0)$ \\
$T_{02}\cong T(-1,-1)$&$T_{12}\cong T(-3,-2)$&$T_{22}\cong R$&$T_{32}\cong T(0,-1)$&$T_{42}\cong T(-2,-2)$ \\
$T_{03}\cong T(-1,0)$&$T_{13}\cong T(-3,-1)$&$T_{23}\cong T(0,1)$&$T_{33}\cong R$&$T_{43}\cong T(-2,-1)$\\ 
$T_{04}\cong T(1,1)$&$T_{14}\cong T(-1,0)$&$T_{24}\cong T(2,2)$&$T_{34}\cong T(2,1)$&$T_{44}\cong R$\\ \hline 
\end{tabular} 
}}
\end{center}

\medskip

%%%%%%%%%%%%%%%%%%%%%%%%%%%%%%%%%%%%%%%%%%%%
\subsubsection{Exchange graph of type 5a}
By combining the above results, we describe the exchange graph $\EG(\TMMG(R))$ for the case of type 5a as follows, 
and it is actually connected. 

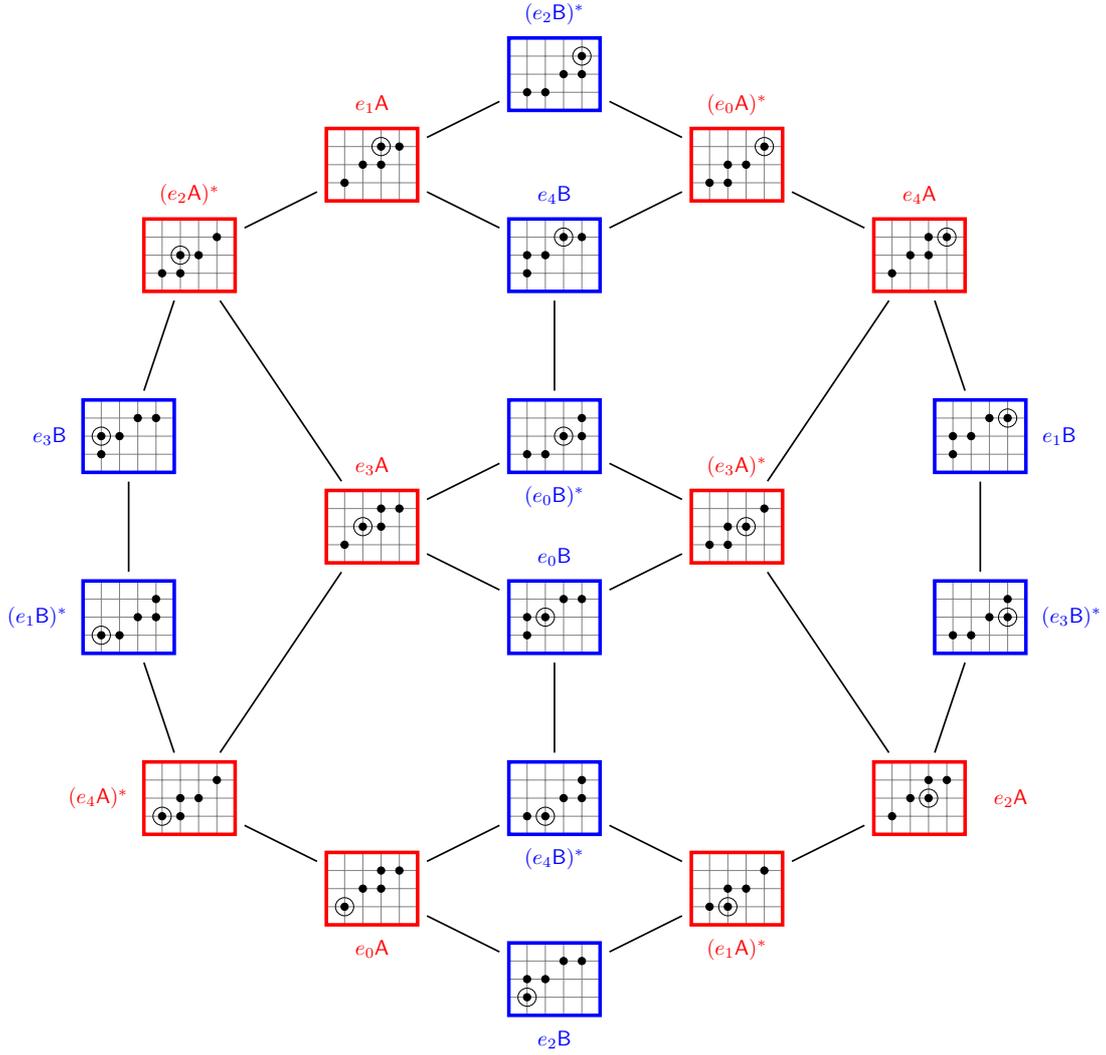
\begin{center}
{\scalebox{0.8}{
%
}};

\draw[line width=0.028cm] (B2+)--(A0+) ; \draw[line width=0.028cm] (B2+)--(A1) ; \draw[line width=0.028cm] (B4)--(A0+) ; 
\draw[line width=0.028cm] (B4)--(A1) ;   \draw[line width=0.028cm] (B4)--(B0+) ; 
\draw[line width=0.028cm] (A0+)--(A4) ; \draw[line width=0.028cm] (A4)--(B1) ; \draw[line width=0.028cm] (B1)--(B3+) ; 
\draw[line width=0.028cm] (A3+)--(A4) ;\draw[line width=0.028cm] (A3+)--(B0+) ; \draw[line width=0.028cm] (A3+)--(B0) ;
\draw[line width=0.028cm] (A3+)--(A2) ; \draw[line width=0.028cm] (B3+)--(A2) ; \draw[line width=0.028cm] (B0)--(B4+) ; 
\draw[line width=0.028cm] (B4+)--(A1+) ; \draw[line width=0.028cm] (A1+)--(B2) ; \draw[line width=0.028cm] (A1+)--(A2) ;
\draw[line width=0.028cm] (A1)--(A2+) ; \draw[line width=0.028cm] (A3)--(A2+) ;  \draw[line width=0.028cm] (A2+)--(B3) ; 
\draw[line width=0.028cm] (A3)--(B0+) ; \draw[line width=0.028cm] (B3)--(B1+) ; \draw[line width=0.028cm] (A4+)--(B1+) ;   
\draw[line width=0.028cm] (A3)--(A4+) ; \draw[line width=0.028cm] (A3)--(B0) ; \draw[line width=0.028cm] (A0)--(B4+) ; 
\draw[line width=0.028cm] (A0)--(B2) ; \draw[line width=0.028cm] (A0)--(A4+) ;
\end{tikzpicture} }}
\end{center}

%%%%%%%%%%%%%%%%%%%%%%%%%%%%%%%%%%%%%%%%%%%%%%%%%%%%%%%
\subsection{Type 5b}
\label{type5b}

%%%%%%%%%%%%%%%%%%%%%
%%% type 5b %%%%%%%%%%%%%
%%%%%%%%%%%%%%%%%%%%%

In this subsection, we consider the reflexive polygon of type 5b. 
Thus, let $R$ be the three dimensional complete local Gorenstein toric singularity defined by the cone $\sigma$: 
\[
\sigma=\mathrm{Cone}\{v_1=(0,1,1), v_2=(-1,0,1), v_3=(-1,-1,1), v_4=(1,-1,1) \}. 
\]
As an element in $\Cl(R)$, we obtain $[I_1]+2[I_4]=0$, $[I_2]-4[I_4]=0$, $[I_3]+3[I_4]=0$.  
Therefore, we have that $\Cl(R)\cong\ZZ$, and each divisorial ideal is represented by $T(0,0,0,d)$ where $d\in\ZZ$. 
There is the unique consistent dimer model (up to right equivalence) written below which give the reflexive polygon of type 5b as the perfect matching polygon. 

\subsubsection{Type 5b-1}

We consider the following consistent dimer model and the associated QP $(Q, W_Q)$, 
and denote by $\sfA$ the complete Jacobian algebra $\calP(Q, W_Q)$. 
We consider splitting MM generators $e_i\sfA$. 

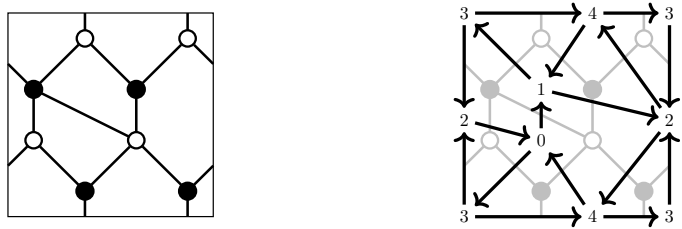
\medskip

\begin{center}
\begin{tikzpicture}
\node (DM) at (0,0) 
{\scalebox{0.675}{
\begin{tikzpicture}
%vertex
\node (B1) at (1,1){$$}; \node (B2) at (3,1){$$}; \node (W1) at (0,2){$$}; \node (W2) at (2,2){$$};
\node (B3) at (0,3){$$}; \node (B4) at (2,3){$$}; \node (W3) at (1,4){$$}; \node (W4) at (3,4){$$};
\draw[line width=0.024cm]  (-0.5,0.5) rectangle (3.5,4.5);

%edge
\draw[line width=0.05cm]  (B1)--(W2)--(B4)--(W3)--(B3)--(W1)--(B1);
\draw[line width=0.05cm]  (W2)--(B2)--(3.5,1.5);\draw[line width=0.05cm] (W2)--(B3);
\draw[line width=0.05cm] (-0.5,1.5)--(W1);\draw[line width=0.05cm] (B4)--(W4)--(3.5,3.5); \draw[line width=0.05cm] (B3)--(-0.5,3.5);
\draw[line width=0.05cm] (B1)--(1,0.5); \draw[line width=0.05cm] (B2)--(3,0.5); 
\draw[line width=0.05cm] (W3)--(1,4.5); \draw[line width=0.05cm] (W4)--(3,4.5); 

%black
\filldraw  [ultra thick, fill=black] (1,1) circle [radius=0.16] ;\filldraw  [ultra thick, fill=black] (3,1) circle [radius=0.16] ;
\filldraw  [ultra thick, fill=black] (0,3) circle [radius=0.16] ;\filldraw  [ultra thick, fill=black] (2,3) circle [radius=0.16] ;
%white
\draw  [line width=1.28pt,fill=white] (0,2) circle [radius=0.16] ;\draw [line width=1.28pt, fill=white] (2,2) circle [radius=0.16] ;
\draw  [line width=1.28pt, fill=white] (1,4) circle [radius=0.16] ;\draw [line width=1.28pt, fill=white] (3,4) circle [radius=0.16] ;
\end{tikzpicture}

} }; 

\node (QV) at (6,0) 
{\scalebox{0.675}{
\begin{tikzpicture}
%vertex
\node (B1) at (1,1){$$}; \node (B2) at (3,1){$$}; \node (W1) at (0,2){$$}; \node (W2) at (2,2){$$};
\node (B3) at (0,3){$$}; \node (B4) at (2,3){$$}; \node (W3) at (1,4){$$}; \node (W4) at (3,4){$$};

%vertex
\node (Q0) at (1,2){$0$}; \node (Q1) at (1,3){$1$}; \node (Q2) at (3.5,2.4){$2$}; \node (Q2a) at (-0.5,2.4){$2$};
\node (Q3a) at (-0.5,0.5){$3$}; \node (Q3b) at (-0.5,4.5){$3$};  \node (Q3c) at (3.5,0.5){$3$};  \node (Q3d) at (3.5,4.5){$3$}; 
\node (Q4a) at (2,0.5){$4$}; \node (Q4b) at (2,4.5){$4$};
%\draw[line width=0.024cm]  (-0.5,0.5) rectangle (3.5,4.5);

%edge
\draw[lightgray, line width=0.05cm]  (B1)--(W2)--(B4)--(W3)--(B3)--(W1)--(B1);
\draw[lightgray, line width=0.05cm]  (W2)--(B2)--(3.5,1.5);\draw[lightgray, line width=0.05cm] (W2)--(B3);
\draw[lightgray, line width=0.05cm] (-0.5,1.5)--(W1);\draw[lightgray, line width=0.05cm] (B4)--(W4)--(3.5,3.5); \draw[lightgray, line width=0.05cm] (B3)--(-0.5,3.5);
\draw[lightgray, line width=0.05cm] (B1)--(1,0.5); \draw[lightgray, line width=0.05cm] (B2)--(3,0.5); 
\draw[lightgray, line width=0.05cm] (W3)--(1,4.5); \draw[lightgray, line width=0.05cm] (W4)--(3,4.5); 

%black
\filldraw  [ultra thick, draw=lightgray, fill=lightgray] (1,1) circle [radius=0.16] ;\filldraw  [ultra thick, draw=lightgray, fill=lightgray] (3,1) circle [radius=0.16] ;
\filldraw  [ultra thick,draw=lightgray, fill=lightgray] (0,3) circle [radius=0.16] ;\filldraw  [ultra thick,draw=lightgray, fill=lightgray] (2,3) circle [radius=0.16] ;
%white
\draw  [line width=1.28pt, draw=lightgray,fill=white] (0,2) circle [radius=0.16] ;\draw  [line width=1.28pt, draw=lightgray, fill=white] (2,2) circle [radius=0.16] ;
\draw  [line width=1.28pt, draw=lightgray, fill=white] (1,4) circle [radius=0.16] ;\draw  [line width=1.28pt, draw=lightgray, fill=white] (3,4) circle [radius=0.16] ;

\draw[->, line width=0.064cm] (Q0)--(Q1); \draw[->, line width=0.064cm] (Q1)--(Q2);\draw[->, line width=0.064cm] (Q2)--(Q4a);
\draw[->, line width=0.064cm] (Q4a)--(Q0); \draw[->, line width=0.064cm] (Q0)--(Q3a);\draw[->, line width=0.064cm] (Q1)--(Q3b);
\draw[->, line width=0.064cm] (Q4b)--(Q1); \draw[->, line width=0.064cm] (Q3a)--(Q4a); \draw[->, line width=0.064cm] (Q3b)--(Q4b);
\draw[->, line width=0.064cm] (Q2)--(Q4b);\draw[->, line width=0.064cm] (Q4a)--(Q3c);\draw[->, line width=0.064cm] (Q4b)--(Q3d);
\draw[->, line width=0.064cm] (Q2a)--(Q0);\draw[->, line width=0.064cm] (Q3d)--(Q2);
\draw[->, line width=0.064cm] (Q3b)--(Q2a);\draw[->, line width=0.064cm] (Q3c)--(Q2); \draw[->, line width=0.064cm] (Q3a)--(Q2a);
\end{tikzpicture}
} }; 
\end{tikzpicture}
\end{center}

\medskip

In the following figure, we fix the perfect matching $\sfP_0$.  
Then, $\sfP_1, \cdots, \sfP_4$ are extremal perfect matchings corresponding to $v_1, \cdots, v_4$. 
Namely, $v_i=(h(\sfP_i, \sfP_0), 1)\in\ZZ^3$.

\medskip

%%%%%perfect_matching%%%%%
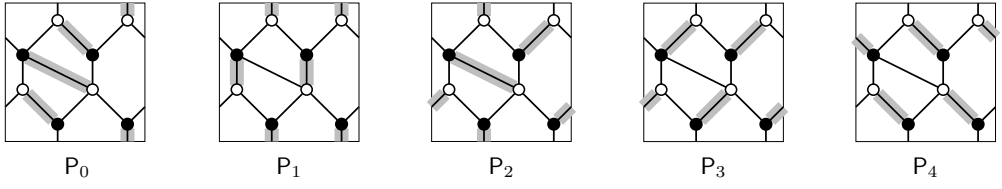
\begin{center}
{\scalebox{0.9}{
\begin{tikzpicture} 
\node at (0,-1.4) {$\sfP_0$};\node at (3.1,-1.4) {$\sfP_1$}; \node at (6.2,-1.4) {$\sfP_2$}; 
\node at (9.3,-1.4) {$\sfP_3$}; \node at (12.4,-1.4) {$\sfP_4$};

\node (PM0) at (0,0) 
{\scalebox{0.51}{\begin{tikzpicture}
%vertex
\node (B1) at (1,1){$$}; \node (B2) at (3,1){$$}; \node (W1) at (0,2){$$}; \node (W2) at (2,2){$$};
\node (B3) at (0,3){$$}; \node (B4) at (2,3){$$}; \node (W3) at (1,4){$$}; \node (W4) at (3,4){$$};

\draw[line width=0.024cm]  (-0.5,0.5) rectangle (3.5,4.5);
%perfect matching
\draw[line width=0.4cm,color=lightgray] (W1)--(B1); \draw[line width=0.4cm,color=lightgray] (W2)--(B3);
\draw[line width=0.4cm,color=lightgray] (B2)--(3,0.5); \draw[line width=0.4cm,color=lightgray] (3,4.5)--(W4); 
\draw[line width=0.4cm,color=lightgray] (W3)--(B4);

%edge
\draw[line width=0.05cm]  (B1)--(W2)--(B4)--(W3)--(B3)--(W1)--(B1);
\draw[line width=0.05cm]  (W2)--(B2)--(3.5,1.5);\draw[line width=0.05cm] (W2)--(B3);
\draw[line width=0.05cm] (-0.5,1.5)--(W1);\draw[line width=0.05cm] (B4)--(W4)--(3.5,3.5); \draw[line width=0.05cm] (B3)--(-0.5,3.5);
\draw[line width=0.05cm] (B1)--(1,0.5); \draw[line width=0.05cm] (B2)--(3,0.5); 
\draw[line width=0.05cm] (W3)--(1,4.5); \draw[line width=0.05cm] (W4)--(3,4.5); 
%black
\filldraw  [ultra thick, fill=black] (1,1) circle [radius=0.16] ;\filldraw  [ultra thick, fill=black] (3,1) circle [radius=0.16] ;
\filldraw  [ultra thick, fill=black] (0,3) circle [radius=0.16] ;\filldraw  [ultra thick, fill=black] (2,3) circle [radius=0.16] ;
%white
\draw  [line width=1.28pt,fill=white] (0,2) circle [radius=0.16] ;\draw [line width=1.28pt, fill=white] (2,2) circle [radius=0.16] ;
\draw  [line width=1.28pt, fill=white] (1,4) circle [radius=0.16] ;\draw [line width=1.28pt, fill=white] (3,4) circle [radius=0.16] ;
\end{tikzpicture} }}; 

\node (PM1) at (3.1,0) 
{\scalebox{0.51}{
\begin{tikzpicture}
%vertex
\node (B1) at (1,1){$$}; \node (B2) at (3,1){$$}; \node (W1) at (0,2){$$}; \node (W2) at (2,2){$$};
\node (B3) at (0,3){$$}; \node (B4) at (2,3){$$}; \node (W3) at (1,4){$$}; \node (W4) at (3,4){$$};

\draw[line width=0.024cm]  (-0.5,0.5) rectangle (3.5,4.5);
%perfect matching
\draw[line width=0.4cm,color=lightgray] (W1)--(B3); \draw[line width=0.4cm,color=lightgray] (W2)--(B4);
\draw[line width=0.4cm,color=lightgray] (B2)--(3,0.5); \draw[line width=0.4cm,color=lightgray] (3,4.5)--(W4); 
\draw[line width=0.4cm,color=lightgray] (B1)--(1,0.5); \draw[line width=0.4cm,color=lightgray] (1,4.5)--(W3); 

%edge
\draw[line width=0.05cm]  (B1)--(W2)--(B4)--(W3)--(B3)--(W1)--(B1);
\draw[line width=0.05cm]  (W2)--(B2)--(3.5,1.5);\draw[line width=0.05cm] (W2)--(B3);
\draw[line width=0.05cm] (-0.5,1.5)--(W1);\draw[line width=0.05cm] (B4)--(W4)--(3.5,3.5); \draw[line width=0.05cm] (B3)--(-0.5,3.5);
\draw[line width=0.05cm] (B1)--(1,0.5); \draw[line width=0.05cm] (B2)--(3,0.5); 
\draw[line width=0.05cm] (W3)--(1,4.5); \draw[line width=0.05cm] (W4)--(3,4.5); 
%black
\filldraw  [ultra thick, fill=black] (1,1) circle [radius=0.16] ;\filldraw  [ultra thick, fill=black] (3,1) circle [radius=0.16] ;
\filldraw  [ultra thick, fill=black] (0,3) circle [radius=0.16] ;\filldraw  [ultra thick, fill=black] (2,3) circle [radius=0.16] ;
%white
\draw  [line width=1.28pt,fill=white] (0,2) circle [radius=0.16] ;\draw [line width=1.28pt, fill=white] (2,2) circle [radius=0.16] ;
\draw  [line width=1.28pt, fill=white] (1,4) circle [radius=0.16] ;\draw [line width=1.28pt, fill=white] (3,4) circle [radius=0.16] ;
\end{tikzpicture} }}; 

\node (PM2) at (6.2,0) 
{\scalebox{0.51}{
\begin{tikzpicture}
%vertex
\node (B1) at (1,1){$$}; \node (B2) at (3,1){$$}; \node (W1) at (0,2){$$}; \node (W2) at (2,2){$$};
\node (B3) at (0,3){$$}; \node (B4) at (2,3){$$}; \node (W3) at (1,4){$$}; \node (W4) at (3,4){$$};

\draw[line width=0.024cm]  (-0.5,0.5) rectangle (3.5,4.5);
%perfect matching
\draw[line width=0.4cm,color=lightgray] (W2)--(B3); \draw[line width=0.4cm,color=lightgray] (W4)--(B4);
\draw[line width=0.4cm,color=lightgray] (B1)--(1,0.5); \draw[line width=0.4cm,color=lightgray] (1,4.5)--(W3); 
\draw[line width=0.4cm,color=lightgray] (B2)--(3.5,1.5); \draw[line width=0.4cm,color=lightgray] (W1)--(-0.5,1.5);

%edge
\draw[line width=0.05cm]  (B1)--(W2)--(B4)--(W3)--(B3)--(W1)--(B1);
\draw[line width=0.05cm]  (W2)--(B2)--(3.5,1.5);\draw[line width=0.05cm] (W2)--(B3);
\draw[line width=0.05cm] (-0.5,1.5)--(W1);\draw[line width=0.05cm] (B4)--(W4)--(3.5,3.5); \draw[line width=0.05cm] (B3)--(-0.5,3.5);
\draw[line width=0.05cm] (B1)--(1,0.5); \draw[line width=0.05cm] (B2)--(3,0.5); 
\draw[line width=0.05cm] (W3)--(1,4.5); \draw[line width=0.05cm] (W4)--(3,4.5); 
%black
\filldraw  [ultra thick, fill=black] (1,1) circle [radius=0.16] ;\filldraw  [ultra thick, fill=black] (3,1) circle [radius=0.16] ;
\filldraw  [ultra thick, fill=black] (0,3) circle [radius=0.16] ;\filldraw  [ultra thick, fill=black] (2,3) circle [radius=0.16] ;
%white
\draw  [line width=1.28pt,fill=white] (0,2) circle [radius=0.16] ;\draw [line width=1.28pt, fill=white] (2,2) circle [radius=0.16] ;
\draw  [line width=1.28pt, fill=white] (1,4) circle [radius=0.16] ;\draw [line width=1.28pt, fill=white] (3,4) circle [radius=0.16] ;
\end{tikzpicture} }} ;  

\node (PM3) at (9.3,0) 
{\scalebox{0.51}{
\begin{tikzpicture}
%vertex
\node (B1) at (1,1){$$}; \node (B2) at (3,1){$$}; \node (W1) at (0,2){$$}; \node (W2) at (2,2){$$};
\node (B3) at (0,3){$$}; \node (B4) at (2,3){$$}; \node (W3) at (1,4){$$}; \node (W4) at (3,4){$$};

\draw[line width=0.024cm]  (-0.5,0.5) rectangle (3.5,4.5);
%perfect matching
\draw[line width=0.4cm,color=lightgray] (W2)--(B1); \draw[line width=0.4cm,color=lightgray] (W3)--(B3);
\draw[line width=0.4cm,color=lightgray] (W4)--(B4); \draw[line width=0.4cm,color=lightgray] (W1)--(-0.5,1.5); \draw[line width=0.4cm,color=lightgray] (B2)--(3.5,1.5);  

%edge
\draw[line width=0.05cm]  (B1)--(W2)--(B4)--(W3)--(B3)--(W1)--(B1);
\draw[line width=0.05cm]  (W2)--(B2)--(3.5,1.5);\draw[line width=0.05cm] (W2)--(B3);
\draw[line width=0.05cm] (-0.5,1.5)--(W1);\draw[line width=0.05cm] (B4)--(W4)--(3.5,3.5); \draw[line width=0.05cm] (B3)--(-0.5,3.5);
\draw[line width=0.05cm] (B1)--(1,0.5); \draw[line width=0.05cm] (B2)--(3,0.5); 
\draw[line width=0.05cm] (W3)--(1,4.5); \draw[line width=0.05cm] (W4)--(3,4.5); 
%black
\filldraw  [ultra thick, fill=black] (1,1) circle [radius=0.16] ;\filldraw  [ultra thick, fill=black] (3,1) circle [radius=0.16] ;
\filldraw  [ultra thick, fill=black] (0,3) circle [radius=0.16] ;\filldraw  [ultra thick, fill=black] (2,3) circle [radius=0.16] ;
%white
\draw  [line width=1.28pt,fill=white] (0,2) circle [radius=0.16] ;\draw [line width=1.28pt, fill=white] (2,2) circle [radius=0.16] ;
\draw  [line width=1.28pt, fill=white] (1,4) circle [radius=0.16] ;\draw [line width=1.28pt, fill=white] (3,4) circle [radius=0.16] ;
\end{tikzpicture} }}; 

\node (PM4) at (12.4,0) 
{\scalebox{0.51}{
\begin{tikzpicture}
%vertex
\node (B1) at (1,1){$$}; \node (B2) at (3,1){$$}; \node (W1) at (0,2){$$}; \node (W2) at (2,2){$$};
\node (B3) at (0,3){$$}; \node (B4) at (2,3){$$}; \node (W3) at (1,4){$$}; \node (W4) at (3,4){$$};

\draw[line width=0.024cm]  (-0.5,0.5) rectangle (3.5,4.5);
%perfect matching
\draw[line width=0.4cm,color=lightgray] (W1)--(B1); \draw[line width=0.4cm,color=lightgray] (W2)--(B2); 
\draw[line width=0.4cm,color=lightgray] (W3)--(B4); \draw[line width=0.4cm,color=lightgray] (W4)--(3.5,3.5);\draw[line width=0.4cm,color=lightgray] (-0.5,3.5)--(B3);   
%edge
\draw[line width=0.05cm]  (B1)--(W2)--(B4)--(W3)--(B3)--(W1)--(B1);
\draw[line width=0.05cm]  (W2)--(B2)--(3.5,1.5);\draw[line width=0.05cm] (W2)--(B3);
\draw[line width=0.05cm] (-0.5,1.5)--(W1);\draw[line width=0.05cm] (B4)--(W4)--(3.5,3.5); \draw[line width=0.05cm] (B3)--(-0.5,3.5);
\draw[line width=0.05cm] (B1)--(1,0.5); \draw[line width=0.05cm] (B2)--(3,0.5); 
\draw[line width=0.05cm] (W3)--(1,4.5); \draw[line width=0.05cm] (W4)--(3,4.5); 
%black
\filldraw  [ultra thick, fill=black] (1,1) circle [radius=0.16] ;\filldraw  [ultra thick, fill=black] (3,1) circle [radius=0.16] ;
\filldraw  [ultra thick, fill=black] (0,3) circle [radius=0.16] ;\filldraw  [ultra thick, fill=black] (2,3) circle [radius=0.16] ;
%white
\draw  [line width=1.28pt,fill=white] (0,2) circle [radius=0.16] ;\draw [line width=1.28pt, fill=white] (2,2) circle [radius=0.16] ;
\draw  [line width=1.28pt, fill=white] (1,4) circle [radius=0.16] ;\draw [line width=1.28pt, fill=white] (3,4) circle [radius=0.16] ;
\end{tikzpicture} }} ;

\end{tikzpicture}
}}
\end{center}

\medskip

Then, we have the following table, and 
$(e_0\sfA)^*\cong e_1\sfA$, $(e_2\sfA)^*\cong e_2\sfA$, $(e_3\sfA)^*\cong e_4\sfA$. 
%For simplicity, we denote by $T(a,b)$ a divisorial ideal $T(0,0,0,0,d)$. 

\medskip

\begin{center}
{\footnotesize{
\begin{tabular}{|l|l|l|l|l|} \hline 
$e_0\sfA$&$e_1\sfA$&$e_2\sfA$&$e_3\sfA$&$e_4\sfA$ \\ \hline
$T_{00}\cong R$&$T_{10}\cong T(0,0,0,-4)$&$T_{20}\cong T(0,0,0,-2)$&$T_{30}\cong T(0,0,0,-1)$&$T_{40}\cong T(0,0,0,-3)$ \\
$T_{01}\cong T(0,0,0,4)$&$T_{11}\cong R$&$T_{21}\cong T(0,0,0,2)$&$T_{31}\cong T(0,0,0,3)$&$T_{41}\cong T(0,0,0,1)$ \\
$T_{02}\cong T(0,0,0,2)$&$T_{12}\cong T(0,0,0,-2)$&$T_{22}\cong R$&$T_{32}\cong T(0,0,0,1)$&$T_{42}\cong T(0,0,0,-1)$ \\
$T_{03}\cong T(0,0,0,1)$&$T_{13}\cong T(0,0,0,-3)$&$T_{23}\cong T(0,0,0,-1)$&$T_{33}\cong R$&$T_{43}\cong T(0,0,0,-2)$\\ 
$T_{04}\cong T(0,0,0,3)$&$T_{14}\cong T(0,0,0,-1)$&$T_{24}\cong T(0,0,0,1)$&$T_{34}\cong T(0,0,0,2)$&$T_{44}\cong R$\\ \hline 
\end{tabular} 
}}
\end{center}

\medskip

\subsubsection{Exchange graph of type 5b} 
By the above results, we describe the exchange graph $\EG(\TMMG(R))$ for the case of type 5b as follows, 
and it is connected. 

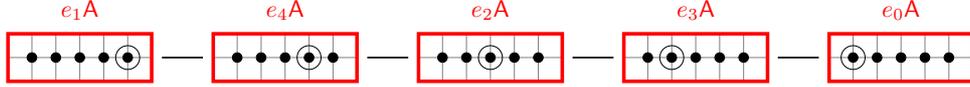
\begin{center}
{\scalebox{0.9}{
\begin{tikzpicture}

\node [red]  at (0,0.7) {$e_1\sfA$}; 
\node (A1) at (0,0)
{\scalebox{0.35}{\begin{tikzpicture}
\draw [step=1,thin, gray] (0,-1) grid (6,1); 
\draw [red, line width=0.15cm] (0,-1) rectangle (6,1);
\draw (5,0) circle [ultra thick, radius=0.5]; 
\filldraw [thick, fill=black] (5,0) circle [radius=0.2]; 
\filldraw  [thick, fill=black] (4,0) circle [radius=0.2] ;
\filldraw  [thick, fill=black] (3,0) circle [radius=0.2] ;
\filldraw  [thick, fill=black] (2,0) circle [radius=0.2] ;
\filldraw  [thick, fill=black] (1,0) circle [radius=0.2] ;
\end{tikzpicture}}};

\node [red]  at (3,0.7) {$e_4\sfA$}; 
\node (A4) at (3,0)
{\scalebox{0.35}{\begin{tikzpicture}
\draw [step=1,thin, gray] (0,-1) grid (6,1); 
\draw [red, line width=0.15cm] (0,-1) rectangle (6,1);
\draw (4,0) circle [ultra thick, radius=0.5]; 
\filldraw [thick, fill=black] (5,0) circle [radius=0.2]; 
\filldraw  [thick, fill=black] (4,0) circle [radius=0.2] ;
\filldraw  [thick, fill=black] (3,0) circle [radius=0.2] ;
\filldraw  [thick, fill=black] (2,0) circle [radius=0.2] ;
\filldraw  [thick, fill=black] (1,0) circle [radius=0.2] ;
\end{tikzpicture}}};

\node [red]  at (6,0.7) {$e_2\sfA$}; 
\node (A2) at (6,0)
{\scalebox{0.35}{\begin{tikzpicture}
\draw [step=1,thin, gray] (0,-1) grid (6,1); 
\draw [red, line width=0.15cm] (0,-1) rectangle (6,1);
\draw (3,0) circle [ultra thick, radius=0.5]; 
\filldraw [thick, fill=black] (5,0) circle [radius=0.2]; 
\filldraw  [thick, fill=black] (4,0) circle [radius=0.2] ;
\filldraw  [thick, fill=black] (3,0) circle [radius=0.2] ;
\filldraw  [thick, fill=black] (2,0) circle [radius=0.2] ;
\filldraw  [thick, fill=black] (1,0) circle [radius=0.2] ;
\end{tikzpicture}}};

\node [red]  at (9,0.7) {$e_3\sfA$}; 
\node (A3) at (9,0)
{\scalebox{0.35}{\begin{tikzpicture}
\draw [step=1,thin, gray] (0,-1) grid (6,1); 
\draw [red, line width=0.15cm] (0,-1) rectangle (6,1);
\draw (2,0) circle [ultra thick, radius=0.5]; 
\filldraw [thick, fill=black] (5,0) circle [radius=0.2]; 
\filldraw  [thick, fill=black] (4,0) circle [radius=0.2] ;
\filldraw  [thick, fill=black] (3,0) circle [radius=0.2] ;
\filldraw  [thick, fill=black] (2,0) circle [radius=0.2] ;
\filldraw  [thick, fill=black] (1,0) circle [radius=0.2] ;
\end{tikzpicture}}};

\node [red]  at (12,0.7) {$e_0\sfA$}; 
\node (A0) at (12,0)
{\scalebox{0.35}{\begin{tikzpicture}
\draw [step=1,thin, gray] (0,-1) grid (6,1); 
\draw [red, line width=0.15cm] (0,-1) rectangle (6,1);
\draw (1,0) circle [ultra thick, radius=0.5]; 
\filldraw [thick, fill=black] (5,0) circle [radius=0.2]; 
\filldraw  [thick, fill=black] (4,0) circle [radius=0.2] ;
\filldraw  [thick, fill=black] (3,0) circle [radius=0.2] ;
\filldraw  [thick, fill=black] (2,0) circle [radius=0.2] ;
\filldraw  [thick, fill=black] (1,0) circle [radius=0.2] ;
\end{tikzpicture}}};

\draw[line width=0.025cm] (A1)--(A4)--(A2)--(A3)--(A0) ;

\end{tikzpicture} }}
\end{center}

%%%%%%%%%%%%%%%%%%%%%%%%%%%%%%%%%%%%%%%%%%%%%%%%%%%%%%%
\subsection{Type 6a}
\label{type6a}

In this subsection, we consider the reflexive polygon of type 6a. 
Thus, let $R$ be the three dimensional complete local Gorenstein toric singularity defined by the cone $\sigma$: 
\[
\sigma=\mathrm{Cone}\{v_1=(1,0,1), v_2=(0,1,1), v_3=(-1,1,1), v_4=(-1,0,1), v_5=(0,-1,1), v_6=(1,-1,1) \}. 
\]
As an element in $\Cl(R)$, we obtain $[I_1]+2[I_2]+2[I_3]+[I_4]=0$, $2[I_1]+3[I_2]+2[I_3]-[I_5]=0$, $2[I_1]+2[I_2]+[I_3]+[I_6]=0$.  
Therefore, we have that $\Cl(R)\cong\ZZ^3$, and each divisorial ideal is represented by $T(a,b,c,0,0,0)$ where $a, b, c\in\ZZ$. 
Also, there are four consistent dimer models (up to right equivalence) written below which give the reflexive polygon of type 6a as the perfect matching polygon.

%%%%%%%%%%%%%%%%%%%%%%%%%%%%%%%%%%%
%%%%%%%%%%%%Type_6a-1%%%%%%%%%%%%%%%%%%
%%%%%%%%%%%%%%%%%%%%%%%%%%%%%%%%%%%
\subsubsection{Type 6a-1}

First, we consider the following consistent dimer model and the associated QP $(Q, W_Q)$. 
For simplicity, we denote by $\sfA$ the complete Jacobian algebra $\calP(Q, W_Q)$. 
Then, we have that $\sfA\cong\End_R(\bigoplus_{j\in Q_0}T_{ij})$ and $e_i\sfA\cong\bigoplus_{j\in Q_0}T_{ij}$ for any $i\in Q_0$ by Theorem~\ref{NCCR1}. 
We will give these splitting MM generators $e_i\sfA\cong\bigoplus_{j\in Q_0}T_{ij}$.

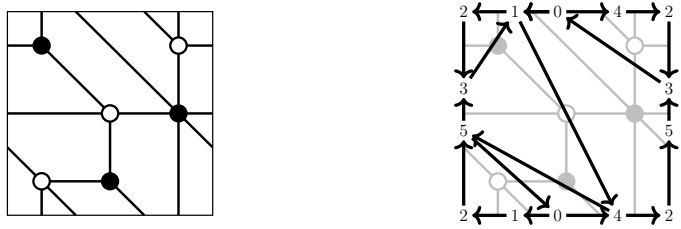
\begin{center}

\end{center}

\medskip

In the following figure, we fix the perfect matching $\sfP_0$.  
Then $\sfP_1, \cdots, \sfP_6$ are extremal perfect matchings corresponding to $v_1, \cdots, v_6$. 
Namely, $v_i=(h(\sfP_i, \sfP_0), 1)\in\ZZ^3$. 

\medskip
%%%%%perfect_matching%%%%%

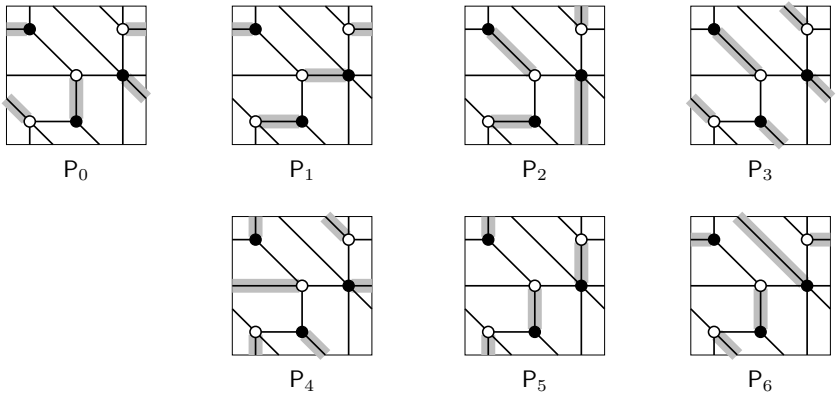
\begin{center}
{\scalebox{0.9}{
%
}}
\end{center}

\medskip

Then, we have the following table. We denote by $T(a,b,c)$ a divisorial ideal $T(a,b,c,0,0,0)$. 

\medskip

\begin{center}
{\small{
\begin{tabular}{|l|l|l|} \hline 
$e_0\sfA$&$e_1\sfA$&$e_2\sfA$\\ \hline
$T_{00}\cong R$&$T_{10}\cong T(2,2,1)$&$T_{20}\cong T(1,1,1)$ \\
$T_{01}\cong T(-2,-2,-1)$&$T_{11}\cong R$&$T_{21}\cong T(-1,-1,0)$ \\
$T_{02}\cong T(-1,-1,-1)$&$T_{12}\cong T(1,1,0)$&$T_{22}\cong R$ \\
$T_{03}\cong T(-2,-3,-2)$&$T_{13}\cong T(0,-1,-1)$&$T_{23}\cong T(-1,-2,-1)$ \\
$T_{04}\cong T(-1,-2,-1)$&$T_{14}\cong T(1,0,0)$&$T_{24}\cong T(0,-1,0)$ \\
$T_{05}\cong T(-1,-1,0)$&$T_{15}\cong T(1,1,1)$&$T_{25}\cong T(0,0,1)$ \\ \hline
$e_3\sfA$&$e_4\sfA$&$e_5\sfA$ \\ \hline
$T_{30}\cong T(2,3,2)$&$T_{40}\cong T(1,2,1)$&$T_{50}\cong T(1,1,0)$ \\
$T_{31}\cong T(0,1,1)$&$T_{41}\cong T(-1,0,0)$&$T_{51}\cong T(-1,-1,-1)$ \\
$T_{32}\cong T(1,2,1)$&$T_{42}\cong T(0,1,0)$&$T_{52}\cong T(0,0,-1)$ \\
$T_{33}\cong R$&$T_{43}\cong T(-1,-1,-1)$&$T_{53}\cong T(-1,-2,-2)$\\ 
$T_{34}\cong T(1,1,1)$&$T_{44}\cong R$&$T_{54}\cong T(0,-1,-1)$\\
$T_{35}\cong T(1,2,2)$&$T_{45}\cong T(0,1,1)$&$T_{55}\cong R$\\ \hline 
\end{tabular} 
}}
\end{center}

\medskip

%%%%%%%%%%%%%%%%%%%%%%%%%%%%%%%%%%%
%%%%%%%%%%%%Type_6a-2%%%%%%%%%%%%%%%%%%
%%%%%%%%%%%%%%%%%%%%%%%%%%%%%%%%%%%
\subsubsection{Type 6a-2}
Next, we consider the following consistent dimer model and the associated QP $(Q, W_Q)$. 
We denote by $\sfB$ the complete Jacobian algebra $\calP(Q, W_Q)$, and consider splitting MM generators $e_i\sfB$.

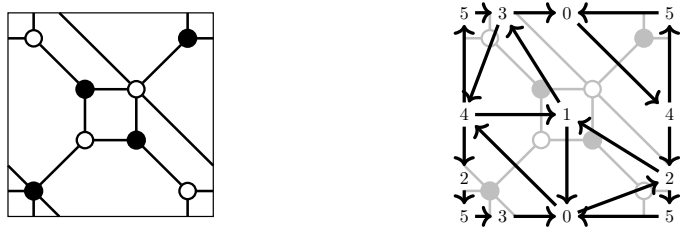
\begin{center}
%
\end{center}

\medskip

In the following figure, we fix the perfect matching $\sfP_0$.  
Then, $\sfP_1, \cdots, \sfP_6$ are extremal perfect matchings corresponding to $v_1, \cdots, v_6$. 

\medskip

%%%%%perfect_matching%%%%%

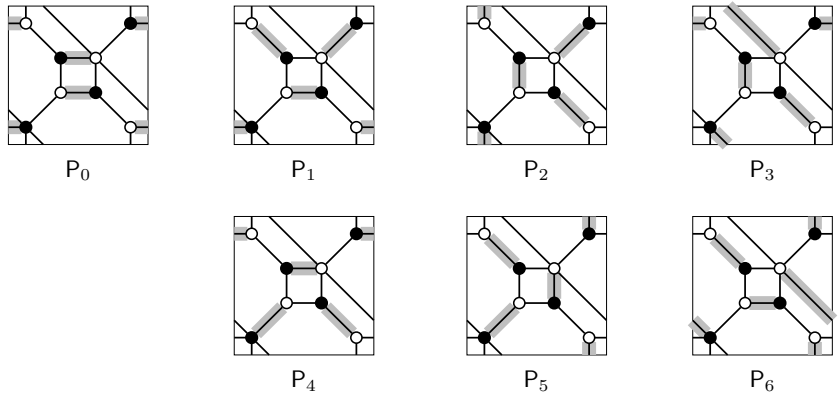
\begin{center}
{\scalebox{0.9}{
%
}}
\end{center}

\medskip

Then, we have the following table. We denote by $T(a,b,c)$ a divisorial ideal $T(a,b,c,0,0,0)$. 

\medskip

\begin{center}
{\small{
\begin{tabular}{|l|l|l|} \hline 
$e_0\sfB$&$e_1\sfB$&$e_2\sfB$ \\ \hline
$T_{00}\cong R$&$T_{10}\cong T(-1,-2,-1)$&$T_{20}\cong T(1,1,1)$ \\
$T_{01}\cong T(1,2,1)$&$T_{11}\cong R$&$T_{21}\cong T(2,3,2)$ \\
$T_{02}\cong T(-1,-1,-1)$&$T_{12}\cong T(-2,-3,-2)$&$T_{22}\cong R$ \\
$T_{03}\cong T(0,0,-1)$&$T_{13}\cong T(-1,-2,-2)$&$T_{23}\cong T(1,1,0)$ \\
$T_{04}\cong T(1,1,0)$&$T_{14}\cong T(0,-1,-1)$&$T_{24}\cong T(2,2,1)$ \\
$T_{05}\cong T(0,-1,-1)$&$T_{15}\cong T(-1,-3,-2)$&$T_{25}\cong T(1,0,0)$ \\ \hline
$e_3\sfB$&$e_4\sfB$&$e_5\sfB$ \\ \hline
$T_{30}\cong T(0,0,1)$&$T_{40}\cong T(-1,-1,0)$&$T_{50}\cong T(0,1,1)$ \\
$T_{31}\cong T(1,2,2)$&$T_{41}\cong T(0,1,1)$&$T_{51}\cong T(1,3,2)$ \\
$T_{32}\cong T(-1,-1,0)$&$T_{42}\cong T(-2,-2,-1)$&$T_{52}\cong T(-1,0,0)$ \\
$T_{33}\cong R$&$T_{43}\cong T(-1,-1,-1)$&$T_{53}\cong T(0,1,0)$\\ 
$T_{34}\cong T(1,1,1)$&$T_{44}\cong R$&$T_{54}\cong T(1,2,1)$\\
$T_{35}\cong T(0,-1,0)$&$T_{45}\cong T(-1,-2,-1)$&$T_{55}\cong R$\\ \hline 
\end{tabular} 
}}
\end{center}

\medskip

%%%%%%%%%%%%%%%%%%%%%%%%%%%%%%%%%%%
%%%%%%%%%%%%Type_6a-2(2)%%%%%%%%%%%%%%%%%%
%%%%%%%%%%%%%%%%%%%%%%%%%%%%%%%%%%%

Let $x$ (resp. $y$) be a $1$-cycle on the two torus $\sfT$ which corresponds to $(1,0)\in\rmH_1(\sfT)$ (resp. $(0,1)\in\rmH_1(\sfT)$) in the above case. 
By replacing these cycles by $y, -(x+y)$ respectively, we have the following dimer model and the associated QP $(Q, W_Q)$. 
We denote by $\sfB^\prime$ the complete Jacobian algebra $\calP(Q, W_Q)$. 
Since this dimer model is right equivalent to the previous one, we have that $\sfB\cong\sfB^\prime$ as an $R$-algebra. 
However, a change of cycles induces an automorphism on $R$, and it gives another description of each splitting MM generator (see Remark~\ref{rem_auto}). 
Thus we give splitting MM generators $e_i\sfB^\prime$.

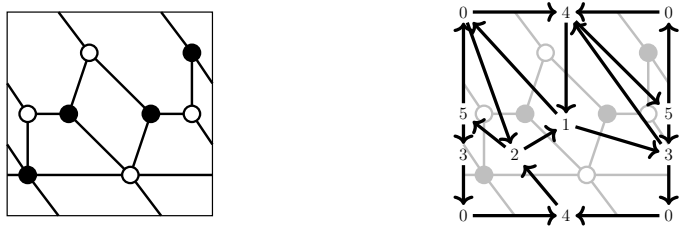
\medskip

\begin{center}
\begin{tikzpicture}
\node (DM) at (0,0) 
{\scalebox{0.54}{
\begin{tikzpicture}
%vertex
\node (B1) at (0.5,1){$$}; \node (B2) at (1.5,2.5){$$}; \node (B3) at (3.5,2.5){$$}; \node (B4) at (4.5,4){$$};  
\node (W1) at (0.5,2.5){$$}; \node (W2) at (2,4){$$}; \node (W3) at (3,1){$$}; \node (W4) at (4.5,2.5){$$}; 
\draw[line width=0.03cm]  (0,0) rectangle (5,5);
%edge
\draw[line width=0.06cm]  (B1)--(W1)--(B2)--(W2)--(B3)--(W3)--(B1) ; 
\draw[line width=0.06cm]  (W3)--(B2); \draw[line width=0.06cm]  (W4)--(B3); \draw[line width=0.06cm]  (W4)--(B4); 
\draw[line width=0.06cm]  (B1)--(0,1); \draw[line width=0.06cm]  (W3)--(5,1); 
\draw[line width=0.06cm]  (B1)--(1.25,0); \draw[line width=0.06cm]  (W2)--(1.25,5); 
\draw[line width=0.06cm]  (B1)--(0,1.75); \draw[line width=0.06cm]  (W4)--(5,1.75);  
\draw[line width=0.06cm]  (B4)--(3.75,5); \draw[line width=0.06cm]  (W3)--(3.75,0); 
\draw[line width=0.06cm]  (B4)--(5,3.25); \draw[line width=0.06cm]  (W1)--(0,3.25);  
%black
\filldraw  [ultra thick, fill=black] (0.5,1) circle [radius=0.2] ; \filldraw  [ultra thick, fill=black] (1.5,2.5) circle [radius=0.2] ;
\filldraw  [ultra thick, fill=black] (3.5,2.5) circle [radius=0.2] ; \filldraw  [ultra thick, fill=black] (4.5,4) circle [radius=0.2] ;
%white
\draw [line width=1.6pt,fill=white] (0.5,2.5) circle [radius=0.2] ; \draw [line width=1.6pt, fill=white] (2,4) circle [radius=0.2] ;
\draw [line width=1.6pt,fill=white] (3,1)circle [radius=0.2] ; \draw [line width=1.6pt,fill=white] (4.5,2.5)circle [radius=0.2] ;
\end{tikzpicture}
} }; 

\node (QV) at (6,0) 
{\scalebox{0.54}{
\begin{tikzpicture}
%vertex
\node (B1) at (0.5,1){$$}; \node (B2) at (1.5,2.5){$$}; \node (B3) at (3.5,2.5){$$}; \node (B4) at (4.5,4){$$};  
\node (W1) at (0.5,2.5){$$}; \node (W2) at (2,4){$$}; \node (W3) at (3,1){$$}; \node (W4) at (4.5,2.5){$$}; 
%\draw[line width=0.03cm]  (0,0) rectangle (5,5);
\node (Q0a) at (0,0){{\Large$0$}}; \node (Q0b) at (5,0){{\Large$0$}}; \node (Q0c) at (5,5){{\Large$0$}}; \node (Q0d) at (0,5){{\Large$0$}}; 
\node (Q1) at (2.5,2.25){{\Large$1$}}; \node (Q2) at (1.25,1.5){{\Large$2$}}; \node (Q3a) at (0,1.5){{\Large$3$}}; \node (Q3b) at (5,1.5){{\Large$3$}};
\node (Q4a) at (2.5,0){{\Large$4$}}; \node (Q4b) at (2.5,5){{\Large$4$}}; 
\node (Q5a) at (0,2.5){{\Large$5$}}; \node (Q5b) at (5,2.5){{\Large$5$}};
%edge
\draw[lightgray, line width=0.06cm]  (B1)--(W1)--(B2)--(W2)--(B3)--(W3)--(B1) ; 
\draw[lightgray, line width=0.06cm]  (W3)--(B2); \draw[lightgray, line width=0.06cm]  (W4)--(B3); \draw[lightgray, line width=0.06cm]  (W4)--(B4); 
\draw[lightgray, line width=0.06cm]  (B1)--(0,1); \draw[lightgray, line width=0.06cm]  (W3)--(5,1); 
\draw[lightgray, line width=0.06cm]  (B1)--(1.25,0); \draw[lightgray, line width=0.06cm]  (W2)--(1.25,5); 
\draw[lightgray, line width=0.06cm]  (B1)--(0,1.75); \draw[lightgray, line width=0.06cm]  (W4)--(5,1.75);  
\draw[lightgray, line width=0.06cm]  (B4)--(3.75,5); \draw[lightgray, line width=0.06cm]  (W3)--(3.75,0); 
\draw[lightgray, line width=0.06cm]  (B4)--(5,3.25); \draw[lightgray, line width=0.06cm]  (W1)--(0,3.25);  
%black
\filldraw  [ultra thick, draw=lightgray, fill=lightgray] (0.5,1) circle [radius=0.2] ; \filldraw  [ultra thick, draw=lightgray, fill=lightgray] (1.5,2.5) circle [radius=0.2] ;
\filldraw  [ultra thick, draw=lightgray, fill=lightgray] (3.5,2.5) circle [radius=0.2] ; \filldraw  [ultra thick, draw=lightgray, fill=lightgray] (4.5,4) circle [radius=0.2] ;
%white
\draw  [line width=1.6pt,draw=lightgray, fill=white] (0.5,2.5) circle [radius=0.2] ; \draw  [line width=1.6pt, draw=lightgray, fill=white] (2,4) circle [radius=0.2] ;
\draw  [line width=1.6pt,draw=lightgray, fill=white] (3,1)circle [radius=0.2] ; \draw  [line width=1.6pt, draw=lightgray, fill=white] (4.5,2.5)circle [radius=0.2] ;
%arrow
\draw[->, line width=0.08cm] (Q0a)--(Q4a); \draw[->, line width=0.08cm] (Q0b)--(Q4a); \draw[->, line width=0.08cm] (Q0c)--(Q4b); 
\draw[->, line width=0.08cm] (Q0d)--(Q4b); \draw[->, line width=0.08cm] (Q0d)--(Q2);  \draw[->, line width=0.08cm] (Q1)--(Q0d); 
\draw[->, line width=0.08cm] (Q1)--(Q3b); \draw[->, line width=0.08cm] (Q2)--(Q1); \draw[->, line width=0.08cm] (Q2)--(Q5a); 
\draw[->, line width=0.08cm] (Q3a)--(Q0a); \draw[->, line width=0.08cm] (Q3b)--(Q0b); \draw[->, line width=0.08cm] (Q3b)--(Q4b); 
\draw[->, line width=0.08cm] (Q4a)--(Q2); \draw[->, line width=0.08cm] (Q4b)--(Q1); \draw[->, line width=0.08cm] (Q4b)--(Q5b); 
\draw[->, line width=0.08cm] (Q5a)--(Q3a); \draw[->, line width=0.08cm] (Q5a)--(Q0d); \draw[->, line width=0.08cm] (Q5b)--(Q3b); 
\draw[->, line width=0.08cm] (Q5b)--(Q0c);
\end{tikzpicture}
} }; 
\end{tikzpicture}
\end{center}

\medskip

In the following figure, we fix the perfect matching $\sfP_0$.  
Then $\sfP_1, \cdots, \sfP_6$ are extremal perfect matchings corresponding to $v_1, \cdots, v_6$. 

\medskip

%%%%%perfect_matching%%%%%

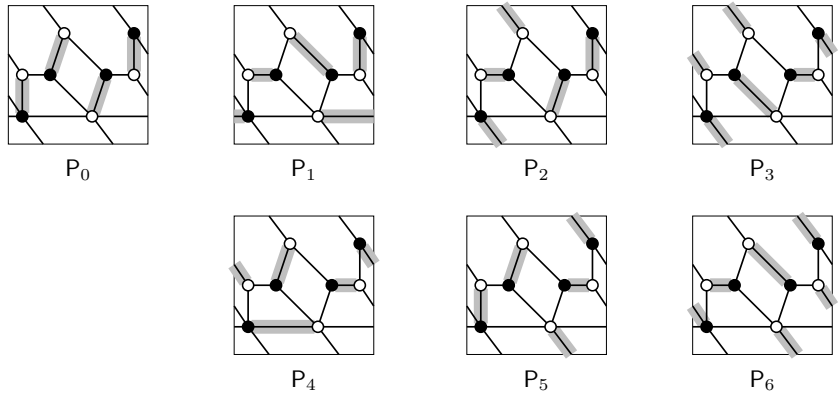
\begin{center}
{\scalebox{0.9}{
\begin{tikzpicture} 
\node at (0,-1.4) {$\sfP_0$};\node at (3.3,-1.4) {$\sfP_1$}; \node at (6.7,-1.4) {$\sfP_2$}; \node at (10,-1.4) {$\sfP_3$}; 
\node at (3.3,-4.5) {$\sfP_4$};\node at (6.7,-4.5) {$\sfP_5$}; \node at (10,-4.5) {$\sfP_6$}; 

\node (PM0) at (0,0) 
{\scalebox{0.408}{
\begin{tikzpicture}
%vertex
\node (B1) at (0.5,1){$$}; \node (B2) at (1.5,2.5){$$}; \node (B3) at (3.5,2.5){$$}; \node (B4) at (4.5,4){$$};  
\node (W1) at (0.5,2.5){$$}; \node (W2) at (2,4){$$}; \node (W3) at (3,1){$$}; \node (W4) at (4.5,2.5){$$}; 
\draw[line width=0.03cm]  (0,0) rectangle (5,5);
%perfectmatching
\draw[line width=0.5cm,color=lightgray] (W1)--(B1); \draw[line width=0.5cm,color=lightgray] (W2)--(B2); 
\draw[line width=0.5cm,color=lightgray] (W3)--(B3); \draw[line width=0.5cm,color=lightgray] (W4)--(B4); 
%edge
\draw[line width=0.06cm]  (B1)--(W1)--(B2)--(W2)--(B3)--(W3)--(B1) ; 
\draw[line width=0.06cm]  (W3)--(B2); \draw[line width=0.06cm]  (W4)--(B3); \draw[line width=0.06cm]  (W4)--(B4); 
\draw[line width=0.06cm]  (B1)--(0,1); \draw[line width=0.06cm]  (W3)--(5,1); 
\draw[line width=0.06cm]  (B1)--(1.25,0); \draw[line width=0.06cm]  (W2)--(1.25,5); 
\draw[line width=0.06cm]  (B1)--(0,1.75); \draw[line width=0.06cm]  (W4)--(5,1.75);  
\draw[line width=0.06cm]  (B4)--(3.75,5); \draw[line width=0.06cm]  (W3)--(3.75,0); 
\draw[line width=0.06cm]  (B4)--(5,3.25); \draw[line width=0.06cm]  (W1)--(0,3.25);  
%black
\filldraw  [ultra thick, fill=black] (0.5,1) circle [radius=0.2] ; \filldraw  [ultra thick, fill=black] (1.5,2.5) circle [radius=0.2] ;
\filldraw  [ultra thick, fill=black] (3.5,2.5) circle [radius=0.2] ; \filldraw  [ultra thick, fill=black] (4.5,4) circle [radius=0.2] ;
%white
\draw [line width=1.6pt,fill=white] (0.5,2.5) circle [radius=0.2] ; \draw [line width=1.6pt, fill=white] (2,4) circle [radius=0.2] ;
\draw [line width=1.6pt,fill=white] (3,1)circle [radius=0.2] ; \draw [line width=1.6pt,fill=white] (4.5,2.5)circle [radius=0.2] ;
\end{tikzpicture}
 }}; 

\node (PM1) at (3.3,0) 
{\scalebox{0.408}{
\begin{tikzpicture}
%vertex
\node (B1) at (0.5,1){$$}; \node (B2) at (1.5,2.5){$$}; \node (B3) at (3.5,2.5){$$}; \node (B4) at (4.5,4){$$};  
\node (W1) at (0.5,2.5){$$}; \node (W2) at (2,4){$$}; \node (W3) at (3,1){$$}; \node (W4) at (4.5,2.5){$$}; 
\draw[line width=0.03cm]  (0,0) rectangle (5,5);
%perfectmatching
\draw[line width=0.5cm,color=lightgray] (W1)--(B2); \draw[line width=0.5cm,color=lightgray] (W2)--(B3); 
\draw[line width=0.5cm,color=lightgray] (W4)--(B4); 
\draw[line width=0.5cm,color=lightgray] (B1)--(0,1); \draw[line width=0.5cm,color=lightgray] (W3)--(5,1);  
%edge
\draw[line width=0.06cm]  (B1)--(W1)--(B2)--(W2)--(B3)--(W3)--(B1) ; 
\draw[line width=0.06cm]  (W3)--(B2); \draw[line width=0.06cm]  (W4)--(B3); \draw[line width=0.06cm]  (W4)--(B4); 
\draw[line width=0.06cm]  (B1)--(0,1); \draw[line width=0.06cm]  (W3)--(5,1); 
\draw[line width=0.06cm]  (B1)--(1.25,0); \draw[line width=0.06cm]  (W2)--(1.25,5); 
\draw[line width=0.06cm]  (B1)--(0,1.75); \draw[line width=0.06cm]  (W4)--(5,1.75);  
\draw[line width=0.06cm]  (B4)--(3.75,5); \draw[line width=0.06cm]  (W3)--(3.75,0); 
\draw[line width=0.06cm]  (B4)--(5,3.25); \draw[line width=0.06cm]  (W1)--(0,3.25);  
%black
\filldraw  [ultra thick, fill=black] (0.5,1) circle [radius=0.2] ; \filldraw  [ultra thick, fill=black] (1.5,2.5) circle [radius=0.2] ;
\filldraw  [ultra thick, fill=black] (3.5,2.5) circle [radius=0.2] ; \filldraw  [ultra thick, fill=black] (4.5,4) circle [radius=0.2] ;
%white
\draw [line width=1.6pt,fill=white] (0.5,2.5) circle [radius=0.2] ; \draw [line width=1.6pt, fill=white] (2,4) circle [radius=0.2] ;
\draw [line width=1.6pt,fill=white] (3,1)circle [radius=0.2] ; \draw [line width=1.6pt,fill=white] (4.5,2.5)circle [radius=0.2] ;
\end{tikzpicture}
 }}; 

\node (PM2) at (6.7,0) 
{\scalebox{0.408}{
\begin{tikzpicture}
%vertex
\node (B1) at (0.5,1){$$}; \node (B2) at (1.5,2.5){$$}; \node (B3) at (3.5,2.5){$$}; \node (B4) at (4.5,4){$$};  
\node (W1) at (0.5,2.5){$$}; \node (W2) at (2,4){$$}; \node (W3) at (3,1){$$}; \node (W4) at (4.5,2.5){$$}; 
\draw[line width=0.03cm]  (0,0) rectangle (5,5);
%perfectmatching
\draw[line width=0.5cm,color=lightgray] (W1)--(B2); \draw[line width=0.5cm,color=lightgray] (W3)--(B3); 
\draw[line width=0.5cm,color=lightgray] (W4)--(B4); 
\draw[line width=0.5cm,color=lightgray] (B1)--(1.25,0); \draw[line width=0.5cm,color=lightgray] (W2)--(1.25,5); 
%edge
\draw[line width=0.06cm]  (B1)--(W1)--(B2)--(W2)--(B3)--(W3)--(B1) ; 
\draw[line width=0.06cm]  (W3)--(B2); \draw[line width=0.06cm]  (W4)--(B3); \draw[line width=0.06cm]  (W4)--(B4); 
\draw[line width=0.06cm]  (B1)--(0,1); \draw[line width=0.06cm]  (W3)--(5,1); 
\draw[line width=0.06cm]  (B1)--(1.25,0); \draw[line width=0.06cm]  (W2)--(1.25,5); 
\draw[line width=0.06cm]  (B1)--(0,1.75); \draw[line width=0.06cm]  (W4)--(5,1.75);  
\draw[line width=0.06cm]  (B4)--(3.75,5); \draw[line width=0.06cm]  (W3)--(3.75,0); 
\draw[line width=0.06cm]  (B4)--(5,3.25); \draw[line width=0.06cm]  (W1)--(0,3.25);  
%black
\filldraw  [ultra thick, fill=black] (0.5,1) circle [radius=0.2] ; \filldraw  [ultra thick, fill=black] (1.5,2.5) circle [radius=0.2] ;
\filldraw  [ultra thick, fill=black] (3.5,2.5) circle [radius=0.2] ; \filldraw  [ultra thick, fill=black] (4.5,4) circle [radius=0.2] ;
%white
\draw [line width=1.6pt,fill=white] (0.5,2.5) circle [radius=0.2] ; \draw [line width=1.6pt, fill=white] (2,4) circle [radius=0.2] ;
\draw [line width=1.6pt,fill=white] (3,1)circle [radius=0.2] ; \draw [line width=1.6pt,fill=white] (4.5,2.5)circle [radius=0.2] ;
\end{tikzpicture}
 }} ;  

\node (PM3) at (10,0) 
{\scalebox{0.408}{
\begin{tikzpicture}
%vertex
\node (B1) at (0.5,1){$$}; \node (B2) at (1.5,2.5){$$}; \node (B3) at (3.5,2.5){$$}; \node (B4) at (4.5,4){$$};  
\node (W1) at (0.5,2.5){$$}; \node (W2) at (2,4){$$}; \node (W3) at (3,1){$$}; \node (W4) at (4.5,2.5){$$}; 
\draw[line width=0.03cm]  (0,0) rectangle (5,5);
%perfectmatching
\draw[line width=0.5cm,color=lightgray] (W3)--(B2); \draw[line width=0.5cm,color=lightgray] (W4)--(B3); 
\draw[line width=0.5cm,color=lightgray] (B1)--(1.25,0); \draw[line width=0.5cm,color=lightgray] (W2)--(1.25,5); 
\draw[line width=0.5cm,color=lightgray] (B4)--(5,3.25); \draw[line width=0.5cm,color=lightgray] (W1)--(0,3.25); 
%edge
\draw[line width=0.06cm]  (B1)--(W1)--(B2)--(W2)--(B3)--(W3)--(B1) ; 
\draw[line width=0.06cm]  (W3)--(B2); \draw[line width=0.06cm]  (W4)--(B3); \draw[line width=0.06cm]  (W4)--(B4); 
\draw[line width=0.06cm]  (B1)--(0,1); \draw[line width=0.06cm]  (W3)--(5,1); 
\draw[line width=0.06cm]  (B1)--(1.25,0); \draw[line width=0.06cm]  (W2)--(1.25,5); 
\draw[line width=0.06cm]  (B1)--(0,1.75); \draw[line width=0.06cm]  (W4)--(5,1.75);  
\draw[line width=0.06cm]  (B4)--(3.75,5); \draw[line width=0.06cm]  (W3)--(3.75,0); 
\draw[line width=0.06cm]  (B4)--(5,3.25); \draw[line width=0.06cm]  (W1)--(0,3.25);  
%black
\filldraw  [ultra thick, fill=black] (0.5,1) circle [radius=0.2] ; \filldraw  [ultra thick, fill=black] (1.5,2.5) circle [radius=0.2] ;
\filldraw  [ultra thick, fill=black] (3.5,2.5) circle [radius=0.2] ; \filldraw  [ultra thick, fill=black] (4.5,4) circle [radius=0.2] ;
%white
\draw [line width=1.6pt,fill=white] (0.5,2.5) circle [radius=0.2] ; \draw [line width=1.6pt, fill=white] (2,4) circle [radius=0.2] ;
\draw [line width=1.6pt,fill=white] (3,1)circle [radius=0.2] ; \draw [line width=1.6pt,fill=white] (4.5,2.5)circle [radius=0.2] ;
\end{tikzpicture}
 }} ; 

\node (PM4) at (3.3,-3.1) 
{\scalebox{0.408}{
\begin{tikzpicture}
%vertex
\node (B1) at (0.5,1){$$}; \node (B2) at (1.5,2.5){$$}; \node (B3) at (3.5,2.5){$$}; \node (B4) at (4.5,4){$$};  
\node (W1) at (0.5,2.5){$$}; \node (W2) at (2,4){$$}; \node (W3) at (3,1){$$}; \node (W4) at (4.5,2.5){$$}; 
\draw[line width=0.03cm]  (0,0) rectangle (5,5); 
%perfectmatching
\draw[line width=0.5cm,color=lightgray] (W3)--(B1); \draw[line width=0.5cm,color=lightgray] (W2)--(B2); 
\draw[line width=0.5cm,color=lightgray] (W4)--(B3); 
\draw[line width=0.5cm,color=lightgray] (B4)--(5,3.25); \draw[line width=0.5cm,color=lightgray] (W1)--(0,3.25);
%edge
\draw[line width=0.06cm]  (B1)--(W1)--(B2)--(W2)--(B3)--(W3)--(B1) ; 
\draw[line width=0.06cm]  (W3)--(B2); \draw[line width=0.06cm]  (W4)--(B3); \draw[line width=0.06cm]  (W4)--(B4); 
\draw[line width=0.06cm]  (B1)--(0,1); \draw[line width=0.06cm]  (W3)--(5,1); 
\draw[line width=0.06cm]  (B1)--(1.25,0); \draw[line width=0.06cm]  (W2)--(1.25,5); 
\draw[line width=0.06cm]  (B1)--(0,1.75); \draw[line width=0.06cm]  (W4)--(5,1.75);  
\draw[line width=0.06cm]  (B4)--(3.75,5); \draw[line width=0.06cm]  (W3)--(3.75,0); 
\draw[line width=0.06cm]  (B4)--(5,3.25); \draw[line width=0.06cm]  (W1)--(0,3.25);  
%black
\filldraw  [ultra thick, fill=black] (0.5,1) circle [radius=0.2] ; \filldraw  [ultra thick, fill=black] (1.5,2.5) circle [radius=0.2] ;
\filldraw  [ultra thick, fill=black] (3.5,2.5) circle [radius=0.2] ; \filldraw  [ultra thick, fill=black] (4.5,4) circle [radius=0.2] ;
%white
\draw [line width=1.6pt,fill=white] (0.5,2.5) circle [radius=0.2] ; \draw [line width=1.6pt, fill=white] (2,4) circle [radius=0.2] ;
\draw [line width=1.6pt,fill=white] (3,1)circle [radius=0.2] ; \draw [line width=1.6pt,fill=white] (4.5,2.5)circle [radius=0.2] ;
\end{tikzpicture}
 }}; 

\node (PM5) at (6.7,-3.1) 
{\scalebox{0.408}{
\begin{tikzpicture}
%vertex
\node (B1) at (0.5,1){$$}; \node (B2) at (1.5,2.5){$$}; \node (B3) at (3.5,2.5){$$}; \node (B4) at (4.5,4){$$};  
\node (W1) at (0.5,2.5){$$}; \node (W2) at (2,4){$$}; \node (W3) at (3,1){$$}; \node (W4) at (4.5,2.5){$$}; 
\draw[line width=0.03cm]  (0,0) rectangle (5,5);
%perfectmatching
\draw[line width=0.5cm,color=lightgray] (W1)--(B1); \draw[line width=0.5cm,color=lightgray] (W2)--(B2); 
\draw[line width=0.5cm,color=lightgray] (W4)--(B3); 
\draw[line width=0.5cm,color=lightgray] (B4)--(3.75,5); \draw[line width=0.5cm,color=lightgray] (W3)--(3.75,0);
%edge
\draw[line width=0.06cm]  (B1)--(W1)--(B2)--(W2)--(B3)--(W3)--(B1) ; 
\draw[line width=0.06cm]  (W3)--(B2); \draw[line width=0.06cm]  (W4)--(B3); \draw[line width=0.06cm]  (W4)--(B4); 
\draw[line width=0.06cm]  (B1)--(0,1); \draw[line width=0.06cm]  (W3)--(5,1); 
\draw[line width=0.06cm]  (B1)--(1.25,0); \draw[line width=0.06cm]  (W2)--(1.25,5); 
\draw[line width=0.06cm]  (B1)--(0,1.75); \draw[line width=0.06cm]  (W4)--(5,1.75);  
\draw[line width=0.06cm]  (B4)--(3.75,5); \draw[line width=0.06cm]  (W3)--(3.75,0); 
\draw[line width=0.06cm]  (B4)--(5,3.25); \draw[line width=0.06cm]  (W1)--(0,3.25);  
%black
\filldraw  [ultra thick, fill=black] (0.5,1) circle [radius=0.2] ; \filldraw  [ultra thick, fill=black] (1.5,2.5) circle [radius=0.2] ;
\filldraw  [ultra thick, fill=black] (3.5,2.5) circle [radius=0.2] ; \filldraw  [ultra thick, fill=black] (4.5,4) circle [radius=0.2] ;
%white
\draw [line width=1.6pt,fill=white] (0.5,2.5) circle [radius=0.2] ; \draw [line width=1.6pt, fill=white] (2,4) circle [radius=0.2] ;
\draw [line width=1.6pt,fill=white] (3,1)circle [radius=0.2] ; \draw [line width=1.6pt,fill=white] (4.5,2.5)circle [radius=0.2] ;
\end{tikzpicture}
 }}; 

\node (PM6) at (10,-3.1) 
{\scalebox{0.408}{
\begin{tikzpicture}
%vertex
\node (B1) at (0.5,1){$$}; \node (B2) at (1.5,2.5){$$}; \node (B3) at (3.5,2.5){$$}; \node (B4) at (4.5,4){$$};  
\node (W1) at (0.5,2.5){$$}; \node (W2) at (2,4){$$}; \node (W3) at (3,1){$$}; \node (W4) at (4.5,2.5){$$}; 
\draw[line width=0.03cm]  (0,0) rectangle (5,5);
%perfectmatching
\draw[line width=0.5cm,color=lightgray] (W1)--(B2); \draw[line width=0.5cm,color=lightgray] (W2)--(B3);  
\draw[line width=0.5cm,color=lightgray] (B4)--(3.75,5); \draw[line width=0.5cm,color=lightgray] (W3)--(3.75,0);
\draw[line width=0.5cm,color=lightgray] (B1)--(0,1.75); \draw[line width=0.5cm,color=lightgray] (W4)--(5,1.75);
%edge
\draw[line width=0.06cm]  (B1)--(W1)--(B2)--(W2)--(B3)--(W3)--(B1) ; 
\draw[line width=0.06cm]  (W3)--(B2); \draw[line width=0.06cm]  (W4)--(B3); \draw[line width=0.06cm]  (W4)--(B4); 
\draw[line width=0.06cm]  (B1)--(0,1); \draw[line width=0.06cm]  (W3)--(5,1); 
\draw[line width=0.06cm]  (B1)--(1.25,0); \draw[line width=0.06cm]  (W2)--(1.25,5); 
\draw[line width=0.06cm]  (B1)--(0,1.75); \draw[line width=0.06cm]  (W4)--(5,1.75);  
\draw[line width=0.06cm]  (B4)--(3.75,5); \draw[line width=0.06cm]  (W3)--(3.75,0); 
\draw[line width=0.06cm]  (B4)--(5,3.25); \draw[line width=0.06cm]  (W1)--(0,3.25);  
%black
\filldraw  [ultra thick, fill=black] (0.5,1) circle [radius=0.2] ; \filldraw  [ultra thick, fill=black] (1.5,2.5) circle [radius=0.2] ;
\filldraw  [ultra thick, fill=black] (3.5,2.5) circle [radius=0.2] ; \filldraw  [ultra thick, fill=black] (4.5,4) circle [radius=0.2] ;
%white
\draw [line width=1.6pt,fill=white] (0.5,2.5) circle [radius=0.2] ; \draw [line width=1.6pt, fill=white] (2,4) circle [radius=0.2] ;
\draw [line width=1.6pt,fill=white] (3,1)circle [radius=0.2] ; \draw [line width=1.6pt,fill=white] (4.5,2.5)circle [radius=0.2] ;
\end{tikzpicture}
 }} ;  

\end{tikzpicture}
}}
\end{center}

\medskip

Then, we have the following table. We denote by $T(a,b,c)$ a divisorial ideal $T(a,b,c,0,0,0)$. 

\medskip

\begin{center}
{\small{
\begin{tabular}{|l|l|l|} \hline 
$e_0\sfB^\prime$&$e_1\sfB^\prime$&$e_2\sfB^\prime$ \\ \hline
$T_{00}\cong R$&$T_{10}\cong T(1,1,0)$&$T_{20}\cong T(1,1,1)$ \\
$T_{01}\cong T(-1,-1,0)$&$T_{11}\cong R$&$T_{21}\cong T(0,0,1)$ \\
$T_{02}\cong T(-1,-1,-1)$&$T_{12}\cong T(0,0,-1)$&$T_{22}\cong R$ \\
$T_{03}\cong T(-1,0,0)$&$T_{13}\cong T(0,1,0)$&$T_{23}\cong T(0,1,1)$ \\
$T_{04}\cong T(0,1,1)$&$T_{14}\cong T(1,2,1)$&$T_{24}\cong T(1,2,2)$ \\
$T_{05}\cong T(1,2,1)$&$T_{15}\cong T(2,3,1)$&$T_{25}\cong T(2,3,2)$ \\ \hline
$e_3\sfB^\prime$&$e_4\sfB^\prime$&$e_5\sfB^\prime$ \\ \hline
$T_{30}\cong T(1,0,0)$&$T_{40}\cong T(0,-1,-1)$&$T_{50}\cong T(-1,-2,-1)$ \\
$T_{31}\cong T(0,-1,0)$&$T_{41}\cong T(-1,-2,-1)$&$T_{51}\cong T(-2,-3,-1)$ \\
$T_{32}\cong T(0,-1,-1)$&$T_{42}\cong T(-1,-2,-2)$&$T_{52}\cong T(-2,-3,-2)$ \\
$T_{33}\cong R$&$T_{43}\cong T(-1,-1,-1)$&$T_{53}\cong T(-2,-2,-1)$\\ 
$T_{34}\cong T(1,1,1)$&$T_{44}\cong R$&$T_{54}\cong T(-1,-1,0)$\\
$T_{35}\cong T(2,2,1)$&$T_{45}\cong T(1,1,0)$&$T_{55}\cong R$\\ \hline 
\end{tabular} 
}}
\end{center}

\medskip

%%%%%%%%%%%%%%%%%%%%%%%%%%%%%%%%%%%
%%%%%%%%%%%%Type_6a-2(3)%%%%%%%%%%%%%%%%%%
%%%%%%%%%%%%%%%%%%%%%%%%%%%%%%%%%%%

Let $x$ (resp. $y$) be a $1$-cycle on the two torus $\sfT$ which corresponds to $(1,0)\in\rmH_1(\sfT)$ (resp. $(0,1)\in\rmH_1(\sfT)$) in the first case of Type 6a-2. 
By replacing these cycles by $x+y,-y$ respectively, we have the following dimer model and the associated QP $(Q, W_Q)$. 
We denote by $\sfB^{\prime\prime}$ the complete Jacobian algebra $\calP(Q, W_Q)$.  
By the same reason as the above case, we also consider splitting MM generators $e_i\sfB^{\prime\prime}$.

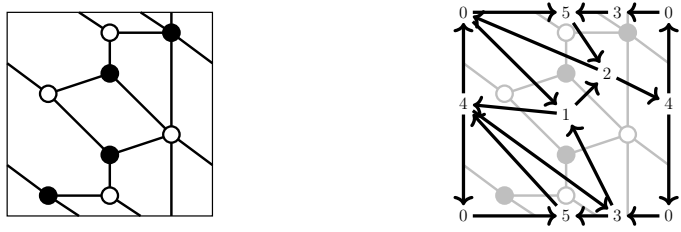
\medskip

\begin{center}
\begin{tikzpicture}
\node (DM) at (0,0) 
{\scalebox{0.54}{
\begin{tikzpicture}
%vertex
\node (B1) at (1,0.5){$$}; \node (B2) at (2.5,1.5){$$}; \node (B3) at (2.5,3.5){$$}; \node (B4) at (4,4.5){$$};  
\node (W1) at (2.5,0.5){$$}; \node (W2) at (4,2){$$}; \node (W3) at (1,3){$$}; \node (W4) at (2.5,4.5){$$}; 
\draw[line width=0.03cm]  (0,0) rectangle (5,5);
%edge
\draw[line width=0.06cm]  (B4)--(W4)--(B3)--(W3)--(B2)--(W2)--(B4) ; 
\draw[line width=0.06cm]  (W2)--(B3); \draw[line width=0.06cm]  (W1)--(B1); \draw[line width=0.06cm]  (W1)--(B2); 
\draw[line width=0.06cm]  (B1)--(0,1.25); \draw[line width=0.06cm]  (W2)--(5,1.25); 
\draw[line width=0.06cm]  (B1)--(1.75,0); \draw[line width=0.06cm]  (W4)--(1.75,5); 
\draw[line width=0.06cm]  (B4)--(3.25,5); \draw[line width=0.06cm]  (W1)--(3.25,0);  
\draw[line width=0.06cm]  (B4)--(4,5); \draw[line width=0.06cm]  (W2)--(4,0); 
\draw[line width=0.06cm]  (B4)--(5,3.75); \draw[line width=0.06cm]  (W3)--(0,3.75);  
%black
\filldraw  [ultra thick, fill=black] (1,0.5) circle [radius=0.2] ; \filldraw  [ultra thick, fill=black] (2.5,1.5) circle [radius=0.2] ;
\filldraw  [ultra thick, fill=black] (2.5,3.5) circle [radius=0.2] ; \filldraw  [ultra thick, fill=black] (4,4.5) circle [radius=0.2] ;
%white
\draw  [line width=1.6pt,fill=white] (2.5,0.5) circle [radius=0.2] ; \draw  [line width=1.6pt, fill=white] (4,2) circle [radius=0.2] ;
\draw  [line width=1.6pt,fill=white] (1,3)circle [radius=0.2] ; \draw  [line width=1.6pt,fill=white] (2.5,4.5)circle [radius=0.2] ;
\end{tikzpicture}
} }; 

\node (QV) at (6,0) 
{\scalebox{0.54}{
\begin{tikzpicture}
%vertex
\node (B1) at (1,0.5){$$}; \node (B2) at (2.5,1.5){$$}; \node (B3) at (2.5,3.5){$$}; \node (B4) at (4,4.5){$$};  
\node (W1) at (2.5,0.5){$$}; \node (W2) at (4,2){$$}; \node (W3) at (1,3){$$}; \node (W4) at (2.5,4.5){$$}; 
%\draw[line width=0.03cm]  (0,0) rectangle (5,5);
\node (Q0a) at (0,0){{\Large$0$}}; \node (Q0b) at (5,0){{\Large$0$}}; \node (Q0c) at (5,5){{\Large$0$}}; \node (Q0d) at (0,5){{\Large$0$}}; 
\node (Q1) at (2.5,2.5){{\Large$1$}}; \node (Q2) at (3.5,3.5){{\Large$2$}}; \node (Q3a) at (3.75,0){{\Large$3$}}; \node (Q3b) at (3.75,5){{\Large$3$}};
\node (Q4a) at (0,2.75){{\Large$4$}}; \node (Q4b) at (5,2.75){{\Large$4$}}; \node (Q5a) at (2.5,0){{\Large$5$}}; \node (Q5b) at (2.5,5){{\Large$5$}};
%edge
\draw[lightgray, line width=0.06cm]  (B4)--(W4)--(B3)--(W3)--(B2)--(W2)--(B4) ; 
\draw[lightgray, line width=0.06cm]  (W2)--(B3); \draw[lightgray, line width=0.06cm]  (W1)--(B1); \draw[lightgray, line width=0.06cm]  (W1)--(B2); 
\draw[lightgray, line width=0.06cm]  (B1)--(0,1.25); \draw[lightgray, line width=0.06cm]  (W2)--(5,1.25); 
\draw[lightgray, line width=0.06cm]  (B1)--(1.75,0); \draw[lightgray, line width=0.06cm]  (W4)--(1.75,5); 
\draw[lightgray, line width=0.06cm]  (B4)--(3.25,5); \draw[lightgray, line width=0.06cm]  (W1)--(3.25,0);  
\draw[lightgray, line width=0.06cm]  (B4)--(4,5); \draw[lightgray, line width=0.06cm]  (W2)--(4,0); 
\draw[lightgray, line width=0.06cm]  (B4)--(5,3.75); \draw[lightgray, line width=0.06cm]  (W3)--(0,3.75);  
%black
\filldraw  [ultra thick, draw=lightgray, fill=lightgray] (1,0.5) circle [radius=0.2] ; \filldraw  [ultra thick, draw=lightgray, fill=lightgray] (2.5,1.5) circle [radius=0.2] ;
\filldraw  [ultra thick, draw=lightgray, fill=lightgray] (2.5,3.5) circle [radius=0.2] ; \filldraw  [ultra thick, draw=lightgray, fill=lightgray] (4,4.5) circle [radius=0.2] ;
%white
\draw  [line width=1.6pt, draw=lightgray, fill=white] (2.5,0.5) circle [radius=0.2] ; \draw  [line width=1.6pt, draw=lightgray, fill=white] (4,2) circle [radius=0.2] ;
\draw  [line width=1.6pt, draw=lightgray, fill=white] (1,3)circle [radius=0.2] ; \draw  [line width=1.6pt, draw=lightgray, fill=white] (2.5,4.5)circle [radius=0.2] ;
%arrow
\draw[->, line width=0.08cm] (Q0a)--(Q5a); \draw[->, line width=0.08cm] (Q0b)--(Q3a); \draw[->, line width=0.08cm] (Q0c)--(Q3b); 
\draw[->, line width=0.08cm] (Q0d)--(Q5b); \draw[->, line width=0.08cm] (Q0d)--(Q1);  \draw[->, line width=0.08cm] (Q1)--(Q2); 
\draw[->, line width=0.08cm] (Q1)--(Q4a); \draw[->, line width=0.08cm] (Q2)--(Q0d); \draw[->, line width=0.08cm] (Q2)--(Q4b); 
\draw[->, line width=0.08cm] (Q3a)--(Q5a); \draw[->, line width=0.08cm] (Q3a)--(Q1); \draw[->, line width=0.08cm] (Q3b)--(Q5b); 
\draw[->, line width=0.08cm] (Q4a)--(Q0a); \draw[->, line width=0.08cm] (Q4a)--(Q0d); \draw[->, line width=0.08cm] (Q4a)--(Q3a); 
\draw[->, line width=0.08cm] (Q4b)--(Q0b); \draw[->, line width=0.08cm] (Q4b)--(Q0c); 
\draw[->, line width=0.08cm] (Q5a)--(Q4a); \draw[->, line width=0.08cm] (Q5b)--(Q2);
\end{tikzpicture}
} }; 
\end{tikzpicture}
\end{center}

\medskip

In the following figure, we fix the perfect matching $\sfP_0$.  
Then $\sfP_1, \cdots, \sfP_6$ are extremal perfect matchings corresponding to $v_1, \cdots, v_6$. 

\medskip

%%%%%perfect_matching%%%%%

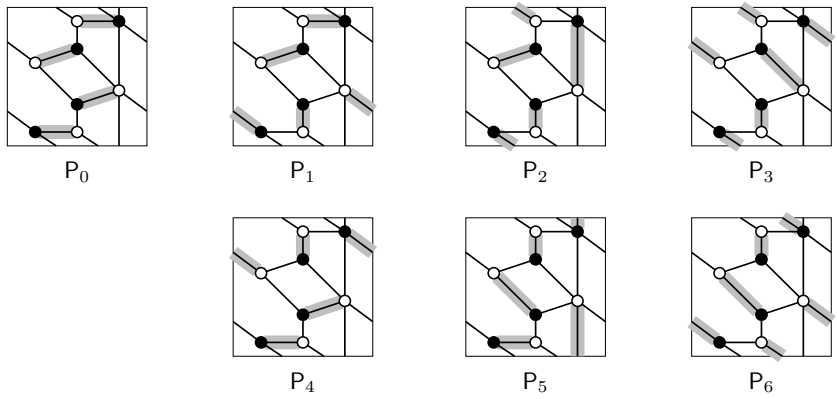
\begin{center}
{\scalebox{0.9}{
\begin{tikzpicture} 
\node at (0,-1.4) {$\sfP_0$};\node at (3.3,-1.4) {$\sfP_1$}; \node at (6.7,-1.4) {$\sfP_2$}; \node at (10,-1.4) {$\sfP_3$}; 
\node at (3.3,-4.5) {$\sfP_4$};\node at (6.7,-4.5) {$\sfP_5$}; \node at (10,-4.5) {$\sfP_6$};

\node (PM0) at (0,0) 
{\scalebox{0.408}{
\begin{tikzpicture}
%vertex
\node (B1) at (1,0.5){$$}; \node (B2) at (2.5,1.5){$$}; \node (B3) at (2.5,3.5){$$}; \node (B4) at (4,4.5){$$};  
\node (W1) at (2.5,0.5){$$}; \node (W2) at (4,2){$$}; \node (W3) at (1,3){$$}; \node (W4) at (2.5,4.5){$$}; 
\draw[line width=0.03cm]  (0,0) rectangle (5,5);
%perfectmatching
\draw[line width=0.5cm,color=lightgray] (W1)--(B1); \draw[line width=0.5cm,color=lightgray] (W2)--(B2); 
\draw[line width=0.5cm,color=lightgray] (W3)--(B3); \draw[line width=0.5cm,color=lightgray] (W4)--(B4);
%edge
\draw[line width=0.06cm]  (B4)--(W4)--(B3)--(W3)--(B2)--(W2)--(B4) ; 
\draw[line width=0.06cm]  (W2)--(B3); \draw[line width=0.06cm]  (W1)--(B1); \draw[line width=0.06cm]  (W1)--(B2); 
\draw[line width=0.06cm]  (B1)--(0,1.25); \draw[line width=0.06cm]  (W2)--(5,1.25); 
\draw[line width=0.06cm]  (B1)--(1.75,0); \draw[line width=0.06cm]  (W4)--(1.75,5); 
\draw[line width=0.06cm]  (B4)--(3.25,5); \draw[line width=0.06cm]  (W1)--(3.25,0);  
\draw[line width=0.06cm]  (B4)--(4,5); \draw[line width=0.06cm]  (W2)--(4,0); 
\draw[line width=0.06cm]  (B4)--(5,3.75); \draw[line width=0.06cm]  (W3)--(0,3.75);  
%black
\filldraw  [ultra thick, fill=black] (1,0.5) circle [radius=0.2] ; \filldraw  [ultra thick, fill=black] (2.5,1.5) circle [radius=0.2] ;
\filldraw  [ultra thick, fill=black] (2.5,3.5) circle [radius=0.2] ; \filldraw  [ultra thick, fill=black] (4,4.5) circle [radius=0.2] ;
%white
\draw  [line width=1.6pt,fill=white] (2.5,0.5) circle [radius=0.2] ; \draw  [line width=1.6pt, fill=white] (4,2) circle [radius=0.2] ;
\draw  [line width=1.6pt,fill=white] (1,3)circle [radius=0.2] ; \draw  [line width=1.6pt,fill=white] (2.5,4.5)circle [radius=0.2] ;
\end{tikzpicture}
 }}; 

\node (PM1) at (3.3,0) 
{\scalebox{0.408}{
\begin{tikzpicture}
%vertex
\node (B1) at (1,0.5){$$}; \node (B2) at (2.5,1.5){$$}; \node (B3) at (2.5,3.5){$$}; \node (B4) at (4,4.5){$$};  
\node (W1) at (2.5,0.5){$$}; \node (W2) at (4,2){$$}; \node (W3) at (1,3){$$}; \node (W4) at (2.5,4.5){$$}; 
\draw[line width=0.03cm]  (0,0) rectangle (5,5);
%perfectmatching
\draw[line width=0.5cm,color=lightgray] (W1)--(B2); 
\draw[line width=0.5cm,color=lightgray] (W3)--(B3); \draw[line width=0.5cm,color=lightgray] (W4)--(B4);
\draw[line width=0.5cm,color=lightgray] (B1)--(0,1.25); \draw[line width=0.5cm,color=lightgray] (W2)--(5,1.25);  
%edge
\draw[line width=0.06cm]  (B4)--(W4)--(B3)--(W3)--(B2)--(W2)--(B4) ; 
\draw[line width=0.06cm]  (W2)--(B3); \draw[line width=0.06cm]  (W1)--(B1); \draw[line width=0.06cm]  (W1)--(B2); 
\draw[line width=0.06cm]  (B1)--(0,1.25); \draw[line width=0.06cm]  (W2)--(5,1.25); 
\draw[line width=0.06cm]  (B1)--(1.75,0); \draw[line width=0.06cm]  (W4)--(1.75,5); 
\draw[line width=0.06cm]  (B4)--(3.25,5); \draw[line width=0.06cm]  (W1)--(3.25,0);  
\draw[line width=0.06cm]  (B4)--(4,5); \draw[line width=0.06cm]  (W2)--(4,0); 
\draw[line width=0.06cm]  (B4)--(5,3.75); \draw[line width=0.06cm]  (W3)--(0,3.75);  
%black
\filldraw  [ultra thick, fill=black] (1,0.5) circle [radius=0.2] ; \filldraw  [ultra thick, fill=black] (2.5,1.5) circle [radius=0.2] ;
\filldraw  [ultra thick, fill=black] (2.5,3.5) circle [radius=0.2] ; \filldraw  [ultra thick, fill=black] (4,4.5) circle [radius=0.2] ;
%white
\draw  [line width=1.6pt,fill=white] (2.5,0.5) circle [radius=0.2] ; \draw  [line width=1.6pt, fill=white] (4,2) circle [radius=0.2] ;
\draw  [line width=1.6pt,fill=white] (1,3)circle [radius=0.2] ; \draw  [line width=1.6pt,fill=white] (2.5,4.5)circle [radius=0.2] ;
\end{tikzpicture}
 }}; 

\node (PM2) at (6.7,0) 
{\scalebox{0.408}{
\begin{tikzpicture}
%vertex
\node (B1) at (1,0.5){$$}; \node (B2) at (2.5,1.5){$$}; \node (B3) at (2.5,3.5){$$}; \node (B4) at (4,4.5){$$};  
\node (W1) at (2.5,0.5){$$}; \node (W2) at (4,2){$$}; \node (W3) at (1,3){$$}; \node (W4) at (2.5,4.5){$$}; 
\draw[line width=0.03cm]  (0,0) rectangle (5,5);
%perfectmatching
\draw[line width=0.5cm,color=lightgray] (W1)--(B2); 
\draw[line width=0.5cm,color=lightgray] (W3)--(B3); \draw[line width=0.5cm,color=lightgray] (W2)--(B4);
\draw[line width=0.5cm,color=lightgray] (B1)--(1.75,0); \draw[line width=0.5cm,color=lightgray] (W4)--(1.75,5);
%edge
\draw[line width=0.06cm]  (B4)--(W4)--(B3)--(W3)--(B2)--(W2)--(B4) ; 
\draw[line width=0.06cm]  (W2)--(B3); \draw[line width=0.06cm]  (W1)--(B1); \draw[line width=0.06cm]  (W1)--(B2); 
\draw[line width=0.06cm]  (B1)--(0,1.25); \draw[line width=0.06cm]  (W2)--(5,1.25); 
\draw[line width=0.06cm]  (B1)--(1.75,0); \draw[line width=0.06cm]  (W4)--(1.75,5); 
\draw[line width=0.06cm]  (B4)--(3.25,5); \draw[line width=0.06cm]  (W1)--(3.25,0);  
\draw[line width=0.06cm]  (B4)--(4,5); \draw[line width=0.06cm]  (W2)--(4,0); 
\draw[line width=0.06cm]  (B4)--(5,3.75); \draw[line width=0.06cm]  (W3)--(0,3.75);  
%black
\filldraw  [ultra thick, fill=black] (1,0.5) circle [radius=0.2] ; \filldraw  [ultra thick, fill=black] (2.5,1.5) circle [radius=0.2] ;
\filldraw  [ultra thick, fill=black] (2.5,3.5) circle [radius=0.2] ; \filldraw  [ultra thick, fill=black] (4,4.5) circle [radius=0.2] ;
%white
\draw  [line width=1.6pt,fill=white] (2.5,0.5) circle [radius=0.2] ; \draw  [line width=1.6pt, fill=white] (4,2) circle [radius=0.2] ;
\draw  [line width=1.6pt,fill=white] (1,3)circle [radius=0.2] ; \draw  [line width=1.6pt,fill=white] (2.5,4.5)circle [radius=0.2] ;
\end{tikzpicture}
 }} ;  

\node (PM3) at (10,0) 
{\scalebox{0.408}{
\begin{tikzpicture}
%vertex
\node (B1) at (1,0.5){$$}; \node (B2) at (2.5,1.5){$$}; \node (B3) at (2.5,3.5){$$}; \node (B4) at (4,4.5){$$};  
\node (W1) at (2.5,0.5){$$}; \node (W2) at (4,2){$$}; \node (W3) at (1,3){$$}; \node (W4) at (2.5,4.5){$$}; 
\draw[line width=0.03cm]  (0,0) rectangle (5,5);
%perfectmatching
\draw[line width=0.5cm,color=lightgray] (W1)--(B2); \draw[line width=0.5cm,color=lightgray] (W2)--(B3); 
\draw[line width=0.5cm,color=lightgray] (B1)--(1.75,0); \draw[line width=0.5cm,color=lightgray] (W4)--(1.75,5); 
\draw[line width=0.5cm,color=lightgray] (B4)--(5,3.75); \draw[line width=0.5cm,color=lightgray] (W3)--(0,3.75);
%edge
\draw[line width=0.06cm]  (B4)--(W4)--(B3)--(W3)--(B2)--(W2)--(B4) ; 
\draw[line width=0.06cm]  (W2)--(B3); \draw[line width=0.06cm]  (W1)--(B1); \draw[line width=0.06cm]  (W1)--(B2); 
\draw[line width=0.06cm]  (B1)--(0,1.25); \draw[line width=0.06cm]  (W2)--(5,1.25); 
\draw[line width=0.06cm]  (B1)--(1.75,0); \draw[line width=0.06cm]  (W4)--(1.75,5); 
\draw[line width=0.06cm]  (B4)--(3.25,5); \draw[line width=0.06cm]  (W1)--(3.25,0);  
\draw[line width=0.06cm]  (B4)--(4,5); \draw[line width=0.06cm]  (W2)--(4,0); 
\draw[line width=0.06cm]  (B4)--(5,3.75); \draw[line width=0.06cm]  (W3)--(0,3.75);  
%black
\filldraw  [ultra thick, fill=black] (1,0.5) circle [radius=0.2] ; \filldraw  [ultra thick, fill=black] (2.5,1.5) circle [radius=0.2] ;
\filldraw  [ultra thick, fill=black] (2.5,3.5) circle [radius=0.2] ; \filldraw  [ultra thick, fill=black] (4,4.5) circle [radius=0.2] ;
%white
\draw  [line width=1.6pt,fill=white] (2.5,0.5) circle [radius=0.2] ; \draw  [line width=1.6pt, fill=white] (4,2) circle [radius=0.2] ;
\draw  [line width=1.6pt,fill=white] (1,3)circle [radius=0.2] ; \draw  [line width=1.6pt,fill=white] (2.5,4.5)circle [radius=0.2] ;
\end{tikzpicture}
 }} ; 

\node (PM4) at (3.3,-3.1) 
{\scalebox{0.408}{
\begin{tikzpicture}
%vertex
\node (B1) at (1,0.5){$$}; \node (B2) at (2.5,1.5){$$}; \node (B3) at (2.5,3.5){$$}; \node (B4) at (4,4.5){$$};  
\node (W1) at (2.5,0.5){$$}; \node (W2) at (4,2){$$}; \node (W3) at (1,3){$$}; \node (W4) at (2.5,4.5){$$}; 
\draw[line width=0.03cm]  (0,0) rectangle (5,5);
%perfectmatching
\draw[line width=0.5cm,color=lightgray] (W1)--(B1); \draw[line width=0.5cm,color=lightgray] (W2)--(B2); \draw[line width=0.5cm,color=lightgray] (W4)--(B3);  
\draw[line width=0.5cm,color=lightgray] (B4)--(5,3.75); \draw[line width=0.5cm,color=lightgray] (W3)--(0,3.75);
%edge
\draw[line width=0.06cm]  (B4)--(W4)--(B3)--(W3)--(B2)--(W2)--(B4) ; 
\draw[line width=0.06cm]  (W2)--(B3); \draw[line width=0.06cm]  (W1)--(B1); \draw[line width=0.06cm]  (W1)--(B2); 
\draw[line width=0.06cm]  (B1)--(0,1.25); \draw[line width=0.06cm]  (W2)--(5,1.25); 
\draw[line width=0.06cm]  (B1)--(1.75,0); \draw[line width=0.06cm]  (W4)--(1.75,5); 
\draw[line width=0.06cm]  (B4)--(3.25,5); \draw[line width=0.06cm]  (W1)--(3.25,0);  
\draw[line width=0.06cm]  (B4)--(4,5); \draw[line width=0.06cm]  (W2)--(4,0); 
\draw[line width=0.06cm]  (B4)--(5,3.75); \draw[line width=0.06cm]  (W3)--(0,3.75);  
%black
\filldraw  [ultra thick, fill=black] (1,0.5) circle [radius=0.2] ; \filldraw  [ultra thick, fill=black] (2.5,1.5) circle [radius=0.2] ;
\filldraw  [ultra thick, fill=black] (2.5,3.5) circle [radius=0.2] ; \filldraw  [ultra thick, fill=black] (4,4.5) circle [radius=0.2] ;
%white
\draw  [line width=1.6pt,fill=white] (2.5,0.5) circle [radius=0.2] ; \draw  [line width=1.6pt, fill=white] (4,2) circle [radius=0.2] ;
\draw  [line width=1.6pt,fill=white] (1,3)circle [radius=0.2] ; \draw  [line width=1.6pt,fill=white] (2.5,4.5)circle [radius=0.2] ;
\end{tikzpicture}
 }}; 

\node (PM5) at (6.7,-3.1) 
{\scalebox{0.408}{
\begin{tikzpicture}
%vertex
\node (B1) at (1,0.5){$$}; \node (B2) at (2.5,1.5){$$}; \node (B3) at (2.5,3.5){$$}; \node (B4) at (4,4.5){$$};  
\node (W1) at (2.5,0.5){$$}; \node (W2) at (4,2){$$}; \node (W3) at (1,3){$$}; \node (W4) at (2.5,4.5){$$}; 
\draw[line width=0.03cm]  (0,0) rectangle (5,5);
%perfectmatching
\draw[line width=0.5cm,color=lightgray] (W1)--(B1); \draw[line width=0.5cm,color=lightgray] (W3)--(B2); \draw[line width=0.5cm,color=lightgray] (W4)--(B3);  
\draw[line width=0.5cm,color=lightgray] (B4)--(4,5); \draw[line width=0.5cm,color=lightgray] (W2)--(4,0);
%edge
\draw[line width=0.06cm]  (B4)--(W4)--(B3)--(W3)--(B2)--(W2)--(B4) ; 
\draw[line width=0.06cm]  (W2)--(B3); \draw[line width=0.06cm]  (W1)--(B1); \draw[line width=0.06cm]  (W1)--(B2); 
\draw[line width=0.06cm]  (B1)--(0,1.25); \draw[line width=0.06cm]  (W2)--(5,1.25); 
\draw[line width=0.06cm]  (B1)--(1.75,0); \draw[line width=0.06cm]  (W4)--(1.75,5); 
\draw[line width=0.06cm]  (B4)--(3.25,5); \draw[line width=0.06cm]  (W1)--(3.25,0);  
\draw[line width=0.06cm]  (B4)--(4,5); \draw[line width=0.06cm]  (W2)--(4,0); 
\draw[line width=0.06cm]  (B4)--(5,3.75); \draw[line width=0.06cm]  (W3)--(0,3.75);  
%black
\filldraw  [ultra thick, fill=black] (1,0.5) circle [radius=0.2] ; \filldraw  [ultra thick, fill=black] (2.5,1.5) circle [radius=0.2] ;
\filldraw  [ultra thick, fill=black] (2.5,3.5) circle [radius=0.2] ; \filldraw  [ultra thick, fill=black] (4,4.5) circle [radius=0.2] ;
%white
\draw  [line width=1.6pt,fill=white] (2.5,0.5) circle [radius=0.2] ; \draw  [line width=1.6pt, fill=white] (4,2) circle [radius=0.2] ;
\draw  [line width=1.6pt,fill=white] (1,3)circle [radius=0.2] ; \draw  [line width=1.6pt,fill=white] (2.5,4.5)circle [radius=0.2] ;
\end{tikzpicture}
 }}; 

\node (PM6) at (10,-3.1) 
{\scalebox{0.408}{
\begin{tikzpicture}
%vertex
\node (B1) at (1,0.5){$$}; \node (B2) at (2.5,1.5){$$}; \node (B3) at (2.5,3.5){$$}; \node (B4) at (4,4.5){$$};  
\node (W1) at (2.5,0.5){$$}; \node (W2) at (4,2){$$}; \node (W3) at (1,3){$$}; \node (W4) at (2.5,4.5){$$}; 
\draw[line width=0.03cm]  (0,0) rectangle (5,5);
%perfectmatching
\draw[line width=0.5cm,color=lightgray] (W3)--(B2); \draw[line width=0.5cm,color=lightgray] (W4)--(B3);  
\draw[line width=0.5cm,color=lightgray] (B1)--(0,1.25); \draw[line width=0.5cm,color=lightgray] (W2)--(5,1.25); 
\draw[line width=0.5cm,color=lightgray] (B4)--(3.25,5); \draw[line width=0.5cm,color=lightgray] (W1)--(3.25,0);
%edge
\draw[line width=0.06cm]  (B4)--(W4)--(B3)--(W3)--(B2)--(W2)--(B4) ; 
\draw[line width=0.06cm]  (W2)--(B3); \draw[line width=0.06cm]  (W1)--(B1); \draw[line width=0.06cm]  (W1)--(B2); 
\draw[line width=0.06cm]  (B1)--(0,1.25); \draw[line width=0.06cm]  (W2)--(5,1.25); 
\draw[line width=0.06cm]  (B1)--(1.75,0); \draw[line width=0.06cm]  (W4)--(1.75,5); 
\draw[line width=0.06cm]  (B4)--(3.25,5); \draw[line width=0.06cm]  (W1)--(3.25,0);  
\draw[line width=0.06cm]  (B4)--(4,5); \draw[line width=0.06cm]  (W2)--(4,0); 
\draw[line width=0.06cm]  (B4)--(5,3.75); \draw[line width=0.06cm]  (W3)--(0,3.75);  
%black
\filldraw  [ultra thick, fill=black] (1,0.5) circle [radius=0.2] ; \filldraw  [ultra thick, fill=black] (2.5,1.5) circle [radius=0.2] ;
\filldraw  [ultra thick, fill=black] (2.5,3.5) circle [radius=0.2] ; \filldraw  [ultra thick, fill=black] (4,4.5) circle [radius=0.2] ;
%white
\draw  [line width=1.6pt,fill=white] (2.5,0.5) circle [radius=0.2] ; \draw  [line width=1.6pt, fill=white] (4,2) circle [radius=0.2] ;
\draw  [line width=1.6pt,fill=white] (1,3)circle [radius=0.2] ; \draw  [line width=1.6pt,fill=white] (2.5,4.5)circle [radius=0.2] ;
\end{tikzpicture}
 }} ;  

\end{tikzpicture}
}}
\end{center}

\medskip

Then, we have the following table. We denote by $T(a,b,c)$ a divisorial ideal $T(a,b,c,0,0,0)$. 

\medskip

\begin{center}
{\small{
\begin{tabular}{|l|l|l|} \hline 
$e_0\sfB^{\prime\prime}$&$e_1\sfB^{\prime\prime}$&$e_2\sfB^{\prime\prime}$ \\ \hline
$T_{00}\cong R$&$T_{10}\cong T(-1,-1,0)$&$T_{20}\cong T(-1,-1,-1)$ \\
$T_{01}\cong T(1,1,0)$&$T_{11}\cong R$&$T_{21}\cong T(0,0,-1)$ \\
$T_{02}\cong T(1,1,1)$&$T_{12}\cong T(0,0,1)$&$T_{22}\cong R$ \\
$T_{03}\cong T(2,3,2)$&$T_{13}\cong T(1,2,2)$&$T_{23}\cong T(1,2,1)$ \\
$T_{04}\cong T(1,2,1)$&$T_{14}\cong T(0,1,1)$&$T_{24}\cong T(0,1,0)$ \\
$T_{05}\cong T(0,1,1)$&$T_{15}\cong T(-1,0,1)$&$T_{25}\cong T(-1,0,0)$ \\ \hline
$e_3\sfB^{\prime\prime}$&$e_4\sfB^{\prime\prime}$&$e_5\sfB^{\prime\prime}$ \\ \hline
$T_{30}\cong T(-2,-3,-2)$&$T_{40}\cong T(-1,-2,-1)$&$T_{50}\cong T(0,-1,-1)$ \\
$T_{31}\cong T(-1,-2,-2)$&$T_{41}\cong T(0,-1,-1)$&$T_{51}\cong T(1,0,-1)$ \\
$T_{32}\cong T(-1,-2,-1)$&$T_{42}\cong T(0,-1,0)$&$T_{52}\cong T(1,0,0)$ \\
$T_{33}\cong R$&$T_{43}\cong T(1,1,1)$&$T_{53}\cong T(2,2,1)$\\ 
$T_{34}\cong T(-1,-1,-1)$&$T_{44}\cong R$&$T_{54}\cong T(1,1,0)$\\
$T_{35}\cong T(-2,-2,-1)$&$T_{45}\cong T(-1,-1,0)$&$T_{55}\cong R$\\ \hline 
\end{tabular} 
}}
\end{center}

\medskip

%%%%%%%%%%%%%%%%%%%%%%%%%%%%%%%%%%%
%%%%%%%%%%%%Type_6a-3%%%%%%%%%%%%%%%%%%
%%%%%%%%%%%%%%%%%%%%%%%%%%%%%%%%%%%

\subsubsection{Type 6a-3}
Next, we consider the following consistent dimer model and the associated QP $(Q, W_Q)$. 
We denote by $\sfC$ the complete Jacobian algebra $\calP(Q, W_Q)$.  
We will give splitting MM generators $e_i\sfC$. 

\newcommand{\figsixathreea}{
%vertex
\node (B1) at (1.5,0.5){$$}; \node (B2) at (3.5,1.5){$$}; \node (B3) at (2.5,3.5){$$}; \node (B4) at (0.5,2.5){$$};  
\node (W1) at (2.5,0.5){$$}; \node (W2) at (3.5,2.5){$$}; \node (W3) at (1.5,3.5){$$}; \node (W4) at (0.5,1.5){$$}; 
\draw[line width=0.024cm]  (0,0) rectangle (4,4);
%edge
\draw[line width=0.05cm]  (B1)--(W1)--(B2)--(W2)--(B3)--(W3)--(B4)--(W4)--(B1) ; 
\draw[line width=0.05cm]  (W1)--(B4); \draw[line width=0.05cm]  (W3)--(B2); 
\draw[line width=0.05cm]  (B1)--(1.5,0); \draw[line width=0.05cm]  (W1)--(2.5,0); 
\draw[line width=0.05cm]  (B2)--(4,1.5); \draw[line width=0.05cm]  (W2)--(4,2.5);
\draw[line width=0.05cm]  (B3)--(2.5,4); \draw[line width=0.05cm]  (W3)--(1.5,4); 
\draw[line width=0.05cm]  (B4)--(0,2.5); \draw[line width=0.05cm]  (W4)--(0,1.5);
%black
\filldraw  [ultra thick, fill=black] (1.5,0.5) circle [radius=0.16] ; \filldraw  [ultra thick, fill=black] (3.5,1.5) circle [radius=0.16] ;
\filldraw  [ultra thick, fill=black] (2.5,3.5) circle [radius=0.16] ; \filldraw  [ultra thick, fill=black] (0.5,2.5) circle [radius=0.16] ;
%white
\draw  [line width=1.28pt,fill=white] (2.5,0.5) circle [radius=0.16] ; \draw  [line width=1.28pt, fill=white] (3.5,2.5) circle [radius=0.16] ;
\draw  [line width=1.28pt,fill=white] (1.5,3.5)circle [radius=0.16] ; \draw  [line width=1.28pt,fill=white] (0.5,1.5)circle [radius=0.16] ;
}

\begin{center}
\begin{tikzpicture}
\node (DM) at (0,0) 
{\scalebox{0.675}{
\begin{tikzpicture}
\figsixathreea
\end{tikzpicture}
} }; 

\node (QV) at (6,0) 
{\scalebox{0.675}{
\begin{tikzpicture}

\node (Q0) at (2,2){$0$}; \node (Q1a) at (2,0){$1$}; \node (Q1b) at (2,4){$1$}; \node (Q2) at (2.75,2.75){$2$}; \node (Q3) at (1.25,1.25){$3$}; 
\node (Q4a) at (0,0){$4$}; \node (Q4b) at (4,0){$4$}; \node (Q4c) at (4,4){$4$}; \node (Q4d) at (0,4){$4$}; 
\node (Q5a) at (0,2){$5$}; \node (Q5b) at (4,2){$5$}; 
%edge
\draw[lightgray, line width=0.05cm]  (B1)--(W1)--(B2)--(W2)--(B3)--(W3)--(B4)--(W4)--(B1) ; 
\draw[lightgray, line width=0.05cm]  (W1)--(B4); \draw[lightgray, line width=0.05cm]  (W3)--(B2); 
\draw[lightgray, line width=0.05cm]  (B1)--(1.5,0); \draw[lightgray, line width=0.05cm]  (W1)--(2.5,0); 
\draw[lightgray, line width=0.05cm]  (B2)--(4,1.5); \draw[lightgray, line width=0.05cm]  (W2)--(4,2.5);
\draw[lightgray, line width=0.05cm]  (B3)--(2.5,4); \draw[lightgray, line width=0.05cm]  (W3)--(1.5,4); 
\draw[lightgray, line width=0.05cm]  (B4)--(0,2.5); \draw[lightgray, line width=0.05cm]  (W4)--(0,1.5);
%black
\filldraw  [ultra thick, draw=lightgray, fill=lightgray] (1.5,0.5) circle [radius=0.16] ; \filldraw  [ultra thick, draw=lightgray, fill=lightgray] (3.5,1.5) circle [radius=0.16] ;
\filldraw  [ultra thick, draw=lightgray, fill=lightgray] (2.5,3.5) circle [radius=0.16] ; \filldraw  [ultra thick, draw=lightgray, fill=lightgray] (0.5,2.5) circle [radius=0.16] ;
%white
\draw  [line width=1.28pt,draw=lightgray, fill=white] (2.5,0.5) circle [radius=0.16] ; \draw  [line width=1.28pt, draw=lightgray, fill=white] (3.5,2.5) circle [radius=0.16] ;
\draw  [line width=1.28pt,draw=lightgray, fill=white] (1.5,3.5)circle [radius=0.16] ; \draw  [line width=1.28pt, draw=lightgray, fill=white] (0.5,1.5)circle [radius=0.16] ;
%arrow
\draw[->, line width=0.064cm] (Q0)--(Q4b); \draw[->, line width=0.064cm] (Q0)--(Q4d); \draw[->, line width=0.064cm] (Q1b)--(Q2);
\draw[->, line width=0.064cm] (Q1a)--(Q3); \draw[->, line width=0.064cm] (Q2)--(Q0); \draw[->, line width=0.064cm] (Q2)--(Q4c); 
\draw[->, line width=0.064cm] (Q3)--(Q0); \draw[->, line width=0.064cm] (Q3)--(Q4a); 
\draw[->, line width=0.064cm] (Q4a)--(Q5a); \draw[->, line width=0.064cm] (Q4b)--(Q5b); 
\draw[->, line width=0.064cm] (Q4d)--(Q5a); \draw[->, line width=0.064cm] (Q4c)--(Q5b); 
\draw[->, line width=0.064cm] (Q4a)--(Q1a); \draw[->, line width=0.064cm] (Q4b)--(Q1a); 
\draw[->, line width=0.064cm] (Q4d)--(Q1b); \draw[->, line width=0.064cm] (Q4c)--(Q1b); 
\draw[->, line width=0.064cm] (Q5a)--(Q3); \draw[->, line width=0.064cm] (Q5b)--(Q2); 
\end{tikzpicture}
} }; 
\end{tikzpicture}
\end{center}

\medskip

In the following figure, we fix the perfect matching $\sfP_0$.  
Then $\sfP_1, \cdots, \sfP_6$ are extremal perfect matchings corresponding to $v_1, \cdots, v_6$. 

\medskip

%%%%%perfect_matching%%%%%

\begin{center}
{\scalebox{0.9}{
\begin{tikzpicture} 
\node at (0,-1.4) {$\sfP_0$};\node at (3.3,-1.4) {$\sfP_1$}; \node at (6.7,-1.4) {$\sfP_2$}; \node at (10,-1.4) {$\sfP_3$}; 
\node at (3.3,-4.5) {$\sfP_4$};\node at (6.7,-4.5) {$\sfP_5$}; \node at (10,-4.5) {$\sfP_6$}; 

\node (PM0) at (0,0) 
{\scalebox{0.51}{
\begin{tikzpicture}
%perfectmatching
\draw[line width=0.4cm,color=lightgray] (W1)--(B1); \draw[line width=0.4cm,color=lightgray] (W3)--(B3); 
\draw[line width=0.4cm,color=lightgray] (W2)--(B2); \draw[line width=0.4cm,color=lightgray] (W4)--(B4); 
%\draw[line width=0.4cm,color=lightgray] (W2)--(4,2.5); \draw[line width=0.4cm,color=lightgray] (B2)--(4,1.5); 
%\draw[line width=0.4cm,color=lightgray] (W4)--(0,1.5); \draw[line width=0.4cm,color=lightgray] (B4)--(0,2.5); 
\figsixathreea
\end{tikzpicture}
 }}; 

\node (PM1) at (3.3,0) 
{\scalebox{0.51}{
\begin{tikzpicture}
%perfectmatching
\draw[line width=0.4cm,color=lightgray] (W3)--(B3); \draw[line width=0.4cm,color=lightgray] (W1)--(B2); 
\draw[line width=0.4cm,color=lightgray] (W4)--(B1);  
\draw[line width=0.4cm,color=lightgray] (W2)--(4,2.5); \draw[line width=0.4cm,color=lightgray] (B4)--(0,2.5); 
\figsixathreea
\end{tikzpicture}
 }}; 

\node (PM2) at (6.7,0) 
{\scalebox{0.51}{
\begin{tikzpicture}
%perfectmatching
\draw[line width=0.4cm,color=lightgray] (W4)--(B4); \draw[line width=0.4cm,color=lightgray] (W1)--(B2); 
\draw[line width=0.4cm,color=lightgray] (W2)--(B3);  
\draw[line width=0.4cm,color=lightgray] (W3)--(1.5,4); \draw[line width=0.4cm,color=lightgray] (B1)--(1.5,0); 
\figsixathreea
\end{tikzpicture}
 }} ;  

\node (PM3) at (10,0) 
{\scalebox{0.51}{
\begin{tikzpicture}
%perfectmatching
\draw[line width=0.4cm,color=lightgray] (W1)--(B4); \draw[line width=0.4cm,color=lightgray] (W2)--(B3); 
\draw[line width=0.4cm,color=lightgray] (W4)--(0,1.5); \draw[line width=0.4cm,color=lightgray] (B2)--(4,1.5); 
\draw[line width=0.4cm,color=lightgray] (W3)--(1.5,4); \draw[line width=0.4cm,color=lightgray] (B1)--(1.5,0); 
\figsixathreea
\end{tikzpicture}
 }} ; 

\node (PM4) at (3.3,-3.1) 
{\scalebox{0.51}{
\begin{tikzpicture}
%perfectmatching
\draw[line width=0.4cm,color=lightgray] (W1)--(B1); \draw[line width=0.4cm,color=lightgray] (W2)--(B3); 
\draw[line width=0.4cm,color=lightgray] (W3)--(B4); 
\draw[line width=0.4cm,color=lightgray] (W4)--(0,1.5); \draw[line width=0.4cm,color=lightgray] (B2)--(4,1.5); 
\figsixathreea
\end{tikzpicture}
 }}; 

\node (PM5) at (6.7,-3.1) 
{\scalebox{0.51}{
\begin{tikzpicture}
%perfectmatching
\draw[line width=0.4cm,color=lightgray] (W4)--(B1); \draw[line width=0.4cm,color=lightgray] (W2)--(B2); 
\draw[line width=0.4cm,color=lightgray] (W3)--(B4); 
\draw[line width=0.4cm,color=lightgray] (W1)--(2.5,0); \draw[line width=0.4cm,color=lightgray] (B3)--(2.5,4);
\figsixathreea
\end{tikzpicture}
 }}; 

\node (PM6) at (10,-3.1) 
{\scalebox{0.51}{
\begin{tikzpicture}
%perfectmatching
\draw[line width=0.4cm,color=lightgray] (W4)--(B1); \draw[line width=0.4cm,color=lightgray] (W3)--(B2);  
\draw[line width=0.4cm,color=lightgray] (W1)--(2.5,0); \draw[line width=0.4cm,color=lightgray] (B3)--(2.5,4); 
\draw[line width=0.4cm,color=lightgray] (W2)--(4,2.5); \draw[line width=0.4cm,color=lightgray] (B4)--(0,2.5); 
\figsixathreea
\end{tikzpicture}
 }} ;  
\end{tikzpicture}
}}
\end{center}

\medskip

Then we have the following table. We denote by $T(a,b,c)$ a divisorial ideal $T(a,b,c,0,0,0)$. 

\medskip

\begin{center}
{\small{
\begin{tabular}{|l|l|l|} \hline 
$e_0\sfC$&$e_1\sfC$&$e_2\sfC$\\ \hline
$T_{00}\cong R$&$T_{10}\cong T(-1,-2,-1)$&$T_{20}\cong T(-2,-2,-1)$ \\
$T_{01}\cong T(1,2,1)$&$T_{11}\cong R$&$T_{21}\cong T(-1,0,0)$ \\
$T_{02}\cong T(2,2,1)$&$T_{12}\cong T(1,0,0)$&$T_{22}\cong R$ \\
$T_{03}\cong T(0,0,-1)$&$T_{13}\cong T(-1,-2,-2)$&$T_{23}\cong T(-2,-2,-2)$ \\
$T_{04}\cong T(1,1,0)$&$T_{14}\cong T(0,-1,-1)$&$T_{24}\cong T(-1,-1,-1)$ \\
$T_{05}\cong T(0,-1,-1)$&$T_{15}\cong T(-1,-3,-2)$&$T_{25}\cong T(-2,-3,-2)$ \\ \hline
$e_3\sfC$&$e_4\sfC$&$e_5\sfC$ \\ \hline
$T_{30}\cong T(0,0,1)$&$T_{40}\cong T(-1,-1,0)$&$T_{50}\cong T(0,1,1)$ \\
$T_{31}\cong T(1,2,2)$&$T_{41}\cong T(0,1,1)$&$T_{51}\cong T(1,3,2)$ \\
$T_{32}\cong T(2,2,2)$&$T_{42}\cong T(1,1,1)$&$T_{52}\cong T(2,3,2)$ \\
$T_{33}\cong R$&$T_{43}\cong T(-1,-1,-1)$&$T_{53}\cong T(0,1,0)$\\
$T_{34}\cong T(1,1,1)$&$T_{44}\cong R$&$T_{54}\cong T(1,2,1)$\\
$T_{35}\cong T(0,-1,0)$&$T_{45}\cong T(-1,-2,-1)$&$T_{55}\cong R$\\ \hline
\end{tabular} 
}}
\end{center}

%%%%%%%%%%%%
%%%%%%%%%%%%
%%%%%%%%%%%%
\iffalse
%%%%%%%%%%%%
%%%%%%%%%%%%
%%%%%%%%%%%%

\medskip

\begin{center}
{\small{
\begin{tabular}{|l|l|l|} \hline 
$(e_0\sfC)^*$&$(e_1\sfC)^*$&$(e_2\sfC)^*$ \\ \hline
$T_{00}^*\cong R$&$T_{10}^*\cong T(1,2,1)$&$T_{20}^*\cong T(2,2,1)$ \\
$T_{01}^*\cong T(-1,-2,-1)$&$T_{11}^*\cong R$&$T_{21}^*\cong T(1,0,0)$ \\
$T_{02}^*\cong T(-2,-2,-1)$&$T_{12}^*\cong T(-1,0,0)$&$T_{22}^*\cong R$ \\
$T_{03}^*\cong T(0,0,1)$&$T_{13}^*\cong T(1,2,2)$&$T_{23}^*\cong T(2,2,2)$ \\
$T_{04}^*\cong T(-1,-1,0)$&$T_{14}^*\cong T(0,1,1)$&$T_{24}^*\cong T(1,1,1)$ \\ 
$T_{05}^*\cong T(0,1,1)$&$T_{15}^*\cong T(1,3,2)$&$T_{25}^*\cong T(2,3,2)$ \\ \hline
$(e_3\sfC)^*$&$(e_4\sfC)^*$&$(e_5\sfC)^*$ \\ \hline
$T_{30}^*\cong T(0,0,-1)$&$T_{40}^*\cong T(1,1,0)$&$T_{50}^*\cong T(0,-1,-1)$ \\
$T_{31}^*\cong T(-1,-2,-2)$&$T_{41}^*\cong T(0,-1,-1)$&$T_{51}^*\cong T(-1,-3,-2)$ \\
$T_{32}^*\cong T(-2,-2,-2)$&$T_{42}^*\cong T(-1,-1,-1)$&$T_{52}^*\cong T(-2,-3,-2)$ \\
$T_{33}^*\cong R$&$T_{43}^*\cong T(1,1,1)$&$T_{53}^*\cong T(0,-1,0)$\\
$T_{34}^*\cong T(-1,-1,-1)$&$T_{44}^*\cong R$&$T_{54}^*\cong T(-1,-2,-1)$\\
$T_{35}^*\cong T(0,1,0)$&$T_{45}^*\cong T(1,2,1)$&$T_{55}^*\cong R$\\ \hline
\end{tabular} 
}}
\end{center}

%%%%%%%%%%%%
%%%%%%%%%%%%
%%%%%%%%%%%%
\fi
%%%%%%%%%%%%
%%%%%%%%%%%%
%%%%%%%%%%%%

\medskip

%%%%%%%%%%%%%%%%%%%%%%%%%%%%%%%%%%%
%%%%%%%%%%%%Type_6a-3(2)%%%%%%%%%%%%%%%%%%
%%%%%%%%%%%%%%%%%%%%%%%%%%%%%%%%%%%
Let $x$ (resp. $y$) be a $1$-cycle on the two torus $\sfT$ which corresponds to $(1,0)\in\rmH_1(\sfT)$ (resp. $(0,1)\in\rmH_1(\sfT)$) in the above case. 
By replacing these cycles by $x, -(x+y)$ respectively, we have the following dimer model and the associated QP $(Q, W_Q)$. 
We denote by $\sfC^\prime$ the complete Jacobian algebra $\calP(Q, W_Q)$. 
Since this dimer model is right equivalent to the previous one, we have $\sfC\cong\sfC^\prime$ as an $R$-algebra. 
However, a change of cycles induces an automorphism on $R$, and it gives another description of each splitting MM generator (see Remark~\ref{rem_auto}). 
Thus we also give splitting MM generators $e_i\sfC^\prime$.

\medskip

\newcommand{\figsixathreeb}{
%vertex
\node (B1) at (2.5,1){$$}; \node (B2) at (4.5,2){$$}; \node (B3) at (0.5,3){$$}; \node (B4) at (2.5,4){$$};  
\node (W1) at (0.5,1){$$}; \node (W2) at (3.5,2){$$}; \node (W3) at (1.5,3){$$}; \node (W4) at (4.5,4){$$}; 
\draw[line width=0.03cm]  (0,0) rectangle (5,5);
%edge
\draw[line width=0.06cm]  (B3)--(W1)--(B1)--(W2)--(B2)--(W4)--(B4)--(W3)--(B3) ; 
\draw[line width=0.06cm]  (W3)--(B1); \draw[line width=0.06cm]  (W2)--(B4); 
\draw[line width=0.06cm]  (W1)--(1.5,0); \draw[line width=0.06cm]  (W1)--(0,1.5); 
\draw[line width=0.06cm]  (B1)--(3.5,0); \draw[line width=0.06cm]  (B2)--(5,1.5);
\draw[line width=0.06cm]  (B3)--(0,3.5); \draw[line width=0.06cm]  (B4)--(1.5,5); 
\draw[line width=0.06cm]  (W4)--(3.5,5); \draw[line width=0.06cm]  (W4)--(5,3.5);
%black
\filldraw  [ultra thick, fill=black] (2.5,1) circle [radius=0.2] ; \filldraw  [ultra thick, fill=black] (4.5,2) circle [radius=0.2] ;
\filldraw  [ultra thick, fill=black] (0.5,3) circle [radius=0.2] ; \filldraw  [ultra thick, fill=black] (2.5,4) circle [radius=0.2] ;
%white
\draw  [line width=1.6pt,fill=white] (0.5,1) circle [radius=0.2] ; \draw  [line width=1.6pt, fill=white] (3.5,2) circle [radius=0.2] ;
\draw  [line width=1.6pt,fill=white] (1.5,3)circle [radius=0.2] ; \draw  [line width=1.6pt, fill=white] (4.5,4)circle [radius=0.2] ;
}

\begin{center}
\begin{tikzpicture}
\node (DM) at (0,0) 
{\scalebox{0.54}{
\begin{tikzpicture}
\figsixathreeb
\end{tikzpicture}
} }; 

\node (QV) at (6,0) 
{\scalebox{0.54}{
\begin{tikzpicture}

\node (Q0a) at (2.5,0){{\Large$0$}}; \node (Q0b) at (2.5,5){{\Large$0$}}; \node (Q1) at (2.5,2.5){{\Large$1$}}; 
\node (Q2) at (1.25,2){{\Large$2$}}; \node (Q3) at (3.75,3){{\Large$3$}}; 
\node (Q4a) at (0,0){{\Large$4$}}; \node (Q4b) at (5,0){{\Large$4$}}; \node (Q4c) at (5,5){{\Large$4$}}; \node (Q4d) at (0,5){{\Large$4$}}; 
\node (Q5a) at (0,2.5){{\Large$5$}}; \node (Q5b) at (5,2.5){{\Large$5$}};
%edge
\draw[lightgray, line width=0.06cm]  (B3)--(W1)--(B1)--(W2)--(B2)--(W4)--(B4)--(W3)--(B3) ; 
\draw[lightgray, line width=0.06cm]  (W3)--(B1); \draw[lightgray, line width=0.06cm]  (W2)--(B4); 
\draw[lightgray, line width=0.06cm]  (W1)--(1.5,0); \draw[lightgray, line width=0.06cm]  (W1)--(0,1.5); 
\draw[lightgray, line width=0.06cm]  (B1)--(3.5,0); \draw[lightgray, line width=0.06cm]  (B2)--(5,1.5);
\draw[lightgray, line width=0.06cm]  (B3)--(0,3.5); \draw[lightgray, line width=0.06cm]  (B4)--(1.5,5); 
\draw[lightgray, line width=0.06cm]  (W4)--(3.5,5); \draw[lightgray, line width=0.06cm]  (W4)--(5,3.5);
%black
\filldraw  [ultra thick, draw=lightgray, fill=lightgray] (2.5,1) circle [radius=0.2] ; \filldraw  [ultra thick, draw=lightgray, fill=lightgray] (4.5,2) circle [radius=0.2] ;
\filldraw  [ultra thick, draw=lightgray, fill=lightgray] (0.5,3) circle [radius=0.2] ; \filldraw  [ultra thick, draw=lightgray, fill=lightgray] (2.5,4) circle [radius=0.2] ;
%white
\draw  [line width=1.6pt,draw=lightgray, fill=white] (0.5,1) circle [radius=0.2] ; \draw  [line width=1.6pt, draw=lightgray, fill=white] (3.5,2) circle [radius=0.2] ;
\draw  [line width=1.6pt,draw=lightgray, fill=white] (1.5,3)circle [radius=0.2] ; \draw  [line width=1.6pt, draw=lightgray, fill=white] (4.5,4)circle [radius=0.2] ;
%arrow
\draw[->, line width=0.08cm] (Q0a)--(Q4a); \draw[->, line width=0.08cm] (Q0a)--(Q4b); 
\draw[->, line width=0.08cm] (Q0b)--(Q4c); \draw[->, line width=0.08cm] (Q0b)--(Q4d);
\draw[->, line width=0.08cm] (Q1)--(Q2); \draw[->, line width=0.08cm] (Q1)--(Q3); \draw[->, line width=0.08cm] (Q2)--(Q0a); 
\draw[->, line width=0.08cm] (Q2)--(Q4d); \draw[->, line width=0.08cm] (Q3)--(Q0b); \draw[->, line width=0.08cm] (Q3)--(Q4b);
\draw[->, line width=0.08cm] (Q4a)--(Q5a); \draw[->, line width=0.08cm] (Q4b)--(Q5b); \draw[->, line width=0.08cm] (Q4c)--(Q5b); 
\draw[->, line width=0.08cm] (Q4d)--(Q5a); \draw[->, line width=0.08cm] (Q4b)--(Q1); \draw[->, line width=0.08cm] (Q4d)--(Q1);
\draw[->, line width=0.08cm] (Q5a)--(Q2); \draw[->, line width=0.08cm] (Q5b)--(Q3);
\end{tikzpicture}
} }; 
\end{tikzpicture}
\end{center}

\medskip

In the following figure, we fix the perfect matching $\sfP_0$.  
Then $\sfP_1, \cdots, \sfP_6$ are extremal perfect matchings corresponding to $v_1, \cdots, v_6$. 

\medskip

%%%%%perfect_matching%%%%%

\begin{center}
{\scalebox{0.9}{
\begin{tikzpicture} 
\node at (0,-1.4) {$\sfP_0$};\node at (3.3,-1.4) {$\sfP_1$}; \node at (6.7,-1.4) {$\sfP_2$}; \node at (10,-1.4) {$\sfP_3$}; 
\node at (3.3,-4.5) {$\sfP_4$};\node at (6.7,-4.5) {$\sfP_5$}; \node at (10,-4.5) {$\sfP_6$};  

\node (PM0) at (0,0) 
{\scalebox{0.408}{
\begin{tikzpicture}
%perfectmatching
\draw[line width=0.5cm,color=lightgray] (W1)--(B3); \draw[line width=0.5cm,color=lightgray] (W3)--(B1); 
\draw[line width=0.5cm,color=lightgray] (W2)--(B4); \draw[line width=0.5cm,color=lightgray] (W4)--(B2);
\figsixathreeb
\end{tikzpicture}
 }}; 

\node (PM1) at (3.3,0) 
{\scalebox{0.408}{
\begin{tikzpicture}
%perfectmatching
\draw[line width=0.5cm,color=lightgray] (W1)--(B1); \draw[line width=0.5cm,color=lightgray] (W2)--(B2); \draw[line width=0.5cm,color=lightgray] (W3)--(B4); 
\draw[line width=0.5cm,color=lightgray] (B3)--(0,3.5); \draw[line width=0.5cm,color=lightgray] (W4)--(5,3.5); 
\figsixathreeb
\end{tikzpicture}
 }}; 

\node (PM2) at (6.7,0) 
{\scalebox{0.408}{
\begin{tikzpicture}
%perfectmatching
\draw[line width=0.5cm,color=lightgray] (W1)--(B3); \draw[line width=0.5cm,color=lightgray] (W2)--(B2); \draw[line width=0.5cm,color=lightgray] (W3)--(B4); 
\draw[line width=0.5cm,color=lightgray] (B1)--(3.5,0); \draw[line width=0.5cm,color=lightgray] (W4)--(3.5,5);
\figsixathreeb
\end{tikzpicture}
 }} ;  

\node (PM3) at (10,0) 
{\scalebox{0.408}{
\begin{tikzpicture}
%perfectmatching
\draw[line width=0.5cm,color=lightgray] (W2)--(B4); \draw[line width=0.5cm,color=lightgray] (W3)--(B3);  
\draw[line width=0.5cm,color=lightgray] (B1)--(3.5,0); \draw[line width=0.5cm,color=lightgray] (W4)--(3.5,5);
\draw[line width=0.5cm,color=lightgray] (B2)--(5,1.5); \draw[line width=0.5cm,color=lightgray] (W1)--(0,1.5);
\figsixathreeb
\end{tikzpicture}
 }} ; 

\node (PM4) at (3.3,-3.1) 
{\scalebox{0.408}{
\begin{tikzpicture}
%perfectmatching
\draw[line width=0.5cm,color=lightgray] (W4)--(B4); \draw[line width=0.5cm,color=lightgray] (W3)--(B3); \draw[line width=0.5cm,color=lightgray] (W2)--(B1); 
\draw[line width=0.5cm,color=lightgray] (B2)--(5,1.5); \draw[line width=0.5cm,color=lightgray] (W1)--(0,1.5);
\figsixathreeb
\end{tikzpicture}
 }}; 

\node (PM5) at (6.7,-3.1) 
{\scalebox{0.408}{
\begin{tikzpicture}
%perfectmatching
\draw[line width=0.5cm,color=lightgray] (W4)--(B2); \draw[line width=0.5cm,color=lightgray] (W3)--(B3); \draw[line width=0.5cm,color=lightgray] (W2)--(B1); 
\draw[line width=0.5cm,color=lightgray] (B4)--(1.5,5); \draw[line width=0.5cm,color=lightgray] (W1)--(1.5,0);
\figsixathreeb
\end{tikzpicture}
 }}; 

\node (PM6) at (10,-3.1) 
{\scalebox{0.408}{
\begin{tikzpicture}
%perfectmatching
\draw[line width=0.5cm,color=lightgray] (W3)--(B1); \draw[line width=0.5cm,color=lightgray] (W2)--(B2); 
\draw[line width=0.5cm,color=lightgray] (B4)--(1.5,5); \draw[line width=0.5cm,color=lightgray] (W1)--(1.5,0); 
\draw[line width=0.5cm,color=lightgray] (B3)--(0,3.5); \draw[line width=0.5cm,color=lightgray] (W4)--(5,3.5);
\figsixathreeb
\end{tikzpicture}
 }} ;  
\end{tikzpicture}
}}
\end{center}

\medskip

Then we have the following table. We denote by $T(a,b,c)$ a divisorial ideal $T(a,b,c,0,0,0)$. 

\medskip

\begin{center}
{\small{
\begin{tabular}{|l|l|l|} \hline 
$e_0\sfC^\prime$&$e_1\sfC^\prime$&$e_2\sfC^\prime$ \\ \hline
$T_{00}\cong R$&$T_{10}\cong T(-1,-2,-1)$&$T_{20}\cong T(1,0,0)$ \\
$T_{01}\cong T(1,2,1)$&$T_{11}\cong R$&$T_{21}\cong T(2,2,1)$ \\
$T_{02}\cong T(-1,0,0)$&$T_{12}\cong T(-2,-2,-1)$&$T_{22}\cong R$ \\
$T_{03}\cong T(1,2,2)$&$T_{13}\cong T(0,0,1)$&$T_{23}\cong T(2,2,2)$ \\
$T_{04}\cong T(0,1,1)$&$T_{14}\cong T(-1,-1,0)$&$T_{24}\cong T(1,1,1)$ \\
$T_{05}\cong T(-1,-1,0)$&$T_{15}\cong T(-2,-3,-1)$&$T_{25}\cong T(0,-1,0)$ \\ \hline
$e_3\sfC^\prime$&$e_4\sfC^\prime$&$e_5\sfC^\prime$ \\ \hline
$T_{30}\cong T(-1,-2,-2)$&$T_{40}\cong T(0,-1,-1)$&$T_{50}\cong T(1,1,0)$ \\
$T_{31}\cong T(0,0,-1)$&$T_{41}\cong T(1,1,0)$&$T_{51}\cong T(2,3,1)$ \\
$T_{32}\cong T(-2,-2,-2)$&$T_{42}\cong T(-1,-1,-1)$&$T_{52}\cong T(0,1,0)$ \\
$T_{33}\cong R$&$T_{43}\cong T(1,1,1)$&$T_{53}\cong T(2,3,2)$\\ 
$T_{34}\cong T(-1,-1,-1)$&$T_{44}\cong R$&$T_{54}\cong T(1,2,1)$\\
$T_{35}\cong T(-2,-3,-2)$&$T_{45}\cong T(-1,-2,-1)$&$T_{55}\cong R$\\ \hline 
\end{tabular} 
}}
\end{center}

%%%%%%%%%%%%
%%%%%%%%%%%%
%%%%%%%%%%%%
\iffalse
%%%%%%%%%%%%
%%%%%%%%%%%%
%%%%%%%%%%%%

\medskip

\begin{center}
{\small{
\begin{tabular}{|l|l|l|} \hline 
$(e_0\sfC^\prime)^*$&$(e_1\sfC^\prime)^*$&$(e_2\sfC^\prime)^*$ \\ \hline
$T_{00}^*\cong R$&$T_{10}^*\cong T(1,2,1)$&$T_{20}^*\cong T(-1,0,0)$ \\
$T_{01}^*\cong T(-1,-2,-1)$&$T_{11}^*\cong R$&$T_{21}^*\cong T(-2,-2,-1)$ \\
$T_{02}^*\cong T(1,0,0)$&$T_{12}^*\cong T(2,2,1)$&$T_{22}^*\cong R$ \\
$T_{03}^*\cong T(-1,-2,-2)$&$T_{13}^*\cong T(0,0,-1)$&$T_{23}^*\cong T(-2,-2,-2)$ \\
$T_{04}^*\cong T(0,-1,-1)$&$T_{14}^*\cong T(1,1,0)$&$T_{24}^*\cong T(-1,-1,-1)$ \\ 
$T_{05}^*\cong T(1,1,0)$&$T_{15}^*\cong T(2,3,1)$&$T_{25}^*\cong T(0,1,0)$ \\ \hline
$(e_3\sfC^\prime)^*$&$(e_4\sfC^\prime)^*$&$(e_5\sfC^\prime)^*$ \\ \hline
$T_{30}^*\cong T(1,2,2)$&$T_{40}^*\cong T(0,1,1)$&$T_{50}^*\cong T(-1,-1,0)$ \\
$T_{31}^*\cong T(0,0,1)$&$T_{41}^*\cong T(-1,-1,0)$&$T_{51}^*\cong T(-2,-3,-1)$ \\
$T_{32}^*\cong T(2,2,2)$&$T_{42}^*\cong T(1,1,1)$&$T_{52}^*\cong T(0,-1,0)$ \\
$T_{33}^*\cong R$&$T_{43}^*\cong T(-1,-1,-1)$&$T_{53}^*\cong T(-2,-3,-2)$\\
$T_{34}^*\cong T(1,1,1)$&$T_{44}^*\cong R$&$T_{54}^*\cong T(-1,-2,-1)$\\
$T_{35}^*\cong T(2,3,2)$&$T_{45}^*\cong T(1,2,1)$&$T_{55}^*\cong R$\\ \hline
\end{tabular} 
}}
\end{center}

%%%%%%%%%%%%
%%%%%%%%%%%%
%%%%%%%%%%%%
\fi
%%%%%%%%%%%%
%%%%%%%%%%%%
%%%%%%%%%%%%

\medskip

Let $x$ (resp. $y$) be a $1$-cycle on the two torus $\sfT$ which corresponds to $(1,0)\in\rmH_1(\sfT)$ (resp. $(0,1)\in\rmH_1(\sfT)$) in the above case. 
By replacing these cycles by $y, x$ respectively, we have the following dimer model and the associated QP $(Q, W_Q)$. 
We denote by $\sfC^{\prime\prime}$ the complete Jacobian algebra $\calP(Q, W_Q)$. 
Since this dimer model is right equivalent to the previous one, we have $\sfC^\prime\cong\sfC^{\prime\prime}$ as an $R$-algebra. 
However, a change of cycles induces an automorphism on $R$, and it gives another description of each splitting MM generator (see Remark~\ref{rem_auto}). 
Thus we also give splitting MM generators $e_i\sfC^{\prime\prime}$.

\newcommand{\figsixathreec}{
%vertex
\node (B1) at (2,0.5){$$}; \node (B2) at (3.5,1.5){$$}; \node (B3) at (0.5,2.5){$$}; \node (B4) at (1.5,3.5){$$};  
\node (W1) at (0.5,0.5){$$}; \node (W2) at (2,1.5){$$}; \node (W3) at (2,2.5){$$}; \node (W4) at (3.5,3){$$}; 
\draw[line width=0.024cm]  (0,0) rectangle (4,4);
%edge
\draw[line width=0.05cm]  (W1)--(B1)--(W2)--(B2)--(W4)--(B4)--(W3)--(B3)--(W1) ; 
\draw[line width=0.05cm]  (W2)--(B3); \draw[line width=0.05cm]  (W3)--(B2); 
\draw[line width=0.05cm]  (W1)--(1,0); \draw[line width=0.05cm]  (W1)--(0,1); 
\draw[line width=0.05cm]  (B1)--(2.5,0); \draw[line width=0.05cm]  (B2)--(4,1); 
\draw[line width=0.05cm]  (B3)--(0,2.75); \draw[line width=0.05cm]  (B4)--(1,4); 
\draw[line width=0.05cm]  (W4)--(2.5,4); \draw[line width=0.05cm]  (W4)--(4,2.75); 
%black
\filldraw  [ultra thick, fill=black] (2,0.5) circle [radius=0.16] ; \filldraw  [ultra thick, fill=black] (3.5,1.5) circle [radius=0.16] ;
\filldraw  [ultra thick, fill=black] (0.5,2.5) circle [radius=0.16] ; \filldraw  [ultra thick, fill=black] (1.5,3.5) circle [radius=0.16] ;
%white
\draw  [line width=1.28pt,fill=white] (0.5,0.5) circle [radius=0.16] ; \draw  [line width=1.28pt, fill=white] (2,1.5) circle [radius=0.16] ;
\draw  [line width=1.28pt,fill=white] (2,2.5)circle [radius=0.16] ; \draw  [line width=1.28pt,fill=white] (3.5,3)circle [radius=0.16] ;
}

\begin{center}
\begin{tikzpicture}
\node (DM) at (0,0) 
{\scalebox{0.675}{
\begin{tikzpicture}
\figsixathreec
\end{tikzpicture}
} };

\node (QV) at (6,0) 
{\scalebox{0.675}{
\begin{tikzpicture}

\node (Q0a) at (0,0){$0$}; \node (Q0b) at (4,0){$0$}; \node (Q0c) at (4,4){$0$}; \node (Q0d) at (0,4){$0$}; 
\node (Q1) at (2,2){$1$}; \node (Q2) at (1.25,1.25){$2$}; \node (Q3) at (2.75,2.65){$3$}; 
\node (Q4a) at (0,2){$4$}; \node (Q4b) at (4,2){$4$}; \node (Q5a) at (2,0){$5$}; \node (Q5b) at (2,4){$5$};
%edge
\draw[lightgray, line width=0.05cm]  (W1)--(B1)--(W2)--(B2)--(W4)--(B4)--(W3)--(B3)--(W1) ; 
\draw[lightgray, line width=0.05cm]  (W2)--(B3); \draw[lightgray, line width=0.05cm]  (W3)--(B2); 
\draw[lightgray, line width=0.05cm]  (W1)--(1,0); \draw[lightgray, line width=0.05cm]  (W1)--(0,1); 
\draw[lightgray, line width=0.05cm]  (B1)--(2.5,0); \draw[lightgray, line width=0.05cm]  (B2)--(4,1); 
\draw[lightgray, line width=0.05cm]  (B3)--(0,2.75); \draw[lightgray, line width=0.05cm]  (B4)--(1,4); 
\draw[lightgray, line width=0.05cm]  (W4)--(2.5,4); \draw[lightgray, line width=0.05cm]  (W4)--(4,2.75); 
%black
\filldraw  [ultra thick, draw=lightgray, fill=lightgray] (2,0.5) circle [radius=0.16] ; \filldraw  [ultra thick, draw=lightgray, fill=lightgray] (3.5,1.5) circle [radius=0.16] ;
\filldraw  [ultra thick, draw=lightgray, fill=lightgray] (0.5,2.5) circle [radius=0.16] ; \filldraw  [ultra thick, draw=lightgray, fill=lightgray] (1.5,3.5) circle [radius=0.16] ;
%white
\draw  [line width=1.28pt, draw=lightgray, fill=white] (0.5,0.5) circle [radius=0.16] ; \draw  [line width=1.28pt, draw=lightgray, fill=white] (2,1.5) circle [radius=0.16] ;
\draw  [line width=1.28pt, draw=lightgray, fill=white] (2,2.5)circle [radius=0.16] ; \draw  [line width=1.28pt, draw=lightgray, fill=white] (3.5,3)circle [radius=0.16] ;
%arrow
\draw[->, line width=0.064cm] (Q0a)--(Q4a); \draw[->, line width=0.064cm] (Q0b)--(Q4b); \draw[->, line width=0.064cm] (Q0c)--(Q4b); 
\draw[->, line width=0.064cm] (Q0d)--(Q4a); \draw[->, line width=0.064cm] (Q0b)--(Q2); \draw[->, line width=0.064cm] (Q0d)--(Q3);
\draw[->, line width=0.064cm] (Q1)--(Q0b); \draw[->, line width=0.064cm] (Q1)--(Q0d); \draw[->, line width=0.064cm] (Q2)--(Q1); 
\draw[->, line width=0.064cm] (Q2)--(Q5a); \draw[->, line width=0.064cm] (Q3)--(Q1); \draw[->, line width=0.064cm] (Q3)--(Q5b); 
\draw[->, line width=0.064cm] (Q4a)--(Q2); \draw[->, line width=0.064cm] (Q4b)--(Q3); \draw[->, line width=0.064cm] (Q5a)--(Q0a); 
\draw[->, line width=0.064cm] (Q5b)--(Q0d); \draw[->, line width=0.064cm] (Q5b)--(Q0c); \draw[->, line width=0.064cm] (Q5a)--(Q0b);
\end{tikzpicture}
} }; 
\end{tikzpicture}
\end{center}

\medskip

In the following figure, we fix the perfect matching $\sfP_0$.  
Then $\sfP_1, \cdots, \sfP_6$ are extremal perfect matchings corresponding to $v_1, \cdots, v_6$. 

\medskip

%%%%%perfect_matching%%%%%

\begin{center}
{\scalebox{0.9}{
\begin{tikzpicture} 
\node at (0,-1.4) {$\sfP_0$};\node at (3.3,-1.4) {$\sfP_1$}; \node at (6.7,-1.4) {$\sfP_2$}; \node at (10,-1.4) {$\sfP_3$}; 
\node at (3.3,-4.5) {$\sfP_4$};\node at (6.7,-4.5) {$\sfP_5$}; \node at (10,-4.5) {$\sfP_6$}; 

\node (PM0) at (0,0) 
{\scalebox{0.51}{
\begin{tikzpicture}
%perfectmatching
\draw[line width=0.4cm,color=lightgray] (W1)--(B1); \draw[line width=0.4cm,color=lightgray] (W2)--(B2);
\draw[line width=0.4cm,color=lightgray] (W3)--(B3); \draw[line width=0.4cm,color=lightgray] (W4)--(B4);
\figsixathreec
\end{tikzpicture}
 }}; 

\node (PM1) at (3.3,0) 
{\scalebox{0.51}{
\begin{tikzpicture}
%perfectmatching
\draw[line width=0.4cm,color=lightgray] (W1)--(B1); \draw[line width=0.4cm,color=lightgray] (W2)--(B2);
\draw[line width=0.4cm,color=lightgray] (W3)--(B4); \draw[line width=0.4cm,color=lightgray] (W4)--(4,2.75); 
\draw[line width=0.4cm,color=lightgray] (B3)--(0,2.75); 
\figsixathreec
\end{tikzpicture}
 }}; 

\node (PM2) at (6.7,0) 
{\scalebox{0.51}{
\begin{tikzpicture}
%perfectmatching
\draw[line width=0.4cm,color=lightgray] (W1)--(B3); \draw[line width=0.4cm,color=lightgray] (W2)--(B2);
\draw[line width=0.4cm,color=lightgray] (W3)--(B4); \draw[line width=0.4cm,color=lightgray] (W4)--(2.5,4); 
\draw[line width=0.4cm,color=lightgray] (B1)--(2.5,0); 
\figsixathreec
\end{tikzpicture}
 }} ;  

\node (PM3) at (10,0) 
{\scalebox{0.51}{
\begin{tikzpicture}
%perfectmatching
\draw[line width=0.4cm,color=lightgray] (W2)--(B3); \draw[line width=0.4cm,color=lightgray] (W3)--(B4);
\draw[line width=0.4cm,color=lightgray] (W1)--(0,1); \draw[line width=0.4cm,color=lightgray] (B2)--(4,1); 
\draw[line width=0.4cm,color=lightgray] (W4)--(2.5,4); \draw[line width=0.4cm,color=lightgray] (B1)--(2.5,0); 
\figsixathreec
\end{tikzpicture}
 }} ; 

\node (PM4) at (3.3,-3.1) 
{\scalebox{0.51}{
\begin{tikzpicture}
%perfectmatching
\draw[line width=0.4cm,color=lightgray] (W2)--(B1); \draw[line width=0.4cm,color=lightgray] (W3)--(B3);
\draw[line width=0.4cm,color=lightgray] (W4)--(B4);  
\draw[line width=0.4cm,color=lightgray] (W1)--(0,1); \draw[line width=0.4cm,color=lightgray] (B2)--(4,1);
\figsixathreec
\end{tikzpicture}
 }}; 

\node (PM5) at (6.7,-3.1) 
{\scalebox{0.51}{
\begin{tikzpicture}
%perfectmatching
\draw[line width=0.4cm,color=lightgray] (W2)--(B1); \draw[line width=0.4cm,color=lightgray] (W4)--(B2);
\draw[line width=0.4cm,color=lightgray] (W3)--(B3);  
\draw[line width=0.4cm,color=lightgray] (W1)--(1,0); \draw[line width=0.4cm,color=lightgray] (B4)--(1,4);
\figsixathreec
\end{tikzpicture}
 }}; 

\node (PM6) at (10,-3.1) 
{\scalebox{0.51}{
\begin{tikzpicture}
%perfectmatching
\draw[line width=0.4cm,color=lightgray] (W2)--(B1); \draw[line width=0.4cm,color=lightgray] (W3)--(B2); 
\draw[line width=0.4cm,color=lightgray] (W1)--(1,0); \draw[line width=0.4cm,color=lightgray] (B4)--(1,4);
\draw[line width=0.4cm,color=lightgray] (W4)--(4,2.75); \draw[line width=0.4cm,color=lightgray] (B3)--(0,2.75);
\figsixathreec
\end{tikzpicture}
 }} ;  

\end{tikzpicture}
}}
\end{center}

\medskip

Then we have the following table. We denote by $T(a,b,c)$ a divisorial ideal $T(a,b,c,0,0,0)$. 

\medskip

\begin{center}
{\small{
\begin{tabular}{|l|l|l|} \hline 
$e_0\sfC^{\prime\prime}$&$e_1\sfC^{\prime\prime}$&$e_2\sfC^{\prime\prime}$\\ \hline
$T_{00}\cong R$&$T_{10}\cong T(1,1,0)$&$T_{20}\cong T(1,1,1)$ \\
$T_{01}\cong T(-1,-1,0)$&$T_{11}\cong R$&$T_{21}\cong T(0,0,1)$ \\
$T_{02}\cong T(-1,-1,-1)$&$T_{12}\cong T(0,0,-1)$&$T_{22}\cong R$ \\
$T_{03}\cong T(1,1,1)$&$T_{13}\cong T(2,2,1)$&$T_{23}\cong T(2,2,2)$ \\
$T_{04}\cong T(-1,-2,-1)$&$T_{14}\cong T(0,-1,-1)$&$T_{24}\cong T(0,-1,0)$ \\
$T_{05}\cong T(0,-1,-1)$&$T_{15}\cong T(1,0,-1)$&$T_{25}\cong T(1,0,0)$ \\ \hline
$e_3\sfC^{\prime\prime}$&$e_4\sfC^{\prime\prime}$&$e_5\sfC^{\prime\prime}$ \\ \hline
$T_{30}\cong T(-1,-1,-1)$&$T_{40}\cong T(1,2,1)$&$T_{50}\cong T(0,1,1)$ \\
$T_{31}\cong T(-2,-2,-1)$&$T_{41}\cong T(0,1,1)$&$T_{51}\cong T(-1,0,1)$ \\
$T_{32}\cong T(-2,-2,-2)$&$T_{42}\cong T(0,1,0)$&$T_{52}\cong T(-1,0,0)$ \\
$T_{33}\cong R$&$T_{43}\cong T(2,3,2)$&$T_{53}\cong T(1,2,2)$\\
$T_{34}\cong T(-2,-3,-2)$&$T_{44}\cong R$&$T_{54}\cong T(-1,-1,0)$\\
$T_{35}\cong T(-1,-2,-2)$&$T_{45}\cong T(1,1,0)$&$T_{55}\cong R$\\ \hline
\end{tabular} 
}}
\end{center}

%%%%%%%%%%%%
%%%%%%%%%%%%
%%%%%%%%%%%%
\iffalse
%%%%%%%%%%%%
%%%%%%%%%%%%
%%%%%%%%%%%%

\medskip

\begin{center}
{\small{
\begin{tabular}{|l|l|l|} \hline 
$(e_0\sfC^{\prime\prime})^*$&$(e_1\sfC^{\prime\prime})^*$&$(e_2\sfC^{\prime\prime})^*$ \\ \hline
$T_{00}^*\cong R$&$T_{10}^*\cong T(-1,-1,0)$&$T_{20}^*\cong T(-1,-1,-1)$ \\
$T_{01}^*\cong T(1,1,0)$&$T_{11}^*\cong R$&$T_{21}^*\cong T(0,0,-1)$ \\
$T_{02}^*\cong T(1,1,1)$&$T_{12}^*\cong T(0,0,1)$&$T_{22}^*\cong R$ \\
$T_{03}^*\cong T(-1,-1,-1)$&$T_{13}^*\cong T(-2,-2,-1)$&$T_{23}^*\cong T(-2,-2,-2)$ \\
$T_{04}^*\cong T(1,2,1)$&$T_{14}^*\cong T(0,1,1)$&$T_{24}^*\cong T(0,1,0)$ \\ 
$T_{05}^*\cong T(0,1,1)$&$T_{15}^*\cong T(-1,0,1)$&$T_{25}^*\cong T(-1,0,0)$ \\ \hline
$(e_3\sfC^{\prime\prime})^*$&$(e_4\sfC^{\prime\prime})^*$&$(e_5\sfC^{\prime\prime})^*$ \\ \hline
$T_{30}^*\cong T(1,1,1)$&$T_{40}^*\cong T(-1,-2,-1)$&$T_{50}^*\cong T(0,-1,-1)$ \\
$T_{31}^*\cong T(2,2,1)$&$T_{41}^*\cong T(0,-1,-1)$&$T_{51}^*\cong T(1,0,-1)$ \\
$T_{32}^*\cong T(2,2,2)$&$T_{42}^*\cong T(0,-1,0)$&$T_{52}^*\cong T(1,0,0)$ \\
$T_{33}^*\cong R$&$T_{43}^*\cong T(-2,-3,-2)$&$T_{53}^*\cong T(-1,-2,-2)$\\
$T_{34}^*\cong T(2,3,2)$&$T_{44}^*\cong R$&$T_{54}^*\cong T(1,1,0)$\\
$T_{35}^*\cong T(1,2,2)$&$T_{45}^*\cong T(-1,-1,0)$&$T_{55}^*\cong R$\\ \hline
\end{tabular} 
}}
\end{center}

%%%%%%%%%%%%
%%%%%%%%%%%%
%%%%%%%%%%%%
\fi
%%%%%%%%%%%%
%%%%%%%%%%%%
%%%%%%%%%%%%

\medskip

%%%%%%%%%%%%%%%%%%%%%%%%%%%%%%%%%%%
%%%%%%%%%%%%Type_6a-4%%%%%%%%%%%%%%%%%%
%%%%%%%%%%%%%%%%%%%%%%%%%%%%%%%%%%%
\subsubsection{Type 6a-4}
Finally, we consider the following consistent dimer model and the associated QP $(Q, W_Q)$. 
We denote by $\sfD$ the complete Jacobian algebra $\calP(Q, W_Q)$, and give splitting MM generators $e_i\sfD$. 

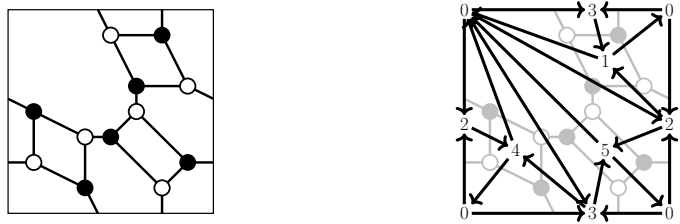
\newcommand{\figsixafour}{
%vertex
\node (B1) at (3,1){$$}; \node (B2) at (1,4){$$}; \node (B3) at (4,3){$$}; \node (B4) at (7,2){$$}; 
\node (B5) at (5,5){$$}; \node (B6) at (6,7){$$}; 
\node (W1) at (1,2){$$}; \node (W2) at (3,3){$$}; \node (W3) at (6,1){$$}; \node (W4) at (5,4){$$}; 
\node (W5) at (7,5){$$}; \node (W6) at (4,7){$$}; 
\draw[line width=0.048cm]  (0,0) rectangle (8,8);
%edge
\draw[line width=0.096cm]  (W6)--(B5)--(W4)--(B3)--(W2)--(B2) ; 
\draw[line width=0.096cm]  (B2)--(W1)--(B1)--(W2); \draw[line width=0.096cm]  (B3)--(W3)--(B4)--(W4); \draw[line width=0.096cm]  (B5)--(W5)--(B6)--(W6); 
\draw[line width=0.096cm]  (B1)--(3.5,0); \draw[line width=0.096cm]  (W3)--(6,0); \draw[line width=0.096cm]  (W1)--(0,2); 
\draw[line width=0.096cm]  (B4)--(8,2); \draw[line width=0.096cm]  (B2)--(0,4.5); \draw[line width=0.096cm]  (W5)--(8,4.5); 
\draw[line width=0.096cm]  (B6)--(6,8); \draw[line width=0.096cm]  (W6)--(3.5,8); 
%black
\filldraw  [ultra thick, fill=black] (3,1) circle [radius=0.3] ; \filldraw  [ultra thick, fill=black] (1,4) circle [radius=0.3] ;
\filldraw  [ultra thick, fill=black] (4,3) circle [radius=0.3] ; \filldraw  [ultra thick, fill=black] (7,2) circle [radius=0.3] ;
\filldraw  [ultra thick, fill=black] (5,5) circle [radius=0.3] ; \filldraw  [ultra thick, fill=black] (6,7) circle [radius=0.3] ;
%white
\draw  [line width=2pt,fill=white] (1,2) circle [radius=0.3] ; \draw  [line width=2pt, fill=white] (3,3) circle [radius=0.3] ;
\draw  [line width=2pt,fill=white] (6,1)circle [radius=0.3] ; \draw  [line width=2pt,fill=white] (5,4)circle [radius=0.3] ; 
\draw  [line width=2pt,fill=white] (7,5)circle [radius=0.3] ; \draw  [line width=2pt,fill=white] (4,7)circle [radius=0.3] ;
}

\begin{center}
\begin{tikzpicture}

\node (DM) at (0,0) 
{\scalebox{0.3375}{
\begin{tikzpicture}
\figsixafour
\end{tikzpicture}
} }; 

\node (QV) at (6,0) 
{\scalebox{0.3375}{
\begin{tikzpicture}

\node (Q0a) at (0,0){{\Huge$0$}}; \node (Q0b) at (8,0){{\Huge$0$}}; \node (Q0c) at (8,8){{\Huge$0$}}; \node (Q0d) at (0,8){{\Huge$0$}}; 
\node (Q1) at (5.5,6){{\Huge$1$}}; \node (Q2a) at (0,3.5){{\Huge$2$}}; \node (Q2b) at (8,3.5){{\Huge$2$}}; 
\node (Q3a) at (5,0){{\Huge$3$}}; \node (Q3b) at (5,8){{\Huge$3$}}; \node (Q4) at (2,2.5){{\Huge$4$}}; \node (Q5) at (5.5,2.5){{\Huge$5$}}; 
%edge
\draw[lightgray, line width=0.096cm]  (W6)--(B5)--(W4)--(B3)--(W2)--(B2) ; 
\draw[lightgray, line width=0.096cm]  (B2)--(W1)--(B1)--(W2); \draw[lightgray, line width=0.096cm]  (B3)--(W3)--(B4)--(W4); 
\draw[lightgray, line width=0.096cm]  (B5)--(W5)--(B6)--(W6); 
\draw[lightgray, line width=0.096cm]  (B1)--(3.5,0); \draw[lightgray, line width=0.096cm]  (W3)--(6,0); \draw[lightgray, line width=0.096cm]  (W1)--(0,2); 
\draw[lightgray, line width=0.096cm]  (B4)--(8,2); \draw[lightgray, line width=0.096cm]  (B2)--(0,4.5); \draw[lightgray, line width=0.096cm]  (W5)--(8,4.5); 
\draw[lightgray, line width=0.096cm]  (B6)--(6,8); \draw[lightgray, line width=0.096cm]  (W6)--(3.5,8); 

%black
\filldraw  [ultra thick, draw=lightgray, fill=lightgray] (3,1) circle [radius=0.3] ; \filldraw  [ultra thick, draw=lightgray, fill=lightgray] (1,4) circle [radius=0.3] ;
\filldraw  [ultra thick, draw=lightgray, fill=lightgray] (4,3) circle [radius=0.3] ; \filldraw  [ultra thick, draw=lightgray, fill=lightgray] (7,2) circle [radius=0.3] ;
\filldraw  [ultra thick, draw=lightgray, fill=lightgray] (5,5) circle [radius=0.3] ; \filldraw  [ultra thick, draw=lightgray, fill=lightgray] (6,7) circle [radius=0.3] ;
%white
\draw  [line width=2pt, draw=lightgray, fill=white] (1,2) circle [radius=0.3] ; \draw  [line width=2pt, draw=lightgray, fill=white] (3,3) circle [radius=0.3] ;
\draw  [line width=2pt, draw=lightgray, fill=white] (6,1) circle [radius=0.3] ; \draw  [line width=2pt, draw=lightgray, fill=white] (5,4) circle [radius=0.3] ; 
\draw  [line width=2pt, draw=lightgray, fill=white] (7,5) circle [radius=0.3] ; \draw  [line width=2pt, draw=lightgray, fill=white] (4,7) circle [radius=0.3] ;
%arrow
\draw[->, line width=0.12cm] (Q0a)--(Q2a); \draw[->, line width=0.12cm] (Q0b)--(Q2b); \draw[->, line width=0.12cm] (Q0a)--(Q3a);
\draw[->, line width=0.12cm] (Q0d)--(Q3b); \draw[->, line width=0.12cm] (Q0d)--(Q2a); \draw[->, line width=0.12cm] (Q0c)--(Q2b); 
\draw[->, line width=0.12cm] (Q1)--(Q0d); \draw[->, line width=0.12cm] (Q1)--(Q0c); \draw[->, line width=0.12cm] (Q0d)--(Q2b);
\draw[->, line width=0.12cm] (Q0b)--(Q3a); \draw[->, line width=0.12cm] (Q0c)--(Q3b); \draw[->, line width=0.12cm] (Q0d)--(Q3a);
\draw[->, line width=0.12cm] (Q2a)--(Q4); \draw[->, line width=0.12cm] (Q2b)--(Q1); \draw[->, line width=0.12cm] (Q2b)--(Q5);
\draw[->, line width=0.12cm] (Q3a)--(Q4); \draw[->, line width=0.12cm] (Q3a)--(Q5); \draw[->, line width=0.12cm] (Q3b)--(Q1);
\draw[->, line width=0.12cm] (Q4)--(Q0a); \draw[->, line width=0.12cm] (Q4)--(Q0d); 
\draw[->, line width=0.12cm] (Q5)--(Q0b); \draw[->, line width=0.12cm] (Q5)--(Q0d);
\end{tikzpicture}
} }; 
\end{tikzpicture}
\end{center}

\medskip

In the following figure, we fix the perfect matching $\sfP_0$.  
Then $\sfP_1, \cdots, \sfP_6$ are extremal perfect matchings corresponding to $v_1, \cdots, v_6$. 

\medskip

%%%%%perfect_matching%%%%%

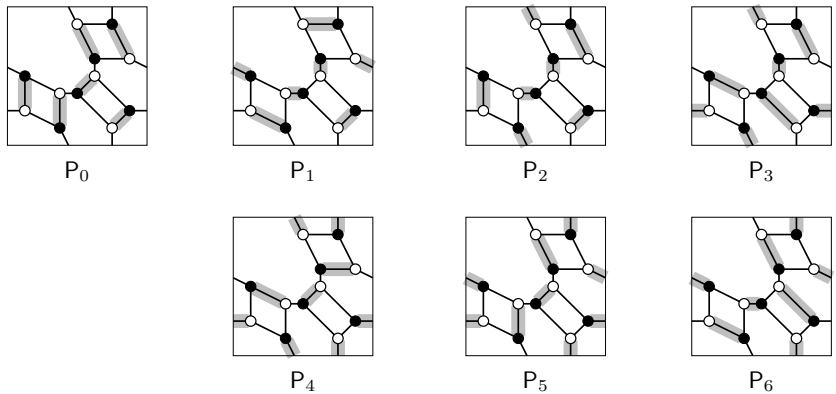
\begin{center}
{\scalebox{0.9}{
\begin{tikzpicture} 
\node at (0,-1.4) {$\sfP_0$};\node at (3.3,-1.4) {$\sfP_1$}; \node at (6.7,-1.4) {$\sfP_2$}; \node at (10,-1.4) {$\sfP_3$}; 
\node at (3.3,-4.5) {$\sfP_4$};\node at (6.7,-4.5) {$\sfP_5$}; \node at (10,-4.5) {$\sfP_6$}; 

\node (PM0) at (0,0) 
{\scalebox{0.255}{
\begin{tikzpicture}
%perfectmatching
\draw[line width=0.8cm,color=lightgray] (W1)--(B2); \draw[line width=0.8cm,color=lightgray] (W2)--(B1);
\draw[line width=0.8cm,color=lightgray] (W3)--(B4); \draw[line width=0.8cm,color=lightgray] (W4)--(B3);
\draw[line width=0.8cm,color=lightgray] (W5)--(B6); \draw[line width=0.8cm,color=lightgray] (W6)--(B5);
\figsixafour
\end{tikzpicture}
 }}; 

\node (PM1) at (3.3,0) 
{\scalebox{0.255}{
\begin{tikzpicture}
%perfectmatching
\draw[line width=0.8cm,color=lightgray] (W1)--(B1); \draw[line width=0.8cm,color=lightgray] (W2)--(B3);
\draw[line width=0.8cm,color=lightgray] (W3)--(B4); \draw[line width=0.8cm,color=lightgray] (W4)--(B5);
\draw[line width=0.8cm,color=lightgray] (W6)--(B6); \draw[line width=0.8cm,color=lightgray] (W5)--(8,4.5); 
\draw[line width=0.8cm,color=lightgray] (B2)--(0,4.5);
\figsixafour
\end{tikzpicture}
 }}; 

\node (PM2) at (6.7,0) 
{\scalebox{0.255}{
\begin{tikzpicture}
%perfectmatching
\draw[line width=0.8cm,color=lightgray] (W1)--(B2); \draw[line width=0.8cm,color=lightgray] (W2)--(B3);
\draw[line width=0.8cm,color=lightgray] (W3)--(B4); \draw[line width=0.8cm,color=lightgray] (W4)--(B5);
\draw[line width=0.8cm,color=lightgray] (W5)--(B6); \draw[line width=0.8cm,color=lightgray] (W6)--(3.5,8); 
\draw[line width=0.8cm,color=lightgray] (B1)--(3.5,0);
\figsixafour
\end{tikzpicture}
 }} ;  

\node (PM3) at (10,0) 
{\scalebox{0.255}{
\begin{tikzpicture}
%perfectmatching
\draw[line width=0.8cm,color=lightgray] (W2)--(B2); \draw[line width=0.8cm,color=lightgray] (W3)--(B3);
\draw[line width=0.8cm,color=lightgray] (W4)--(B5); \draw[line width=0.8cm,color=lightgray] (W5)--(B6); 
\draw[line width=0.8cm,color=lightgray] (W6)--(3.5,8); \draw[line width=0.8cm,color=lightgray] (B1)--(3.5,0);
\draw[line width=0.8cm,color=lightgray] (W1)--(0,2); \draw[line width=0.8cm,color=lightgray] (B4)--(8,2);
\figsixafour
\end{tikzpicture}
 }} ; 

\node (PM4) at (3.3,-3.1) 
{\scalebox{0.255}{
\begin{tikzpicture}
%perfectmatching
\draw[line width=0.8cm,color=lightgray] (W2)--(B2); \draw[line width=0.8cm,color=lightgray] (W4)--(B3);
\draw[line width=0.8cm,color=lightgray] (W5)--(B5);  
\draw[line width=0.8cm,color=lightgray] (W6)--(3.5,8); \draw[line width=0.8cm,color=lightgray] (B1)--(3.5,0);
\draw[line width=0.8cm,color=lightgray] (W1)--(0,2); \draw[line width=0.8cm,color=lightgray] (B4)--(8,2);
\draw[line width=0.8cm,color=lightgray] (W3)--(6,0); \draw[line width=0.8cm,color=lightgray] (B6)--(6,8);
\figsixafour
\end{tikzpicture}
 }}; 

\node (PM5) at (6.7,-3.1) 
{\scalebox{0.255}{
\begin{tikzpicture}
%perfectmatching
\draw[line width=0.8cm,color=lightgray] (W2)--(B1); \draw[line width=0.8cm,color=lightgray] (W4)--(B3);
\draw[line width=0.8cm,color=lightgray] (W6)--(B5);  
\draw[line width=0.8cm,color=lightgray] (W5)--(8,4.5); \draw[line width=0.8cm,color=lightgray] (B2)--(0,4.5);
\draw[line width=0.8cm,color=lightgray] (W1)--(0,2); \draw[line width=0.8cm,color=lightgray] (B4)--(8,2);
\draw[line width=0.8cm,color=lightgray] (W3)--(6,0); \draw[line width=0.8cm,color=lightgray] (B6)--(6,8);
\figsixafour
\end{tikzpicture}
 }}; 

\node (PM6) at (10,-3.1) 
{\scalebox{0.255}{
\begin{tikzpicture}
%perfectmatching
\draw[line width=0.8cm,color=lightgray] (W1)--(B1); \draw[line width=0.8cm,color=lightgray] (W2)--(B3);
\draw[line width=0.8cm,color=lightgray] (W4)--(B4); \draw[line width=0.8cm,color=lightgray] (W6)--(B5); 
\draw[line width=0.8cm,color=lightgray] (W5)--(8,4.5); \draw[line width=0.8cm,color=lightgray] (B2)--(0,4.5);
\draw[line width=0.8cm,color=lightgray] (W3)--(6,0); \draw[line width=0.8cm,color=lightgray] (B6)--(6,8);
\figsixafour
\end{tikzpicture}
 }} ;  
\end{tikzpicture}
}}
\end{center}

\medskip

Then we have the following table. We denote by $T(a,b,c)$ a divisorial ideal $T(a,b,c,0,0,0)$. 

\medskip

\begin{center}
{\small{
\begin{tabular}{|l|l|l|} \hline 
$e_0\sfD$&$e_1\sfD$&$e_2\sfD$\\ \hline
$T_{00}\cong R$&$T_{10}\cong T(0,1,1)$&$T_{20}\cong T(-1,-1,-1)$ \\
$T_{01}\cong T(0,-1,-1)$&$T_{11}\cong R$&$T_{21}\cong T(-1,-2,-2)$ \\
$T_{02}\cong T(1,1,1)$&$T_{12}\cong T(1,2,2)$&$T_{22}\cong R$ \\
$T_{03}\cong T(-1,-1,-1)$&$T_{13}\cong T(-1,0,0)$&$T_{23}\cong T(-2,-2,-2)$ \\
$T_{04}\cong T(1,2,1)$&$T_{14}\cong T(1,3,2)$&$T_{24}\cong T(0,1,0)$ \\
$T_{05}\cong T(-1,-1,0)$&$T_{15}\cong T(-1,0,1)$&$T_{25}\cong T(-2,-2,-1)$ \\ \hline
$e_3\sfD$&$e_4\sfD$&$e_5\sfD$ \\ \hline
$T_{30}\cong T(1,1,1)$&$T_{40}\cong T(-1,-2,-1)$&$T_{50}\cong T(1,1,0)$ \\
$T_{31}\cong T(1,0,0)$&$T_{41}\cong T(-1,-3,-2)$&$T_{51}\cong T(1,0,-1)$ \\
$T_{32}\cong T(2,2,2)$&$T_{42}\cong T(0,-1,0)$&$T_{52}\cong T(2,2,1)$ \\
$T_{33}\cong R$&$T_{43}\cong T(-2,-3,-2)$&$T_{53}\cong T(0,0,-1)$\\
$T_{34}\cong T(2,3,2)$&$T_{44}\cong R$&$T_{54}\cong T(2,3,1)$\\
$T_{35}\cong T(0,0,1)$&$T_{45}\cong T(-2,-3,-1)$&$T_{55}\cong R$\\ \hline
\end{tabular} 
}}
\end{center}

%%%%%%%%%%%%
%%%%%%%%%%%%
%%%%%%%%%%%%
\iffalse
%%%%%%%%%%%%
%%%%%%%%%%%%
%%%%%%%%%%%%

\medskip

\begin{center}
{\small{
\begin{tabular}{|l|l|l|} \hline 
$(e_0\sfD)^*$&$(e_1\sfD)^*$&$(e_2\sfD)^*$ \\ \hline
$T_{00}^*\cong R$&$T_{10}^*\cong T(0,-1,-1)$&$T_{20}^*\cong T(1,1,1)$ \\
$T_{01}^*\cong T(0,1,1)$&$T_{11}^*\cong R$&$T_{21}^*\cong T(1,2,2)$ \\
$T_{02}^*\cong T(-1,-1,-1)$&$T_{12}^*\cong T(-1,-2,-2)$&$T_{22}^*\cong R$ \\
$T_{03}^*\cong T(1,1,1)$&$T_{13}^*\cong T(1,0,0)$&$T_{23}^*\cong T(2,2,2)$ \\
$T_{04}^*\cong T(-1,-2,-1)$&$T_{14}^*\cong T(-1,-3,-2)$&$T_{24}^*\cong T(0,-1,0)$ \\ 
$T_{05}^*\cong T(1,1,0)$&$T_{15}^*\cong T(1,0,-1)$&$T_{25}^*\cong T(2,2,1)$ \\ \hline
$(e_3\sfD)^*$&$(e_4\sfD)^*$&$(e_5\sfD)^*$ \\ \hline
$T_{30}^*\cong T(-1,-1,-1)$&$T_{40}^*\cong T(1,2,1)$&$T_{50}^*\cong T(-1,-1,0)$ \\
$T_{31}^*\cong T(-1,0,0)$&$T_{41}^*\cong T(1,3,2)$&$T_{51}^*\cong T(-1,0,1)$ \\
$T_{32}^*\cong T(-2,-2,-2)$&$T_{42}^*\cong T(0,1,0)$&$T_{52}^*\cong T(-2,-2,-1)$ \\
$T_{33}^*\cong R$&$T_{43}^*\cong T(2,3,2)$&$T_{53}^*\cong T(0,0,1)$\\
$T_{34}^*\cong T(-2,-3,-2)$&$T_{44}^*\cong R$&$T_{54}^*\cong T(-2,-3,-1)$\\
$T_{35}^*\cong T(0,0,-1)$&$T_{45}^*\cong T(2,3,1)$&$T_{55}^*\cong R$\\ \hline
\end{tabular} 
}}
\end{center}

%%%%%%%%%%%%
%%%%%%%%%%%%
%%%%%%%%%%%%
\fi
%%%%%%%%%%%%
%%%%%%%%%%%%
%%%%%%%%%%%%

\medskip

%%%%%%%%%%%%%%%%%%%%%%%%%%%%%%%%%%%%%%%%%%%%
%%%%%%%%%%%%%%%%%%%%%%%%%%%%%%%%%%%%%%%%%%%%
\newpage
\subsubsection{ \textup{Exchange graph of type 6a}} 
By combining the above results, we describe the exchange graph $\EG(\TMMG(R))$ for the case of type 6a. 
To make the graph clear, we divide it into two parts. 
Identifying common vertices surrounded by dashed squares in the following figures, we obtain $\EG(\TMMG(R))$ and it is actually connected. 

%%%%%%%%%%%%%%%%%%%%%
%%% type 6a exchange graph 1%%%%%%%%%%%%%
%%%%%%%%%%%%%%%%%%%%%
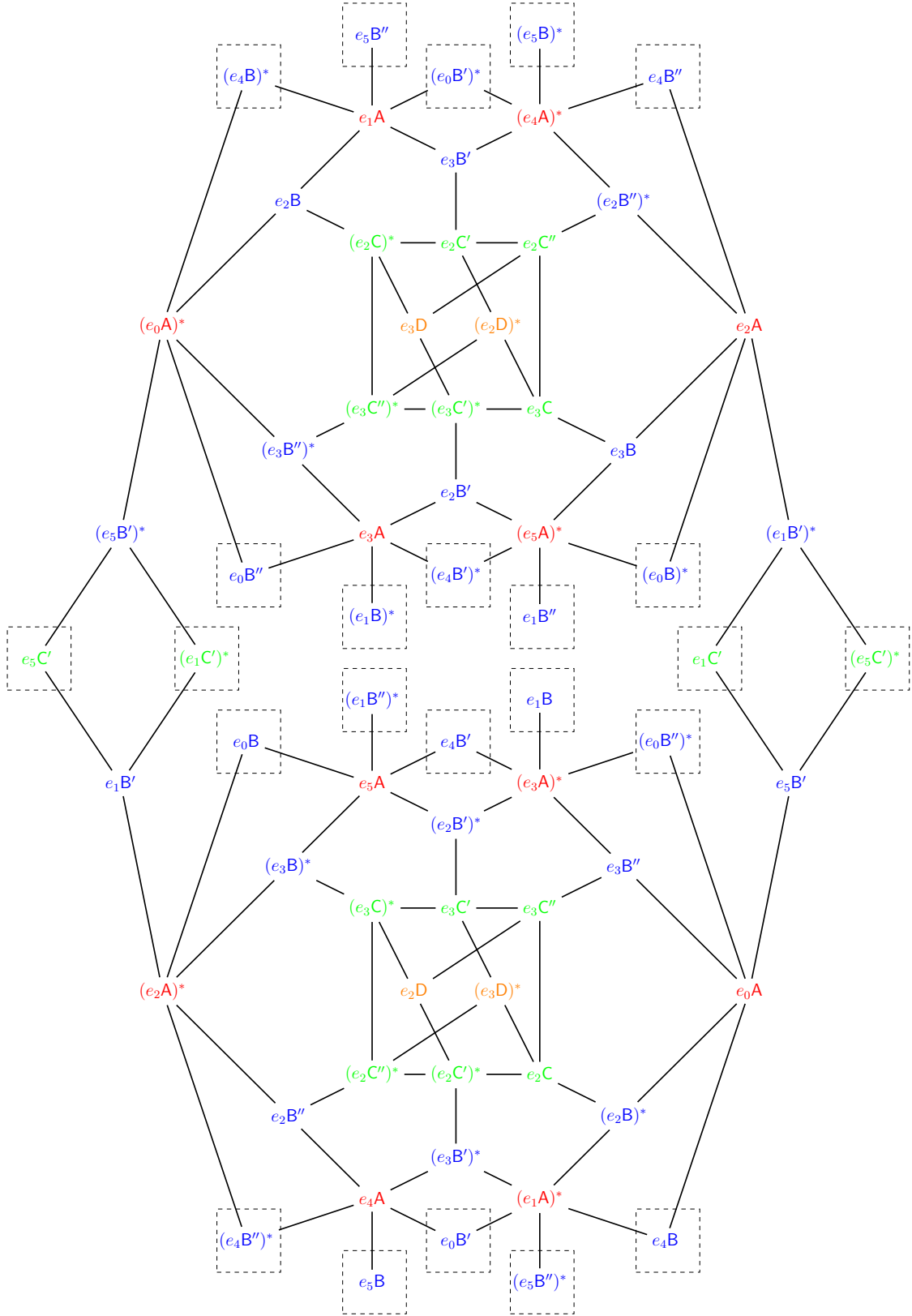
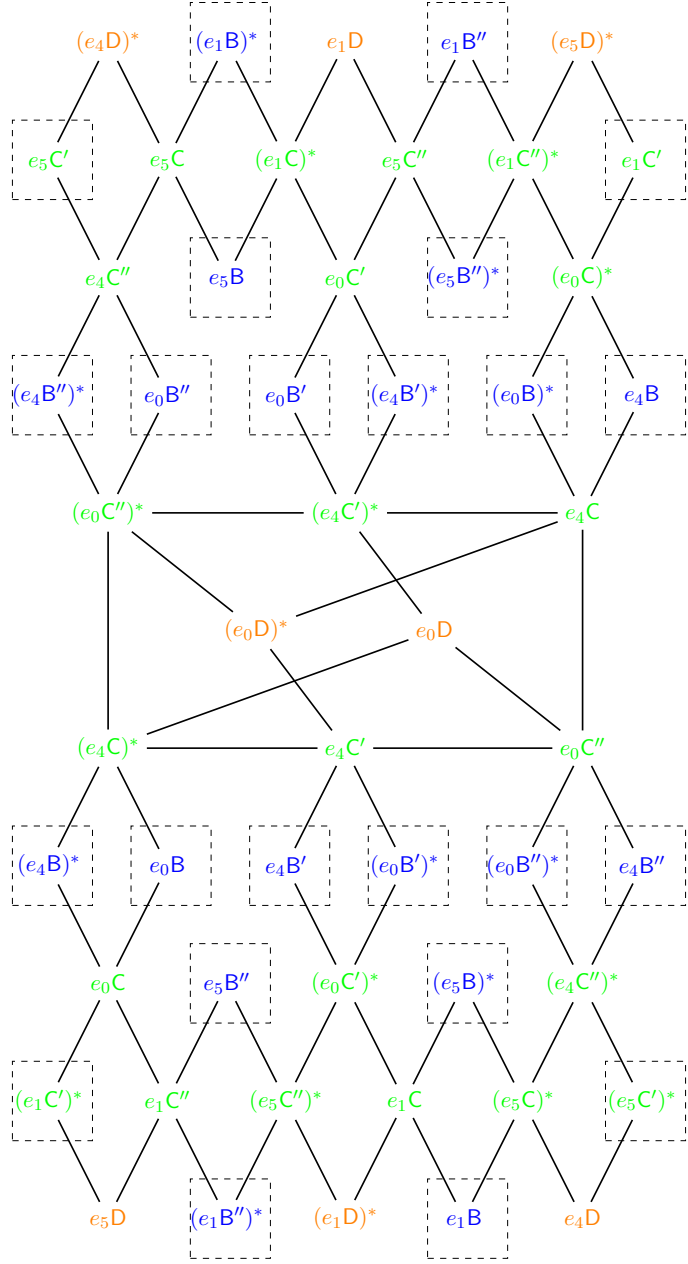
\begin{center}
\newcommand{\myrect} [1] [] {
\begin{tikzpicture}
\draw [dashed]  (0,0) rectangle (1.52,1.52) ;
\end{tikzpicture}}

{\scalebox{0.68}{
\begin{tikzpicture}
\node (A0) [red] at (7,-8) {{\Large$e_0\sfA$}}; 
\node (A1) [red] at (-2,13) {{\Large$e_1\sfA$}}; 
\node (A2) [red] at (7,8) {{\Large$e_2\sfA$}}; 
\node (A3) [red] at (-2,3) {{\Large$e_3\sfA$}}; 
\node (A4) [red] at (-2,-13) {{\Large$e_4\sfA$}}; 
\node (A5) [red] at (-2,-3) {{\Large$e_5\sfA$}}; 
\node (A0+) [red] at (-7,8) {{\Large$(e_0\sfA)^*$}}; 
\node (A1+) [red] at (2,-13) {{\Large$(e_1\sfA)^*$}}; 
\node (A2+) [red] at (-7,-8) {{\Large$(e_2\sfA)^*$}}; 
\node (A3+) [red] at (2,-3) {{\Large$(e_3\sfA)^*$}}; 
\node (A4+) [red] at (2,13) {{\Large$(e_4\sfA)^*$}}; 
\node (A5+) [red] at (2,3) {{\Large$(e_5\sfA)^*$}}; 

\node (B0) [blue] at (-5,-2) {{\Large$e_0\sfB$}}; 
\node (B1) [blue] at (2,-1) {{\Large$e_1\sfB$}}; 
\node (B2) [blue] at (-4,11) {{\Large$e_2\sfB$}}; 
\node (B3) [blue] at (4,5) {{\Large$e_3\sfB$}}; 
\node (B4) [blue] at (5,-14) {{\Large$e_4\sfB$}}; 
\node (B5) [blue] at (-2,-15) {{\Large$e_5\sfB$}}; 
\node (B0+) [blue] at (5,2) {{\Large$(e_0\sfB)^*$}}; 
\node (B1+) [blue] at (-2,1) {{\Large$(e_1\sfB)^*$}}; 
\node (B2+) [blue] at (4,-11) {{\Large$(e_2\sfB)^*$}}; 
\node (B3+) [blue] at (-4,-5) {{\Large$(e_3\sfB)^*$}}; 
\node (B4+) [blue] at (-5,14) {{\Large$(e_4\sfB)^*$}}; 
\node (B5+) [blue] at (2,15) {{\Large$(e_5\sfB)^*$}}; 

\node (BB0) [blue] at (0,-14) {{\Large$e_0\sfB^\prime$}}; 
\node (BB1) [blue] at (-8,-3) {{\Large$e_1\sfB^\prime$}}; 
\node (BB2) [blue] at (0,4) {{\Large$e_2\sfB^\prime$}}; 
\node (BB3) [blue] at (0,12) {{\Large$e_3\sfB^\prime$}}; 
\node (BB4) [blue] at (0,-2) {{\Large$e_4\sfB^\prime$}}; 
\node (BB5) [blue] at (8,-3) {{\Large$e_5\sfB^\prime$}}; 
\node (BB0+) [blue] at (0,14) {{\Large$(e_0\sfB^\prime)^*$}}; 
\node (BB1+) [blue] at (8,3) {{\Large$(e_1\sfB^\prime)^*$}}; 
\node (BB2+) [blue] at (0,-4) {{\Large$(e_2\sfB^\prime)^*$}}; 
\node (BB3+) [blue] at (0,-12) {{\Large$(e_3\sfB^\prime)^*$}}; 
\node (BB4+) [blue] at (0,2) {{\Large$(e_4\sfB^\prime)^*$}}; 
\node (BB5+) [blue] at (-8,3) {{\Large$(e_5\sfB^\prime)^*$}}; 

\node (BBB0) [blue] at (-5,2) {{\Large$e_0\sfB^{\prime\prime}$}}; 
\node (BBB1) [blue] at (2,1) {{\Large$e_1\sfB^{\prime\prime}$}}; 
\node (BBB2) [blue] at (-4,-11) {{\Large$e_2\sfB^{\prime\prime}$}}; 
\node (BBB3) [blue] at (4,-5) {{\Large$e_3\sfB^{\prime\prime}$}}; 
\node (BBB4) [blue] at (5,14) {{\Large$e_4\sfB^{\prime\prime}$}}; 
\node (BBB5) [blue] at (-2,15) {{\Large$e_5\sfB^{\prime\prime}$}}; 
\node (BBB0+) [blue] at (5,-2) {{\Large$(e_0\sfB^{\prime\prime})^*$}}; 
\node (BBB1+) [blue] at (-2,-1) {{\Large$(e_1\sfB^{\prime\prime})^*$}}; 
\node (BBB2+) [blue] at (4,11) {{\Large$(e_2\sfB^{\prime\prime})^*$}}; 
\node (BBB3+) [blue] at (-4,5) {{\Large$(e_3\sfB^{\prime\prime})^*$}}; 
\node (BBB4+) [blue] at (-5,-14) {{\Large$(e_4\sfB^{\prime\prime})^*$}}; 
\node (BBB5+) [blue] at (2,-15) {{\Large$(e_5\sfB^{\prime\prime})^*$}}; 

%\node (C0) [green] at (6,-9) {$e_0\sfC$}; 
%\node (C1) [green] at (5,0) {$e_1\sfC$}; 
\node (C2) [green] at (2,-10) {{\Large$e_2\sfC$}}; 
\node (C3) [green] at (2,6) {{\Large$e_3\sfC$}}; 
%\node (C4) [green] at (6,9) {$e_4\sfC$}; 
%\node (C5) [green] at (-5,0) {$e_5\sfC$}; 
%\node (C0+) [green] at (6,-9) {$(e_0\sfC)^*$}; 
%\node (C1+) [green] at (5,0) {$(e_1\sfC)^*$}; 
\node (C2+) [green] at (-2,10) {{\Large$(e_2\sfC)^*$}}; 
\node (C3+) [green] at (-2,-6) {{\Large$(e_3\sfC)^*$}}; 
%\node (C4+) [green] at (6,9) {$(e_4\sfC)^*$}; 
%\node (C5+) [green] at (-5,0) {$(e_5\sfC)^*$}; 

%\node (CC0) [green] at (6,-9) {$e_0\sfC^\prime$}; 
\node (CC1) [green] at (6,0) {{\Large$e_1\sfC^\prime$}}; 
\node (CC2) [green] at (0,10) {{\Large$e_2\sfC^\prime$}}; 
\node (CC3) [green] at (0,-6) {{\Large$e_3\sfC^\prime$}}; 
%\node (CC4) [green] at (6,9) {$e_4\sfC^\prime$}; 
\node (CC5) [green] at (-10,0) {{\Large$e_5\sfC^\prime$}}; 
%\node (CC0+) [green] at (6,-9) {$(e_0\sfC^\prime)^*$}; 
\node (CC1+) [green] at (-6,0) {{\Large$(e_1\sfC^\prime)^*$}}; 
\node (CC2+) [green] at (0,-10) {{\Large$(e_2\sfC^\prime)^*$}}; 
\node (CC3+) [green] at (0,6) {{\Large$(e_3\sfC^\prime)^*$}}; 
%\node (CC4+) [green] at (6,9) {$(e_4\sfC^\prime)^*$}; 
\node (CC5+) [green] at (10,0) {{\Large$(e_5\sfC^\prime)^*$}}; 

%\node (CCC0) [green]  at (-6,-9) {$e_0\sfC^{\prime\prime}$}; 
%\node (CCC1) [green]  at (-8,0) {$e_1\sfC^{\prime\prime}$}; 
\node (CCC2) [green]  at (2,10) {{\Large$e_2\sfC^{\prime\prime}$}};
\node (CCC3) [green]  at (2,-6) {{\Large$e_3\sfC^{\prime\prime}$}}; 
%\node (CCC4) [green]  at (-6,9) {$e_4\sfC^{\prime\prime}$}; 
%\node (CCC5) [green]  at (-6,9) {$e_5\sfC^{\prime\prime}$};
%\node (CCC0+) [green]  at (-6,-9) {$(e_0\sfC^{\prime\prime})^*$}; 
%\node (CCC1+) [green]  at (-8,0) {$(e_1\sfD)^*$}; 
\node (CCC2+) [green]  at (-2,-10) {{\Large$(e_2\sfC^{\prime\prime})^*$}};
\node (CCC3+) [green]  at (-2,6) {{\Large$(e_3\sfC^{\prime\prime})^*$}}; 
%\node (CCC4+) [green]  at (-6,9) {$(e_4\sfC^{\prime\prime})^*$}; 
%\node (CCC5+) [green]  at (-6,9) {$(e_5\sfC^{\prime\prime})^*$};

%\node (E0) [orange] at (0,-9) {$e_0\sfD$}; 
%\node (E1) [orange] at (-2,0) {$e_1\sfD$}; 
\node (E2) [orange] at (-1,-8) {{\Large$e_2\sfD$}}; 
\node (E3) [orange] at (-1,8) {{\Large$e_3\sfD$}}; 
%\node (E4) [orange] at (0,9) {$e_4\sfD$};
%\node (E5) [orange] at (-9,0) {$e_5\sfD$}; 
%\node (E0+) [orange] at (0,-9) {$(e_0\sfD)^*$}; 
%\node (E1+) [orange] at (-2,0) {$(e_1\sfD)^*$}; 
\node (E2+) [orange] at (1,8) {{\Large$(e_2\sfD)^*$}}; 
\node (E3+) [orange] at (1,-8) {{\Large$(e_3\sfD)^*$}}; 
%\node (E4+) [orange] at (0,9) {$(e_4\sfD)^*$};
%\node (E5+) [orange] at (-9,0) {$(e_5\sfD)^*$}; 

\draw[line width=0.032cm] (C3)--(CC3+); \draw[line width=0.032cm] (CC3+)--(CCC3+); \draw[line width=0.032cm] (CCC3+)--(C2+); 
\draw[line width=0.032cm] (C2+)--(CC2); \draw[line width=0.032cm] (CC2)--(CCC2); \draw[line width=0.032cm] (CCC2)--(C3); 
\draw[line width=0.032cm] (E2+)--(C3); \draw[line width=0.032cm] (E2+)--(CCC3+); \draw[line width=0.032cm] (E2+)--(CC2); 
\draw[line width=0.032cm] (E3)--(CCC2); \draw[line width=0.032cm] (E3)--(CC3+); \draw[line width=0.032cm] (E3)--(C2+);
\draw[line width=0.032cm] (C3+)--(CC3); \draw[line width=0.032cm] (CC3)--(CCC3); \draw[line width=0.032cm] (CCC3)--(C2); 
\draw[line width=0.032cm] (C2)--(CC2+); \draw[line width=0.032cm] (CC2+)--(CCC2+); \draw[line width=0.032cm] (CCC2+)--(C3+); 
\draw[line width=0.032cm] (E2)--(C3+); \draw[line width=0.032cm] (E2)--(CCC3); \draw[line width=0.032cm] (E2)--(CC2+); 
\draw[line width=0.032cm] (E3+)--(CCC2+); \draw[line width=0.032cm] (E3+)--(CC3); \draw[line width=0.032cm] (E3+)--(C2);
\draw[line width=0.032cm] (A5+)--(B3); \draw[line width=0.032cm] (A5+)--(B0+); \draw[line width=0.032cm] (A5+)--(BB2); 
\draw[line width=0.032cm] (A5+)--(BB4+); \draw[line width=0.032cm] (A5+)--(BBB1);
\draw[line width=0.032cm] (A3)--(BBB0); \draw[line width=0.032cm] (A3)--(BBB3+); \draw[line width=0.032cm] (A3)--(BB2); 
\draw[line width=0.032cm] (A3)--(BB4+); \draw[line width=0.032cm] (B3)--(C3); \draw[line width=0.032cm] (BB2)--(CC3+); \draw[line width=0.032cm] (BBB3+)--(CCC3+); 
\draw[line width=0.032cm] (A2)--(B3); \draw[line width=0.032cm] (A2)--(B0+); \draw[line width=0.032cm] (A0+)--(BBB0); 
\draw[line width=0.032cm] (A0+)--(BBB3+); \draw[line width=0.032cm] (A3)--(B1+);
\draw[line width=0.032cm] (A5)--(B3+); \draw[line width=0.032cm] (A5)--(B0); \draw[line width=0.032cm] (A5)--(BB2+); 
\draw[line width=0.032cm] (A5)--(BB4); \draw[line width=0.032cm] (A5)--(BBB1+);
\draw[line width=0.032cm] (A3+)--(BBB0+); \draw[line width=0.032cm] (A3+)--(BBB3); \draw[line width=0.032cm] (A3+)--(BB2+); 
\draw[line width=0.032cm] (A3+)--(BB4); \draw[line width=0.032cm] (B3+)--(C3+); \draw[line width=0.032cm] (BB2+)--(CC3); \draw[line width=0.032cm] (BBB3)--(CCC3); 
\draw[line width=0.032cm] (A2+)--(B3+); \draw[line width=0.032cm] (A2+)--(B0); \draw[line width=0.032cm] (A0)--(BBB0+); 
\draw[line width=0.032cm] (A0)--(BBB3); \draw[line width=0.032cm] (A3+)--(B1);
\draw[line width=0.032cm] (A4+)--(BBB2+); \draw[line width=0.032cm] (A4+)--(BBB4); \draw[line width=0.032cm] (A4+)--(BB3); 
\draw[line width=0.032cm] (A4+)--(BB0+); \draw[line width=0.032cm] (A1)--(B2); \draw[line width=0.032cm] (A1)--(B4+); 
\draw[line width=0.032cm] (A1)--(BB3); \draw[line width=0.032cm] (A1)--(BB0+);
\draw[line width=0.032cm] (BBB2+)--(CCC2); \draw[line width=0.032cm] (BB3)--(CC2); \draw[line width=0.032cm] (B2)--(C2+); 
\draw[line width=0.032cm] (A2)--(BBB2+); \draw[line width=0.032cm] (A2)--(BBB4); 
\draw[line width=0.032cm] (A0+)--(B2); \draw[line width=0.032cm] (A0+)--(B4+); \draw[line width=0.032cm] (A4+)--(B5+); 
\draw[line width=0.032cm] (A1)--(BBB5); \draw[line width=0.032cm] (A4)--(BBB2); \draw[line width=0.032cm] (A4)--(BBB4+); 
\draw[line width=0.032cm] (A4)--(BB3+); \draw[line width=0.032cm] (A4)--(BB0); 
\draw[line width=0.032cm] (A1+)--(B2+); \draw[line width=0.032cm] (A1+)--(B4); \draw[line width=0.032cm] (A1+)--(BB3+); 
\draw[line width=0.032cm] (A1+)--(BB0); \draw[line width=0.032cm] (BBB2)--(CCC2+); \draw[line width=0.032cm] (BB3+)--(CC2+); 
\draw[line width=0.032cm] (B2+)--(C2); \draw[line width=0.032cm] (A2+)--(BBB2); \draw[line width=0.032cm] (A2+)--(BBB4+); 
\draw[line width=0.032cm] (A0)--(B2+); \draw[line width=0.032cm] (A0)--(B4); \draw[line width=0.032cm] (A4)--(B5); \draw[line width=0.032cm] (A1+)--(BBB5+);
\draw[line width=0.032cm] (CC1)--(BB1+); \draw[line width=0.032cm] (CC1)--(BB5); \draw[line width=0.032cm] (CC5+)--(BB1+); 
\draw[line width=0.032cm] (CC5+)--(BB5); \draw[line width=0.032cm] (BB1+)--(A2); \draw[line width=0.032cm] (BB5)--(A0);
\draw[line width=0.032cm] (CC1+)--(BB1); \draw[line width=0.032cm] (CC1+)--(BB5+); \draw[line width=0.032cm] (CC5)--(BB1); 
\draw[line width=0.032cm] (CC5)--(BB5+); \draw[line width=0.032cm] (BB1)--(A2+); \draw[line width=0.032cm] (BB5+)--(A0+);

%\filldraw [draw=lightgray,fill=lightgray]  (-5.5,-2.5) rectangle (-4.5,-1.5) ;

%\node (B0) [blue] at (-5,-2) {{\Large$e_0\sfB$}}; 
\node at (-5,-2) {\myrect};
%\node (B1) [blue] at (2,-1) {{\Large$e_1\sfB$}}; 
\node at (2,-1) {\myrect};
%\node (B4) [blue] at (5,-14) {{\Large$e_4\sfB$}}; 
\node at (5,-14) {\myrect};
%\node (B5) [blue] at (-2,-15) {{\Large$e_5\sfB$}}; 
\node at (-2,-15) {\myrect};
%\node (B0+) [blue] at (5,2) {{\Large$(e_0\sfB)^*$}}; 
\node at (5,2) {\myrect};
%\node (B1+) [blue] at (-2,1) {{\Large$(e_1\sfB)^*$}}; 
\node at (-2,1) {\myrect};
%\node (B4+) [blue] at (-5,14) {{\Large$(e_4\sfB)^*$}}; 
\node at (-5,14) {\myrect};
%\node (B5+) [blue] at (2,15) {{\Large$(e_5\sfB)^*$}}; 
\node at (2,15) {\myrect};
%\node (BB0) [blue] at (0,-14) {{\Large$e_0\sfB^\prime$}}; 
\node at (0,-14) {\myrect};
%\node (BB4) [blue] at (0,-2) {{\Large$e_4\sfB^\prime$}}; 
\node at (0,-2) {\myrect};
%\node (BB0+) [blue] at (0,14) {{\Large$(e_0\sfB^\prime)^*$}}; 
\node at (0,14) {\myrect};
%\node (BB4+) [blue] at (0,2) {{\Large$(e_4\sfB^\prime)^*$}}; 
\node at (0,2) {\myrect};
%\node (BBB0) [blue] at (-5,2) {{\Large$e_0\sfB^{\prime\prime}$}}; 
\node at (-5,2) {\myrect};
%\node (BBB1) [blue] at (2,1) {{\Large$e_1\sfB^{\prime\prime}$}}; 
\node at (2,1) {\myrect};
%\node (BBB4) [blue] at (5,14) {{\Large$e_4\sfB^{\prime\prime}$}}; 
\node at (5,14) {\myrect};
%\node (BBB5) [blue] at (-2,15) {{\Large$e_5\sfB^{\prime\prime}$}}; 
\node at (-2,15) {\myrect};
%\node (BBB0+) [blue] at (5,-2) {{\Large$(e_0\sfB^{\prime\prime})^*$}}; 
\node at (5,-2) {\myrect};
%\node (BBB4+) [blue] at (-5,-14) {{\Large$(e_4\sfB^{\prime\prime})^*$}}; 
\node at (-5,-14) {\myrect};
%\node (BBB5+) [blue] at (2,-15) {{\Large$(e_5\sfB^{\prime\prime})^*$}}; 
\node at (2,-15) {\myrect};
%\node (BBB1+) [blue] at (-2,-1) {{\Large$(e_1\sfB^{\prime\prime})^*$}}; 
\node at (-2,-1) {\myrect};

%\node (CC1) [green] at (6,0) {{\Large$e_1\sfC^\prime$}}; 
\node at (6,0) {\myrect};
%\node (CC5) [green] at (-10,0) {{\Large$e_5\sfC^\prime$}}; 
\node at (-10,0) {\myrect};
%\node (CC1+) [green] at (-6,0) {{\Large$(e_1\sfC^\prime)^*$}}; 
\node at (-6,0) {\myrect};
%\node (CC5+) [green] at (10,0) {{\Large$(e_5\sfC^\prime)^*$}}; 
\node at (10,0) {\myrect};

\end{tikzpicture} }}
\end{center}

%%%%%%%%%%%%%%%%%%%%%
%%% type 6a exchange graph 2%%%%%%%%%%%%%
%%%%%%%%%%%%%%%%%%%%%
\begin{center}
\newcommand{\myrect} [1] [] {
\begin{tikzpicture}
\draw [dashed]  (0,0) rectangle (1.35,1.35) ;
\end{tikzpicture}}
{\scalebox{0.78}{
\begin{tikzpicture}

\node (B0) [blue] at (-3,-4) {\large{$e_0\sfB$}}; 
\node (B1) [blue] at (2,-10) {\large{$e_1\sfB$}}; 
\node (B4) [blue] at (5,4) {\large{$e_4\sfB$}}; 
\node (B5) [blue] at (-2,6) {\large{$e_5\sfB$}}; 
\node (B0+) [blue] at (3,4) {\large{$(e_0\sfB)^*$}}; 
\node (B1+) [blue] at (-2,10) {\large{$(e_1\sfB)^*$}}; 
\node (B4+) [blue] at (-5,-4) {\large{$(e_4\sfB)^*$}}; 
\node (B5+) [blue] at (2,-6) {\large{$(e_5\sfB)^*$}}; 

\node (BB0) [blue] at (-1,4) {\large{$e_0\sfB^\prime$}}; 
\node (BB4) [blue] at (-1,-4) {\large{$e_4\sfB^\prime$}}; 
\node (BB0+) [blue] at (1,-4) {\large{$(e_0\sfB^\prime)^*$}}; 
\node (BB4+) [blue] at (1,4) {\large{$(e_4\sfB^\prime)^*$}}; 

\node (BBB0) [blue] at (-3,4) {\large{$e_0\sfB^{\prime\prime}$}}; 
\node (BBB1) [blue] at (2,10) {\large{$e_1\sfB^{\prime\prime}$}}; 
\node (BBB4) [blue] at (5,-4) {\large{$e_4\sfB^{\prime\prime}$}}; 
\node (BBB5) [blue] at (-2,-6) {\large{$e_5\sfB^{\prime\prime}$}}; 
\node (BBB0+) [blue] at (3,-4) {\large{$(e_0\sfB^{\prime\prime})^*$}}; 
\node (BBB1+) [blue] at (-2,-10) {\large{$(e_1\sfB^{\prime\prime})^*$}}; 
\node (BBB4+) [blue] at (-5,4) {\large{$(e_4\sfB^{\prime\prime})^*$}}; 
\node (BBB5+) [blue] at (2,6) {\large{$(e_5\sfB^{\prime\prime})^*$}}; 

\node (CC1) [green] at (5,8) {\large{$e_1\sfC^\prime$}}; 
\node (CC5) [green] at (-5,8) {\large{$e_5\sfC^\prime$}}; 
\node (CC1+) [green] at (-5,-8) {\large{$(e_1\sfC^\prime)^*$}}; 
\node (CC5+) [green] at (5,-8) {\large{$(e_5\sfC^\prime)^*$}}; 
%%%%%%%%%%%%%%%%%%%%%%%%%%%%%%%%%%%%%%%%

\node (C0) [green] at (-4,-6) {\large{$e_0\sfC$}}; 
\node (C1) [green] at (1,-8) {\large{$e_1\sfC$}}; 
\node (C4) [green] at (4,2) {\large{$e_4\sfC$}}; 
\node (C5) [green] at (-3,8) {\large{$e_5\sfC$}}; 
\node (C0+) [green] at (4,6) {\large{$(e_0\sfC)^*$}}; 
\node (C1+) [green] at (-1,8) {\large{$(e_1\sfC)^*$}}; 
\node (C4+) [green] at (-4,-2) {\large{$(e_4\sfC)^*$}}; 
\node (C5+) [green] at (3,-8) {\large{$(e_5\sfC)^*$}}; 

\node (CC0) [green] at (0,6) {\large{$e_0\sfC^\prime$}}; 
\node (CC4) [green] at (0,-2) {\large{$e_4\sfC^\prime$}}; 
\node (CC0+) [green] at (0,-6) {\large{$(e_0\sfC^\prime)^*$}}; 
\node (CC4+) [green] at (0,2) {\large{$(e_4\sfC^\prime)^*$}}; 

\node (CCC0) [green]  at (4,-2) {\large{$e_0\sfC^{\prime\prime}$}}; 
\node (CCC1) [green]  at (-3,-8) {\large{$e_1\sfC^{\prime\prime}$}};  
\node (CCC4) [green]  at (-4,6) {\large{$e_4\sfC^{\prime\prime}$}}; 
\node (CCC5) [green]  at (1,8) {\large{$e_5\sfC^{\prime\prime}$}};
\node (CCC0+) [green]  at (-4,2) {\large{$(e_0\sfC^{\prime\prime})^*$}}; 
\node (CCC1+) [green]  at (3,8) {\large{$(e_1\sfC^{\prime\prime})^*$}}; 
\node (CCC4+) [green]  at (4,-6) {\large{$(e_4\sfC^{\prime\prime})^*$}}; 
\node (CCC5+) [green]  at (-1,-8) {\large{$(e_5\sfC^{\prime\prime})^*$}};

\node (E0) [orange] at (1.5,0) {\large{$e_0\sfD$}}; 
\node (E1) [orange] at (0,10) {\large{$e_1\sfD$}}; 
\node (E4) [orange] at (4,-10) {\large{$e_4\sfD$}};
\node (E5) [orange] at (-4,-10) {\large{$e_5\sfD$}}; 
\node (E0+) [orange] at (-1.5,0) {\large{$(e_0\sfD)^*$}}; 
\node (E1+) [orange] at (0,-10) {\large{$(e_1\sfD)^*$}}; 
\node (E4+) [orange] at (-4,10) {\large{$(e_4\sfD)^*$}};
\node (E5+) [orange] at (4,10) {\large{$(e_5\sfD)^*$}}; 

\draw[line width=0.028cm] (E0)--(CCC0); \draw[line width=0.028cm] (E0)--(C4+); \draw[line width=0.028cm] (E0)--(CC4+); 
\draw[line width=0.028cm] (E0+)--(CCC0+); \draw[line width=0.028cm] (E0+)--(C4); \draw[line width=0.028cm] (E0+)--(CC4);
\draw[line width=0.028cm] (C4)--(CCC0); \draw[line width=0.028cm] (C4+)--(CCC0+); \draw[line width=0.028cm] (C4)--(CC4+); 
\draw[line width=0.028cm] (CC4+)--(CCC0+); \draw[line width=0.028cm] (C4+)--(CC4); \draw[line width=0.028cm] (CC4)--(CCC0);
\draw[line width=0.028cm] (CCC0+)--(BBB4+); \draw[line width=0.028cm] (CCC0+)--(BBB0); \draw[line width=0.028cm] (CC4+)--(BB0); 
\draw[line width=0.028cm] (CC4+)--(BB4+); \draw[line width=0.028cm] (C4)--(B0+); \draw[line width=0.028cm] (C4)--(B4); 
\draw[line width=0.028cm] (CCC4)--(BBB4+); \draw[line width=0.028cm] (CCC4)--(BBB0); \draw[line width=0.028cm] (CC0)--(BB0); 
\draw[line width=0.028cm] (CC0)--(BB4+); \draw[line width=0.028cm] (C0+)--(B0+); \draw[line width=0.028cm] (C0+)--(B4);
\draw[line width=0.028cm] (CCC0)--(BBB4); \draw[line width=0.028cm] (CCC0)--(BBB0+); \draw[line width=0.028cm] (CC4)--(BB0+); 
\draw[line width=0.028cm] (CC4)--(BB4); \draw[line width=0.028cm] (C4+)--(B0); \draw[line width=0.028cm] (C4+)--(B4+); 
\draw[line width=0.028cm] (CCC4+)--(BBB4); \draw[line width=0.028cm] (CCC4+)--(BBB0+); \draw[line width=0.028cm] (CC0+)--(BB0+); 
\draw[line width=0.028cm] (CC0+)--(BB4); \draw[line width=0.028cm] (C0)--(B0); \draw[line width=0.028cm] (C0)--(B4+);
\draw[line width=0.028cm] (CC1)--(C0+); \draw[line width=0.028cm] (CC1)--(E5+); \draw[line width=0.028cm] (CCC1+)--(C0+); 
\draw[line width=0.028cm] (CCC1+)--(E5+); \draw[line width=0.028cm] (CCC1+)--(BBB1); \draw[line width=0.028cm] (CCC1+)--(BBB5+);
\draw[line width=0.028cm] (CCC5)--(BBB1); \draw[line width=0.028cm] (CCC5)--(BBB5+); \draw[line width=0.028cm] (CCC5)--(E1); 
\draw[line width=0.028cm] (CCC5)--(CC0); \draw[line width=0.028cm] (C1+)--(E1); \draw[line width=0.028cm] (C1+)--(CC0); 
\draw[line width=0.028cm] (C1+)--(B1+); \draw[line width=0.028cm] (C1+)--(B5); 
\draw[line width=0.028cm] (C5)--(B1+); \draw[line width=0.028cm] (C5)--(B5); \draw[line width=0.028cm] (C5)--(E4+); 
\draw[line width=0.028cm] (C5)--(CCC4); \draw[line width=0.028cm] (CC5)--(E4+); \draw[line width=0.028cm] (CC5)--(CCC4);
\draw[line width=0.028cm] (CC1+)--(C0); \draw[line width=0.028cm] (CC1+)--(E5); \draw[line width=0.028cm] (CCC1)--(C0); 
\draw[line width=0.028cm] (CCC1)--(E5); \draw[line width=0.028cm] (CCC1)--(BBB1+); \draw[line width=0.028cm] (CCC1)--(BBB5);
\draw[line width=0.028cm] (CCC5+)--(BBB1+); \draw[line width=0.028cm] (CCC5+)--(BBB5); \draw[line width=0.028cm] (CCC5+)--(E1+); 
\draw[line width=0.028cm] (CCC5+)--(CC0+); \draw[line width=0.028cm] (C1)--(E1+); \draw[line width=0.028cm] (C1)--(CC0+); 
\draw[line width=0.028cm] (C1)--(B1); \draw[line width=0.028cm] (C1)--(B5+); 
\draw[line width=0.028cm] (C5+)--(B1); \draw[line width=0.028cm] (C5+)--(B5+); \draw[line width=0.028cm] (C5+)--(E4); 
\draw[line width=0.028cm] (C5+)--(CCC4+); \draw[line width=0.028cm] (CC5+)--(E4); \draw[line width=0.028cm] (CC5+)--(CCC4+);

%\node (B0) [blue] at (-3,-4) {$e_0\sfB$}; 
\node at (-3,-4) {\myrect};
%\node (B1) [blue] at (2,-10) {$e_1\sfB$}; 
\node at (2,-10) {\myrect};
%\node (B4) [blue] at (5,4) {$e_4\sfB$}; 
\node at (5,4) {\myrect};
%\node (B5) [blue] at (-2,6) {$e_5\sfB$};
\node at (-2,6) {\myrect};
%\node (B0+) [blue] at (3,4) {$(e_0\sfB)^*$}; 
\node at (3,4) {\myrect};
%\node (B1+) [blue] at (-2,10) {$(e_1\sfB)^*$}; 
\node at (-2,10) {\myrect};
%\node (B4+) [blue] at (-5,-4) {$(e_4\sfB)^*$}; 
\node at (-5,-4) {\myrect};
%\node (B5+) [blue] at (2,-6) {$(e_5\sfB)^*$}; 
\node at (2,-6) {\myrect};
%\node (BB0) [blue] at (-1,4) {$e_0\sfB^\prime$}; 
\node at (-1,4) {\myrect};
%\node (BB4) [blue] at (-1,-4) {$e_4\sfB^\prime$}; 
\node at (-1,-4) {\myrect};
%\node (BB0+) [blue] at (1,-4) {$(e_0\sfB^\prime)^*$}; 
\node at (1,-4) {\myrect};
%\node (BB4+) [blue] at (1,4) {$(e_4\sfB^\prime)^*$}; 
\node at (1,4) {\myrect};
%\node (BBB0) [blue] at (-3,4) {$e_0\sfB^{\prime\prime}$}; 
\node at (-3,4) {\myrect};
%\node (BBB1) [blue] at (2,10) {$e_1\sfB^{\prime\prime}$}; 
\node at (2,10) {\myrect};
%\node (BBB4) [blue] at (5,-4) {$e_4\sfB^{\prime\prime}$}; 
\node at (5,-4) {\myrect};
%\node (BBB5) [blue] at (-2,-6) {$e_5\sfB^{\prime\prime}$}; 
\node at (-2,-6) {\myrect};
%\node (BBB0+) [blue] at (3,-4) {$(e_0\sfB^{\prime\prime})^*$}; 
\node at (3,-4) {\myrect};
%\node (BBB1+) [blue] at (-2,-10) {$(e_1\sfB^{\prime\prime})^*$};
\node at (-2,-10) {\myrect};
%\node (BBB4+) [blue] at (-5,4) {$(e_4\sfB^{\prime\prime})^*$}; 
\node at (-5,4) {\myrect};
%\node (BBB5+) [blue] at (2,6) {$(e_5\sfB^{\prime\prime})^*$}; 
\node at (2,6) {\myrect};
%\node (CC1) [green] at (5,8) {$e_1\sfC^\prime$}; 
\node at (5,8) {\myrect};
%\node (CC5) [green] at (-5,8) {$e_5\sfC^\prime$}; 
\node at (-5,8) {\myrect};
%\node (CC1+) [green] at (-5,-8) {$(e_1\sfC^\prime)^*$}; 
\node at (-5,-8) {\myrect};
%\node (CC5+) [green] at (5,-8) {$(e_5\sfC^\prime)^*$}; 
\node at (5,-8) {\myrect};

\end{tikzpicture} }}
\end{center}
 
%%%%%%%%%%%%%%%%%%%%%%%%%%%%%%%%%%%%%%%%%%%%%%%%%%%%%%%
\subsection{Type 6b}
\label{type6b}

In this subsection, we consider the reflexive polygon of type 6b. 
Thus, let $R$ be the three dimensional complete local Gorenstein toric singularity defined by the cone $\sigma$: 
\[
\sigma=\mathrm{Cone}\{v_1=(1,0,1), v_2=(0,1,1), v_3=(-1,1,1), v_4=(-1,-1,1), v_5=(0,-1,1) \}. 
\]
As an element in $\Cl(R)$, we obtain $[I_1]+2[I_4]+2[I_5]=0$, $[I_2]-4[I_4]-3[I_5]=0$, $[I_3]+3[I_4]+2[I_5]=0$.  
Therefore, we have that $\Cl(R)\cong\ZZ^2$, and each divisorial ideal is represented by $T(0,0,0,d,e)$ where $d, e\in\ZZ$. 
Also, there are three consistent dimer models (up to right equivalence) written below which give the reflexive polygon of type 6b as the perfect matching polygon. 

%%%%%%%%%%%%%%%%%%%%%%%%%%%%%%%%%%%
%%%%%%%%%%%%Type_6b-1%%%%%%%%%%%%%%%%%%
%%%%%%%%%%%%%%%%%%%%%%%%%%%%%%%%%%%

\subsubsection{Type 6b-1}

First, we consider the following consistent dimer model and the associated QP $(Q, W_Q)$. 
For simplicity, we denote by $\sfA$ the complete Jacobian algebra $\calP(Q, W_Q)$. 
Then, we have that $\sfA\cong\End_R(\bigoplus_{j\in Q_0}T_{ij})$ and $e_i\sfA\cong\bigoplus_{j\in Q_0}T_{ij}$ for any $i\in Q_0$ by Theorem~\ref{NCCR1}. 
We will give these splitting MM generators $\bigoplus_{j\in Q_0}T_{ij}$. 

\medskip

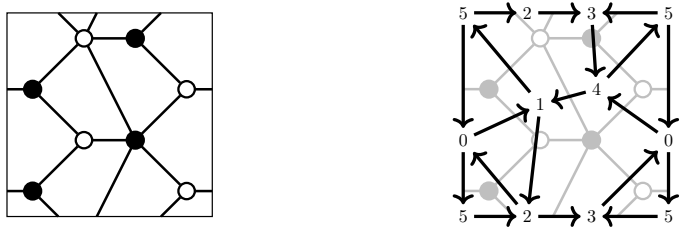
\begin{center}

\end{center}

\medskip

We fix the perfect matching $\sfP_0$, 
then the following $\sfP_1, \cdots, \sfP_5$ are extremal perfect matchings corresponding to $v_1, \cdots, v_5$. 
Namely, $v_i=(h(\sfP_i, \sfP_0), 1)\in\ZZ^3$. 

\medskip

%%%%%perfect_matching%%%%%

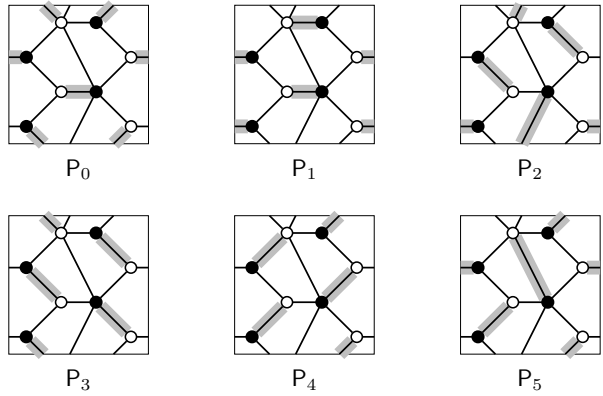
\begin{center}
{\scalebox{0.9}{
%
}}
\end{center}

\medskip

Then, we have the following table. For simplicity, we denote by $T(d,e)$ a divisorial ideal $T(0,0,0,d,e)$. 
We remark that $(e_0\sfA)^*\cong e_5\sfA$, $(e_1\sfA)^*\cong e_4\sfA$ and $(e_2\sfA)^*\cong e_3\sfA$. 

\medskip

\begin{center}
{\footnotesize{
\begin{tabular}{|l|l|l|l|l|l|} \hline 
$e_0\sfA$&$e_1\sfA$&$e_2\sfA$&$e_3\sfA$&$e_4\sfA$ &$e_5\sfA$\\ \hline
$T_{00}\cong R$&$T_{10}\cong T(-1,-1)$&$T_{20}\cong T(1,1)$&$T_{30}\cong T(-3,-2)$&$T_{40}\cong T(-1,0)$&$T_{50}\cong T(-2,-1)$ \\
$T_{01}\cong T(1,1)$&$T_{11}\cong R$&$T_{21}\cong T(2,2)$&$T_{31}\cong T(-2,-1)$&$T_{41}\cong T(0,1)$&$T_{51}\cong T(-1,0)$ \\
$T_{02}\cong T(-1,-1)$&$T_{12}\cong T(-2,-2)$&$T_{22}\cong R$&$T_{32}\cong T(-4,-3)$&$T_{42}\cong T(-2,-1)$&$T_{52}\cong T(-3,-2)$ \\
$T_{03}\cong T(3,2)$&$T_{13}\cong T(2,1)$&$T_{23}\cong T(4,3)$&$T_{33}\cong R$&$T_{43}\cong T(2,2)$&$T_{53}\cong T(1,1)$ \\
$T_{04}\cong T(1,0)$&$T_{14}\cong T(0,-1)$&$T_{24}\cong T(2,1)$&$T_{34}\cong T(-2,-2)$&$T_{44}\cong R$&$T_{54}\cong T(-1,-1)$ \\
$T_{05}\cong T(2,1)$&$T_{15}\cong T(1,0)$&$T_{25}\cong T(3,2)$&$T_{35}\cong T(-1,-1)$&$T_{45}\cong T(1,1)$&$T_{55}\cong R$ \\ \hline 
\end{tabular} 
}}
\end{center}

\medskip

%%%%%%%%%%%%%%%%%%%%%%%%%%%%%%%%%%%
%%%%%%%%%%%%Type_6b-2%%%%%%%%%%%%%%%%%%
%%%%%%%%%%%%%%%%%%%%%%%%%%%%%%%%%%%

\subsubsection{Type 6b-2}

Next, we consider the following consistent dimer model and the associated QP $(Q, W_Q)$. 
We denote by $\sfB$ the complete Jacobian algebra $\calP(Q, W_Q)$, and consider splitting MM generators $e_i\sfB$. 

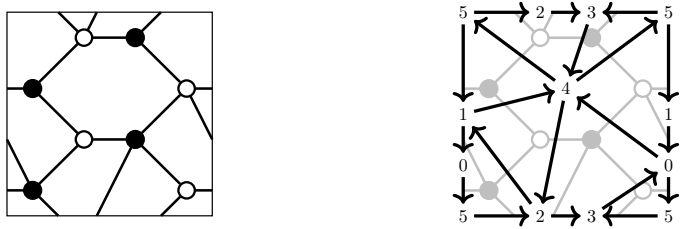
\medskip

\begin{center}
%
\end{center}

\medskip

We fix the perfect matching $\sfP_0$, 
then the following $\sfP_1, \cdots, \sfP_5$ are extremal perfect matchings corresponding to $v_1, \cdots, v_5$. 

\medskip

%%%%%perfect_matching%%%%%

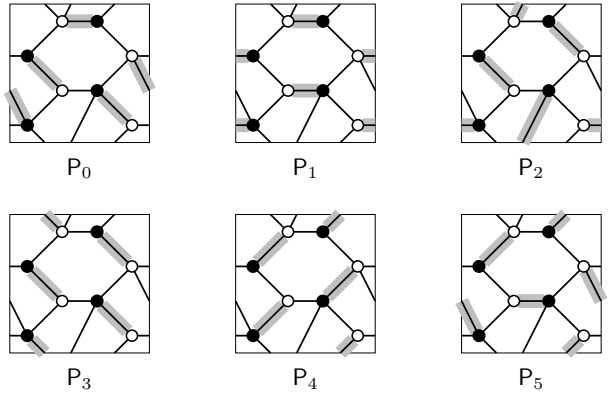
\begin{center}
{\scalebox{0.9}{
%
}}
\end{center}

\medskip

Then, we have the following tables. For simplicity, we denote by $T(d,e)$ a divisorial ideal $T(0,0,0,d,e)$. 

\medskip

\begin{center}
{\footnotesize{
\begin{tabular}{|l|l|l|l|l|l|} \hline 
$e_0\sfB$&$e_1\sfB$&$e_2\sfB$&$e_3\sfB$&$e_4\sfB$ &$e_5\sfB$\\ \hline
$T_{00}\cong R$&$T_{10}\cong T(0,1)$&$T_{20}\cong T(1,1)$&$T_{30}\cong T(-3,-2)$&$T_{40}\cong T(-1,0)$&$T_{50}\cong T(-2,-1)$ \\
$T_{01}\cong T(0,-1)$&$T_{11}\cong R$&$T_{21}\cong T(1,0)$&$T_{31}\cong T(-3,-3)$&$T_{41}\cong T(-1,-1)$&$T_{51}\cong T(-2,-2)$ \\
$T_{02}\cong T(-1,-1)$&$T_{12}\cong T(-1,0)$&$T_{22}\cong R$&$T_{32}\cong T(-4,-3)$&$T_{42}\cong T(-2,-1)$&$T_{52}\cong T(-3,-2)$ \\
$T_{03}\cong T(3,2)$&$T_{13}\cong T(3,3)$&$T_{23}\cong T(4,3)$&$T_{33}\cong R$&$T_{43}\cong T(2,2)$&$T_{53}\cong T(1,1)$ \\
$T_{04}\cong T(1,0)$&$T_{14}\cong T(1,1)$&$T_{24}\cong T(2,1)$&$T_{34}\cong T(-2,-2)$&$T_{44}\cong R$&$T_{54}\cong T(-1,-1)$ \\
$T_{05}\cong T(2,1)$&$T_{15}\cong T(2,2)$&$T_{25}\cong T(3,2)$&$T_{35}\cong T(-1,-1)$&$T_{45}\cong T(1,1)$&$T_{55}\cong R$ \\ \hline 
\end{tabular} 
}}
\end{center}

\medskip

%%%%%%%%%%%%%%%%%%%%%%%%%%%%%%%%%%%
%%%%%%%%%%%%Type_6b-3%%%%%%%%%%%%%%%%%%
%%%%%%%%%%%%%%%%%%%%%%%%%%%%%%%%%%%

\subsubsection{Type 6b-3}

Next, we consider the following consistent dimer model and the associated QP $(Q, W_Q)$. 
We denote by $\sfC$ the complete Jacobian algebra $\calP(Q, W_Q)$, and consider splitting MM generators $e_i\sfC$. 

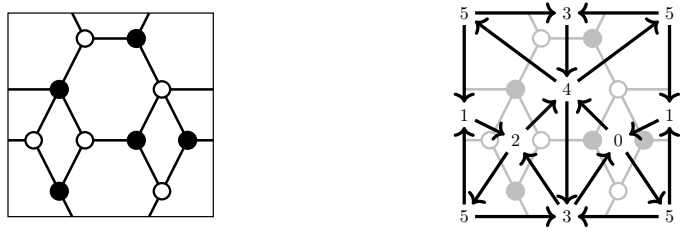
\medskip

\begin{center}
%
\end{center}

\medskip

We fix the perfect matching $\sfP_0$, 
then the following $\sfP_1, \cdots, \sfP_5$ are extremal perfect matchings corresponding to $v_1, \cdots, v_5$. 
Namely, $v_i=(h(\sfP_i, \sfP_0), 1)\in\ZZ^3$. 

\medskip

%%%%%perfect_matching%%%%%
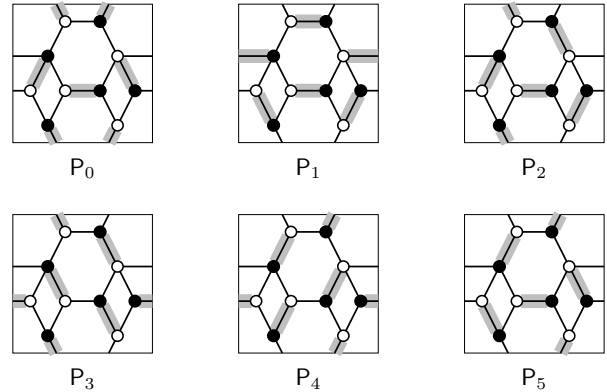
\begin{center}
{\scalebox{0.9}{
%
}}
\end{center}

\medskip

Then, we have the following tables. For simplicity, we denote by $T(d,e)$ a divisorial ideal $T(0,0,0,d,e)$. 

\medskip

\begin{center}
{\footnotesize{
\begin{tabular}{|l|l|l|l|l|l|} \hline 
$e_0\sfC$&$e_1\sfC$&$e_2\sfC$&$e_3\sfC$&$e_4\sfC$ &$e_5\sfC$\\ \hline
$T_{00}\cong R$&$T_{10}\cong T(0,1)$&$T_{20}\cong T(-4,-2)$&$T_{30}\cong T(-3,-2)$&$T_{40}\cong T(-1,0)$&$T_{50}\cong T(-2,-1)$ \\
$T_{01}\cong T(0,-1)$&$T_{11}\cong R$&$T_{21}\cong T(-4,-3)$&$T_{31}\cong T(-3,-3)$&$T_{41}\cong T(-1,-1)$&$T_{51}\cong T(-2,-2)$ \\
$T_{02}\cong T(4,2)$&$T_{12}\cong T(4,3)$&$T_{22}\cong R$&$T_{32}\cong T(1,0)$&$T_{42}\cong T(3,2)$&$T_{52}\cong T(2,1)$ \\
$T_{03}\cong T(3,2)$&$T_{13}\cong T(3,3)$&$T_{23}\cong T(-1,0)$&$T_{33}\cong R$&$T_{43}\cong T(2,2)$&$T_{53}\cong T(1,1)$ \\
$T_{04}\cong T(1,0)$&$T_{14}\cong T(1,1)$&$T_{24}\cong T(-3,-2)$&$T_{34}\cong T(-2,-2)$&$T_{44}\cong R$&$T_{54}\cong T(-1,-1)$ \\
$T_{05}\cong T(2,1)$&$T_{15}\cong T(2,2)$&$T_{25}\cong T(-2,-1)$&$T_{35}\cong T(-1,-1)$&$T_{45}\cong T(1,1)$&$T_{55}\cong R$ \\ \hline 
\end{tabular} 
}}
\end{center}

\medskip

%%%%%%%%%%%%%%%%%%%%%%%%%%%%%%%%%%%%%%%%%%%%

\subsubsection{Exchange graph of type 6b}
By combining the above results, we describe the exchange graph $\EG(\TMMG(R))$ for the case of type 6b as follows, 
and it is actually connected. 

%%%%%%%%%%%%%%%%%%%%%
%%% type 6b exchange graph %%%%%%%%%%%%%
%%%%%%%%%%%%%%%%%%%%%
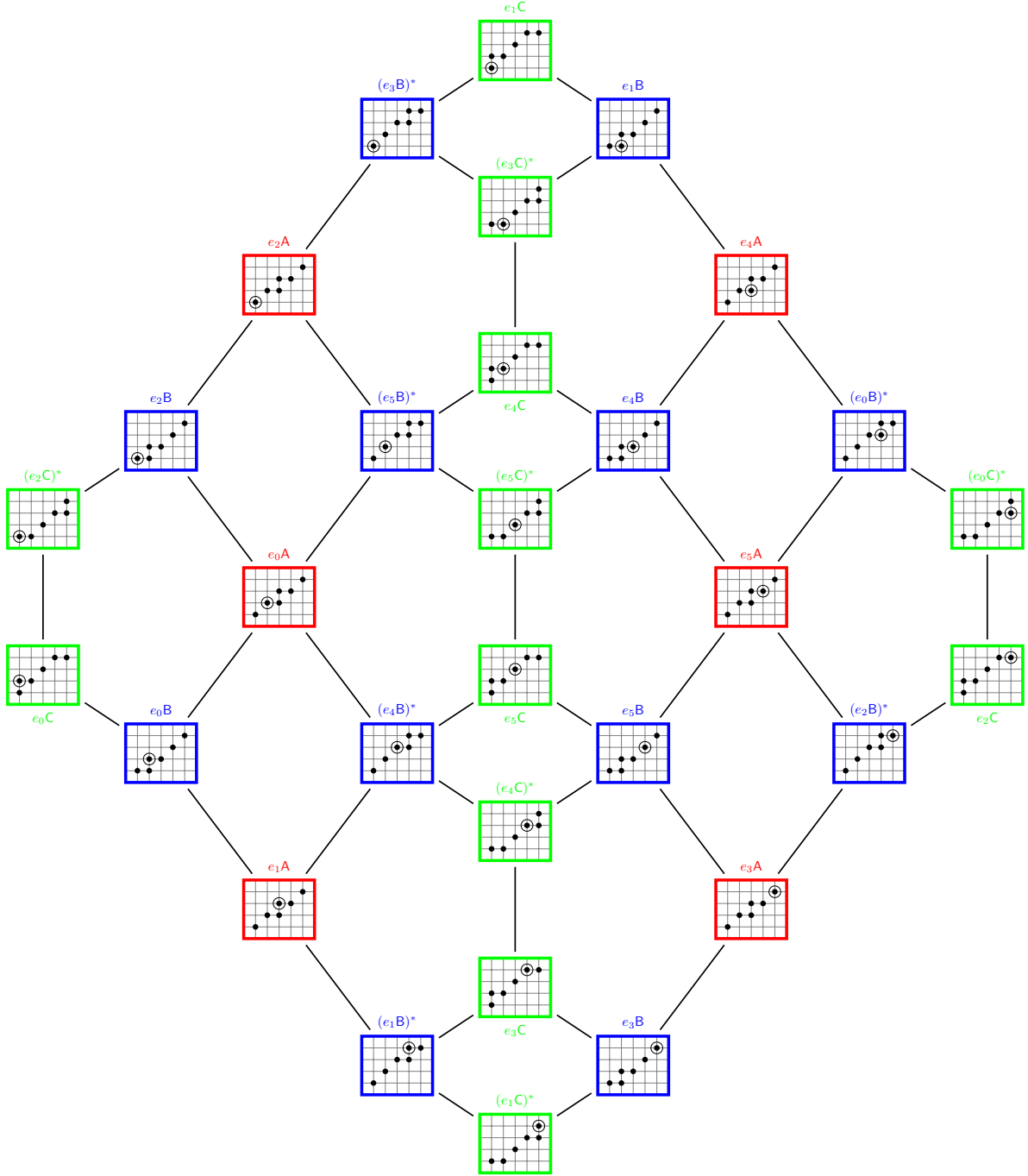
\begin{center}
{\scalebox{0.6}{
%
}};

\draw[line width=0.0375cm] (A5)--(B5); \draw[line width=0.0375cm] (A5)--(B4); \draw[line width=0.0375cm] (A5)--(B0+); 
\draw[line width=0.0375cm] (A5)--(B2+); \draw[line width=0.0375cm] (A1)--(B0); \draw[line width=0.0375cm] (A1)--(B4+); 
\draw[line width=0.0375cm] (A2)--(B2); \draw[line width=0.0375cm] (A2)--(B5+); 
\draw[line width=0.0375cm] (A3)--(B2+); \draw[line width=0.0375cm] (A3)--(B5); \draw[line width=0.0375cm] (A4)--(B0+); 
\draw[line width=0.0375cm] (A4)--(B4);\draw[line width=0.0375cm] (C5)--(C5+); \draw[line width=0.0375cm] (C5+)--(B4); 
\draw[line width=0.0375cm] (C5+)--(B5+); \draw[line width=0.0375cm] (B5)--(C5);  \draw[line width=0.0375cm] (B4+)--(C5); 
\draw[line width=0.0375cm] (B0+)--(C0+); \draw[line width=0.0375cm] (B2+)--(C2); \draw[line width=0.0375cm] (C0+)--(C2);   
\draw[line width=0.0375cm] (B0)--(C0); \draw[line width=0.0375cm] (B2)--(C2+); \draw[line width=0.0375cm] (C0)--(C2+); 
\draw[line width=0.0375cm] (A0)--(B5+); \draw[line width=0.0375cm] (A0)--(B4+); \draw[line width=0.0375cm] (A0)--(B0); 
\draw[line width=0.0375cm] (A0)--(B2); \draw[line width=0.0375cm] (C4)--(B5+); \draw[line width=0.0375cm] (C4)--(B4); 
\draw[line width=0.0375cm] (C4)--(C3+); \draw[line width=0.0375cm] (A2)--(B3+); \draw[line width=0.0375cm] (A4)--(B1); 
\draw[line width=0.0375cm] (B3+)--(C3+); \draw[line width=0.0375cm] (B1)--(C3+);
\draw[line width=0.0375cm] (C4+)--(B5); \draw[line width=0.0375cm] (C4+)--(B4+); \draw[line width=0.0375cm] (C4+)--(C3); 
\draw[line width=0.0375cm] (A3)--(B3); \draw[line width=0.0375cm] (A1)--(B1+); \draw[line width=0.0375cm] (B3)--(C3); 
\draw[line width=0.0375cm] (B1+)--(C3); \draw[line width=0.0375cm] (B3+)--(C1); \draw[line width=0.0375cm] (B1)--(C1); 
\draw[line width=0.0375cm] (B3)--(C1+); \draw[line width=0.0375cm] (B1+)--(C1+); 
\end{tikzpicture} }}
\end{center}

%%%%%%%%%%%%%%%%%%%%%%%%%%%%%%%%%%%%%%%%%%%%%%%%%%%%%%%
\subsection{Type 6c}

In this subsection, we consider the reflexive polygon of type 6c. 
Thus, let $R$ be the three dimensional complete local Gorenstein toric singularity defined by the cone $\sigma$: 
\[
\sigma=\mathrm{Cone}\{v_1=(1,0,1), v_2=(0,1,1), v_3=(-2,-1,1), v_4=(0,-1,1) \}. 
\]
As an element in $\Cl(R)$, we obtain $[I_1]=2[I_3]$, $[I_2]-[I_3]-[I_4]=0$, $4[I_3]+2[I_4]=0$.  
Therefore, we have that $\Cl(R)\cong\ZZ\times \ZZ/2\ZZ$, and each divisorial ideal is represented by $T(0,0,c,d)$ where $c\in\ZZ, d\in\ZZ/2\ZZ$. 
Also, there are two consistent dimer models (up to right equivalence) written below which give the reflexive polygon of type 6c as the perfect matching polygon. 

\medskip

%%%%%%%%%%%%%%%%%%%%%%%%%%%%%%%%%%%
%%%%%%%%%%%%Type_6c-1%%%%%%%%%%%%%%%%%%
%%%%%%%%%%%%%%%%%%%%%%%%%%%%%%%%%%%

\subsubsection{Type 6c-1}

First, we consider the following consistent dimer model and the associated QP $(Q, W_Q)$. 
For simplicity, we denote by $\sfA$ the complete Jacobian algebra $\calP(Q, W_Q)$, and we will give splitting MM generators $e_i\sfA\cong\bigoplus_{j\in Q_0}T_{ij}$. 

\medskip

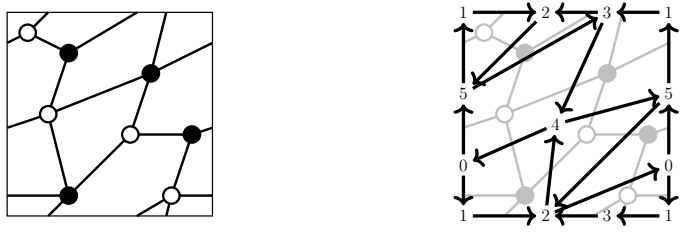
\begin{center}

\end{center}

\medskip

We fix the perfect matching $\sfP_0$, 
then the following $\sfP_1, \cdots, \sfP_4$ are extremal perfect matchings corresponding to $v_1, \cdots, v_4$. 
Namely, $v_i=(h(\sfP_i, \sfP_0), 1)\in\ZZ^3$. 

\medskip

%%%%%perfect_matching%%%%%

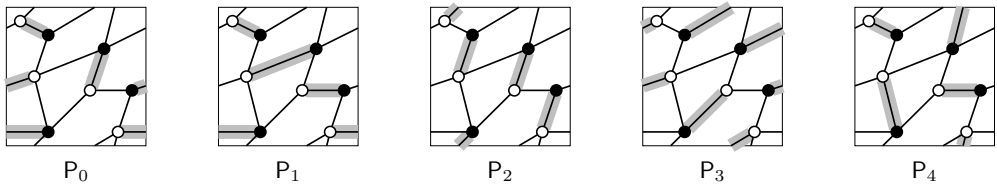
\begin{center}
{\scalebox{0.9}{
%
}}
\end{center}

\medskip

Then, we have the following table. For simplicity, we denote by $T(c,d)$ a divisorial ideal $T(0,0,c,d)$. 
We remark that $e_0\sfA\cong e_3\sfA$, $e_1\sfA\cong e_4\sfA$, $e_2\sfA\cong e_5\sfA$, $(e_0\sfA)^*\cong e_1\sfA$ and $(e_2\sfA)^*\cong e_2\sfA$. 

\medskip

\begin{center}
{\small{
\begin{tabular}{|l|l|l|l|l|l|} \hline 
$e_0\sfA$&$e_1\sfA$&$e_2\sfA$&$e_3\sfA$&$e_4\sfA$ &$e_5\sfA$\\ \hline
$T_{00}\cong R$&$T_{10}\cong T(-2,0)$&$T_{20}\cong T(1,1)$&$T_{30}\cong T(2,1)$&$T_{40}\cong T(0,1)$&$T_{50}\cong T(-1,0)$ \\
$T_{01}\cong T(2,0)$&$T_{11}\cong R$&$T_{21}\cong T(3,1)$&$T_{31}\cong T(4,1)$&$T_{41}\cong T(2,1)$&$T_{51}\cong T(1,0)$ \\
$T_{02}\cong T(3,1)$&$T_{12}\cong T(1,1)$&$T_{22}\cong R$&$T_{32}\cong T(1,0)$&$T_{42}\cong T(-1,0)$&$T_{52}\cong T(2,1)$ \\
$T_{03}\cong T(2,1)$&$T_{13}\cong T(0,1)$&$T_{23}\cong T(-1,0)$&$T_{33}\cong R$&$T_{43}\cong T(-2,0)$&$T_{53}\cong T(1,1)$ \\
$T_{04}\cong T(4,1)$&$T_{14}\cong T(2,1)$&$T_{24}\cong T(1,0)$&$T_{34}\cong T(2,0)$&$T_{44}\cong R$&$T_{54}\cong T(3,1)$ \\
$T_{05}\cong T(1,0)$&$T_{15}\cong T(-1,0)$&$T_{25}\cong T(2,1)$&$T_{35}\cong T(3,1)$&$T_{45}\cong T(1,1)$&$T_{55}\cong R$ \\ \hline 
\end{tabular} 
}}
\end{center}

\medskip

%%%%%%%%%%%%%%%%%%%%%%%%%%%%%%%%%%%
%%%%%%%%%%%%Type_6c-2%%%%%%%%%%%%%%%%%%
%%%%%%%%%%%%%%%%%%%%%%%%%%%%%%%%%%%

\subsubsection{Type 6c-2}

Next, we consider the following consistent dimer model and the associated QP $(Q, W_Q)$. 
We denote by $\sfB$ the complete Jacobian algebra $\calP(Q, W_Q)$, and consider splitting MM generators $e_i\sfB$. 

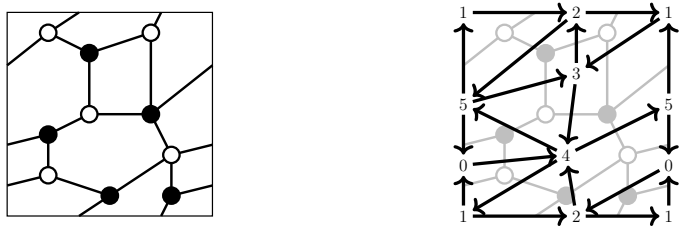
\medskip

\begin{center}
%
\end{center}

\medskip

We fix the perfect matching $\sfP_0$, 
then the following $\sfP_1, \cdots, \sfP_4$ are extremal perfect matchings corresponding to $v_1, \cdots, v_4$. 

\medskip

%%%%%perfect_matching%%%%%
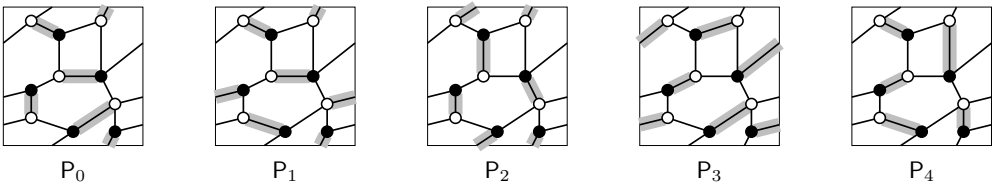
\begin{center}
{\scalebox{0.9}{
%
}}
\end{center}

\medskip

Then, we have the following table. For simplicity, we denote by $T(c,d)$ a divisorial ideal $T(0,0,c,d)$. 
We remark that $(e_0\sfB)^*\cong e_3\sfB$, $(e_1\sfB)^*\cong e_2\sfB$ and $(e_4\sfB)^*\cong e_5\sfB$. 

\medskip

\begin{center}
{\small{
\begin{tabular}{|l|l|l|l|l|l|} \hline 
$e_0\sfB$&$e_1\sfB$&$e_2\sfB$&$e_3\sfB$&$e_4\sfB$ &$e_5\sfB$\\ \hline
$T_{00}\cong R$&$T_{10}\cong T(1,0)$&$T_{20}\cong T(4,1)$&$T_{30}\cong T(5,1)$&$T_{40}\cong T(3,1)$&$T_{50}\cong T(2,0)$ \\
$T_{01}\cong T(-1,0)$&$T_{11}\cong R$&$T_{21}\cong T(3,1)$&$T_{31}\cong T(4,1)$&$T_{41}\cong T(2,1)$&$T_{51}\cong T(1,0)$ \\
$T_{02}\cong T(0,1)$&$T_{12}\cong T(1,1)$&$T_{22}\cong R$&$T_{32}\cong T(1,0)$&$T_{42}\cong T(-1,0)$&$T_{52}\cong T(2,1)$ \\
$T_{03}\cong T(-1,1)$&$T_{13}\cong T(0,1)$&$T_{23}\cong T(-1,0)$&$T_{33}\cong R$&$T_{43}\cong T(-2,0)$&$T_{53}\cong T(1,1)$ \\
$T_{04}\cong T(1,1)$&$T_{14}\cong T(2,1)$&$T_{24}\cong T(1,0)$&$T_{34}\cong T(2,0)$&$T_{44}\cong R$&$T_{54}\cong T(3,1)$ \\
$T_{05}\cong T(-2,0)$&$T_{15}\cong T(-1,0)$&$T_{25}\cong T(2,1)$&$T_{35}\cong T(3,1)$&$T_{45}\cong T(1,1)$&$T_{55}\cong R$ \\ \hline 
\end{tabular} 
}}
\end{center}

\medskip

%%%%%%%%%%%%%%%%%%%%%%%%%%%%%%%%%%%%%%%%%%%%
\subsubsection{Exchange graph of type 6c}
By combining the above results, we describe the exchange graph $\EG(\TMMG(R))$ for the case of type 6c as follows, 
and it is actually connected.

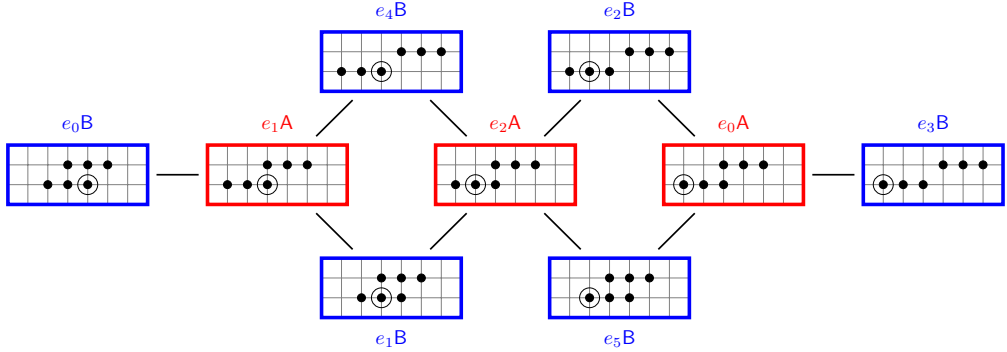
\begin{center}
{\scalebox{0.75}{
%
}};

\draw[line width=0.03cm] (A2)--(B5) ; \draw[line width=0.03cm] (A2)--(B2) ; \draw[line width=0.03cm] (A2)--(B4) ; \draw[line width=0.03cm] (A2)--(B1) ;
\draw[line width=0.03cm] (A0)--(B5) ; \draw[line width=0.03cm] (A0)--(B2) ; \draw[line width=0.03cm] (A0)--(B3) ;
\draw[line width=0.03cm] (A1)--(B4) ; \draw[line width=0.03cm] (A1)--(B1) ; \draw[line width=0.03cm] (A1)--(B0) ;

\end{tikzpicture} }} 
\end{center}

%%%%%%%%%%%%%%%%%%%%%%%%%%%%%%%%%%%%%%%%%%%%%%%%%%%%%%%
\subsection{Type 7a}
\label{type7a}

In this subsection, we consider the reflexive polygon of type 7a. 
Thus, let $R$ be the three dimensional complete local Gorenstein toric singularity defined by the cone $\sigma$: 
\[
\sigma=\mathrm{Cone}\{v_1=(1,0,1), v_2=(0,1,1), v_3=(-1,1,1), v_4=(-1,-1,1), v_5=(1,-1,1) \}. 
\]
As an element in $\Cl(R)$, we obtain $[I_1]+2[I_4]+2[I_5]=0$, $[I_2]-4[I_4]-2[I_5]=0$, $[I_3]+3[I_4]+[I_5]=0$. 
Therefore, we have that $\Cl(R)\cong\ZZ^2$, and each divisorial ideal is represented by $T(0,0,0,d,e)$ where $d, e\in\ZZ$. 
Also, there are four consistent dimer models (up to right equivalence) written below which give the reflexive polygon of type 7a as the perfect matching polygon. 

\medskip

\subsubsection{Type 7a-1}
First, we consider the following consistent dimer model and the associated QP $(Q, W_Q)$. 
We denote by $\sfA$ the complete Jacobian algebra $\calP(Q, W_Q)$, 
and give splitting MM generators $e_i\sfA$. 

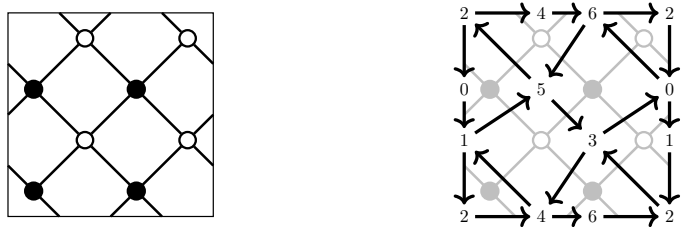
\medskip

%%%%%%%%%%%%%%%%%%%%%%%%%%%%%%%%%%%
%%%%%%%%%%%%Type_7a-1%%%%%%%%%%%%%%%%%%
%%%%%%%%%%%%%%%%%%%%%%%%%%%%%%%%%%%

\begin{center}

\end{center}

\medskip

In the following figure, we fix the perfect matching $\sfP_0$.  
Then $\sfP_1, \cdots, \sfP_5$ are extremal perfect matchings corresponding to $v_1, \cdots, v_5$. 
Namely, $v_i=(h(\sfP_i, \sfP_0), 1)\in\ZZ^3$.

\medskip

%%%%%perfect_matching%%%%%

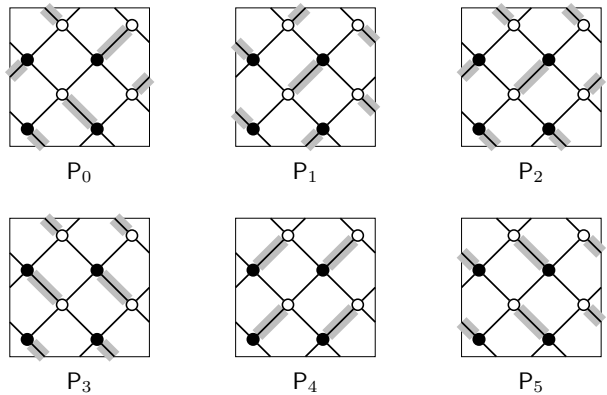
\begin{center}
{\scalebox{0.9}{
%
}}
\end{center}

\medskip

Then, we have the following table. For simplicity, we denote by $T(d,e)$ a divisorial ideal $T(0,0,0,d,e)$. 
We remark that $(e_0\sfA)^*\cong e_1\sfA$, $(e_2\sfA)^*\cong e_2\sfA$, $(e_3\sfA)^*\cong e_5\sfA$, and $(e_4\sfA)^*\cong e_6\sfA$. 

\medskip

\begin{center}
{\small{
\begin{tabular}{|l|l|l|l|} \hline 
$e_0\sfA$&$e_1\sfA$&$e_2\sfA$&$e_3\sfA$\\ \hline
$T_{00}\cong R$&$T_{10}\cong T(-4,-2)$&$T_{20}\cong T(-2,-1)$&$T_{30}\cong T(-3,-1)$ \\
$T_{01}\cong T(4,2)$&$T_{11}\cong R$&$T_{21}\cong T(2,1)$&$T_{31}\cong T(1,1)$ \\
$T_{02}\cong T(2,1)$&$T_{12}\cong T(-2,-1)$&$T_{22}\cong R$&$T_{32}\cong T(-1,0)$\\
$T_{03}\cong T(3,1)$&$T_{13}\cong T(-1,-1)$&$T_{23}\cong T(1,0)$&$T_{33}\cong R$ \\
$T_{04}\cong T(3,2)$&$T_{14}\cong T(-1,0)$&$T_{24}\cong T(1,1)$&$T_{34}\cong T(0,1)$ \\
$T_{05}\cong T(1,1)$&$T_{15}\cong T(-3,-1)$&$T_{25}\cong T(-1,0)$&$T_{35}\cong T(-2,0)$ \\ 
$T_{06}\cong T(1,0)$&$T_{16}\cong T(-3,-2)$&$T_{26}\cong T(-1,-1)$&$T_{36}\cong T(-2,-1)$ \\ \hline
$e_4\sfA$ &$e_5\sfA$&$e_6\sfA$\\ \cline{1-3}
$T_{40}\cong T(-3,-2)$&$T_{50}\cong T(-1,-1)$&$T_{60}\cong T(-1,0)$ \\
$T_{41}\cong T(1,0)$&$T_{51}\cong T(3,1)$&$T_{61}\cong T(3,2)$ \\
$T_{42}\cong T(-1,-1)$&$T_{52}\cong T(1,0)$&$T_{62}\cong T(1,1)$ \\
$T_{43}\cong T(0,-1)$&$T_{53}\cong T(2,0)$&$T_{63}\cong T(2,1)$ \\
$T_{44}\cong R$&$T_{54}\cong T(2,1)$&$T_{64}\cong T(2,2)$ \\
$T_{45}\cong T(-2,-1)$&$T_{55}\cong R$&$T_{65}\cong T(0,1)$ \\ 
$T_{46}\cong T(-2,-2)$&$T_{56}\cong T(0,-1)$&$T_{66}\cong R$ \\ \cline{1-3}
\end{tabular} 
}}
\end{center}

\subsubsection{Type 7a-2}

Next, we consider the following consistent dimer model and the associated QP $(Q, W_Q)$. 
We denote by $\sfB$ the complete Jacobian algebra $\calP(Q, W_Q)$, and consider splitting MM generators $e_i\sfB$. 

%%%%%%%%%%%%%%%%%%%%%%%%%%%%%%%%%%%
%%%%%%%%%%%%Type_7a-2%%%%%%%%%%%%%%%%%%
%%%%%%%%%%%%%%%%%%%%%%%%%%%%%%%%%%%
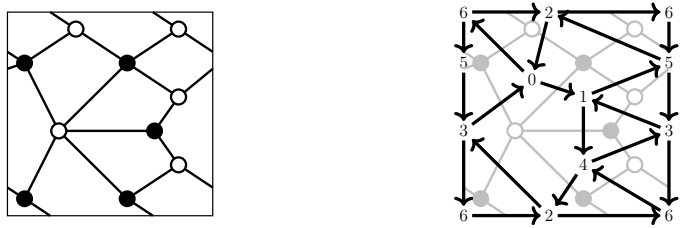
\begin{center}
%
\end{center}

\medskip

In the following figure, we fix the perfect matching $\sfP_0$.  
Then $\sfP_1, \cdots, \sfP_5$ are extremal perfect matchings corresponding to $v_1, \cdots, v_5$. 

\medskip

%%%%%perfect_matching%%%%%

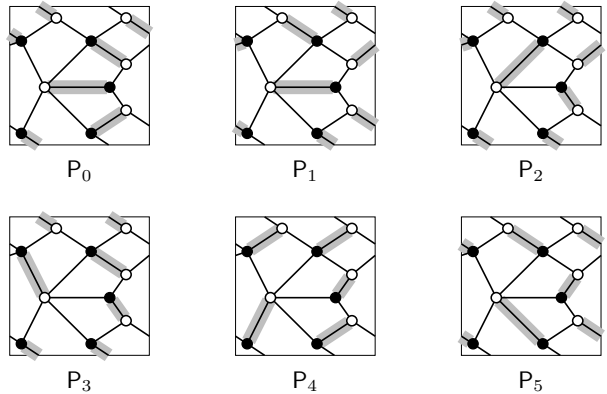
\begin{center}
{\scalebox{0.9}{
%
}}
\end{center}

\medskip

Then, we have the following tables. For simplicity, we denote by $T(d,e)$ a divisorial ideal $T(0,0,0,d,e)$. 

\medskip

\begin{center}
{\small{
\begin{tabular}{|l|l|l|l|} \hline 
$e_0\sfB$&$e_1\sfB$&$e_2\sfB$&$e_3\sfB$\\ \hline
$T_{00}\cong R$&$T_{10}\cong T(-4,-2)$&$T_{20}\cong T(-2,-1)$&$T_{30}\cong T(-3,-1)$ \\
$T_{01}\cong T(4,2)$&$T_{11}\cong R$&$T_{21}\cong T(2,1)$&$T_{31}\cong T(1,1)$ \\
$T_{02}\cong T(2,1)$&$T_{12}\cong T(-2,-1)$&$T_{22}\cong R$&$T_{32}\cong T(-1,0)$ \\
$T_{03}\cong T(3,1)$&$T_{13}\cong T(-1,-1)$&$T_{23}\cong T(1,0)$&$T_{33}\cong R$ \\
$T_{04}\cong T(2,0)$&$T_{14}\cong T(-2,-2)$&$T_{24}\cong T(0,-1)$&$T_{34}\cong T(-1,-1)$ \\
$T_{05}\cong T(1,1)$&$T_{15}\cong T(-3,-1)$&$T_{25}\cong T(-1,0)$&$T_{35}\cong T(-2,0)$ \\ 
$T_{06}\cong T(1,0)$&$T_{16}\cong T(-3,-2)$&$T_{26}\cong T(-1,-1)$&$T_{36}\cong T(-2,-1)$ \\ \hline
$e_4\sfB$&$e_5\sfB$&$e_6\sfB$\\ \cline{1-3}
$T_{40}\cong T(-2,0)$&$T_{50}\cong T(-1,-1)$&$T_{60}\cong T(-1,0)$ \\
$T_{41}\cong T(2,2)$&$T_{51}\cong T(3,1)$&$T_{61}\cong T(3,2)$ \\
$T_{42}\cong T(0,1)$&$T_{52}\cong T(1,0)$&$T_{62}\cong T(1,1)$ \\
$T_{43}\cong T(1,1)$&$T_{53}\cong T(2,0)$&$T_{63}\cong T(2,1)$ \\
$T_{44}\cong R$&$T_{54}\cong T(1,-1)$&$T_{64}\cong T(1,0)$ \\
$T_{45}\cong T(-1,1)$&$T_{55}\cong R$&$T_{65}\cong T(0,1)$ \\ 
$T_{46}\cong T(-1,0)$&$T_{56}\cong T(0,-1)$&$T_{66}\cong R$ \\ \cline{1-3}
\end{tabular} 
}}
\end{center}

\medskip

Let $x$ (resp. $y$) be a $1$-cycle on the two torus $\sfT$ which corresponds to $(1,0)\in\rmH_1(\sfT)$ (resp. $(0,1)\in\rmH_1(\sfT)$) in the above case. 
By replacing these cycles by $y, x$ respectively, we have the following dimer model and the associated QP $(Q, W_Q)$. 
We denote by $\sfB^\prime$ the complete Jacobian algebra $\calP(Q, W_Q)$. 
Since this dimer model is right equivalent to the previous one, we have that $\sfB\cong\sfB^\prime$ as an $R$-algebra. 
However, a change of cycles induces an automorphism on $R$, and it gives another description of each splitting MM generator (see Remark~\ref{rem_auto}). 
Thus we also consider splitting MM generators $e_i\sfB^\prime$. 

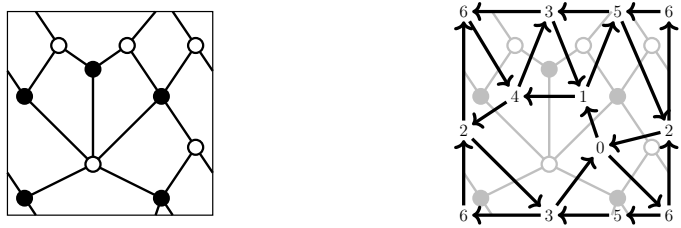
\medskip 

\begin{center}

\end{center}

\medskip

In the following figure, we fix the perfect matching $\sfP_0$.  
Then $\sfP_1, \cdots, \sfP_5$ are extremal perfect matchings corresponding to $v_1, \cdots, v_5$. 

\medskip

%%%%%perfect_matching%%%%%

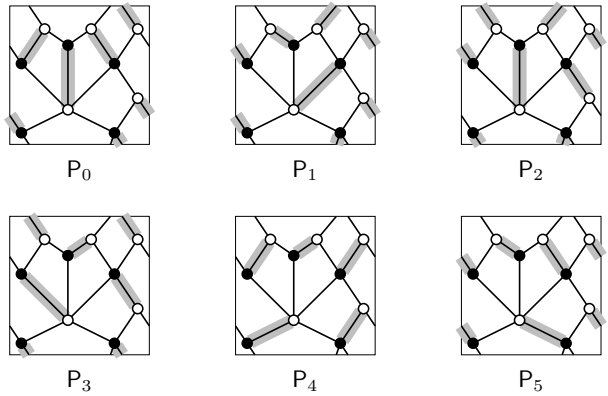
\begin{center}
{\scalebox{0.9}{
%
}}
\end{center}

\medskip

Then, we have the following tables. For simplicity, we denote by $T(d,e)$ a divisorial ideal $T(0,0,0,d,e)$. 

\medskip

\begin{center}
{\small{
\begin{tabular}{|l|l|l|l|} \hline 
$e_0\sfB^\prime$&$e_1\sfB^\prime$&$e_2\sfB^\prime$&$e_3\sfB^\prime$\\ \hline
$T_{00}\cong R$&$T_{10}\cong T(2,2)$&$T_{20}\cong T(1,1)$&$T_{30}\cong T(0,1)$ \\
$T_{01}\cong T(-2,-2)$&$T_{11}\cong R$&$T_{21}\cong T(-1,-1)$&$T_{31}\cong T(-2,-1)$ \\
$T_{02}\cong T(-1,-1)$&$T_{12}\cong T(1,1)$&$T_{22}\cong R$&$T_{32}\cong T(-1,0)$ \\
$T_{03}\cong T(0,-1)$&$T_{13}\cong T(2,1)$&$T_{23}\cong T(1,0)$&$T_{33}\cong R$ \\
$T_{04}\cong T(2,0)$&$T_{14}\cong T(4,2)$&$T_{24}\cong T(3,1)$&$T_{34}\cong T(2,1)$ \\
$T_{05}\cong T(-2,-1)$&$T_{15}\cong T(0,1)$&$T_{25}\cong T(-1,0)$&$T_{35}\cong T(-2,0)$ \\ 
$T_{06}\cong T(1,0)$&$T_{16}\cong T(3,2)$&$T_{26}\cong T(2,1)$&$T_{36}\cong T(1,1)$ \\ \hline
$e_4\sfB^\prime$&$e_5\sfB^\prime$&$e_6\sfB^\prime$\\ \cline{1-3}
$T_{40}\cong T(-2,0)$&$T_{50}\cong T(2,1)$&$T_{60}\cong T(-1,0)$ \\
$T_{41}\cong T(-4,-2)$&$T_{51}\cong T(0,-1)$&$T_{61}\cong T(-3,-2)$ \\
$T_{42}\cong T(-3,-1)$&$T_{52}\cong T(1,0)$&$T_{62}\cong T(-2,-1)$ \\
$T_{43}\cong T(-2,-1)$&$T_{53}\cong T(2,0)$&$T_{63}\cong T(-1,-1)$ \\
$T_{44}\cong R$&$T_{54}\cong T(4,1)$&$T_{64}\cong T(1,0)$ \\
$T_{45}\cong T(-4,-1)$&$T_{55}\cong R$&$T_{65}\cong T(-3,-1)$ \\ 
$T_{46}\cong T(-1,0)$&$T_{56}\cong T(3,1)$&$T_{66}\cong R$ \\ \cline{1-3}
\end{tabular} 
}}
\end{center}

\medskip

\subsubsection{Type 7a-3}

Next, we consider the following consistent dimer model and the associated QP $(Q, W_Q)$. 
We denote by $\sfC$ the complete Jacobian algebra $\calP(Q, W_Q)$, and consider splitting MM generators $e_i\sfC$. 

%%%%%%%%%%%%%%%%%%%%%%%%%%%%%%%%%%%
%%%%%%%%%%%%Type_7a-3%%%%%%%%%%%%%%%%%%
%%%%%%%%%%%%%%%%%%%%%%%%%%%%%%%%%%%
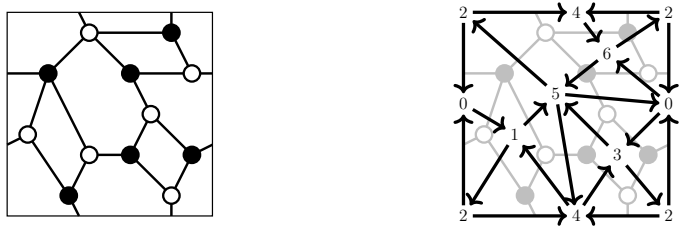
\begin{center}
%
\end{center}

\medskip

In the following figure, we fix the perfect matching $\sfP_0$.  
Then $\sfP_1, \cdots, \sfP_5$ are extremal perfect matchings corresponding to $v_1, \cdots, v_5$. 

\medskip

%%%%%perfect_matching%%%%%

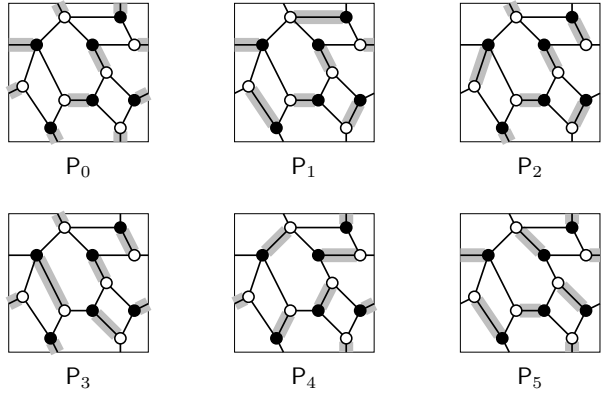
\begin{center}
{\scalebox{0.9}{
%
}}
\end{center}

\medskip

Then, we have the following tables. For simplicity, we denote by $T(d,e)$ a divisorial ideal $T(0,0,0,d,e)$. 

\medskip

\begin{center}
{\small{
\begin{tabular}{|l|l|l|l|} \hline 
$e_0\sfC$&$e_1\sfC$&$e_2\sfC$&$e_3\sfC$\\ \hline
$T_{00}\cong R$&$T_{10}\cong T(-4,-2)$&$T_{20}\cong T(-2,-1)$&$T_{30}\cong T(0,-1)$ \\
$T_{01}\cong T(4,2)$&$T_{11}\cong R$&$T_{21}\cong T(2,1)$&$T_{31}\cong T(4,1)$ \\
$T_{02}\cong T(2,1)$&$T_{12}\cong T(-2,-1)$&$T_{22}\cong R$&$T_{32}\cong T(2,0)$ \\
$T_{03}\cong T(0,1)$&$T_{13}\cong T(-4,-1)$&$T_{23}\cong T(-2,0)$&$T_{33}\cong R$ \\
$T_{04}\cong T(3,2)$&$T_{14}\cong T(-1,0)$&$T_{24}\cong T(1,1)$&$T_{34}\cong T(3,1)$ \\
$T_{05}\cong T(1,1)$&$T_{15}\cong T(-3,-1)$&$T_{25}\cong T(-1,0)$&$T_{35}\cong T(1,0)$ \\ 
$T_{06}\cong T(1,0)$&$T_{16}\cong T(-3,-2)$&$T_{26}\cong T(-1,-1)$&$T_{36}\cong T(1,-1)$ \\ \hline
$e_4\sfC$&$e_5\sfC$&$e_6\sfC$\\ \cline{1-3}
$T_{40}\cong T(-3,-2)$&$T_{50}\cong T(-1,-1)$&$T_{60}\cong T(-1,0)$ \\
$T_{41}\cong T(1,0)$&$T_{51}\cong T(3,1)$&$T_{61}\cong T(3,2)$ \\
$T_{42}\cong T(-1,-1)$&$T_{52}\cong T(1,0)$&$T_{62}\cong T(1,1)$ \\
$T_{43}\cong T(-3,-1)$&$T_{53}\cong T(-1,0)$&$T_{63}\cong T(-1,1)$ \\
$T_{44}\cong R$&$T_{54}\cong T(2,1)$&$T_{64}\cong T(2,2)$ \\
$T_{45}\cong T(-2,-1)$&$T_{55}\cong R$&$T_{65}\cong T(0,1)$ \\ 
$T_{46}\cong T(-2,-2)$&$T_{56}\cong T(0,-1)$&$T_{66}\cong R$ \\ \cline{1-3}
\end{tabular} 
}}
\end{center}

\medskip

%%%%%%%%%%%%%%%%%%%%%%%%%%%%%%%%%%%%%%%%%%%%
%\newpage
\subsubsection{Exchange graph of type 7a}
By combining the above results, we describe the exchange graph $\EG(\TMMG(R))$ for the case of type 7a as follows, 
and it is actually connected. 

%\newpage
%%%%%%%%%%%%%%%%%%%%%
%%% type 7a exchange graph %%%%%%%%%%%%%
%%%%%%%%%%%%%%%%%%%%%

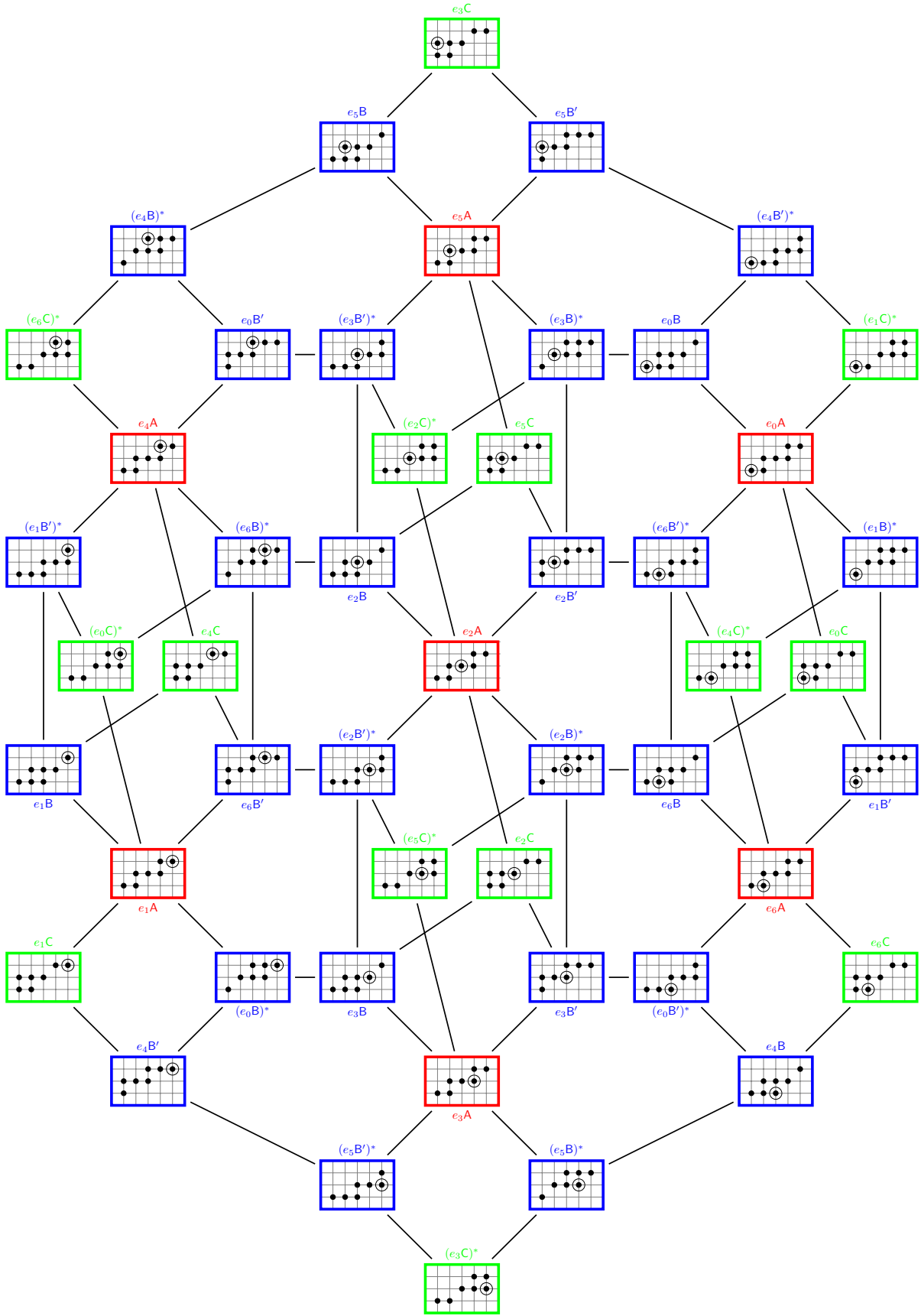
\begin{center}
{\scalebox{0.6}{
%
}};

\draw[line width=0.0375cm] (A2)--(B2); \draw[line width=0.0375cm] (A2)--(B2+); \draw[line width=0.0375cm] (A2)--(C2); 
\draw[line width=0.0375cm] (A2)--(C2+); \draw[line width=0.0375cm] (A2)--(D2); \draw[line width=0.0375cm] (A2)--(D2+); 
\draw[line width=0.0375cm] (B2)--(D3+); \draw[line width=0.0375cm] (B2)--(C5); \draw[line width=0.0375cm] (C2+)--(D3+); 
\draw[line width=0.0375cm] (C2+)--(B3+); \draw[line width=0.0375cm] (D2)--(C5); \draw[line width=0.0375cm] (D2)--(B3+); 
\draw[line width=0.0375cm] (D3+)--(A5); \draw[line width=0.0375cm] (A5)--(C5); \draw[line width=0.0375cm] (A5)--(B3+);
\draw[line width=0.0375cm] (B2+)--(C5+); \draw[line width=0.0375cm] (B2+)--(D3); \draw[line width=0.0375cm] (D2+)--(C5+); 
\draw[line width=0.0375cm] (D2+)--(B3); \draw[line width=0.0375cm] (C2)--(D3); \draw[line width=0.0375cm] (C2)--(B3);
\draw[line width=0.0375cm] (A3)--(C5+);  \draw[line width=0.0375cm] (A3)--(D3); \draw[line width=0.0375cm] (A3)--(B3); 
\draw[line width=0.0375cm] (A3)--(B5+); \draw[line width=0.0375cm] (A3)--(D5+); \draw[line width=0.0375cm] (B5+)--(C3+); 
\draw[line width=0.0375cm] (D5+)--(C3+); \draw[line width=0.0375cm] (A5)--(B5); \draw[line width=0.0375cm] (A5)--(D5); 
\draw[line width=0.0375cm] (C3)--(B5); \draw[line width=0.0375cm] (C3)--(D5); \draw[line width=0.0375cm] (D6)--(D2+); 
\draw[line width=0.0375cm] (D3)--(D0+); \draw[line width=0.0375cm] (D0)--(D3+); \draw[line width=0.0375cm] (D2)--(D6+); 
\draw[line width=0.0375cm] (A6)--(C6); \draw[line width=0.0375cm] (A6)--(D0+); \draw[line width=0.0375cm] (B4)--(C6); 
\draw[line width=0.0375cm] (B4)--(D0+); \draw[line width=0.0375cm] (B4)--(B5+); \draw[line width=0.0375cm] (B4+)--(B5); 
\draw[line width=0.0375cm] (B4+)--(D0); \draw[line width=0.0375cm] (B4+)--(C6+);
\draw[line width=0.0375cm] (A6)--(B6); \draw[line width=0.0375cm] (B6)--(D6+); \draw[line width=0.0375cm] (A4)--(B6+); 
\draw[line width=0.0375cm] (B6+)--(D6); \draw[line width=0.0375cm] (A4)--(C6+); \draw[line width=0.0375cm] (A4)--(D0);
\draw[line width=0.0375cm] (A6)--(D1); \draw[line width=0.0375cm] (A6)--(C4+);  \draw[line width=0.0375cm] (C4+)--(D6+); 
\draw[line width=0.0375cm] (C4+)--(B1+); \draw[line width=0.0375cm] (D1)--(B1+); \draw[line width=0.0375cm] (D1)--(C0); 
\draw[line width=0.0375cm] (A0)--(B1+); \draw[line width=0.0375cm] (A0)--(D6+); \draw[line width=0.0375cm] (A0)--(C0); 
\draw[line width=0.0375cm] (A0)--(C1+); \draw[line width=0.0375cm] (A0)--(B0); \draw[line width=0.0375cm] (C1+)--(D4+); 
\draw[line width=0.0375cm] (B0)--(D4+); \draw[line width=0.0375cm] (D5)--(D4+); 
\draw[line width=0.0375cm] (A4)--(D1+); \draw[line width=0.0375cm] (A4)--(C4); \draw[line width=0.0375cm] (C4)--(D6); 
\draw[line width=0.0375cm] (C4)--(B1); \draw[line width=0.0375cm] (D1+)--(B1);  \draw[line width=0.0375cm] (D1+)--(C0+); 
\draw[line width=0.0375cm] (A1)--(B1); \draw[line width=0.0375cm] (A1)--(D6); \draw[line width=0.0375cm] (A1)--(C0+); 
\draw[line width=0.0375cm] (A1)--(C1); \draw[line width=0.0375cm] (A1)--(B0+); \draw[line width=0.0375cm] (C1)--(D4); 
\draw[line width=0.0375cm] (B0+)--(D4); \draw[line width=0.0375cm] (D5+)--(D4); 
\draw[line width=0.0375cm] (B0)--(B3+); \draw[line width=0.0375cm] (B6)--(B2+); \draw[line width=0.0375cm] (B6)--(C0); 
\draw[line width=0.0375cm] (B0+)--(B3); \draw[line width=0.0375cm] (B6+)--(B2); \draw[line width=0.0375cm] (B6+)--(C0+);

\end{tikzpicture} }}
\end{center}

%%%%%%%%%%%%%%%%%%%%%%%%%%%%%%%%%%%%%%%%%%%%%%%%%%%%%%%
\subsection{Type 7b}
\label{type7b}

In this subsection, we consider the reflexive polygon of type 7b. 
Thus, let $R$ be the three dimensional complete local Gorenstein toric singularity defined by the cone $\sigma$: 
\[
\sigma=\mathrm{Cone}\{v_1=(1,0,1), v_2=(0,1,1), v_3=(-2,-1,1), v_4=(1,-1,1) \}. 
\]
As an element in $\Cl(R)$, we obtain $[I_1]-6[I_3]=0$, $[I_2]+3[I_3]=0$, $[I_4]+4[I_3]=0$. 
Therefore, we have that $\Cl(R)\cong\ZZ$, and each divisorial ideal is represented by $T(0,0,c,0)$ where $c\in\ZZ$. 
There is the unique consistent dimer model (up to right equivalence) written below which give the reflexive polygon of type 7b as the perfect matching polygon. 

\medskip

\subsubsection{Type 7b-1}
We consider the following consistent dimer model and the associated QP $(Q, W_Q)$. 
We denote by $\sfA$ the complete Jacobian algebra $\calP(Q, W_Q)$, and consider splitting MM generators $e_i\sfA$. 

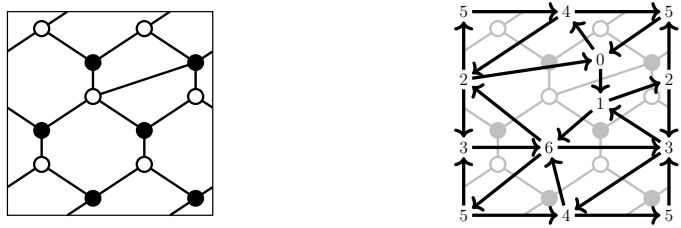
\medskip

\begin{center}

\end{center}

\medskip

We fix the perfect matching $\sfP_0$, 
then the following $\sfP_1, \cdots, \sfP_4$ are extremal perfect matchings corresponding to $v_1, \cdots, v_4$. 
Namely, $v_i=(h(\sfP_i, \sfP_0), 1)\in\ZZ^3$. 

\medskip

%%%%%perfect_matching%%%%%
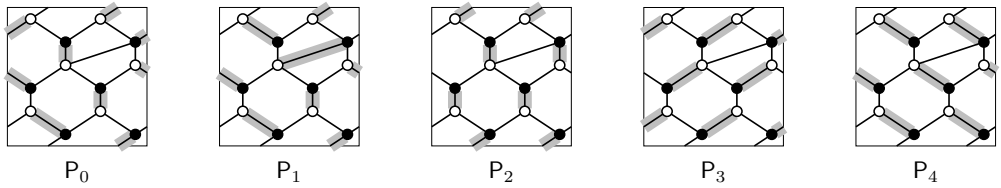
\begin{center}
{\scalebox{0.9}{
%
}}
\end{center}

\medskip

Then, we have the following table. For simplicity, we denote by $T(c)$ a divisorial ideal $T(0,0,c,0)$. 
We remark that $e_0\sfA\cong(e_1\sfA)^*$, $e_2\sfA\cong(e_2\sfA)^*$, $e_3\sfA\cong(e_4\sfA)^*$, and $e_5\sfA\cong(e_6\sfA)^*$. 
\medskip

\begin{center}
{\footnotesize{
\begin{tabular}{|l|l|l|l|l|l|l|} \hline 
$e_0\sfA$&$e_1\sfA$&$e_2\sfA$&$e_3\sfA$&$e_4\sfA$ &$e_5\sfA$&$e_6\sfA$\\ \hline
$T_{00}\cong R$&$T_{10}\cong T(-6)$&$T_{20}\cong T(-3)$&$T_{30}\cong T(-5)$&$T_{40}\cong T(-1)$&$T_{50}\cong T(-4)$&$T_{60}\cong T(-2)$ \\
$T_{01}\cong T(6)$&$T_{11}\cong R$&$T_{21}\cong T(3)$&$T_{31}\cong T(1)$&$T_{41}\cong T(5)$&$T_{51}\cong T(2)$&$T_{61}\cong T(4)$ \\
$T_{02}\cong T(3)$&$T_{12}\cong T(-3)$&$T_{22}\cong R$&$T_{32}\cong T(-2)$&$T_{42}\cong T(2)$&$T_{52}\cong T(-1)$&$T_{62}\cong T(1)$ \\
$T_{03}\cong T(5)$&$T_{13}\cong T(-1)$&$T_{23}\cong T(2)$&$T_{33}\cong R$&$T_{43}\cong T(4)$&$T_{53}\cong T(1)$&$T_{63}\cong T(3)$ \\
$T_{04}\cong T(1)$&$T_{14}\cong T(-5)$&$T_{24}\cong T(-2)$&$T_{34}\cong T(-4)$&$T_{44}\cong R$&$T_{54}\cong T(-3)$&$T_{64}\cong T(-1)$ \\
$T_{05}\cong T(4)$&$T_{15}\cong T(-2)$&$T_{25}\cong T(1)$&$T_{35}\cong T(-1)$&$T_{45}\cong T(3)$&$T_{55}\cong R$&$T_{65}\cong T(2)$ \\ 
$T_{06}\cong T(2)$&$T_{16}\cong T(-4)$&$T_{26}\cong T(-1)$&$T_{36}\cong T(-3)$&$T_{46}\cong T(1)$&$T_{56}\cong T(-2)$&$T_{66}\cong R$ \\ \hline 
\end{tabular} 
}}
\end{center}

\medskip

\subsubsection{Exchange graph of type 7b}
By the above computation, we describe the exchange graph $\EG(\TMMG(R))$ for the case of type 7b as follows, 
and it is actually connected. 

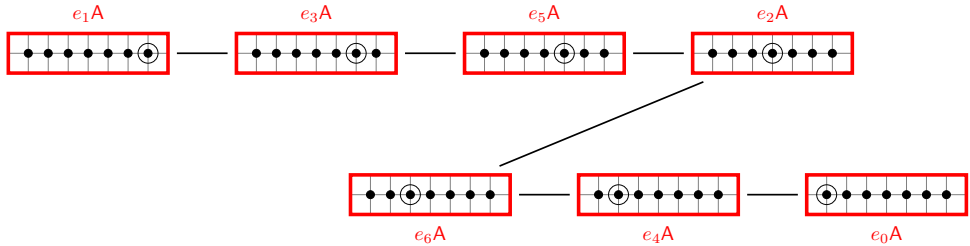
\begin{center}
{\scalebox{0.75}{
%
}};

\draw[line width=0.03cm] (A1)--(A3); \draw[line width=0.03cm] (A3)--(A5); \draw[line width=0.03cm] (A5)--(A2);
\draw[line width=0.03cm] (A2)--(A6); \draw[line width=0.03cm] (A6)--(A4); \draw[line width=0.03cm] (A4)--(A0);

\end{tikzpicture} }}
\end{center}

%%%%%%%%%%%%%%%%%%%%%%%%%%%%%%%%%%%%%%%%%%%%%%%%%%%%%%%
\subsection{Type 8a}
\label{type8a}
In this subsection, we consider the reflexive polygon of type 8a. 
Thus, let $R$ be the three dimensional complete local Gorenstein toric singularity defined by the cone $\sigma$: 
\[
\sigma=\mathrm{Cone}\{v_1=(1,1,1), v_2=(-1,1,1), v_3=(-1,-1,1), v_4=(1,-1,1) \}. 
\]
As an element in $\Cl(R)$, we obtain $2[I_1]=-2[I_2]=2[I_3]=-2[I_4]$, and $[I_4]=[I_1]+[I_2]-[I_3]$.  
Therefore, we have that $\Cl(R)\cong\ZZ\times(\ZZ/2\ZZ)^2$, 
and each divisorial ideal is represented by $T(a,b,c,0)$ where $a\in\ZZ, \, b,c\in\ZZ/2\ZZ$. 
Also, there are four consistent dimer models (up to right equivalence) written below which give the reflexive polygon of type 8a as the perfect matching polygon. 

\medskip

\subsubsection{Type 8a-1}
We consider the following consistent dimer model and the associated QP $(Q, W_Q)$. 
We denote by $\sfA$ the complete Jacobian algebra $\calP(Q, W_Q)$, and consider splitting MM generators $e_i\sfA$. 

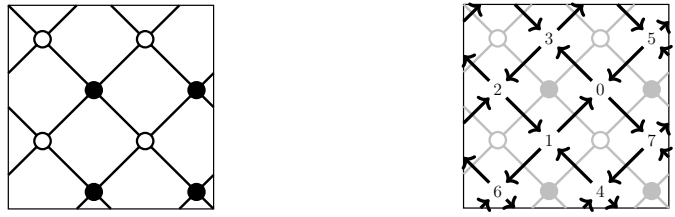
\medskip

%%%%%%%%%%%%%%%%%%%%%%%%%%%%%%%%%%%
%%%%%%%%%%%%Type_8a-1%%%%%%%%%%%%%%%%%%
%%%%%%%%%%%%%%%%%%%%%%%%%%%%%%%%%%%
\begin{center}

\end{center}

\medskip

We fix the perfect matching $\sfP_0$, 
then the following $\sfP_1, \cdots, \sfP_4$ are extremal perfect matchings corresponding to $v_1, \cdots, v_4$. 
Namely, $v_i=(h(\sfP_i, \sfP_0), 1)\in\ZZ^3$. 

\medskip

%%%%%perfect_matching%%%%%
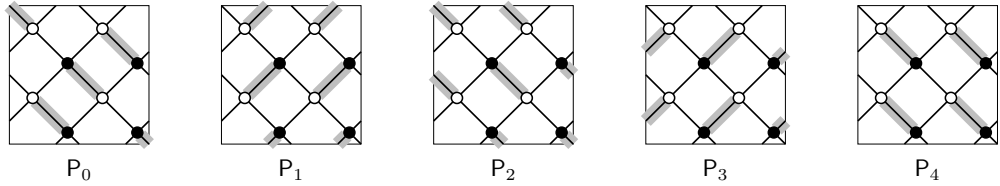
\begin{center}
{\scalebox{0.9}{
%
}}
\end{center}

\medskip

Then, we have the following, and $e_0\sfA\cong e_2\sfA\cong e_4\sfA\cong e_6\sfA$, $e_1\sfA\cong e_3\sfA\cong e_5\sfA\cong e_7\sfA$, 
and $(e_0\sfA)^*\cong e_1\sfA$. 

\medskip

\begin{center}
{\small{
\begin{tabular}{|l|l|l|l|} \hline 
$e_0\sfA$&$e_1\sfA$&$e_2\sfA$&$e_3\sfA$ \\ \hline
$T_{00}\cong R$&$T_{10}\cong T(0,1,0,0)$&$T_{20}\cong T(1,1,0,0)$&$T_{30}\cong T(-2,0,1,0)$ \\
$T_{01}\cong T(2,1,0,0)$&$T_{11}\cong R$&$T_{21}\cong T(1,0,0,0)$&$T_{31}\cong T(0,1,1,0)$ \\
$T_{02}\cong T(1,1,0,0)$&$T_{12}\cong T(-1,0,0,0)$&$T_{22}\cong R$&$T_{32}\cong T(-1,1,1,0)$ \\
$T_{03}\cong T(0,0,1,0)$&$T_{13}\cong T(0,1,1,0)$&$T_{23}\cong T(1,1,1,0)$&$T_{33}\cong R$\\ 
$T_{04}\cong T(0,1,1,0)$&$T_{14}\cong T(-2,0,1,0)$&$T_{24}\cong T(-1,0,1,0)$&$T_{34}\cong T(0,1,0,0)$ \\
$T_{05}\cong T(1,1,1,0)$&$T_{15}\cong T(-1,0,1,0)$&$T_{25}\cong T(0,0,1,0)$&$T_{35}\cong T(1,1,0,0)$ \\
$T_{06}\cong T(-1,0,1,0)$&$T_{16}\cong T(-1,1,1,0)$&$T_{26}\cong T(0,1,1,0)$&$T_{36}\cong T(-1,0,0,0)$ \\
$T_{07}\cong T(1,0,0,0)$&$T_{17}\cong T(1,1,0,0)$&$T_{27}\cong T(2,1,0,0)$&$T_{37}\cong T(-1,0,1,0)$ \\ \hline 
$e_4\sfA$&$e_5\sfA$&$e_6\sfA$&$e_7\sfA$ \\ \hline
$T_{40}\cong T(0,1,1,0)$&$T_{50}\cong T(-1,1,1,0)$&$T_{60}\cong T(-1,0,1,0)$&$T_{70}\cong T(-1,0,0,0)$ \\
$T_{41}\cong T(0,0,1,0)$&$T_{51}\cong T(-1,0,1,0)$&$T_{61}\cong T(1,1,1,0)$&$T_{71}\cong T(1,1,0,0)$ \\
$T_{42}\cong T(-1,0,1,0)$&$T_{52}\cong T(-2,0,1,0)$&$T_{62}\cong T(0,1,1,0)$&$T_{72}\cong T(0,1,0,0)$ \\
$T_{43}\cong T(2,1,0,0)$&$T_{53}\cong T(1,1,0,0)$&$T_{63}\cong T(1,0,0,0)$&$T_{73}\cong T(-1,0,1,0)$\\ 
$T_{44}\cong R$&$T_{54}\cong T(-1,0,0,0)$&$T_{64}\cong T(1,1,0,0)$&$T_{74}\cong T(-1,1,1,0)$ \\
$T_{45}\cong T(1,0,0,0)$&$T_{55}\cong R$&$T_{65}\cong T(2,1,0,0)$&$T_{75}\cong T(0,1,1,0)$ \\
$T_{46}\cong T(1,1,0,0)$&$T_{56}\cong T(0,1,0,0)$&$T_{66}\cong R$&$T_{76}\cong T(-2,0,1,0)$ \\
$T_{47}\cong T(1,1,1,0)$&$T_{57}\cong T(0,1,1,0)$&$T_{67}\cong T(0,0,1,0)$&$T_{77}\cong R$ \\ \hline 
\end{tabular}
}} 
\end{center}

\medskip

\subsubsection{Type 8a-2}
Next, we consider the following consistent dimer model and the associated QP $(Q, W_Q)$. 
We denote by $\sfB$ the complete Jacobian algebra $\calP(Q, W_Q)$, and consider splitting MM generators $e_i\sfB$. 

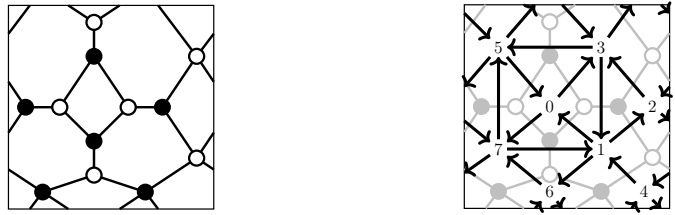
\medskip

%%%%%%%%%%%%%%%%%%%%%%%%%%%%%%%%%%%
%%%%%%%%%%%%Type_8a-2%%%%%%%%%%%%%%%%%%
%%%%%%%%%%%%%%%%%%%%%%%%%%%%%%%%%%%

\begin{center}
%
\end{center}

\medskip

We fix the perfect matching $\sfP_0$, 
then the following $\sfP_1, \cdots, \sfP_4$ are extremal perfect matchings corresponding to $v_1, \cdots, v_4$. 

\medskip

%%%%%perfect_matching%%%%%
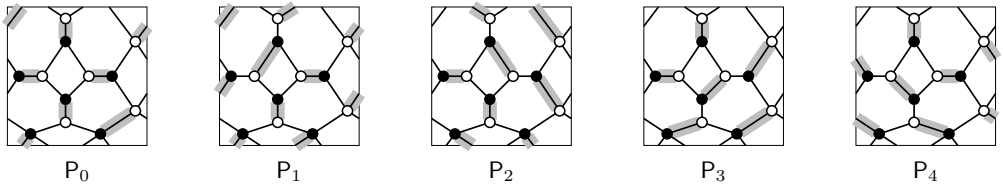
\begin{center}
{\scalebox{0.9}{
%
}}
\end{center}

\medskip

Then, we have the following table. 

\medskip

\begin{center}
{\small{
\begin{tabular}{|l|l|l|l|} \hline 
$e_0\sfB$&$e_1\sfB$&$e_2\sfB$&$e_3\sfB$ \\ \hline
$T_{00}\cong R$&$T_{10}\cong T(0,0,1,0)$&$T_{20}\cong T(2,1,1,0)$&$T_{30}\cong T(2,1,0,0)$ \\
$T_{01}\cong T(-2,0,1,0)$&$T_{11}\cong R$&$T_{21}\cong T(2,1,0,0)$&$T_{31}\cong T(0,1,1,0)$ \\
$T_{02}\cong T(-2,1,1,0)$&$T_{12}\cong T(0,1,0,0)$&$T_{22}\cong R$&$T_{32}\cong T(-2,0,1,0)$ \\
$T_{03}\cong T(0,1,0,0)$&$T_{13}\cong T(0,1,1,0)$&$T_{23}\cong T(0,0,1,0)$&$T_{33}\cong R$\\ 
$T_{04}\cong T(-2,0,0,0)$&$T_{14}\cong T(-2,0,1,0)$&$T_{24}\cong T(0,1,1,0)$&$T_{34}\cong T(0,1,0,0)$ \\
$T_{05}\cong T(-1,0,0,0)$&$T_{15}\cong T(-1,0,1,0)$&$T_{25}\cong T(1,1,1,0)$&$T_{35}\cong T(1,1,0,0)$ \\
$T_{06}\cong T(-1,1,0,0)$&$T_{16}\cong T(-1,1,1,0)$&$T_{26}\cong T(-1,0,1,0)$&$T_{36}\cong T(-1,0,0,0)$ \\
$T_{07}\cong T(-1,1,1,0)$&$T_{17}\cong T(1,1,0,0)$&$T_{27}\cong T(1,0,0,0)$&$T_{37}\cong T(-1,0,1,0)$ \\ \hline 
$e_4\sfB$&$e_5\sfB$&$e_6\sfB$&$e_7\sfB$ \\ \hline
$T_{40}\cong T(2,0,0,0)$&$T_{50}\cong T(1,0,0,0)$&$T_{60}\cong T(3,1,0,0)$&$T_{70}\cong T(1,1,1,0)$ \\
$T_{41}\cong T(0,0,1,0)$&$T_{51}\cong T(-1,0,1,0)$&$T_{61}\cong T(1,1,1,0)$&$T_{71}\cong T(1,1,0,0)$ \\
$T_{42}\cong T(0,1,1,0)$&$T_{52}\cong T(-1,1,1,0)$&$T_{62}\cong T(-1,0,1,0)$&$T_{72}\cong T(-1,0,0,0)$ \\
$T_{43}\cong T(2,1,0,0)$&$T_{53}\cong T(1,1,0,0)$&$T_{63}\cong T(1,0,0,0)$&$T_{73}\cong T(-1,0,1,0)$\\ 
$T_{44}\cong R$&$T_{54}\cong T(-1,0,0,0)$&$T_{64}\cong T(1,1,0,0)$&$T_{74}\cong T(-1,1,1,0)$ \\
$T_{45}\cong T(1,0,0,0)$&$T_{55}\cong R$&$T_{65}\cong T(2,1,0,0)$&$T_{75}\cong T(0,1,1,0)$ \\
$T_{46}\cong T(1,1,0,0)$&$T_{56}\cong T(0,1,0,0)$&$T_{66}\cong R$&$T_{76}\cong T(-2,0,1,0)$ \\
$T_{47}\cong T(1,1,1,0)$&$T_{57}\cong T(0,1,1,0)$&$T_{67}\cong T(0,0,1,0)$&$T_{77}\cong R$ \\ \hline
\end{tabular}
}} 
\end{center}

\medskip

\subsubsection{Type 8a-3}
Next, we consider the following consistent dimer model and the associated QP $(Q, W_Q)$. 
We denote by $\sfC$ the complete Jacobian algebra $\calP(Q, W_Q)$, and consider splitting MM generators $e_i\sfC$. 

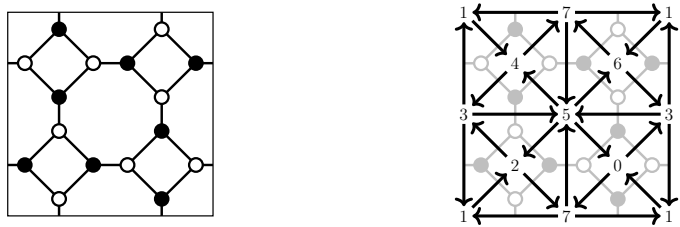
\medskip 

%%%%%%%%%%%%%%%%%%%%%%%%%%%%%%%%%%%
%%%%%%%%%%%%Type_8a-3%%%%%%%%%%%%%%%%%%
%%%%%%%%%%%%%%%%%%%%%%%%%%%%%%%%%%%

\begin{center}
%
\end{center}

\medskip

We fix the perfect matching $\sfP_0$, 
then the following $\sfP_1, \cdots, \sfP_4$ are extremal perfect matchings corresponding to $v_1, \cdots, v_4$. 

\medskip

%%%%%perfect_matching%%%%%
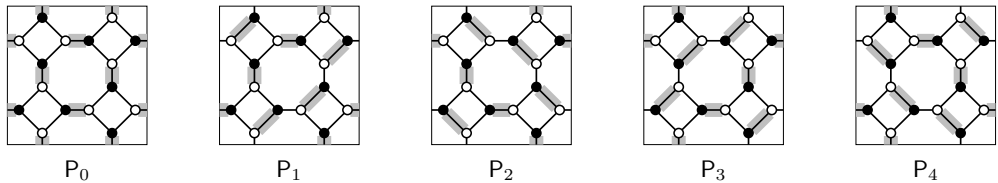
\begin{center}
{\scalebox{0.9}{
%
}}
\end{center}

\medskip

Then, we have the following table, and $e_0\sfC\cong e_4\sfC$, $e_1\sfC\cong e_5\sfC$, $e_2\sfC\cong e_6\sfC$, 
$e_3\sfC\cong e_7\sfC$, $(e_0\sfC)^*\cong e_2\sfC$, and $(e_1\sfC)^*\cong e_3\sfC$. 

\medskip

\begin{center}
{\small{
\begin{tabular}{|l|l|l|l|} \hline 
$e_0\sfC$&$e_1\sfC$&$e_2\sfC$&$e_3\sfC$ \\ \hline
$T_{00}\cong R$&$T_{10}\cong T(0,0,1,0)$&$T_{20}\cong T(2,1,1,0)$&$T_{30}\cong T(2,1,0,0)$ \\
$T_{01}\cong T(-2,0,1,0)$&$T_{11}\cong R$&$T_{21}\cong T(2,1,0,0)$&$T_{31}\cong T(0,1,1,0)$ \\
$T_{02}\cong T(-2,1,1,0)$&$T_{12}\cong T(0,1,0,0)$&$T_{22}\cong R$&$T_{32}\cong T(-2,0,1,0)$ \\
$T_{03}\cong T(0,1,0,0)$&$T_{13}\cong T(0,1,1,0)$&$T_{23}\cong T(0,0,1,0)$&$T_{33}\cong R$\\ 
$T_{04}\cong T(-1,0,1,0)$&$T_{14}\cong T(1,0,0,0)$&$T_{24}\cong T(3,1,0,0)$&$T_{34}\cong T(1,1,1,0)$ \\
$T_{05}\cong T(-1,0,0,0)$&$T_{15}\cong T(-1,0,1,0)$&$T_{25}\cong T(1,1,1,0)$&$T_{35}\cong T(1,1,0,0)$ \\
$T_{06}\cong T(-1,1,0,0)$&$T_{16}\cong T(-1,1,1,0)$&$T_{26}\cong T(-1,0,1,0)$&$T_{36}\cong T(-1,0,0,0)$ \\
$T_{07}\cong T(-1,1,1,0)$&$T_{17}\cong T(1,1,0,0)$&$T_{27}\cong T(1,0,0,0)$&$T_{37}\cong T(-1,0,1,0)$ \\ \hline 
$e_4\sfC$&$e_5\sfC$&$e_6\sfC$&$e_7\sfC$ \\ \hline
$T_{40}\cong T(-1,0,1,0)$&$T_{50}\cong T(1,0,0,0)$&$T_{60}\cong T(3,1,0,0)$&$T_{70}\cong T(1,1,1,0)$ \\
$T_{41}\cong T(-1,0,0,0)$&$T_{51}\cong T(-1,0,1,0)$&$T_{61}\cong T(1,1,1,0)$&$T_{71}\cong T(1,1,0,0)$ \\
$T_{42}\cong T(-1,1,0,0)$&$T_{52}\cong T(-1,1,1,0)$&$T_{62}\cong T(-1,0,1,0)$&$T_{72}\cong T(-1,0,0,0)$ \\
$T_{43}\cong T(-1,1,1,0)$&$T_{53}\cong T(1,1,0,0)$&$T_{63}\cong T(1,0,0,0)$&$T_{73}\cong T(-1,0,1,0)$\\ 
$T_{44}\cong R$&$T_{54}\cong T(0,0,1,0)$&$T_{64}\cong T(2,1,1,0)$&$T_{74}\cong T(2,1,0,0)$ \\
$T_{45}\cong T(-2,0,1,0)$&$T_{55}\cong R$&$T_{65}\cong T(2,1,0,0)$&$T_{75}\cong T(0,1,1,0)$ \\
$T_{46}\cong T(-2,1,1,0)$&$T_{56}\cong T(0,1,0,0)$&$T_{66}\cong R$&$T_{76}\cong T(-2,0,1,0)$ \\
$T_{47}\cong T(0,1,0,0)$&$T_{57}\cong T(0,1,1,0)$&$T_{67}\cong T(0,0,1,0)$&$T_{77}\cong R$ \\ \hline 
\end{tabular}
}} 
\end{center}

\medskip

\subsubsection{Type 8a-4}
Next, we consider the following consistent dimer model and the associated QP $(Q, W_Q)$. 
We denote by $\sfD$ the complete Jacobian algebra $\calP(Q, W_Q)$, and consider splitting MM generators $e_i\sfD$. 

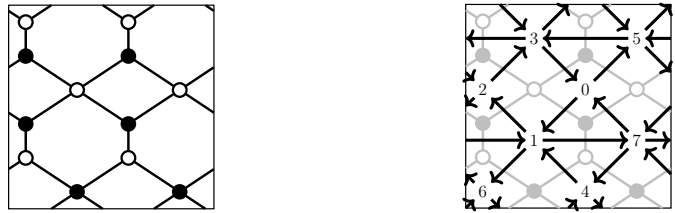
\medskip
%%%%%%%%%%%%%%%%%%%%%%%%%%%%%%%%%%%
%%%%%%%%%%%%Type_8a-4%%%%%%%%%%%%%%%%%%
%%%%%%%%%%%%%%%%%%%%%%%%%%%%%%%%%%%

\begin{center}
%
\end{center}

\medskip

We fix the perfect matching $\sfP_0$, 
then the following $\sfP_1, \cdots, \sfP_4$ are extremal perfect matchings corresponding to $v_1, \cdots, v_4$. 

\medskip

%%%%%perfect_matching%%%%%
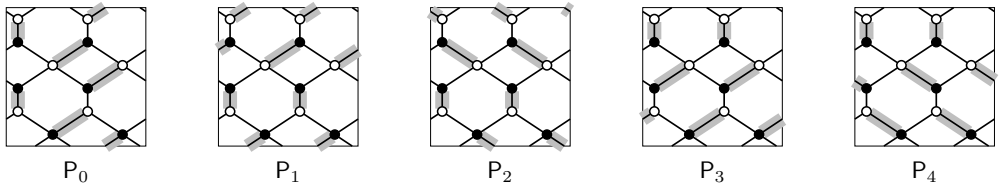
\begin{center}
{\scalebox{0.9}{
%
}}
\end{center}

\medskip

Then, we have the following table, and $e_0\sfD\cong e_2\sfD$, $e_1\sfD\cong e_7\sfD$, $e_3\sfD\cong e_5\sfD$, 
$e_4\sfD\cong e_6\sfD$, $(e_0\sfD)^*\cong e_4\sfD$, $(e_1\sfD)^*\cong e_1\sfD$, and $(e_3\sfD)^*\cong e_3\sfD$. 

\medskip

\begin{center}
{\small{
\begin{tabular}{|l|l|l|l|} \hline 
$e_0\sfD$&$e_1\sfD$&$e_2\sfD$&$e_3\sfD$ \\ \hline
$T_{00}\cong R$&$T_{10}\cong T(1,1,1,0)$&$T_{20}\cong T(1,1,0,0)$&$T_{30}\cong T(1,0,0,0)$ \\
$T_{01}\cong T(-1,1,1,0)$&$T_{11}\cong R$&$T_{21}\cong T(-2,0,1,0)$&$T_{31}\cong T(0,1,1,0)$ \\
$T_{02}\cong T(1,1,0,0)$&$T_{12}\cong T(0,0,1,0)$&$T_{22}\cong R$&$T_{32}\cong T(2,1,0,0)$ \\
$T_{03}\cong T(-1,0,0,0)$&$T_{13}\cong T(0,1,1,0)$&$T_{23}\cong T(0,1,0,0)$&$T_{33}\cong R$\\ 
$T_{04}\cong T(-1,1,0,0)$&$T_{14}\cong T(-2,0,1,0)$&$T_{24}\cong T(-2,0,0,0)$&$T_{34}\cong T(0,1,0,0)$ \\
$T_{05}\cong T(0,1,0,0)$&$T_{15}\cong T(-1,0,1,0)$&$T_{25}\cong T(-1,0,0,0)$&$T_{35}\cong T(1,1,0,0)$ \\
$T_{06}\cong T(-2,0,0,0)$&$T_{16}\cong T(-1,1,1,0)$&$T_{26}\cong T(-1,1,0,0)$&$T_{36}\cong T(-1,0,0,0)$ \\
$T_{07}\cong T(-2,0,1,0)$&$T_{17}\cong T(1,1,0,0)$&$T_{27}\cong T(-1,1,1,0)$&$T_{37}\cong T(-1,0,1,0)$ \\ \hline 
$e_4\sfD$&$e_5\sfD$&$e_6\sfD$&$e_7\sfD$ \\ \hline
$T_{40}\cong T(3,1,0,0)$&$T_{50}\cong T(2,1,0,0)$&$T_{60}\cong T(2,0,0,0)$&$T_{70}\cong T(0,0,1,0)$ \\
$T_{41}\cong T(0,0,1,0)$&$T_{51}\cong T(-1,0,1,0)$&$T_{61}\cong T(1,1,1,0)$&$T_{71}\cong T(1,1,0,0)$ \\
$T_{42}\cong T(2,0,0,0)$&$T_{52}\cong T(1,0,0,0)$&$T_{62}\cong T(3,1,0,0)$&$T_{72}\cong T(1,1,1,0)$ \\
$T_{43}\cong T(2,1,0,0)$&$T_{53}\cong T(1,1,0,0)$&$T_{63}\cong T(1,0,0,0)$&$T_{73}\cong T(-1,0,1,0)$\\ 
$T_{44}\cong R$&$T_{54}\cong T(-1,0,0,0)$&$T_{64}\cong T(1,1,0,0)$&$T_{74}\cong T(-1,1,1,0)$ \\
$T_{45}\cong T(1,0,0,0)$&$T_{55}\cong R$&$T_{65}\cong T(2,1,0,0)$&$T_{75}\cong T(0,1,1,0)$ \\
$T_{46}\cong T(1,1,0,0)$&$T_{56}\cong T(0,1,0,0)$&$T_{66}\cong R$&$T_{76}\cong T(-2,0,1,0)$ \\
$T_{47}\cong T(1,1,1,0)$&$T_{57}\cong T(0,1,1,0)$&$T_{67}\cong T(0,0,1,0)$&$T_{77}\cong R$ \\ \hline 
\end{tabular}
}} 
\end{center}

\medskip

Let $x$ (resp. $y$) be a $1$-cycle on the two torus $\sfT$ which corresponds to $(1,0)\in\rmH_1(\sfT)$ (resp. $(0,1)\in\rmH_1(\sfT)$) in the above case. 
By replacing these cycles by $y, -x$ respectively, we have the following dimer model and the associated QP $(Q, W_Q)$. 
We denote by $\sfD^\prime$ the complete Jacobian algebra $\calP(Q, W_Q)$. 
Since this dimer model is right equivalent to the previous one, we have that $\sfD\cong\sfD^\prime$ as an $R$-algebra. 
However, a change of cycles induces an automorphism on $R$, and it gives another description of each splitting MM generator (see Remark~\ref{rem_auto}). 
Thus we consider splitting MM generators $e_i\sfD^\prime$.

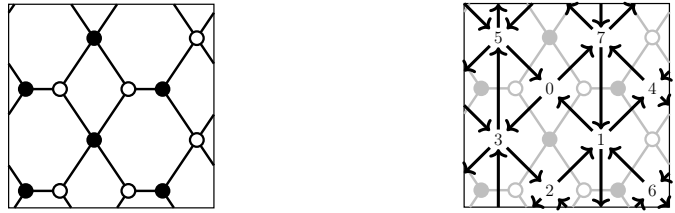
\medskip

%%%%%%%%%%%%%%%%%%%%%%%%%%%%%%%%%%%
%%%%%%%%%%%%Type_8a-4+%%%%%%%%%%%%%%%%%%
%%%%%%%%%%%%%%%%%%%%%%%%%%%%%%%%%%%

\begin{center}

\end{center}

\medskip

We fix the perfect matching $\sfP_0$. 
Then, extremal perfect matchings corresponding to $v_1, \cdots, v_4$ are following. 

\medskip

%%%%%perfect_matching%%%%%
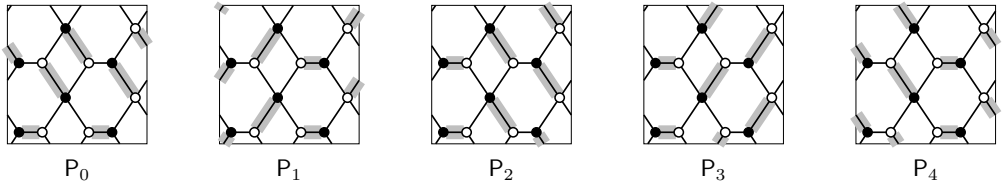
\begin{center}
{\scalebox{0.9}{
%
}}
\end{center}

\medskip

Then, we have the following table, and $e_0\sfD^\prime\cong e_2\sfD^\prime$, $e_1\sfD^\prime\cong e_7\sfD^\prime$, $e_3\sfD^\prime\cong e_5\sfD^\prime$, 
$e_4\sfD^\prime\cong e_6\sfD^\prime$, $(e_0\sfD^\prime)^*\cong e_4\sfD^\prime$, $(e_1\sfD^\prime)^*\cong e_1\sfD^\prime$, and $(e_3\sfD^\prime)^*\cong e_3\sfD^\prime$. 

\medskip

\begin{center}
{\small{
\begin{tabular}{|l|l|l|l|} \hline 
$e_0\sfD^\prime$&$e_1\sfD^\prime$&$e_2\sfD^\prime$&$e_3\sfD^\prime$ \\ \hline
$T_{00}\cong R$&$T_{10}\cong T(0,0,1,0)$&$T_{20}\cong T(0,1,1,0)$&$T_{30}\cong T(1,1,1,0)$ \\
$T_{01}\cong T(-2,0,1,0)$&$T_{11}\cong R$&$T_{21}\cong T(0,1,0,0)$&$T_{31}\cong T(1,1,0,0)$ \\
$T_{02}\cong T(0,1,1,0)$&$T_{12}\cong T(2,1,0,0)$&$T_{22}\cong R$&$T_{32}\cong T(1,0,0,0)$ \\
$T_{03}\cong T(-1,1,1,0)$&$T_{13}\cong T(1,1,0,0)$&$T_{23}\cong T(-1,0,0,0)$&$T_{33}\cong R$\\ 
$T_{04}\cong T(-2,1,1,0)$&$T_{14}\cong T(0,1,0,0)$&$T_{24}\cong T(-2,0,0,0)$&$T_{34}\cong T(-1,0,0,0)$ \\
$T_{05}\cong T(-1,0,0,0)$&$T_{15}\cong T(-1,0,1,0)$&$T_{25}\cong T(-1,1,1,0)$&$T_{35}\cong T(0,1,1,0)$ \\
$T_{06}\cong T(-2,0,0,0)$&$T_{16}\cong T(-2,0,1,0)$&$T_{26}\cong T(-2,1,1,0)$&$T_{36}\cong T(-1,1,1,0)$ \\
$T_{07}\cong T(0,1,0,0)$&$T_{17}\cong T(0,1,1,0)$&$T_{27}\cong T(-2,0,1,0)$&$T_{37}\cong T(-1,0,1,0)$ \\ \hline 
$e_4\sfD^\prime$&$e_5\sfD^\prime$&$e_6\sfD^\prime$&$e_7\sfD^\prime$ \\ \hline
$T_{40}\cong T(2,1,1,0)$&$T_{50}\cong T(1,0,0,0)$&$T_{60}\cong T(2,0,0,0)$&$T_{70}\cong T(2,1,0,0)$ \\
$T_{41}\cong T(2,1,0,0)$&$T_{51}\cong T(-1,0,1,0)$&$T_{61}\cong T(0,0,1,0)$&$T_{71}\cong T(0,1,1,0)$ \\
$T_{42}\cong T(2,0,0,0)$&$T_{52}\cong T(1,1,1,0)$&$T_{62}\cong T(2,1,1,0)$&$T_{72}\cong T(0,0,1,0)$ \\
$T_{43}\cong T(1,0,0,0)$&$T_{53}\cong T(0,1,1,0)$&$T_{63}\cong T(1,1,1,0)$&$T_{73}\cong T(-1,0,1,0)$\\ 
$T_{44}\cong R$&$T_{54}\cong T(-1,1,1,0)$&$T_{64}\cong T(0,1,1,0)$&$T_{74}\cong T(-2,0,1,0)$ \\
$T_{45}\cong T(1,1,1,0)$&$T_{55}\cong R$&$T_{65}\cong T(1,0,0,0)$&$T_{75}\cong T(1,1,0,0)$ \\
$T_{46}\cong T(0,1,1,0)$&$T_{56}\cong T(-1,0,0,0)$&$T_{66}\cong R$&$T_{76}\cong T(0,1,0,0)$ \\
$T_{47}\cong T(0,0,1,0)$&$T_{57}\cong T(1,1,0,0)$&$T_{67}\cong T(2,1,0,0)$&$T_{77}\cong R$ \\ \hline
\end{tabular}
}} 
\end{center}

\medskip

%%%%%%%%%%%%%%%%%%%%%%%%%%%%%%%%%%%%%%%%%%%%
%\newpage

\subsubsection{Exchange graph of type 8a}
By combining the above results, we describe the exchange graph $\EG(\TMMG(R))$ for the case of type 8a as follows, 
and it is actually connected.
%%%%%%%%%%%%%%%%%%%%%
%%% type 8a exchange graph %%%%%%%%%%%%%
%%%%%%%%%%%%%%%%%%%%%

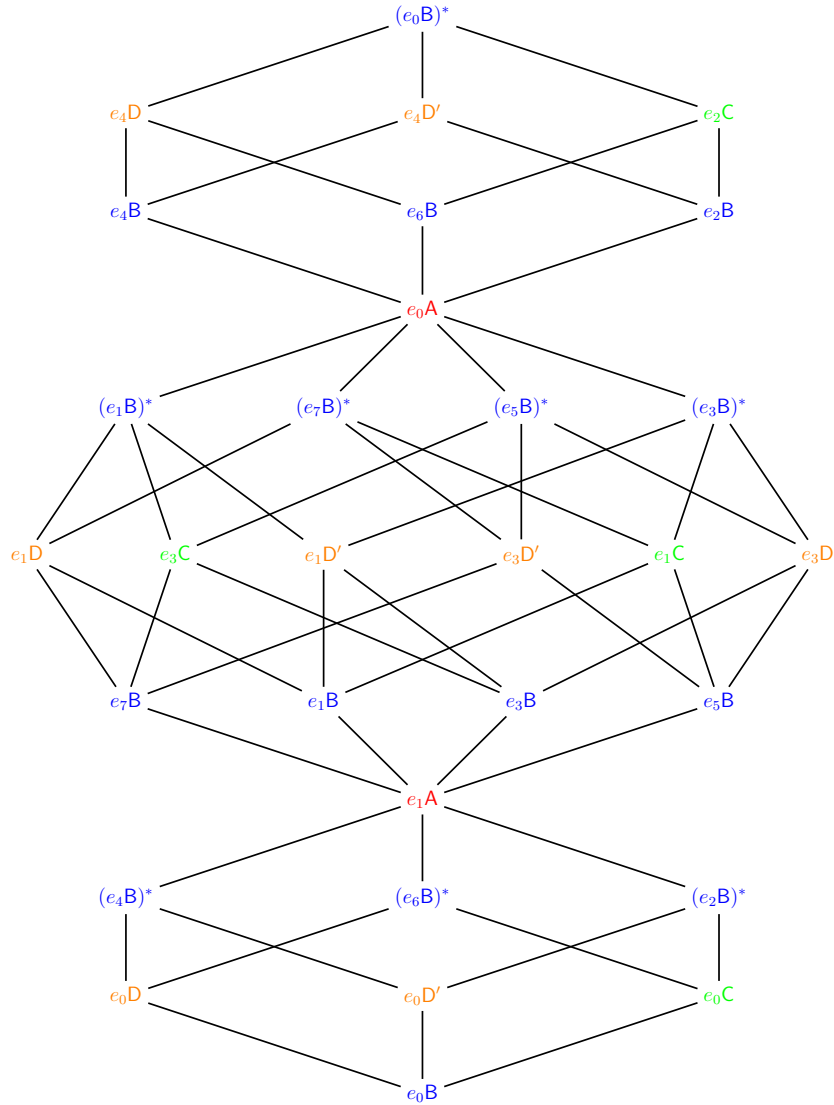
\begin{center}
{\scalebox{0.65}{
%
 }}
\end{center}

%%%%%%%%%%%%%%%%%%%%%%%%%%%%%%%%%%%%%%%%%%%%%%%%%%%%%%%
\subsection{Type 8b}
\label{type8b}
In this subsection, we consider the reflexive polygon of type 8b. 
Thus, let $R$ be the three dimensional complete local Gorenstein toric singularity defined by the cone $\sigma$: 
\[
\sigma=\mathrm{Cone}\{v_1=(0,1,1), v_2=(-1,1,1), v_3=(-1,-1,1), v_4=(2,-1,1) \}. 
\]
As an element in $\Cl(R)$, we obtain $[I_1]-2[I_3]+[I_4]=0$, $[I_2]+[I_3]-2[I_4]=0$, $2[I_3]+2[I_4]=0$.  
Therefore, we have that $\Cl(R)\cong\ZZ\times\ZZ/2\ZZ$, and each divisorial ideal is represented by $T(0,0,c,d)$ where $c\in\ZZ, d\in\ZZ/2\ZZ$. 
Also, there are two consistent dimer models (up to right equivalence) written below which give the reflexive polygon of type 8b as the perfect matching polygon. 

\medskip

\subsubsection{Type 8b-1}
We consider the following consistent dimer model and the associated QP $(Q, W_Q)$. 
We denote by $\sfA$ the complete Jacobian algebra $\calP(Q, W_Q)$, and consider splitting MM generators $e_i\sfA$. 

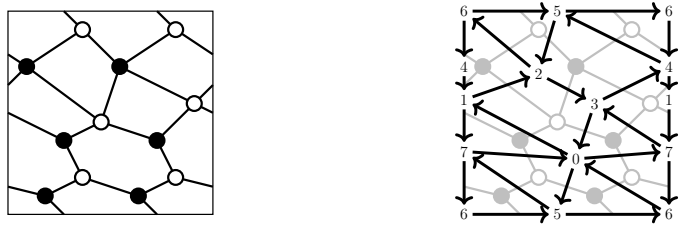
\medskip

%%%%%%%%%%%%%%%%%%%%%%%%%%%%%%%%%%%
%%%%%%%%%%%%Type_8b-1%%%%%%%%%%%%%%%%%%
%%%%%%%%%%%%%%%%%%%%%%%%%%%%%%%%%%%

\begin{center}

\end{center}

\medskip

We fix the perfect matching $\sfP_0$, 
then the following $\sfP_1, \cdots, \sfP_4$ are extremal perfect matchings corresponding to $v_1, \cdots, v_4$. 
Namely, $v_i=(h(\sfP_i, \sfP_0), 1)\in\ZZ^3$. 

\medskip

%%%%%perfect_matching%%%%%
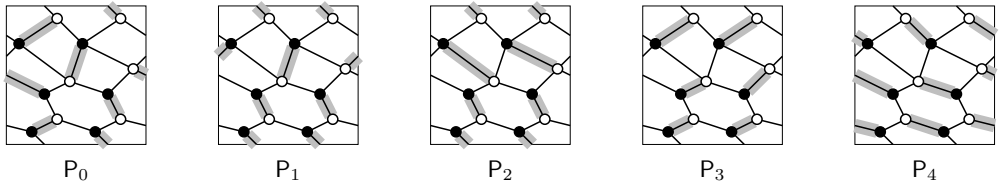
\begin{center}
{\scalebox{0.9}{
%
}}
\end{center}

\medskip

Then, we have the following table. 
%For simplicity, we denote by $T(c,d)$ a divisorial ideal $T(0,0,c,d)$. 
We remark that $e_0\sfA\cong e_7\sfA\cong (e_5\sfA)^*\cong (e_6\sfA)^*$ and $e_1\sfA\cong e_3\sfA\cong (e_2\sfA)^*\cong (e_4\sfA)^*$. 

\medskip

\begin{center}
{\small{
\begin{tabular}{|l|l|l|l|} \hline 
$e_0\sfA$&$e_1\sfA$&$e_2\sfA$&$e_3\sfA$ \\ \hline
$T_{00}\cong R$&$T_{10}\cong T(0,0,-1,0)$&$T_{20}\cong T(0,0,2,0)$&$T_{30}\cong T(0,0,0,1)$ \\
$T_{01}\cong T(0,0,1,0)$&$T_{11}\cong R$&$T_{21}\cong T(0,0,3,0)$&$T_{31}\cong T(0,0,1,1)$ \\
$T_{02}\cong T(0,0,-2,0)$&$T_{12}\cong T(0,0,-3,0)$&$T_{22}\cong R$&$T_{32}\cong T(0,0,-2,1)$ \\
$T_{03}\cong T(0,0,2,1)$&$T_{13}\cong T(0,0,1,1)$&$T_{23}\cong T(0,0,4,1)$&$T_{33}\cong R$ \\
$T_{04}\cong T(0,0,-1,1)$&$T_{14}\cong T(0,0,-2,1)$&$T_{24}\cong T(0,0,1,1)$&$T_{34}\cong T(0,0,-3,0)$ \\
$T_{05}\cong T(0,0,0,1)$&$T_{15}\cong T(0,0,-1,1)$&$T_{25}\cong T(0,0,2,1)$&$T_{35}\cong T(0,0,-2,0)$ \\ 
$T_{06}\cong T(0,0,-1,0)$&$T_{16}\cong T(0,0,-2,0)$&$T_{26}\cong T(0,0,1,0)$&$T_{36}\cong T(0,0,-1,1)$ \\ 
$T_{07}\cong T(0,0,1,1)$&$T_{17}\cong T(0,0,0,1)$&$T_{27}\cong T(0,0,3,1)$&$T_{37}\cong T(0,0,-1,0)$ \\ \hline
$e_4\sfA$&$e_5\sfA$&$e_6\sfA$&$e_7\sfA$ \\ \hline
$T_{40}\cong T(0,0,3,1)$&$T_{50}\cong T(0,0,2,1)$&$T_{60}\cong T(0,0,1,0)$&$T_{70}\cong T(0,0,1,1)$ \\
$T_{41}\cong T(0,0,4,1)$&$T_{51}\cong T(0,0,3,1)$&$T_{61}\cong T(0,0,2,0)$&$T_{71}\cong T(0,0,2,1)$ \\
$T_{42}\cong T(0,0,1,1)$&$T_{52}\cong T(0,0,0,1)$&$T_{62}\cong T(0,0,-1,0)$&$T_{72}\cong T(0,0,-1,1)$ \\
$T_{43}\cong T(0,0,3,0)$&$T_{53}\cong T(0,0,2,0)$&$T_{63}\cong T(0,0,3,1)$&$T_{73}\cong T(0,0,1,0)$ \\
$T_{44}\cong R$&$T_{54}\cong T(0,0,-1,0)$&$T_{64}\cong T(0,0,0,1)$&$T_{74}\cong T(0,0,-2,0)$ \\
$T_{45}\cong T(0,0,1,0)$&$T_{55}\cong R$&$T_{65}\cong T(0,0,1,1)$&$T_{75}\cong T(0,0,-1,0)$ \\ 
$T_{46}\cong T(0,0,2,1)$&$T_{56}\cong T(0,0,1,1)$&$T_{66}\cong R$&$T_{76}\cong T(0,0,0,1)$ \\ 
$T_{47}\cong T(0,0,2,0)$&$T_{57}\cong T(0,0,1,0)$&$T_{67}\cong T(0,0,2,1)$&$T_{77}\cong R$ \\ \hline
\end{tabular} 
}}
\end{center}

\medskip

\subsubsection{Type 8b-2}
Next, we consider the following consistent dimer model and the associated QP $(Q, W_Q)$. 
We denote by $\sfB$ the complete Jacobian algebra $\calP(Q, W_Q)$, and consider splitting MM generators $e_i\sfB$. 

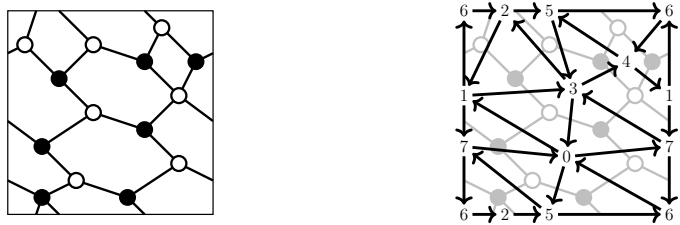
\medskip 
%%%%%%%%%%%%%%%%%%%%%%%%%%%%%%%%%%%
%%%%%%%%%%%%Type_8b-2%%%%%%%%%%%%%%%%%%
%%%%%%%%%%%%%%%%%%%%%%%%%%%%%%%%%%%
\begin{center}
%
\end{center}

\medskip

We fix the perfect matching $\sfP_0$, 
then the following $\sfP_1, \cdots, \sfP_4$ are extremal perfect matchings corresponding to $v_1, \cdots, v_4$. 

\medskip

%%%%%perfect_matching%%%%%
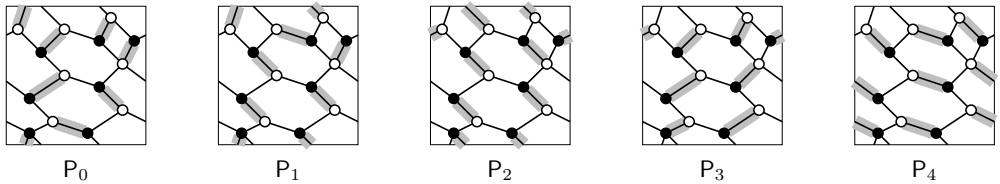
\begin{center}
{\scalebox{0.9}{
%
}}
\end{center}

\medskip

Then, we have the following table. 
%For simplicity, we denote by $T(c,d)$ a divisorial ideal $T(0,0,c,d)$. 
We remark that $e_0\sfB\cong (e_0\sfB)^*$, $e_1\sfB\cong (e_6\sfB)^*$, $e_2\sfB\cong (e_4\sfB)^*$, $e_3\sfB\cong (e_5\sfB)^*$, and $e_7\sfB\cong (e_7\sfB)^*$. 

\medskip

\begin{center}
{\small{
\begin{tabular}{|l|l|l|l|} \hline 
$e_0\sfB$&$e_1\sfB$&$e_2\sfB$&$e_3\sfB$ \\ \hline
$T_{00}\cong R$&$T_{10}\cong T(0,0,-1,0)$&$T_{20}\cong T(0,0,-1,1)$&$T_{30}\cong T(0,0,0,1)$ \\
$T_{01}\cong T(0,0,1,0)$&$T_{11}\cong R$&$T_{21}\cong T(0,0,0,1)$&$T_{31}\cong T(0,0,1,1)$ \\
$T_{02}\cong T(0,0,3,1)$&$T_{12}\cong T(0,0,2,1)$&$T_{22}\cong R$&$T_{32}\cong T(0,0,1,0)$ \\
$T_{03}\cong T(0,0,2,1)$&$T_{13}\cong T(0,0,1,1)$&$T_{23}\cong T(0,0,-1,0)$&$T_{33}\cong R$ \\
$T_{04}\cong T(0,0,-1,1)$&$T_{14}\cong T(0,0,-2,1)$&$T_{24}\cong T(0,0,-4,0)$&$T_{34}\cong T(0,0,-3,0)$ \\
$T_{05}\cong T(0,0,0,1)$&$T_{15}\cong T(0,0,-1,1)$&$T_{25}\cong T(0,0,-3,0)$&$T_{35}\cong T(0,0,-2,0)$ \\ 
$T_{06}\cong T(0,0,-1,0)$&$T_{16}\cong T(0,0,-2,0)$&$T_{26}\cong T(0,0,-2,1)$&$T_{36}\cong T(0,0,-1,1)$ \\ 
$T_{07}\cong T(0,0,1,1)$&$T_{17}\cong T(0,0,0,1)$&$T_{27}\cong T(0,0,-2,0)$&$T_{37}\cong T(0,0,-1,0)$ \\ \hline
$e_4\sfB$&$e_5\sfB$&$e_6\sfB$&$e_7\sfB$ \\ \hline
$T_{40}\cong T(0,0,3,1)$&$T_{50}\cong T(0,0,2,1)$&$T_{60}\cong T(0,0,1,0)$&$T_{70}\cong T(0,0,1,1)$ \\
$T_{41}\cong T(0,0,4,1)$&$T_{51}\cong T(0,0,3,1)$&$T_{61}\cong T(0,0,2,0)$&$T_{71}\cong T(0,0,2,1)$ \\
$T_{42}\cong T(0,0,4,0)$&$T_{52}\cong T(0,0,3,0)$&$T_{62}\cong T(0,0,4,1)$&$T_{72}\cong T(0,0,2,0)$ \\
$T_{43}\cong T(0,0,3,0)$&$T_{53}\cong T(0,0,2,0)$&$T_{63}\cong T(0,0,3,1)$&$T_{73}\cong T(0,0,1,0)$ \\
$T_{44}\cong R$&$T_{54}\cong T(0,0,-1,0)$&$T_{64}\cong T(0,0,0,1)$&$T_{74}\cong T(0,0,-2,0)$ \\
$T_{45}\cong T(0,0,1,0)$&$T_{55}\cong R$&$T_{65}\cong T(0,0,1,1)$&$T_{75}\cong T(0,0,-1,0)$ \\ 
$T_{46}\cong T(0,0,2,1)$&$T_{56}\cong T(0,0,1,1)$&$T_{66}\cong R$&$T_{76}\cong T(0,0,0,1)$ \\ 
$T_{47}\cong T(0,0,2,0)$&$T_{57}\cong T(0,0,1,0)$&$T_{67}\cong T(0,0,2,1)$&$T_{77}\cong R$ \\ \hline
\end{tabular} 
}}
\end{center}

\medskip

%%%%%%%%%%%%%%%%%%%%%%%%%%%%%%%%%%%%%%%%%%%%
\subsubsection{Exchange graph of type 8b}
By combining the above results, we describe the exchange graph $\EG(\TMMG(R))$ for the case of type 8b as follows, 
and it is actually connected.

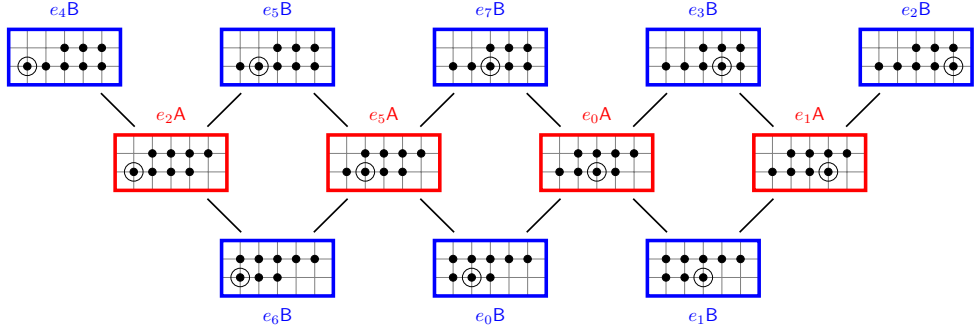
\begin{center}
{\scalebox{0.7}
{%
}};

\draw[line width=0.032cm] (A0)--(B3); \draw[line width=0.032cm] (A0)--(B1); \draw[line width=0.032cm] (A0)--(B7); \draw[line width=0.032cm] (A0)--(B0);
\draw[line width=0.032cm] (A5)--(B7); \draw[line width=0.032cm] (A5)--(B0); \draw[line width=0.032cm] (A5)--(B5); \draw[line width=0.032cm] (A5)--(B6);
\draw[line width=0.032cm] (A1)--(B3); \draw[line width=0.032cm] (A1)--(B1); \draw[line width=0.032cm] (A1)--(B2); \draw[line width=0.032cm] (A2)--(B5);
\draw[line width=0.032cm] (A2)--(B6); \draw[line width=0.032cm] (A2)--(B4);

\end{tikzpicture}}}
\end{center}

%\fi
%%%%%%%%%%%%%%%%%%%%%%%%%%%%%%%%%%%%%%%%%%%%%%%%%%%%%%%
%%%%%%%%%%%%%%%%%%%%%%%%%%%%%%%%%%%%%%%%%%%%%%%%%%%%%%%
\section{Mutations of splitting MM modules} 
\label{sec_toric_mutation_module}

We end this paper by discussing mutations of splitting MM modules. 
Recall that we denote the set of isomorphism classes of basic splitting MM $R$-modules by $\TMM(R)$, 
and denote the exchange graph of $\TMM(R)$ by $\EG(\TMM(R))$. 
We start our discussion in a slightly general situation. 

\begin{lemma}
\label{lem_for_MM}
Let $R$ be a complete local $3$-sCY normal domain, and $I$ is a rank one reflexive $R$-module.
\begin{enumerate}[\rm(1)]
\item For $M\in\TMM(R)$, we have that $(M\otimes_RI)^{**}\in\TMM(R)$. 
\item Suppose that $M, N\in\TMM(R)$, and they are connected in $\EG(\TMM(R))$. Then $(M\otimes_RI)^{**}$ and $(N\otimes_RI)^{**}$ 
are also connected in $\EG(\TMM(R))$. 
\item For $M\in\TMM(R)$, we have that $M\cong(N\otimes_RJ)^{**}$ for some $N\in\TMMG(R)$ and some rank one reflexive module $J$ 
with $J\in\add_RM$. 
\end{enumerate}

\end{lemma}

\begin{proof}
(1) For $M^\prime\coloneqq(M\otimes_RI)^{**}$, we have that $\End_R(M)\cong\End_R(M^\prime)$. 
Indeed, localizing a natural morphism $\End_R(M)\rightarrow\End_R(M^\prime)$ between reflexive $R$-modules at any prime ideal $\fkp$ with ${\rm ht}\,\fkp=1$, 
we have a morphism between free $R_\fkp$-modules with the same rank, and it is isomorphic in particular. 
This induces an isomorphism of the original one (see e.g., \cite[Lemma~5.11]{LW}). 
Thus, we have the assertion. 

(2) We may assume that $M$ and $N$ are connected by a single mutation of splitting MM modules. 
Then, $M$ and $N$ have the same direct summands except one component 
(see the last part of subsection~\ref{subsec_MMmutation}), and this relationship holds for $(M\otimes_RI)^{**}$ and $(N\otimes_RI)^{**}$. 

(3) See \cite[Lemma~5.22]{IW3}. 
\end{proof}

By using these results, we have the next proposition. 

\begin{proposition}
\label{prop_MMmod}
Let $R$ be a complete local $3$-sCY normal domain. Assume that $\Cl(R)$ is generated by $[I_1],\cdots,[I_t]$ and $\EG(\TMMG(R))$ is connected. 
If there is a splitting MM generator $N_j$ such that $(N_j\otimes_RI_j)^{**}$ is connected to $\EG(\TMMG(R))$ for every $j=1,\cdots,t$, 
then $\EG(\TMM(R))$ is connected. 
\end{proposition}

\begin{proof}
In this proof, for $M,M^\prime\in\TMM(R)$, we denote $M\sim M^\prime$ if they are connected in $\EG(\TMM(R))$.  

Let $M$ be a basic splitting MM module. By Lemma~\ref{lem_for_MM}(3), we can describe it as $M\cong(N\otimes_RI)^{**}$ 
where $N\in\TMMG(R)$ and $[I]\in\Cl(R)$. 
Since $\Cl(R)$ is generated by $[I_1],\cdots,[I_t]$, we have that $I\cong(I_1^{\otimes n_1}\otimes\cdots\otimes I_t^{\otimes n_t})^{**}$ for some $n_1,\cdots,n_t\in\ZZ_{\ge0}$. 
Since $\EG(\TMMG(R))$ is connected, we have that $N\sim N_j$ for any $j$. 
Thus, we have that $(N\otimes_R I_j)^{**}\sim(N_j\otimes_R I_j)^{**}$ for any $j$ by Lemma~\ref{lem_for_MM}(2). 
By the assumption, we have that $N_j\sim (N_j\otimes_R I_j)^{**}$, and therefore we conclude $N\sim (N\otimes_R I_j)^{**}$ for any $j$. 
Therefore, if $N\sim (N\otimes I_1^{\otimes m_1}\otimes\cdots\otimes I_j^{\otimes m_j}\otimes\cdots\otimes I_t^{\otimes m_t})^{**}$, then 
we have that $N\sim (N\otimes_R I_j)^{**}\sim (N\otimes I_1^{\otimes m_1}\otimes\cdots\otimes I_j^{\otimes m_j+1}\otimes\cdots\otimes I_t^{\otimes m_t})^{**}$. 
By an inductive argument, we have that $M\sim N$. Thus, we can connect $M$ to $\EG(\TMMG(R))$. 
\end{proof}
 
Finally, we have the connectedness of the exchange graph of splitting MM modules for three dimensional Gorenstein toric singularities associated with reflexive polygons. 

\begin{corollary}
\label{mutation_MMmod}
Let $R$ be a three dimensional complete local Gorenstein toric singularity associated with a reflexive polygon. 
Then, the exchange graph $\EG(\TMM(R))$ of splitting MM modules is connected. 
\end{corollary}

\begin{proof}
We can see that the conditions in Proposition~\ref{prop_MMmod} hold for this toric singularity. 
Indeed, the connectedness of $\EG(\TMMG(R))$ holds by Theorem~\ref{main_reflexive}. 
Also, by checking $\EG(\TMMG(R))$ described in Section~\ref{mutation_reflexive} for each case, 
we can easily see that the existence of a splitting MM generator as in Proposition~\ref{prop_MMmod} for each generator of $\Cl(R)$. 
For example, in the case of type 4a (see subsection~\ref{type4a}), $\Cl(R)$ is generated by $I_1\coloneqq T(1,0,0,0)$ and $I_2\coloneqq T(0,1,0,0)$, 
and we see that $(e_0\sfA\otimes_RI_1)^{**}\cong e_1\sfA$, $(e_0\sfB\otimes_RI_2)^{**}\cong e_1\sfB$. 
Thus, we have the assertion. 
\end{proof}

We finish this paper with the following natural question. 

\begin{question}
Is the exchange graph $\EG(\TMM(R))$ of splitting MM modules connected for any three dimensional Gorenstein toric singularity?
\end{question}

%%%%%%%%%%%%%%%%%%%%%%%%%%%%%%%%%%%%%%%%%%%%%%%%%%%%%%%%%%%%%%%%%
\subsection*{Acknowledgements}
The author would like to thank Professor Osamu Iyama for suggesting to study this subject, 
and for valuable discussions. 
The author also would like to thank Akihiro Higashitani for helpful comments about reflexive polygons, 
and thank Professor Richard Eager for explaining about a quiver gauge theory and Seiberg duality. 
The author would like to thank the anonymous referee for valuable comments on the paper. 

This work is supported by Grant-in-Aid for JSPS Fellows No. 26-422. 

%%%%%%%%%%%%%%%%%%%%%%%%%%%%%%%%%%%%%%%%%%%%%%%%%%%%%%%%%%%%%%%%%

%%%%%---reference---%%%%%

\end{document}